\documentclass[12pt]{amsart} 
\usepackage[all]{xy}
\usepackage{graphicx}

 \newtheorem{theorem}{Theorem}
 \newtheorem{lemma}[theorem]{Lemma}
  \newtheorem{cor}[theorem]{Corollary}
 \newtheorem{definition}[theorem]{Definition}
  \newtheorem{prop}[theorem]{Proposition}
 \newtheorem{remark}[theorem]{Remark}
\numberwithin{theorem}{section}

\usepackage{amssymb}
\usepackage{graphics,epsfig}

\begin{document}

\newcommand{\C}[1]{\mbox{$C_#1(\widetilde X _n / StN_n)$}}
\newcommand{\autfn}{\mbox{${\text {Aut}} (F_n)$}}
\newcommand{\xt}{\mbox{$\widetilde X$}}
\newcommand{\cp}{\mbox{$\mathcal{P}$}}
\newcommand{\cb}{\mbox{$\mathcal{B}$}}
\newcommand{\cs}{\mbox{$\mathcal{S}$}}
\newcommand{\cf}{\mbox{$\mathcal{F}$}}
\newcommand{\cc}{\mbox{$\mathcal{C}$}}
\newcommand{\ct}{\mbox{$\mathcal{T}$}}
\newcommand{\cd}{\mbox{$\mathcal{D}$}}
\newcommand{\sigt}{\mbox{$\tilde{\sigma}$}}
\newcommand{\wt}{\mbox{$\widetilde W$}}
\newcommand{\tp}{\mbox{$\widetilde{\cp}$}}
\newcommand{\lda}{\mbox{$\langle\!\langle$}}
\newcommand{\rda}{\mbox{$\rangle\!\rangle$}}
\newcommand{\cq}{\mbox{$\mathcal{Q}$}}

\centerline{\large{Complexes with free actions by the group of automorphisms of a free group}}
\bigskip\bigskip
\centerline{Jeff Kiralis}
\bigskip\bigskip


\section{Introduction}  We define two finite-dimensional complexes on which the group $\autfn$ acts freely,
and study their homotopy groups. The definitions of these complexes were motivated by low-dimensional pseudo-isotopy
theory.  So not surprisingly, the complexes are related  to  a number of spaces arising in algebraic K-theory, and
to another in low-dimensional topology.

The complexes are denoted by $X_n$ and $X_n^{alg}$.  
The cells of $ X_n$
  correspond bijectively to certain pointed graphs which are marked by homotopy
equivalences with $\bigvee_{i=1}^n S^1$,
the wedge of $n$ circles.  Each of these marked graphs is required 
to have   some additional structure which we
call a coloring.  The action of $\autfn$ on $X_n$ is defined by precomposing the markings with maps
$\bigvee_{i=1}^n S^1\to \bigvee_{i=1}^n S^1$ corresponding to elements of $\autfn$.
The complex $X_n^{alg}$ is defined  by 
algebraically mimicing the definition of $X_n$.
There is an equivarient map $X_n\to X_n^{alg}$ which, after stabalizing, may be a homotopy equivalence.
 
Forgetting the colorings defines an equivarient  map from $X_n$ to the sphere complex $\mathcal S_n$. 
This last complex was introduced by Hatcher  in \cite{H}. 
Forgetting basepoints as well gives a map from  $X_n$ to
Outer space. Outer space was defined by Culer and Vogtmann in \cite{CV},
and denoted by $X_n$  there. The colorings provide just
enough structure to eliminate the fixed points of $\mathcal
S_n$, which $\autfn$ acts on, but not freely. From the
point of view of pseudo-isotopy theory, the cells  of $
X_n$ may be thought of as parametrized families  of certain
handle-body structures of the 4-disk, or equivalently, of smooth
functions on $D^4$.  The colorings of the graphs prevent certain
higher singularities from occuring among these
functions so that the homotopy groups of $X_n$ resemble a $K$-theory.
There is no such restriction on the singularities  with $\mathcal
S_n$, which may be thought of as a sort of low-dimensional Whitehead
space, with a little added structure.

The space $X=\lim_{n\to\infty} X_n$ (and probably $X^{alg}=\lim_{n\to\infty} X^{alg}_n$ as well) should perhaps be thought of as a pre-K-theory, with the corresponding K-theory indicated in the remarks at the end of \S\ref{Xn}.  The homotopy groups $\pi_i(X_n)$ depend on $n$, even for $i>\!\!>n$, but seem to be more amenable to computation using elementary techniques than the probably stable homotopy
groups of the corresponding K-theory.

The homotopy groups of $X$ and
$X^{alg}$ map to the algebraic $K$-groups of the ring $\mathbb Z$, 
essentially by abelianizing
$F_n$, or, with pseudo-isotopy theory in mind, by suspending the families of smooth functions mentioned above. 
More specifically, there is a
map $h\colon X^{alg}\to V$, where $V\equiv V(\mathbb Z)$ is the space 
with $\pi_{n-1}(V)=K_n(\mathbb Z)$ defined by Volodin,  as in \cite{S}. 
%

The complex $X_n^{alg}$ is, in some ways, similar to the unstable Volodin space $V_n$, but with 
$\autfn$ in place of ${\rm GL}_n(\mathbb Z)$.  A  difference is that,  while ${\tilde V}/{\rm St}(\mathbb Z)$ is acylic,
$\widetilde X^{alg}/{\rm StN}$ is not. 
 (Here the tildes denote universal covers of a component of $V$ and $X$.
Also,  the group $\rm StN$,  introduced by Gersten in \cite{G}, and also studied in \cite{KN},
is a nonabelian version of ${\rm St}({\mathbb Z})$, the Steinberg group of $\mathbb Z$.)
 Indeed, the homology group $H_2(\xt^{alg}/{\rm StN})$ is nontrivial,
as is $H_2(\xt/{\rm StN})$. The elements of these groups detect elements of 
$\pi_2(X)$ 
which map onto the cokernel of the map $ \pi _3^S \to K_3(\mathbb Z)$.

There is also a map $f$ which makes the diagram 
$$\xymatrix{\Omega B {\rm StN}^+             \ar[d]^g \ar[r]^-f & \xt^{\text{alg}} \ar[d]^{\tilde h} \\
          \Omega   B {\rm St}({\mathbb Z})^+   \ar[r]^{\quad\sim}      & \tilde V} $$
 homotopy commute,
where $\tilde h$  is induced by $h\colon X_{alg}\to V$. 
The lower horizontal map is a homotopy equivalence coming from the one between $V$ and $\Omega B {\rm GL(\mathbb Z)}^+$
mentioned above. 

The main results, summarized in Theorem \ref{thm}, are the exactness of the sequence 
$$0\longrightarrow  H_3({\rm StN}) \longrightarrow  H_2(\tilde X)_{StN}    
{\buildrel q_*\over \longrightarrow} H_2(\tilde X/{\rm StN})    \longrightarrow 0 $$
and a computation of its terms.  In particular, the group $H_3({\rm StN})$
 is cyclic of order 24.  While this is analogous to Lee and Szarbo's result in \cite{LS} that the group $H_3({\rm St}(\mathbb Z))
\approx  K_3(\mathbb Z)$  is cyclic of order 48, our methods are completely different from  theirs.  
We show in \S\ref{ECT} and the following two that the group $H_3(\tilde X_n/StN_n)$ is trivial, 
which  gives the injectivily of the map  $H_3({\rm StN}) \longrightarrow H_2(\tilde X)_{StN} $.  
The homomorphism $q_*$ in the sequence is easily computed using the isomorphism $\Theta$
defined in \S\ref{H2xmodstn},  where
  $H_2(\tilde X_n/StN_n)$ is also computed.
The contractibility of the sphere complex is used in \S\ref{defP} to show that every 2-cycle in 
$\tilde X$ is homologous to one of a relatively simple type.  Such cycles are worked with in  
\S\ref{pictures} and \S\ref{lastcomp} to compute
$H_2(\tilde X)_{StN} $.
Along the way, we give many explicit descriptions of generators of $ K_3(\mathbb Z)$ as 2-cycles 
in $\xt$ (or $\tilde V$), and  show in an elementary way that they have order 48.
We also point out in Theorem \ref{last_thm} that these generators appear as nontrivial elements of the homology of certain  upper triangular subgroups of Aut($F_n$).

My main motivation for writing this paper was to compute $H_3({\rm StN})=H_3({\rm Aut}(F_n))=
\pi_3(B{\rm Aut}(F_n)^+)$, for large $n$.  But Galatius has recently announced that the homology 
groups of ${\rm Aut}(F_n)$ and the symmetric group are the same stably in {\em all} dimensions.
The spaces $\xt$ and $\xt/StN$ may still have some  interest.


\section{Preliminaries}

\subsection{Partitioned graphs} \quad In this paper, the word graph  means a  connected 1-dimensional $CW$-complex with finitely may vertices.  A graph $G$ will often be specified combinatorially by two sets $E(G)$ and $V(G)$, the edges and 
 vertices of $G$, respectively, along with two functions $d_0,d_1\colon E(G)\to V(G)$ which  determine
  the attaching maps of the 1-cells $E(G)$ of $G$.
The  rank of a graph is, by definition, the rank of its fundamental group.  
We often identify cycles in a graph $G$, as well as subgraphs of $G$ 
  (with each component larger than a vertex), with the subset of $E(G)$
  they comprise.
  
  A partition of a graph $G$ of rank $n$ is a 
partition of the set $E(G)$ of edges of $G$ into $n$ nonempty sets $h_1,\ldots,h_n$ such that any nontrivial cycle of $G$ contains one of the $h_i$.  The sets $h_i$ are called the edge sets of the
partitioned graph $G$. If $e$ is the only  element of an $h_i$, then $e$ is called a singleton edge,
  and $h_i$ a singleton edge set. 
  We often write $(G,\{h_i\},n)$, or just $(G,\{h_i\})$, for the partitioned graph $G$ with edge sets $h_i,\dots ,h_n$.
%

  If $G$ is a graph with edge $e$, let $G_e$ denote the graph obtained from $G$ by collapsing
$e$ to a point.  That is, $G_e$ is the graph with $E(G_e)=E(G) - \{e\}$ and $V(G_e)=V(G)/d_0(e)\sim d_1(e)$, and with attaching maps those of $G$, followed by the quotient map $V(G)\to V(G_e)$.

\begin{lemma} \label{bdspar}  Let $(G,\{h_i\},n)$ be a partioned graph.  If $e \in h_k$ where
$|h_k|>1$, then
 the graph $G_e$ is partitioned by the edge sets $h_1,\ldots,h_{k-1},h_k',h_{k+1},\ldots, h_n$ where $h_k'=h_k -\{e\}$,
  and the quotient map $f\colon G\to G_e$ is a homotopy equivalence. 
\end{lemma}
 
 We say that the partition of $G_e$ in the lemma is induced from the partition of  $(G,\{h_i\})$.

 \begin{proof}[Proof of  \ref{bdspar}]
Since $|h_k|>1$, the edge $e$ cannot itself form a cycle of length one.  For 
such a cycle would not contain any of the edgesets. Thus the collapsing map $f\colon G\to G_e$ is a
homotopy equivalence.

Any cycle $C_e$ of $G_e$ corresponds, via $f$, to a
unique cycle $C$ of $G$.  Since $h_i \subset C$ for some $i$ we have 
either $h_i\subset C_e$, if $i\ne k$, or $h_k'\subset C_e$ if $i=k$.
\end{proof}

 Quotient maps $f\colon G\to G_e$ between partitioned graphs 
 where the cardinality of the edge set containing $e$ is $> 1$, and where $G_e$ has the 
 induced partition, 
  are called blowdowns, as are compositions of such maps.   If $f \colon G\to G'$
is a blowdown, we say that $G'$ is obtained from $G$ by blowing down the edges of $G$ which are 
collapsed by $f$.

An edge basis of a partitioned graph $(G,\{h_i\},n)$ 
is, by definition,  any set of $n$ edges $\{e_1,\ldots ,e_n\}$ where $e_i\in h_i$ for each $i$.  
If $\{e_i\}$ is an edge basis for $(G,\{h_i\})$, we always assume, without comment, that 
$e_i\in h_i$ for each $i$.  

The complement
of any edge basis is a maximal tree.  More precisely
we have:

\begin{lemma} \label{maxtree}
Assume that $E=\{e_1,\ldots ,e_n\}$ is an edge basis for the partitioned graph $(G,\{h_i\})$.  Then
the subgraph $T^E$ of $G$ consisting of all the edges of $G$ except for those in $E$ is a maximal
tree in $G$.
\end{lemma}

\begin{proof}
The subgraph $T=T^E$ has no cycles since any cycle in $G$ contains an edge set $h_i$ for some $i$.

By Lemma \ref{bdspar}, all the edges of $T$ can be blown down giving a graph we call $R$.  Since $R$ 
and $G$ are homotopy equivalent, we have
$$\chi(R)=\chi(G)=\text {rank\ of}\ H_0(G) - \text {rank\ of}\ H_1(G)= 1-n.$$
Since $R$ has $n$ edges, it follows that it has just one vertex.  Therefore $T$ is connected and
contains all the vertices of $G$.  So $T$ is maximal.
\end{proof}
   
Assume that  $E=\{e_i\}$ is an edge basis for the partitioned graph $(G,\{h_i\})$, and let $T^E$ be 
the corresponding maximal tree as defined in the lemma.
Let $G_i$ be the subgraph of $G$ formed by adding the single edge $e_i$ to  $T^E$.
Denote by $C_i^E$, or just $C_i$, the unique (shortest) cycle in $G_i$.  Removing any edge in $h_i$ from $G_i$ gives,
by  Lemma \ref{maxtree}, a tree.  It follows that $h_i\subset C_i$.  Note that $C_i$ contains 
none of the other edge sets $h_j$ for $j\ne i$ since $C_i$  contains none of the edges $e_j$ for $j\ne i$.

The $C_i$ are called basic cycles (associated with the edge basis $E$). 
If a vertex $*$ of $G$ is chosen as a basepoint, then the basic cycles
 determine a basis for the free group $\pi_1(G,*)$.

Some basic facts about basic cycles are given in the next two lemmas.

\begin{lemma} \label{conn}
If $C_j$ and $C_k$ are two basic cycles associated with the edge basis $E=\{e_i\}$, then $C_j \cap C_k$ 
is connected.
\end{lemma}

\begin{proof}
Let $T$ be the maximal tree associated with the edge basis $E$ and let $x,y\in T$.  Denote by $[x,y]$
the smallest connected set in $T$ which contains both $x$ and $y$.   Note that if $u,v\in [x,y]$,
then $[u,v]\subset [x,y]$.  So if $u,v\in C_j \cap C_k= (C_j \cap T) \cap (C_k \cap T)$, then $[u,v]
\subset C_j \cap C_k$ since each set $C_i \cap T$ is of the form $[x_i,y_i]$ (where $x_i$ and $y_i$
are the vertices incident with the edge $e_i$).
 \end{proof}

\begin{lemma} \label{lemma1}
If $C_i$ and $C_j$ are basic cycles associated to the edge basis $E$ of the partitioned graph $(G,\{h_i\})$,
then at least one of the sets of edges $C_i\cap h_j$ or $C_j\cap h_i$ is empty.
\end{lemma}

\begin{proof}
Suppose that $e_j' \in C_i\cap h_j$ and $ e_i'\in C_j\cap h_i$.  Then both of $e_i'$ and $e_j'$ are
contained in the connected set $C_i\cap C_j$.  So the complement of the edge basis 
$E\cup \{e_i',e_j'\} -\{e_i,e_j\}$ contains a cycle, namely  
  $C_i\cup C_j- C_i\cap C_j$.  This contradicts Lemma \ref{maxtree}.
\end{proof}

\subsection{Oriented, Partitioned Graphs}
An edge basis is said to be oriented if each edge of it is oriented.  If $E=\{e_i\}$ is an oriented edge
basis of a partitioned graph $(G,\{h_i\})$, then the orientation of each $e_i$ determines an orientation 
of the corresponding basic cycle $C_i^E$, and hence of all the edges in $\{h_i\}$.  Thus, 
  since $E(G)= \coprod h_i$,   $E$ induces an
orientation of every edge of $G$.  A partitioned graph is said to be oriented if all of its edges 
are oriented by some oriented edge basis in this way.  

If a partitioned graph $G$ has rank $n$ then there 
are just $2^n$  possible orientations of $G$ as the next lemma shows.

\begin{lemma}\label{or}
Assume that $G$ is a partitioned graph which is oriented by the oriented edge basis $E$.  
If the orientation of $G$ is used to orient any edge basis $F$ of
$G$, then the orientations of $G$ induced by $F$ and $E$ agree.
\end{lemma} 

\begin{proof}
If the basic cycles $C_k^E$ and $C_k^F$ are not equal for some $k$, then $E\cap h_j\ne F\cap h_j$ 
for some $j\ne k$.  For such a $j$,  replace $e_j\in
E\cap h_j$ with 
$e_j'\in F\cap h_j$ to obtain a new edge basis $E'$. If $C_k^E\ne C_k^{E'}$, then $e_j'$ is contained in 
$C_k^E\cap C_j^E$.  In this case  $C_k^E\cap C_k^{E'}=C_k^E - (C_k^E\cap C_j^E)$ is connected since $C_k^E\cap C_j^E$
is by Lemma \ref{conn}.  Since $h_k\subset C_k^E\cap C_k^{E'}$, the orientations of the edges in $h_1$
determined by $E$ and $E'$ agree.  Continuing in this way shows that the orientations 
determined by $E$ and $F$ of all the edges in $h_k$, 
and also in the other $h_i$,  agree.
\end{proof}
  
A consequence of the previous lemma is that if $f\colon G\to G'$  is  a blowdown of partitioned
 graphs and $G$ is oriented, then the orientation which each edge of $G'$ inherits from the corresponding
 edge of $G$ together give an orientation of $G'$.  When $G'$ has this orientation, we say that $f$ is orientation
 preserving.

\subsection{Colored Graphs}
A pointed graph is a graph with one of its vertices designated as a basepoint.
  Basepoints are denoted by $*$.
A colored graph is a pointed, oriented, partitioned graph $G$ together with a bijection between the set ${\bf
n}=\{1,\ldots,n\}$ and the set $\{h_1,\ldots,h_n\}$ of edge sets of $G$.  Unless otherwise mentioned, we will always
assume that, for each $i\in {\bf n}$, this bijection takes $i$ to $h_i$.  We often think of the elements  of the set
${\bf n}$ as distinct colors and of drawing the colored graph $G$ so that all the edges of $G$ in $h_i$ have
color $i$.  Thus the set ${\bf n}$ is sometimes called the set of colors  of $G$, and the bijection from
$\bf n$ to the $h_i$, a coloring. 
  
 Colored graphs (really just partitioned graphs) are often drawn using the smooth edge convention.  Namely, two edges
  are in the same edgeset if and only if they  meet in a smooth curve at some vertex. 
  For example, partitions of the graphs  in      
Figure  \ref{reducedcores}, except for the last two, 
are indicated in this way.  Alternatively, edges in the same edgeset of a colored graph are drawn  with the same type of line, eg.\ solid, dotted, or finely dotted. Both conventions are often used on the same
graph.
  
  For convenience, we require that each vertex of a colored graph has valence at least three, except for
  the basepoint, which may be bivalent.

Two colored graphs $G$ and  $G'$ are said to be equivalent if there is a pointed cellular homeomorphism
$f\colon G\to G'$ which preserves colors and orientions.  
Such a homeomorphism  is called
an equivalence of colored graphs.  
\begin{prop}\label{equiv} 
Any two equivalences $G\to G'$ induce the same bijections $E(G)\to E(G')$ and  $V(G)\to V(G')$.
  \end{prop}
  
  The proof uses the next two lemmas.  These are preceded by a definition. 
  
  If $G$ is a graph, the number of half-edges incident with a vertex $v$ of $G$ is the integer
  $|\{e\in E(G)\colon d_0(e)=v   \}|+|\{e\in E(G)\colon d_1(e)=v   \}|$.
  
  \begin{lemma}\label{notthree}
  The number of half-edges of the same color incident with any vertex $v$ of a colored graph $G$ is
  at most two.
  \end{lemma}
  
  \begin{proof} If three half-edges of the same color $k$  were incident with $v$, then some subcycle
  of the basic cycle $C_k$ (with respect to any edge basis) would not contain any edge set of $G$.
    \end{proof}
  
  \begin{lemma}\label{v} If two equivalences $G\to G'$ of colored graphs agree on a vertex $v$ of $G$,
  they agree on all edges of $G$ incident with $v$.
  \end{lemma}
  
  This follows from Lemma \ref{notthree} since equivalences preserve colorings and orientations.
  
   \begin{proof}[Proof of  \ref{equiv}]
   Let $v\in V(G)$.  Choose an edge path $p$  in $G$ from $*$ to $v$.  Since equivalences preserve $*$,
   it follows from repeated use of Lemma \ref{v} that equivalences agree on the edges and, by continuity,
   on the (nonbasepoint) vertices of $p$.  Thus any two equivalences induce the same bijection
   $V(G)\to V(G')$, and so by Lemma \ref{v} the same bijection $E(G)\to E(G')$.
   
      \end{proof}

\subsection{The space $X_n$} \label{Xn}
Let $R_n$ be a fixed colored graph of rank $n$ having just one vertex. A marked, colored graph $(G,\mu)$, or just $G$, 
is, by definition, a colored graph $G$ together with a pointed homotopy equivalence $\mu\colon R_n\to G$.  
The map $\mu$ is called a marking (of $G$).  Two marked,
colored graphs $(G,\mu)$ and $(G',\mu')$ are equivalent if there is an equivalence $f\colon G\to G'$ of colored graphs
and a homotopy between $f\circ \mu$ and $\mu'$ which keeps the basepoint of $R_n$ fixed at all
times. 
We say that $(G',\mu')$ is a blowdown of $(G,\mu)$ if there is a blowdown $h\colon G\to G'$ of the underlying
partitioned graphs which preserves basepoints, colors and orientations, and if $(G',\mu')$ is equivalent to $(G',h\circ \mu)$.

Marked, colored graphs might be called decorated graphs, so we use 
  ${\rm b}\mathcal D_n$  to denote the poset of equivalence classes of marked, colored graphs of rank $n$ ordered by 
blowdowns.  Specifically, two elements $G_1$ and $G_2$ of ${\rm b}\mathcal D_n$ satisfy $G_1\ge G_2$ if $G_2$ is a blowdown 
of $G_1$.  Here and in the sequel we use the same symbol to  denote both a marked, colored graph and its equivalence class in ${\rm b}\mathcal D_n$.
A poset map ${\rm b}\mathcal D_n \to {\rm b}\mathcal D_{n+1}$ is defined by 
   sending $\mu\colon R_n\to G$ to $\mu\vee id \colon R_n\vee R_1\to G\vee R_1$ where the color of $R_1$ is
   $n+1$.
   The limit of the posets ${\rm b}\mathcal D_n$ using these maps is itself
a poset which we denote by ${\rm b}\mathcal D$.

\begin{definition}\label{def}\rm
The space $X_n$ is the geometric realization of the poset ${\rm b}\mathcal D_n$, and $X$ is the direct limit of the 
$X_n$ as $n \to \infty$.
\end{definition}

We now work towards another definition of $X_n$ which is ultimately given in Proposition \ref{altdef} below.

\begin{lemma}
The minimal elements of the poset ${\rm b}\mathcal D_n$ are in bijective correspondence with the elements of
the group $Aut(\pi_1(R_n))$.
\end{lemma}

\begin{proof}
Let $(R, \mu)$ be any minimal element of ${\rm b}\mathcal D_n$ so that the graph $R$ has just one vertex. 
Using the unique (up to homotopy) equivalence $R_n\to R$ of colored graphs, 
the marking $\mu$ of $R$  may be viewed as a self homotopy equivalence of $R_n$.  These are
in bijective correspondence with the elements of $Aut(\pi_1(R_n))$ since $R_n$ is an
Eilenburg-MacClane space of type $K(F_n,1)$.
\end{proof}

We now fix an isomorphism from the free group $F_n$ with basis $\{x_1,\dots,x_n\}$ to 
$\pi_1(R_n)\equiv \pi_1(R_n,*)$.  It takes 
$x_i$ to the (class of the) loop which wraps once around the single edge in $h_i$ in an 
orientation preserving manner.  This isomorphism, along with the previous lemma, gives a bijection 
\begin{equation}\label{bij}
\{\text {minimal\ elements\ of\ }{\rm b} {\cd}_n\}\longleftrightarrow  
\text{Aut}(F_n).
\end{equation}

Each marked, colored graph $G$ determines a full poset ${\rm b}G$ of ${\rm b}\mathcal D$ containing $G$ and all of 
its blowdowns.  As an example, if $G$ has just two vertices, 
with, say, $h_1=\{a,b\}$ the only edgeset with $> 1$ element,
  then ${\rm b}G$ has three elements: $G$, $G_a$ and $G_b$.  The
realization of ${\rm b}G$ is the barycentric subdivision of a 1-simplex, with the 
barycenter corresponding to $G$. 
 The other two vertices $G_a$ and $G_b$, when identified with elements of $Aut(F_n)$,
differ by a Whitehead automorphism as we now explain. 

We first review the definition of Whitehead automorphisms of the  free group $F_n$ with basis $\{x_1,\dots,x_n\}$.  If $a\in L=\{x_1^{\pm 1},\ldots, x_n^{\pm 1}\}$
and $a\notin \{x_k,x_k^{-1}\}$, then the Whitehead automorphism $\langle a,x_k\rangle$ takes $a$ to 
$a x_k$ and fixes those $x_i$ different from $a$ and $a^{-1}$.  In general, if $A\subset L$ with 
$A\cap \{x_k,x_k^{-1}\} = \varnothing$, then the Whitehead automorphism
 $\langle A,x_k\rangle = \prod_{a\in A} \langle a,x_k\rangle$.
We often write $\langle A,x_k^{-1}\rangle$ for $\langle A,x_k\rangle^{-1}$.

 Now assume that the two vertices of $G$ are the basepoint $*$ and $v$, and
that $G$ has edge sets $h_1,\ldots, h_n$. Also assume that 
 the basic cycle $C_1$ is orientated so that $init(b)=*$. Define a subset $A$ of the
letters $L$ of $F_n$ by 
$$A=\{x_i:i\ne 1\ {\rm and\ } term(e_i)=v\ {\rm where}\ e_i\in h_i\}\quad$$ $$\quad\quad \quad\cup\, \{x_i^{-1}: i\ne 1\
{\rm and}\ init(e_i)=v\ {\rm where}\ e_i\in h_i\}.$$

Writing $fg$ for the composite first $f$, then $g$, we have the following.

\begin{lemma}\label{eew}    If $G_a$ and $G_b$ are identified with automorphisms of $F_n$ as in (\ref{bij}) above,   then $G_b\langle A,x_1^{-1}\rangle=G_a$.
\end{lemma}

The omitted proof is not difficult, cf.\ 3.1.1 of \cite{CV}. 
  
  We will denote the edgepath in $|{\rm b}G|\subset X$ from $G_b$ to $G_a$ corresponding to
  the automorphism $\langle A,x_1\rangle$ by $b\wedge a$.  The symbol $b\wedge a$ could be
 read as (blow)up $b$, (blow)down $a$. 

Given $G\in \mathcal {\rm b}{\mathcal D}_n$ with edgesets $h_1,\dots , h_n$, let $P_i$ denote the poset of all nonempty subsets
of $h_i$ ordered by inclusion and let $P_1 \times \ldots \times P_n$ be the product poset where ordered pairs
$(s_i)$ and $(s_i')$ satisfy $(s_i) \ge (s_i')$ if $s_i\supseteq s_i'$ for all $i$.  Let  $f\colon P_1 \times \ldots
\times P_n \to {\rm b}G$ be the poset map which takes
$s=(s_i)\in P_1 \times \ldots \times P_n$ to the marked colored graph $G^s$ obtained from $G$ by collapsing the 
edges $\cup_i(h_i-s_i)$ of $G$.  Thus the edges of $G^s$ correspond to the elements of the set $\cup_is_i$.

\begin{lemma} \label{A}
The map $f\colon P_1 \times \ldots \times P_n \to {\rm b}G$ is a poset isomorphism.  
\end{lemma}
\begin{proof}  Surjectivity is clear.  
To show that $f$
is injective, we first assume that $s$ and $t$ are distinct minimal elements of $ P_1 \times \ldots \times
P_n$ with  say $s=(\{e_i\})$ and $t=(\{e_i'\})$.  
Let $j_1,\ldots ,j_m$ be the elements of the set $\{i\in {\bf n}:e_i\ne e_i'\}$ listed in any order.
Then the path in ${\rm b}G$ from $f(s)$ to $f(t)$   formed by 
 the sequence of length two edgepaths $e_{j_1}' \wedge e_{j_1},  e_{j_2}' \wedge e_{j_2},\dots ,
 e_{j_m}' \wedge e_{j_m}$
gives, by repeated use of Lemma \ref{W},  that
$f(t)=f(s) \prod_{i=1}^m\langle A_{j_i},x_{j_i}^{\pm 1}\rangle$. Here the sets  $A_{j_i}$ 
 are as in the  Lemma, 
 and the choices of the powers of each $x_{j_i}$ are dictated by the orientation of $G$.

Now $\prod_{i=1}^{m-1}\langle A_{j_i},x_{j_i}^{\pm 1}\rangle \ne \langle A_{j_m},x_{j_m}^{\mp 1}\rangle$
since, if $x_k^{\pm 1}\in A_{j_m}$, the automorphism $\langle A_{j_m},x_{j_m}^{\pm 1}\rangle$ takes $x_k$ to a word
containing at least one of $x_{j_m}$ or $x_{j_m}^{-1}$ whereas the image of $x_k$ under 
$\prod_{i=1}^{m-1}\langle A_{j_i},x_{j_i}^{\pm 1}\rangle$  contains neither $x_{j_m}$ nor  $x_{j_m}^{-1}$.  Therefore 
$\prod_{i=1}^{m}\langle A_{j_i},x_{j_i}^{\pm 1}\rangle$ is not the identity and so $f(t)\ne f(s)$.

Now assume that $f(s)=f(t)$ for distinct, nonminimal $s,t\in P_1 \times \ldots \times P_n$.  Choose distinct minimal 
elements $s'=(\{e_i\})$ and $t'=(\{e_i'\})$ such that $s\ge s'$ and $t\ge t'$, and such that the unique
equivalence $f(s)\to f(t)$ of colored graphs 
takes each $e_i$ to $e_i'$. This last condition implies
that 
$f(s')=f(t')$ which contradicts the first part of the proof since $s'\ne t'$.
\end{proof}

	If $\mathcal P$ is any poset, we let $|\mathcal P|$ denote the geometric relization of
$\mathcal P$. Since $| P_1\times\cdots\times  P_n|$ is homeomorphic to
$| P_1| \times\cdots\times | P_n|$, and since $| P_i|$ is simplicially
isomorphic to the first barycentric subdivision of $\sigma^{\|h_i\|}$, the 
standard simplex of dimension $\| h_i\| =|h_i| -1$,
the previous lemma gives the following.

\begin{prop}
\label{prop:cellofG} Let $G\in {\rm b}{\mathcal D}_n$.  Then the space $| {\rm b}G |$ is homeomorphic to the  disk
$\sigma^{\|h_i\|} \times \cdots \times \sigma^{\|h_n\|}$ of dimension $n$.
\end{prop}

The following is a consequence of Lemma \ref{A} and Proposition \ref{prop:cellofG}.

\begin{prop} \label{altdef}
The space $X_n$ is homeomorphic to the $n$-dimensional $CW$-complex having one $k$-cell $| {\rm b}G |$
corresponding to each $G\in {\rm b}{\mathcal D}_n      $ with $k+1$ vertices.  The attaching maps of the cells are as 
indicated by \ref{A} and \ref{prop:cellofG}.
\end{prop}

\begin{proof}
For each $G\in {\rm b}\cd_n$, there is a natural inclusion $|{\rm b}G| \hookrightarrow X_n$.
These inclusions take the interiors of the disks $|{\rm b}G|$ to pairwise disjoint sets which 
cover $X_n$.
\end{proof}

We will make use of both this $CW$-structure and the simplicial structure of $X_n$.  Any mention of the cells 
of $X_n$ or its skelata indicates that we have the $CW$-structure in mind.  The simplicial structure
is inferred by reference to simplices of $X_n$.
   
   We conclude this section with some remarks concerning the space $X_n$. Most of these will not be used in the
sequel.

According to Proposition \ref{prop:cellofG}, the cells of $X_n$ of dimension $k$ have different ``shapes",
 one  for each partition of $k$ as a sum of positive integers.  It may  be helpful to think of $X_n$ 
as the realization of a ``generalized simplicial complex".  We say generalized here since the buildings blocks in degree $k$--one 
for each
marked colored graphs with $k+1$ vertices--are not just the usual $k$-simplices, but all the shapes corresponding to the
partitions  of $k$.  Of course,  taking the ``barycentric  subdivision" of these shapes gives  a simplicial  set. Its realization
is the description of $X_n$ given in definition  \ref{def}.
   
Since the shapes are products of simplices, perhaps $X_n$ is the realization of a multisimplicial set
   in some natural way. It seems that only 
a subcomplex is, as we now explain.
Let ${\rm b}\cd_n'$ be the full subposet of ${\rm b}\cd_n$ consisting of the marked  colored graphs of rank $n$  whose  edgesets can
be canonically linearly ordered. Note that every edgeset of a colored graph is cyclically ordered in a canonical way
by any edgepath containing it.  If, for instance, such an edgepath passes  through the basepoint, then the edges
of the edgeset are canonically linearly ordered as  well.  

Let $Y=Y_{\bullet\dots\bullet}$ be the $n$-simplicial set 
with $Y_{i_1\dots i_n}$  the marked  colored graphs in ${\rm b}\cd_n'$ with $1+i_j$
edges of color $j$, for $1\le j\le n$.  The
$m^{\rm th}$ face map $Y_{i_1\dots i_j\dots i_n}  \longrightarrow Y_{i_1\dots i_j-1\dots i_n}$ is defined by blowing down the $m^{\rm th}$ edge of the edgeset $h_j$ of each $G\in Y_{i_1\dots i_j\dots i_n} $.  
To define the $m^{\rm th}$ edge, we identify the linearly ordered edges in each $h_j$ with the set $\{0,1, \dots ,i_j\}$ in
the obvious way. Degeneracies can be  defined by introducing bivalent vertices.
The realization of $Y$ describes a large CW-subcomplex of $X_n$.  For instance, $|Y|$ contains the entire 2-skeleton of $X_n$ 
except for those 2-cells  corresponding  to theta graphs.
Theta  graphs are those colored graphs equivalent to the colored graph shown  on the left of Figure \ref{theta}
for some choice of colors $x$ and $y$.
The image of the composite $|Y|\hookrightarrow X_n\to V$ agrees with that of the map $ X_n\to V$ but representatives of certain elements of $\pi_2(V)$, namely exotic elements
of $K_3(\mathbb Z)$, lift to $X_n$ but not to $|Y|$.  Here $V=V({\rm GL}(\mathbb Z))$ is the space defined by 
Volodin in \cite{V} and discussed in \S\ref{maps} below.

Theta graphs are also noteworthy for the following reason.  Colored graphs and partition-preserving
maps between them form a Waldhausen category $\mathcal C$ 
(as discussed in \S1.3 of \cite{W})
with weak equivalences those morphims which 
induce bijections on the set of colors and with injections as cofibrations.
If every colored graph $G$ admitted a filtration $*\rightarrowtail G_1\rightarrowtail \dots\rightarrowtail G_n=G$ by cofibrations with all subquotients
$G_i/G_{i-1}$ roses, then essentially Waldhausen's additivity theorem would give a homotopy equivalence 
${\rm pS}_\bullet \mathcal C \simeq \Omega^{\infty} S^{\infty}$, where p denotes the weak equivalences just 
mentioned.
But theta graphs do not admit such 
filtrations.  In fact, $\Omega^{\infty} S^{\infty}$ splits off 
${\rm pS}_\bullet \mathcal C$ with theta graphs
contributing to nontrivial homology of the complementary factor.
 
We note that the realization of the simplicial set  $[n]\mapsto Y_{0\dots 0\,n}$ is the disjoint  union  of  
simplicial copies of $\mathbb{R}^{2n-2} $.  This is discussed in \S\ref{RninXk} below from a slightly 
 different viewpoint.
The free abelian group $\mathbb Z^{2n-2}\subset \autfn$ generated by $\langle x_1,x_n\rangle, \dots ,
\langle x_{n-1},x_n\rangle$ and $\langle x_1^{-1},x_n\rangle, \dots ,
\langle x_{n-1}^{-1},x_n\rangle$ acts freely on each copy.
  Also, the simplicial map $|Y_{0\dots 0 \bullet}| \hookrightarrow X_n \to  S\mathbb A_n$ is injective  on  each copy.  In fact, it is injective  on considerably larger subcomplexes of $X_n$.  Here $S\mathbb A_n$ is a version of the 
 sphere complex defined, for instance, in \cite{HV}.

Lemma \ref{eew} implies that each (oriented) 1-cell of $X_n$ corresponds to a Whitehead automorphism.
According to Nielson's theorem, cf.\ Cor.\ 3.6 of \cite{KN},
these generate the index two subgroup $Aut_+(F_n)$ of $Aut(F_n)$ consisting of all automorphisms which abelianize to 
elements of $GL_n({\mathbb Z})$ with determinent 1. Thus the space $X_n$ has two connected components.  The 1-skeleton
of the component containing the identity is the Caley graph of $Aut_+(F_n)$ with the Whitehead automorphisms as 
generators. 

Let $G\in {\rm b}\cd$.  The induced $CW$-structure on $|{\rm b}G|$ is denoted by $|G|$, and is just the product $CW$-structure
given by Proposition \ref{prop:cellofG} (assuming the $CW$-structure of the standard simplex $\sigma^{||h_i||}$ is just its 
simplicial structure).  Elaborating on the middle paragraph of the proof of Lemma \ref{A} shows that
the 1-skeleton of $|G|$ is a geodesically convex subset of the Caley graph of $Aut_+(F_n)$ mentioned above.
This, along with the next two sections, suggests an algebraic definition  of the space $X_n$ which does not
use colored graphs.

\subsection{The Partial Ordering determined by a Colored Graph} \label{cgordering}
Each colored graph $(G,\{h_i\})$ of rank $n$ defines an ordering of the set ${\bf n}=\{1,\dots ,n\}$  of colors
as follows.  If $b,r\in {\bf n}$, then $b>_G r$ or just $b> r$  if $b\ne r$ and if there is some
edge basis $E$ of $G$ such that $C_b^E\cap h_r\ne \varnothing$.

Note that $b>r$ if and only if  there
exists an edge basis
$E=\{e_1,\dots ,e_n\}$ and one additional edge $e_r'\in h_r$ such that the edge $e_b\in E$ joins the two
vertices of the graph obtained by blowing down all the edges of $G$ except for those in $E\cup \{e_r'\}$. 

 The
following two lemmas will show that the ordering $>$ is antisymmetric and transitive.
   Thus $\ge$ is a partial ordering of the set ${\bf n}$.


\begin{lemma} \label{antisymmetric}
If $b>r$, then $C_r^E\cap h_b=\varnothing$ for every edge basis $E$.
\end{lemma}
It follows from this that both $b>r$ and $r>b$ cannot occur.
\begin{proof}[Proof of  \ref{antisymmetric}]
Let $\mathcal E$ be the set of edgebases defined by $E\in\mathcal E$ if $C_b^E\cap h_r\ne \varnothing$.
Note that if $E\in \mathcal E$, then $C_r^E\cap h_b= \varnothing$ by Lemma \ref{lemma1}.

Let $F=\{f_1,\ldots, f_n\}$ be any edge basis with $f_i\in h_i$. We will show that there is an 
$E\in \mathcal E$ such that $C_r^E\cap F=\{f_r\}$.  Then $C_r^E=C_r^F$ and so $C_r^F\cap h_b =
C_r^E\cap h_b=\varnothing$, which will complete the proof.

So assume that $f_i\in C_r^E\cap F$ where $i\ne r$ and $E\in \mathcal E$.  If $f_i\in C_b^E$, then
$C_i^E\cap h_b=\varnothing$ by Lemma \ref{lemma1}.  So letting $E_1$ be
the edge basis
$E\cup \{f_i\}-\{e_i\}$, we have  $C_b^{E_1}=(C_b^E\cup C_i^E) - (C_b^E\cap C_i^E)$. Therefore 
$E_1\in \mathcal E$ since $C_i^E\cap h_r=\varnothing$ by Lemma \ref{lemma1} once more.  And if $f_i\notin C_b^E$,
then $C_b^{E_1}=C_b^E$ again giving $E_1\in \mathcal E$.  Note that this shows that if $E\in \mathcal E$
and $|C_r^E\cap F|>1$, then there is an $E_1\in \mathcal E$ such that $|E_1\cap F|>|E\cap F|$.
So now either $|C_r^{E_1}\cap F|= 1$ which would complete the proof or, if not, there is, as just noted, an 
$E_2\in \mathcal E$ with $|E_2\cap F|>|E_1\cap F|$.  Repeating this as needed eventually gives an $E\in 
\mathcal E$ with $|C_r^E\cap F|=1$ since $C_r^F\cap F=\{f_r\}$.
\end{proof}

\begin{lemma}\label{trans}
If $b>r$ and $r>g$, then $C_b^E\cap h_r\ne \varnothing$ and $C_b^E\cap h_g\ne \varnothing$ for some edge basis 
$E$.
\end{lemma}
So the relation $>_G$ is transitive.
\begin{proof}[Proof of \ref{trans}]
Assume $r>g$ and let $F=\{f_i\}$ be an edge basis with $f_i\in h_i$ such that $C_r^F\cap h_g\ne \varnothing$.
Also assume $b>r$. Then as in the proof of Lemma \ref{antisymmetric}, there is an edge basis $E$ with 
$C_b^E\cap h_r\ne \varnothing$  and $C_r^E\cap F=\{f_r\}$.  If $C_b^E\cap h_g\ne \varnothing$ the lemma is 
proved, so assume $C_b^E\cap h_g= \varnothing$.  Pick an edge $e_r'\in C_b^E\cap h_r$, and let
$E'=E\cup\{e_r'\}-\{e_r\}$.  Then 
$C_b^{E'}=(C_b^E\cup C_r^E)-(C_b^E\cap C_r^E)$, which contains an edge of $h_g$ since
$C_r^E=C_r^F$ does, and $C_b^E$ does not.  Also $e_r\in C_b^{E'}$ since $e_r\notin
C_b^E$.
\end{proof}

\subsection{Triangular subgroups}\label{tri}
Let $G$ be a colored graph.  The previous section shows that the Whitehead automorphisms 
corresponding to the 1-cells of $|G|$ are 
contained in a triangular subgroup
$T_n^{\sigma}$ of $Aut(F_n)$ which we are about to define.

	Let $\sigma$ be any partial ordering of the set {\bf n}.  In particular, $\sigma=\sigma (G)$ could be the 
order determined by $G$.
Let  $S_n^{\sigma}$
be the subgroup of $\autfn$ generated by the set $\{ \langle
x_i,x_j\rangle\colon i \stackrel{\sigma}{>}j\}  \cup \{ \langle
x_i^{-1},x_j\rangle \colon i \stackrel{\sigma}{>}j\}$, and let $C_n$ be the subgroup
generated by the set
$\{\langle x_i,x_i^{-1};x_j\rangle \colon i\ne j\}$. 
Elements of $C_n$ include automorphisms which conjugate any one basis element
by any other and leave all other basis elements fixed.  The triangular
subgroup
$T_n^{\sigma}$ mentioned above is, by definition, the one generated by
$S_n^{\sigma} \cup C_n$.  Elements of
$T_n^{\sigma}$ are referred to as ($\sigma$-)triangular automorphisms, while
those of
$S_n^{\sigma}$ are called strictly triangular.

\begin{lemma}  The group $T_n^{\sigma}$ of triangular automorphisms is torsion
free.
\end{lemma}
\begin{proof} Let $a\colon \autfn \to \text {GL}_n (\mathbb Z)$ be the
natural homomorphism.  Let $K$ be the kernel of $a$ restricted to
$T_n^{\sigma}$ so that the sequence $$1\to K \to T_n^{\sigma}\to
a(T_n^{\sigma}) \to 1 $$ is exact.  $K$ is contained in the torsion free
subgroup $\text {IA}(F_n)$ of $\autfn$. Also, $a(T_n^{\sigma})$ is a torsion
free subgroup of $\text {GL}_n (\mathbb Z)$. If $\sigma$ is the standard
ordering, this subgroup consists of all upper triangular matrices with
diagonal entries equal to 1, and, for other $\sigma$, is isomorphic to it. Thus
both $K$ and $a(T_n^{\sigma})$ are torsion free and so $T_n^{\sigma}$ is also.
\end{proof}

\subsection{The action of $Aut(F_n)$ on $X_n$}
Let $f\in \autfn$ and $G=(G,\mu) \in {\rm b}\mathcal{D}_n$.  Assume $f$ corresponds to the
minimal element $f'\colon R_n\to R_n$ of ${\rm b}\mathcal{D}_n$.  Denote by
$f\cdot G$ or      
$f G$ the element $(G,f'\mu)$ of
${\rm b}\mathcal{D}_n$. This defines a left action of
$\autfn$ on
${\rm b}\mathcal{D}_n$ by poset isomorphisms.  Thus there is a simplicial left action
of $\autfn$ on
$X_n$.  This action is also cellular. Moreover, we have the following. 

\begin{prop} \label{free}
The action of $\autfn$ on $X_n$ is free.
\end{prop}
\begin{proof} Identify $X_n ^{(0)}$, the 0-skeleton of $X_n$, with \autfn.  Then the action of any $f\in
\autfn$ on $X_n^{(0)}$ amounts to multiplication by $f$ on the left.  So
$\autfn$ clearly acts freely  on $X_n^{(0)}$.  Therefore  suppose $f\in
\autfn$ fixes an interior point of some cell $\sigma$ of $X_n$.  Then $f$
permutes the 0-cells of $\sigma$ and so $f^n= id $ for some $n$.  Let $v\in
\autfn$ be any 0-cell of $\sigma$.  Since $fv$ is also a 0-cell of $\sigma$,
there is a triangular automorphism $t$ such that $fv=vt$.  Therefore
$t^n=v^{-1} f^n v=id$ which, by the preceding lemma, implies that $t=id$. So
$f=id$ also.
\end{proof}

\subsection{The space $\xt_n$} \label{remarks}
Let $X_{*,n}$ be the component of $X_n$ containing the vertex corresponding via the bijection (\ref{bij})
to the identity in $\text{Aut}(F_n)$.
We often work with the universal cover $\xt_n$ of the space $X_{*,n}$. 
The covering $\xt_n$ can be defined using marked colored graphs just as $X_n$ is, except that the markings
involve elements of the group $StN_n$, rather than those of $\autfn$. 
 Details follow, beginning with basic facts about the  group $StN_n$.

Let $B^{\pm}=\{x_1,\ldots ,x_n\}\cup \{x_1^{-1},\ldots ,x_n^{-1}\}$ where $\{x_1,\ldots ,x_n\}$
is a basis for $F_n$.  The group $StN_n$ is
defined by the presentation having one generator $\lda x,y\rda$ for all $x,y\in B^{\pm}$ with $x\ne y^{\pm 1}$, and the
relations 
\begin{description}
\item[{(R0)}] $\lda x,y^{-1}\rda = \lda x,y\rda^{-1}$
\item[(R1)] $[\lda x,y\rda , \lda x',y'\rda ]=1$ if $y^{\pm 1}\ne x'\ne x\ne (y')^{\pm 1}$
\item[{(R2)}] $[\lda x,y\rda,\lda y,z\rda ]=\lda x,z\rda$, if $x\ne z^{\pm 1}$
\item[{(R3)}] $\lda y,x\rda \lda x^{-1},y\rda \lda  y^{-1},x^{-1}\rda  =\lda y^{-1},x^{-1}\rda \lda x,y^{-1}\rda \lda
y,x\rda $.
\end{description}
   
   We sometimes will use the notation $L_{ij}$ for $\lda x_i^{-1}, x_j^{-1} \rda$ and 
$R_{ij}$ for $\lda x_i, x_j^{-1} \rda$.

Two elements of $StN_n$ which we will often use are 
$w_{ij}^L=  \lda x_i,x_j\rda \lda x_j ^{-1},x_i\rda \lda  x_i^{-1}, x_j ^{-1}\rda = R_{ij}^{-1} L_{ji}^{-1} L_{ij}$
and $w_{ij}^R= \lda x_i^{-1}, x_j ^{-1}\rda \lda x_j,x_i^{-1}\rda  \lda x_i, x_j \rda  = 
 L_{ij} R_{ji} R_{ij}^{-1}$
  
   There is a homomorphism $\phi\colon StN_n\to \autfn$ defined by $\phi(\lda x,y\rda)=\langle x,y\rangle$.  Gersten showed
in \cite{G} that the kernel of $\phi$, which we denote by $KN_{n,2}$, has order two.

With the following lemma, many of the facts about marked, colored graphs discussed above carry over to 
$StN$-marked, colored graphs.  Such graphs are defined after its proof.

\begin{lemma} \label{triinj}
Let $T^{\sigma}(F_n)$ be any triangular subgroup of $\autfn$. Then
there is an injection $T^{\sigma}(F_n)\to StN_n$ induced by $\langle x_i,x_j\rangle \mapsto 
\langle\!\langle x_i,x_j\rangle\!\rangle$, 
$\langle  x_i^{-1},x_j \rangle \mapsto \langle\!\langle x_i^{-1},x_j  \rangle\!\rangle$, and     
$\langle  \{x_i,x_i^{- 1}\},x_j \rangle \mapsto \langle\!\langle x_i^{-1},x_j  \rangle\!\rangle 
\langle\!\langle x_i,x_j  \rangle\!\rangle$.      
\end{lemma}

\begin{proof}
   Let $\langle G_{TF}\rangle$ be the free group on the set $S_n^{\sigma} \cup C_n$  of defining generators of $T_n^{\sigma}$ as in \S\ref{tri}, and let $\langle G_{T{\mathbb Z}}\rangle$ be the 
   free group on the image of this set in ${\rm GL}_n({\mathbb Z})$.
   We take the horizontal map $\langle G_{TF}\rangle \to  StN_n    $  in the diagram
  
  $$\xymatrix@R=14pt@C=-6pt{
  {R_F} \ar[dr] \ar[dd] \ar[rr]  & &{\mathbb Z}_2 \ar[dr]^{\approx} \ar'[d] [dd]
\\
&R_{\mathbb Z}\ar[dd]\ar[rr] & & {\mathbb Z}_2\ar[dd]
\\
{\langle G_{TF}\rangle}\ar[dd] \ar[dr] \ar'[r][rr] &  &  StN_n \ar'[d]_{\phi}[dd] \ar[dr]
\\
&{\langle G_{T{\mathbb Z}}\rangle}\ar[dd] \ar[rr] & &St_n({\mathbb Z})\ar[dd]
\\
T^{\sigma}_n\ar[dr] \ar'[r][rr] & & {\rm Aut}(F_n)\ar[dr]
\\
&T^{\sigma}_n({\mathbb Z})\ar[rr] & &  {\rm GL}_n({\mathbb Z})
}$$ 
  to be induced as in the statement of the lemma.  The horizontal map 
  $\langle G_{T{\mathbb Z}}\rangle\to St_n({\mathbb Z})$ is induced similarly.  Each of the 
  four columns of the diagram is part of a short exact sequence. The groups $R_F$ and 
  $R_{\mathbb Z}$ are defined as the kernels of the relevant maps, while one of the $\mathbb Z_2$'s
  follows from the computation of $K_2(\mathbb Z)$, the other from the nonabelian analog.
   It follows from a Lemma in \cite{V} that the homomorphism $R_{\mathbb Z}\to {\mathbb Z}_2$ in the top
horizontal square of the diagram
  is trivial.  Since this square commutes, the homomorphism $R_F\to {\mathbb Z}_2$ in it is also trivial.
Thus the injective map $T^{\sigma}(F_n)\approx \langle G_{TF}\rangle /R_F\to \autfn$ factors through the group $StN_n$.
\end{proof}
A StN-marked, colored graph (of rank $n$) is, by definition, a marked, colored graph $G$ (of rank $n$) together with a 
compatible choice of lifts of the 0-cells of $|G|$ from $\autfn$ to $StN_n$.
Recall that $|G|$ denotes the $CW$-structure of $|{\rm b}G|$ discussed just after the bijection (\ref{bij})
  near the beginning of section \ref{Xn}.  
As usual, we are identifying the elements of $|G|^{(0)}$ with elements of $\autfn$ using this bijection. 
A lift of $v\in |G|^{(0)}\subset \autfn$ then means one of the two preimages of $v$ under the map
$\phi\colon StN_n\to \autfn$.
Lifts are compatible if whenever
$v$ and $w$, respectively, are in $|G|^{(0)}$, their lifts $\tilde v$ and $\tilde w$, respectively,  
satisfy $\tilde v \prod_i
\langle\! \langle A_i,a_i\rangle \!\rangle^{\pm 1}= \tilde w$ in the group $StN_n$, where 
$\prod_i \langle A_i,a_i\rangle^{\pm 1}$
corresponds a the shortest path in $|G|^{(1)}$ from $v$ to $w$.  Here we are identifying (oriented) 1-cells  
of $|G|$ with Whitehead automorphisms in $T_n^{\sigma (G)}$ and then, using Lemma \ref{triinj}, with elements of $StN_n$.
This lemma and section \ref{cgordering} 
 shows that a compatible choice of lifts of the elements of  $|G|^{(0)}$ always exists and is determined by a choice of 
lift of any one of the elements.

Two $StN$-marked, colored graphs are equivalent if, when their StN-markings are converted to just markings 
using the homomorphism $\phi\colon StN_n\to \autfn$,
 they are equal as elements of ${\rm b}\cd$, and if their StN-markings 
agree.  Blowdowns of StN-marked, colored graphs inherit $StN$-markings in a natural way. Let
$\widetilde{{\rm b} {\cd}}_n$ be the poset of equivalence classes of $StN$-marked, colored graphs ordered by blowdowns.
Note that since the kernel of the map $\phi\colon StN_{n} \to Aut(F_{n})$ has 
order two, there is a two-to-one poset map $\phi_{\cp}\colon 
\widetilde{{\rm b} {\cd}}_n\to{\rm b} {\cd}_n$.

Let $\xt_n$ be the geometric realization of the poset $\widetilde{{\rm b} {\cd}}_n$.  The space $\xt_n$ has a $CW$-structure analogous 
to that of $X_n$.  Specifically,  
 each $\tilde{G} \in\widetilde{{\rm b} {\cd}}_n$ determines a full subposet ${\rm b}\tilde G$ of 
$\widetilde{{\rm b} {\cd}}_n$ consisting of $\tilde{G}$ and all of its blowdowns.  The proof of Lemma \ref{A},
 along with Lemma \ref{triinj}, implies that 
  $\phi _{\cp}|\colon {\rm b}\tilde{G}\to {\rm b}G$ is a
poset isomorphism.  Thus, the proof of Proposition \ref{altdef} essentially shows that the indicated
 $CW$-structure of $\xt_n$ consists of one $k$-cell $|{\rm b}\tilde{G}|$ for each $\tilde{G}
\in \widetilde{{\rm b} {\cd}}_n$ with $k+1$ vertices.  
Just as with the space $X_n$, mention of the cells of $\xt_n$ or any of its skelta indicates that
we are working with the $CW$-structure of $\xt_n$, while mention of its simplices or vertices indicates   
the simplicial structure.

The poset map $\phi_{\cp}\colon 
\widetilde{{\rm b} {\cd}}_n\to{\rm b} {\cd}_n$ induces a simplicial and cellular map 
$\phi_X\colon \xt_n \to X_n$ of spaces with image $X_{*,n}$.  These maps are both denoted 
by $\phi$, with subscripts that are usually omitted, since, when suitably restricted, 
they both agree with the homomorphism $\phi\colon StN_n\to \autfn$ 
via the bijections
$$\quad \{\text {minimal\ elements\ of\ }\widetilde{{\rm b} {\cd}}_n\}\leftrightarrow\xt_n^{(0)} \leftrightarrow StN_n\quad \text{and}$$ $$
\quad \{\text {minimal\ elements\ of\ } {\rm b} {\cd}_n\} \leftrightarrow X_n^{(0)} \leftrightarrow \autfn .\quad$$

If $\tilde v_1$ and $\tilde v_2$ are vertices of $\xt_n$ which satisfy $\phi(\tilde v_1)=v=
\phi(\tilde v_2)$ then the closed stars of $\tilde v_1$ and $\tilde v_2$ in $\xt_n$ are 
disjoint since their vertices are.  Also, each closed star is mapped homeomorphically
by $\phi\colon\xt_n\to X_n$ to the closed star of $v$ in $X_n$.  Thus, since $X_n$ is covered by the open stars of its
vertices, the map $\phi\colon\xt_n\to X_n$ is a covering of its image $X_n^{id}$.  
The next lemma shows that it is the universal covering.

\begin{lemma}\label{sc}
The space $\xt_n$ is simply-connected for $n\ge 3$.
\end{lemma}

\begin{proof}
For each (oriented) 1-cell $e$ of $X_n$ which corresponds to a Whitehead automorphism of the form $\langle A,x_i\rangle$ 
  where $|A|=1$, identify both (oriented) lifts of $e$ in $\xt_n$ with the generator 
  $\langle\!\langle A,x_i\rangle\!\rangle$ 
  of the group $StN_n$.  
  Using these identifications, the defining 
  relations of $StN_n$ can be realized by 2-cells of $\xt_n$.  
  Figure \ref{wLeqwR} shows how to do this for the relations (R3).
  Each of the three 2-cells in the figure
  correspond  to the colored graph drawn within it.  Edges of color $j$ are dotted, and those of color
  $i$ solid.  Arrows near the boundary of each 2-cell indicate the boundaries of these cells.
  For instance, the left-hand 1-cell of the boundary of the square cell in the middle corresponds to
  the $\theta$-shaped graph with its edge $c$ blown down.
  
  The various relations (R1) can be realized similarly using the colored graphs numbered 3, 4, 5, and 6
  in  Figure \ref{redsq} below. 
  In case $y=y'$ in relation
  (R1), then a pair of 2-cells, both corresponding to either the colored graph numbered 5 or the one 
  numbered 6, are needed. Otherwise, a single 2-cell corresponding to one of the colored graphs 
  3 or 4 will do.  Finally, the relations (R2) are realized by pairs of 2-cells, with one cell of each pair 
  corresponding to the colored graph numbered 2 in Figure \ref{redsq}, and the other to the one 
  numbered 5.
  
  \begin{figure}[tb]
\centerline{\mbox{\includegraphics*{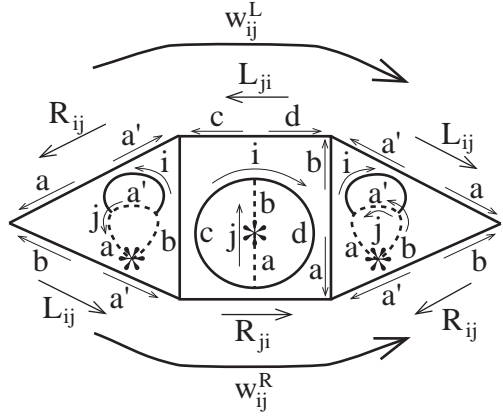}}}
\caption{Realizing the relation $w_{ij}^R=w_{ij}^L$ in $\xt_n$.}
\label{wLeqwR}
\end{figure}
	
To complete the proof, we show that any edgeloop of $\xt_n$ can be deformed to one whose edges were identified in the first paragraph with the defining generators of $StN_n$. 
For this let $\tilde e$ be any 1-cell of $\xt_n$.  
Assume that
$\tilde e$ maps to the edge of $X_n$ which corresponds to the Whitehead automorphism $\langle A,x_i\rangle$.
If $|A|=m>1$, then the edgepath $\tilde e$ can be deformed, rel endpoints, to an edgepath whose edges correspond to the generators $\langle\!\langle a,x_i\rangle\!\rangle$, $a\in A$ of $StN_n$.  
In case   $|A|=m=2$ a cell corresponding to one of the colored graphs numbered 5 or 6 in Figure 
\ref{redsq}  facilitate such a deformation.
  In general, $m-1$ such cells are needed.
\end{proof}

We mention that, since the kernel of the homomorphism $\phi\colon StN_n\to \autfn$ is $KN_{2,n}$, the component
$X_n^{id}$ is homeomorphic to the quotient of $\xt_n$ by the action of the $KN_{2,n}$. 
So it follows from \ref{sc} that $\pi_1(X_n^{id})\approx KN_{2,n}$.

\subsection{The space $X^{\text {alg}}$}\label{maps}
In this section we define a space $X_n^{\text {alg}}$ which is, to some (undetermined) extent, an
algebraic version of the space $X_n$.  We will also define maps $X_n\to X_n^{\text {alg}}$ and 
$X_n^{\text {alg}}\to V(GL_n(\mathbb Z))$. Here $V(GL_n(\mathbb Z))$ is a space 
(denoted by $V(GL_n({\mathbb Z}),\{T_n^{\sigma}({\mathbb Z})\})$ in \cite{S}) whose limit $V(GL({\mathbb Z}))
=\lim_{n\to \infty} V(GL_n({\mathbb Z}))$ has homotopy groups the $K$-groups of $\mathbb Z$. The space
$V(GL({\mathbb Z}))$ was first defined by Volodin in \cite{V} and Hatcher, but we follow Suslin's
description in \cite{S} of Volodin's work.  Indeed much of this section is taken from or is
similar to $\S$1 of \cite{S}.

The definition of $X_n^{\text {alg}}$ is similar to that of $V(GL_n({\mathbb Z}))$, and both are essentially
special instances of the following definitions.  Let $G$ be a group and $\{G_i\}$ a collection
of subgroups of $G$.  
The simplices of the simplicial complex $V(G,\{G_i\})$ are, by definition, certain subsets of $G$.
Specifically, $\sigma=\{g_0,\ldots, g_p\}$
is a $p$-simplex of $V(G,\{G_i\})$ if the corresponding set $\{g_ig_j^{-1} : g_i, g_j\in \sigma\}$ 
of products is contained in one of the subgoups $G_i$.
We also use $V(G,\{G_i\})$ to denote the  geometric realization of the simplicial complex $V(G,\{G_i\})$.
  A space and simplicial complex $V_{\text {op}}(G,\{G_i\})$ is defined just as
$V(G,\{G_i\})$ is, but with $g_i^{-1}g_j$ in place of $g_ig_j^{-1}$.

Now let $T_n^{\sigma}$ be a triangular subgroup of $\autfn$. Such subgroups were defined in section \ref{tri}.
As in \cite{S}, let $T_n^{\sigma}({\mathbb Z})$ denote 
the image of $T_n^{\sigma}$ under the abelianizing homomorphism $S\colon\autfn\to \text{GL}_n({\mathbb Z})$.  Thus, if 
$\sigma$ is the standard ordering of the set ${\bf n}$, the subgroup $T_n^{\sigma}({\mathbb Z})$
of  $\text {GL}_n (\mathbb Z)$ consists of all upper triangular matrices with ones on the diagonal.
  Volodin defined $V(GL_n({\mathbb Z}))$
as $V(GL_n({\mathbb Z}),\{T_n^{\sigma}({\mathbb Z})\})$ where $\sigma$ ranges over all 
partial orders of the set ${\bf
n}$. Similarly, we define the space $X_n^{\text {alg}}$  by
 $$X_n^{\text {alg}} = V_{op}(\autfn,\{T_n^{\sigma}\}),$$ again with $\sigma$
ranging over all partials orders of ${\bf n}$.  The reason for the 
subscript $op$ is that we are using the convention that the composite $fg$ means first $f$, then $g$, if $f,g\in \autfn$
and, the opposite order, first $g$, then $f$, if $f,g\in GL_n({\mathbb Z})$.

The abelianizing map $S$ induces a simplicial map $X_n^{\text {alg}}\to V(GL_n({\mathbb Z}))$.
Also, letting 
$$V(St_n({\mathbb Z})) = V(St_n({\mathbb Z}),\{T_n^{\sigma}({\mathbb Z})\})\quad \text {and} \quad
\xt_n^{\text {alg}} = V_{op}(StN_n,\{T_n^{\sigma}\}),$$  a simplicial map
$\xt_n^{\text {alg}}\to V(St_n({\mathbb Z}))$ is induced by the homomorphism $ StN_n\to St_n({\mathbb Z})$.
Here we are using Lemma \ref{triinj} to identify $T_n^{\sigma}$ with its image in $StN_n$.

The group $StN_n$ acts freely on $\xt_n^{\text {alg}}$ by left multiplication, cf.\ Proposition \ref{free}. The quotient 
$\xt_n^{\text {alg}}/KN_{2,n}$ is homeomorphic to $X_n^{\text {alg}}$.  Since the defining relations for the group $StN$
are triangular, the space $\xt_n^{\text {alg}}$ is
simply connected.  So it is the universal cover of $X_{*,n}^{\text {alg}}$. 
Here $X_{*,n}^{\text {alg}}= V_{op}(\text{Aut}_+(F_n),\{T_n^{\sigma}\})$ is the component of $X_n^{\text {alg}}$ containing 
the 0-simplex $\{id\}$.

Letting $n\to \infty$ we obtain stable versions of the above spaces, complexes, groups, actions and maps.
Stable objects are typically denoted by deleting the subscript $n$ from the notation for the corresponding unstable
object.

The following Proposition is the main reason for introducing the complex $X_n^{\text {alg}}$.
  We conjecture that it's true with $\xt$ in place of $\xt^{\text {alg}}$.

\begin{prop}\label{trivact}
The group $StN$ acts trivially on both the homology and the homotopy groups of $\xt^{\text {alg}}$.
\end{prop}

The proof will make use of cartesian products so, following Suslin \cite{S}, we introduce a simplicial set $\wt_n$, with
limit $\wt$, whose geometric realization is homotopy equivalent to $\xt_n^{\text {alg}}$.  The simplicial set
$\wt_n$ has a $p$-simplex $(g_0,\ldots ,g_p)$
for each {\it sequence} $g_0,\ldots ,g_p$ of elements of $StN_n$ such that all $g_i^{-1}g_j$ are in one of 
the triangular subgroups $T_n^{\sigma}$.  The $i$th face, respectively degeneracy, of the simplex 
$(g_0,\ldots ,g_p)$ is defined as usual by deleting, respectively repeating, $g_i$.

\begin{lemma}\label{sus}
If $x\in StN_{n+1}$, then the maps $u_n\colon \wt_n\to \wt_{n+1}$  and $x\cdot u_n$ are homotopic.
\end{lemma}

This shows that Proposition  \ref{trivact} is true if, in its statement, $X^{\text {alg}}$ is replaced with $\wt$.
Since the obvious $StN$-equivarient map $\wt\to \xt^{\text {alg}}$ is a homotopy equivalence, the Proposition is
true as stated.

\begin{proof}[Proof of Lemma \ref{sus}]
Because of defining relations R2, the group $StN_{n+1}$ is generated by elements of the form 
$\langle\!\langle x_i,x_{n+1}\rangle\!\rangle$, $\langle\!\langle x_i^{-1},x_{n+1}\rangle\!\rangle$,
$\langle\!\langle x_{n+1},x_i\rangle\!\rangle$, and $\langle\!\langle x_{n+1}^{-1},x_i\rangle\!\rangle$.
We first assume that $x\in StN_{n+1}$ has one of these forms and define  a homotopy between 
$u_n$ and $x\cdot u_n$ by taking the simplex 
$$(\underbrace{0,\ldots ,0}_s,\underbrace{1,\ldots ,1}_t)\times (\alpha_1,\ldots , \alpha_{s+t})$$
of $I\times \wt_n$ to the simplex 
$$\tau=(\alpha_1,\ldots,\alpha_s,x\alpha_{s+1},\ldots, x\alpha_{s+t})$$
of $\wt_{n+1}$.  To see that $\tau$ is a simplex of 
$\wt_{n+1}$, we first mention that one can check, by using the defining relations for $StN_n$ and considering a
number of cases (or by just using marked, colored graphs) that if 
$x=\langle\!\langle x_i^{\pm 1},x_{n+1}\rangle\!\rangle^{\pm 1}$ and $\alpha\in StN_n$, then $\alpha^{-1}x\alpha$
is in the subgroup of $StN_{n+1}$ generated by 
$$\{\langle\!\langle x_i,x_{n+1}\rangle\!\rangle : i\in {\bf n}\} \cup 
\{\langle\!\langle x_i^{-1},x_{n+1}\rangle\!\rangle : i\in {\bf n}\}\cup
\{\langle\!\langle x_{n+1},x_i\rangle\!\rangle
\langle\!\langle x^{-1}_{n+1},x_i\rangle\!\rangle : i\in {\bf n}\}.$$
And if $x=\langle\!\langle x_{n+1}^{\pm 1},x_i\rangle\!\rangle^{\pm 1}$, then 
$\alpha^{-1}x\alpha$
is in the subgroup generated by 
$$\{\langle\!\langle x_{n+1},x_i\rangle\!\rangle : i\in {\bf n}\} \cup 
\{\langle\!\langle x_{n+1}^{-1},x_i\rangle\!\rangle : i\in {\bf n}\}.$$

Now let $\sigma$ be the ordering of ${\bf n}$ such that all the $\alpha_i^{-1}\alpha_j$ are in $T_n^{\sigma}$.
Writing $\alpha_i^{-1}x\alpha_j$ as $(\alpha_i^{-1}x\alpha_i)(\alpha_i^{-1}\alpha_j)$, we see that $\tau$ 
is a simplex of $\wt_{n+1}$ where the relevant ordering of the set ${\bf n+1}$ is the ordering $\sigma$ with $n+1$ added
as either a maximal or minimal element.

	Now assume $y\in StN_{n+1}$ is of one of the forms mentioned above.  Let $h$ denote the image of the homotopy
between $u_n$ and $x\cdot u_n$.  Then  the subcomplex $y\cdot h$ gives a homotopy between $y\cdot u_n$ and
$y\cdot (x\cdot u_n)=yx\cdot u_n$.  Since $y\cdot u_n$ is homotopic to $u_n$ by the first part of the 
proof, the lemma is true for $yx\in StN_{n+1}$. Continuing in this way proves the lemma for arbitrary elements of
$StN_{n+1}$.
\end{proof}

\subsection{The map $\mu\colon X\to X^{\text {alg}}$} \label{xtoxalg}  By the (first) barycentric subdivision poset ${\rm Simp}(C)$ of a
simplicial complex $C$,  we mean the poset of simplices of 
$C$ ordered by inclusion.  The realization of ${\rm Simp}(C)$ is often called the barycentric subdivision of $C$ and
is homeomorphic to the realization of $C$. We define a poset map 
  $\mu_n'\colon {\rm b}\cd_n\to {\rm Simp}(X_n^{\text {alg}})$ by
$\mu_n'(G)=|G|^{(0)}$.  The poset ${\rm b}\cd_n$ was defined near the  beginning of section \ref{Xn}.
 Recall that  $|G|^{(0)}$  is the subset of $\autfn$ corresponding to all of the minimal elements of
the poset
${\rm b}G$.  It follows from section \ref{cgordering}, where the ordering determined by $G$ was discussed,
 that $|G|^{(0)}$ is a simplex of 
$X_n^{\text {alg}}$.  Also, $\mu_n'$ preserves orders because if $G_1<G_2$ in ${\rm b}\cd_n$, then ${\rm b}G_1$ is a subposet of 
${\rm b}G_2$.

Taking geometric realizations gives a map $\mu_n\colon X_n\to X_n^{\text {alg}}$ which is $\autfn$-equivarient.  
Letting $n\to\infty$ gives an AUT-equivarient map $\mu\colon X\to X^{\text {alg}}$.  Similarly we define
maps $\tilde \mu_n\colon\xt_n\to \xt_n^{\text {alg}}$ and  $\tilde \mu\colon\xt\to \xt^{\text {alg}}$ which are
$StN_n$- and StN-equivarient.

To help get a feeling for these maps, we mention that $\mu_n$ takes the 0-cells of $X_n$ bijectively to the vertices of 
$X_n^{\text{alg}}$, and takes the set of 1-cells of $X_n$ injectively to the set of edges of  $X_n^{\text {alg}}$. Also, if 
$s$ is a 2-cell of $X_n$ having the shape of the square $\sigma^1\times \sigma^1$, then $\mu_n(s)$ lies 
halfway between the two possible triangulations of $s$ (which don't introduce any new vertices) in the 
following sense.  The interior of $s$ is taken by $\mu_n$ to the interior of a 3-simplex $s'$ of 
 $X_n^{\text {alg}}$ whose four faces correspond to the two possible triangulations of $s$.  
 The four 1-cells forming the boundary of $s$ are taken to the four edges of $s'$ different from the
two diagonals of $s$.


\subsection{Maps between various K-theories and Pre-K-theories} 

The following is an application of Proposition \ref{trivact}.
\begin{theorem}\label{diag1}  There is a map $f\colon \Omega B {\rm StN}^+\to \xt^{\text {alg}}$, defined up to homotopy,
 which makes the diagram 
$$\xymatrix{\Omega B StN^+             \ar[d]^g \ar[r]^-f & \xt^{\text{alg}} \ar[d]^h \\
          \Omega B St({\mathbb Z})^+  \ar[r]^{\sim}      &  V(St({\mathbb Z})) }  $$
 homotopy commute.  The map $g$ in the diagram is induced by the homomorphism $StN\to St(\mathbb Z)$, and 
$h$  was defined in section \ref{maps}.
The lower horizontal map is the  homotopy equivalence of \cite{S}.
\end{theorem}

 Even though maps from algebraically to geometrically defined spaces are rare, we conjecture that there is a  map $\hat f$ such that the diagram

$$\xymatrix{
& \xt \ar[d]^{\tilde \mu}
\\
\Omega B StN^+  \ar[r]^-f \ar[ur]^{\hat f}& \xt^{\text{alg}}  
}  $$ homotopy commutes.

\begin{proof} [Proof of Theorem \ref{diag1}] We will show that there is a homotopy commutative diagram
$$\xymatrix@R=14pt@C=-20pt{&&{\Omega BStN^+} \ar[dd] \ar[dr]^g 
\\
&&&{\Omega BSt({\mathbb Z})^+} \ar[dd]^{\sim}
\\
{\xt^{\text {alg}}} \ar[dr]^h \ar[dd] \ar[rr]^{\sim}  & &F \ar[dr] \ar'[d] [dd]
\\
&V(St({\mathbb Z}))\ar[dd]\ar[rr]^(.4){\sim} & & F_1\ar[dd]
\\
{\xt^{\text {alg}} /StN}\ar[dd] \ar[dr] \ar'[r][rr] &  &  \xt^{\text {alg}} /StN^+ \ar'[d][dd] \ar[dr]
\\
&V(St({\mathbb Z}))/St({\mathbb Z})\ar[dd] \ar[rr] & &V(St({\mathbb Z}))/St({\mathbb Z})^+\ar[dd]
\\
BStN\ar[dr] \ar'[r][rr] & & BStN^+\ar[dr]
\\
&BSt({\mathbb Z})\ar[rr] & & BSt({\mathbb Z})^+
}$$
with the top two horizontal arrows homotopy equivalences.  
The theorem will then follow from the square on the upper right.

The diagram consists of maps of (homotopy) fibrations, with the fibrations drawn vertically.
  The covering spaces $\xt^{\text {alg}} \to \xt^{\text {alg}}/StN$
and $V(St({\mathbb Z})) \to V(St({\mathbb Z}))/St({\mathbb Z})$ along with there classifying maps 
give the two fibratations on the left.  The map between these is induced by $h$.

Applying the Quillen plus construction to the four spaces in the square on  the lower left  of the diagram gives the 
square on the lower right.  The four maps forming the edges of this last square are
induced by the universal property of the  plus construction so that the lower cube in the diagram commutes.

The spaces $F$ and $F_1$ are the homotopy fibers of the maps below them, and the loop spaces above $F$ and $F_1$
are the homotopy fibers of the maps below them in the usual way.  By Proposition \ref{trivact}, the group 
$StN\approx \pi_1(BStN)$ acts trivialy  on the homology groups of the space  $\xt^{\text {alg}}$.  Since the plus
construction does not alter homology, it follows from the spectral sequence comparison theorem that
the induced map  $\xt^{\text {alg}}\to F$ of fibers is a homology equivalence.  Since 
$\pi_2(StN^+)\approx H_2(StN)\approx 0$ and $\pi_1(\xt^{\text {alg}}/StN^+)\approx 0$, the space $F$ is
simply-connected, as is $\xt^{\text {alg}}$.  Thus by Whitehead's theorem the map $\xt^{\text {alg}}\to F$
is a homotopy equivalence.

	Similarly, the map $V(St({\mathbb Z}))\to F_1$ is a homotopy equivalence.  For this,
 one needs that the group $St({\mathbb Z})$ acts trivialy on the homology groups of $V(St({\mathbb Z}))$.
  This follows from (1.3) in \cite{S}.  Also, the complex $V(St({\mathbb Z}))$ is simply-connected since the 
defining relations for the group $St({\mathbb Z})$ are triangular.  And  $F_1$ is simply-connected by an
argument like  that, sketched above, showing  $F$ is.

Finally, the CW-complex $V(St({\mathbb Z}))/St({\mathbb Z})^+$ is simply-connected and, according to \cite{S}, acyclic.
So it is contratible.  Therefore the map $\Omega BSt({\mathbb Z})^+\to F_1$ is a homotopy equivalence.
\end{proof}

Applying the homotopy group functors to the diagram in Theorem \ref{diag1} and using standard K-theoretic identifications
gives the commutative diagram in Corollary \ref{diag2} below--once the following definitions have been made.

For each $i\ge 1$, let 
$$KN_i=\pi_i(B {\rm AUT}^+),$$  where ${\rm AUT}=\lim_{n\to \infty}Aut(F_n)$.
(Groups $KN_i(G)$  were defined in \cite{KN}, where $G$ is any group.
  But we will only use the groups $KN_i(1)$ below,
 where 1 is the trivial group, and abbreviate these by $KN_i$.)
This definition is analogous to Quillen's definition $K_i({\mathbb Z})\equiv \pi_i(B{\rm GL}({\mathbb Z})^+)$ of 
the K-groups of ${\mathbb Z}$, but with automorphisms of free groups in place of automorphisms of free
abelian groups.

	The space $B{\rm AUT}^+$ may  also be thought of as a second approximation to Waldhausen's space 
$A(*)\equiv \lim_{m,n\to \infty}{\rm Aut}(\vee_n \,S^m)$, the topological K-theory of the one-point space $*$.
Here $Aut(X)$, where $X$ is a based space, is the topological monoid of homotopy equivalences of $X$.
Indeed, a first approximation to $A(*)$ is 
$$(B\lim_{n\to \infty}Aut(\vee_n \,S^0))^+ =(B\lim_{n\to \infty} \Sigma_n)^+ = B\Sigma^+
 $$  or a compoment of $\Omega^{\infty}S^{\infty}$,
with the  second approximation being
 $$B\lim_{n\to \infty}Aut(\vee_n \,S^1)^+ \simeq B\lim_{n\to \infty} {\rm Aut}(F_n)^+=B{\rm AUT}^+.$$

For each $i\ge 1$, we also define a group $KN_i^V$ by 
$$KN_i^V=\pi_{i-1}(X^{alg}).$$

\begin{cor}\label{diag2} For each $i\ge 3$, the diagram
$$\xymatrix{
KN_i\ar[dr]^{g_*} \ar[r]^{f_*} & KN_i^{V}\ar[d]^{h_*}
\\
                               & K_i({\mathbb Z})
}$$
induced by the one in Theorem \ref{diag1} commutes.
\end{cor}

\begin{proof} We sketch the necessary K-theoretic identifications.
According to Prop. 11.2.6 of \cite{L}
  there is a homotopy fibration
$$B\,K_2({\mathbb Z}) \to B{\rm St}({\mathbb Z})^+\to B{\rm E}({\mathbb Z})^+.$$  Similarly 
$$BKN_2 \to B{\rm StN}^+\to B{\rm AUT}_*^+$$ is a homotopy fibration.  Here ${\rm AUT}_*$ is the index two subgroup of ${\rm
AUT}$ consisting of all automorphisms which abelianize to elements of ${\rm GL}({\mathbb Z})$ with positive determinant.

Since $ B{\rm E}({\mathbb Z})^+$ is the universal cover of $ B{\rm GL}({\mathbb Z})^+$ and 
$B{\rm AUT}_*^+$  the universal cover of $B{\rm AUT}^+$, the long exact homotopy sequences of the four fibrations
just mentioned give 
$$\pi_i(B{\rm StN}^+)\approx \pi_i(B{\rm AUT}^+)\equiv KN_i, \quad {\rm and}\quad
\pi_i(B{\rm St}({\mathbb Z})^+)\approx \pi_i(B{\rm GL}({\mathbb Z})^+)\equiv K_i({\mathbb Z})
$$ provided $i\ge 3$.  Also, 
$$\pi_i(\xt^{alg})=\pi_i(X^{alg}) \equiv KN_{i+1}^V$$
for $i\ge 2$ since $\xt^{alg}$ is a cover of a component of $X^{alg}$.
\end{proof}

\begin{remark} \rm The prism
$$\xymatrix@R=14pt{{\Omega(B\,KN_2)}\ar[dd] \ar@/_/[drr] \ar[r] &{\mathbb Z}_2 \ar[dd]  \ar[dr]
\\
&&{\Omega( BK_2}({\mathbb Z}) )\ar[dd]
\\
\Omega(B{\rm StN}^+) \ar@/_/[drr]_g \ar[r]^>>>>>f &\xt^{alg}\ar[dr]
\\
&&\Omega (B {\rm St}({\mathbb Z})^+)
\\
}$$ homotopy commutes where the vertical maps come from the inclusions of the fibers in the fibrations mentioned
in the proof just above.  The lower triangle is the diagram in Theorem \ref{diag1} without the space 
$V(St({\mathbb Z})$.  The induced maps of the homotopy fibers of the vertical maps give the homotopy
commutative diagram
 
$$\xymatrix{{\Omega^2(B{\rm AUT}_*^+)}\ar[r] \ar[dr]  &{\Omega X_*^{alg}} \ar[d]  
\\
&{\Omega^2( B{\rm E}}({\mathbb Z})^+ )
}$$

Applying $\pi_0$ to this gives a commutative diagram like the one in Corollary \ref{diag2} with $i=2$.  There 
is also such a commutative diagram with $i=1$.
\end{remark}

\section{Eliminating cubical terms}\label{ECT}

This section and the next two are devoted to the proof of the following.

\begin{theorem}\label{X/S}  
The homology group $H_3(\xt_n/StN_n)$ is trivial provided $n\ge7$.
\end{theorem}

The proof is broken up into a number of lemmas, with each  showing that all the elements of the
third homology group of some subcomplex of $\xt_n/StN_n$ have representatives in a smaller
subcomplex.  All of these lemmas, except for \ref{CC11}, are true if $n=6$, and many are true for lesser $n$.
Nevertheless, for convenience, we assume throughout this section and the next that $n\ge 7$.

We first introduce some terminology. 
Let $G$ be a colored or partitioned graph.  An edge of $G$ is called a singleton edge if it is the only 
element of the edgeset containing it.  The core of $G$, or $core(G)$, is the 
subgraph of $G$ consisting of all the  nonsingleton edges of $G$, together with either the coloring or
partition of these inherited from that of $G$.  
The graph $core(G)$ may contain bivalent vertices other than the basepoint.  
We often identify the vertices of $G$ with those of $core(G)$.

Let $\C{k}$ denote the group of cellular $k$-chains in $\xt_n/StN_n$
and $Z_3(\xt_n/StN_n)$ the subgroup of 3-cycles.
If $G$ is a colored graph with 
$k+1$ vertices, and a choice of orientation of the corresponding cell in $\xt_n/StN_n$ has been made,
 then $gen(G)$ denotes the corresponding generator of $\C{k}$.  
This choice of orientation is oftened referred to as a cell orientation for $G$.
Cell orientations are discussed more after the proof of Lemma \ref{creduced} below.
Conversely, if $x$ is a generator of $\C{k}$
corresponding to some $k$-cell of $\xt_n/StN_n$, then $cg(x)$ denotes the colored graph corresponding to the $k$-cell.
The underlying partitioned graph of a colored graph $G$ is denoted by $par(G)$. Usually $par(cg(x))$ is 
shortened to just $par(x)$.

A non-basepoint vertex $v$ of a colored graph $G$ is called an $m$-valent, $k$-color vertex 
if, in $core(G)$, the set of edges meeting $v$ has a subset consisting of $m$ edges
of exactly $k$ different colors.  Let $\cf_{4,2}$ be the subgroup of $Z_3(\xt_n/StN_n)$ generated by elements
corresponding to colored graphs with no 4-valent, 2-color vertices.
A $d$-cell of $\xt_n/StN_n$ corresponding to the colored graph $G$ with edge sets $h_1,\dots, h_k$
is referred to as an $n_1,\dots,n_k$-cell if there is some permutation $\sigma$ of the set $\bf k$ with
$|h_{\sigma(i)}|-1=n_i$ for each $i$. (In this case, $\sum_{i=1}^k n_i=d$.)

A cubical (resp.~prismatic, simplicial) generator of $\C{3}$ is a generator which 
corresponds to a $1,1,1$-cell (resp.~1,2-cell, 0,3-cell) of $\xt_n/StN_n$.

\begin{figure}[tb]
\centerline{\mbox{\includegraphics*{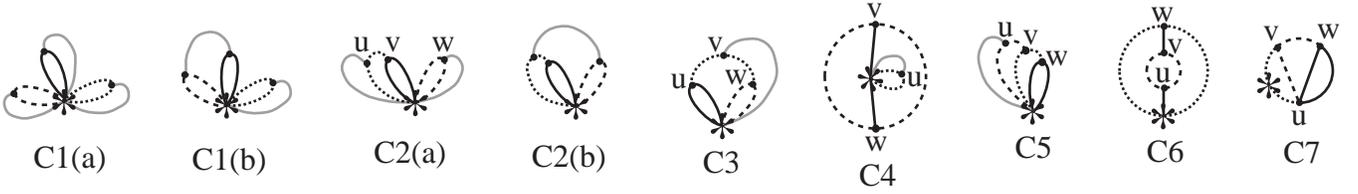}}}
\caption{All reduced, partitioned, cubical graphs.}
\label{reducedcores}
\end{figure}

We say that a colored or partitioned graph $G$ has core of type Cn 
if the partitioned graph $core(G)$ is 
isomorphic to the core of the partitioned graph labeled Cn (or Cn(a) if $n=1$ or $2$) in figure
\ref{reducedcores}.
A generator $x$ of $\C{3}$ and the corresponding colored graph $cg(x)$ are both said to be core-reduced
if $cg(x)$ has core of one of the types C1 through C7.  Let $\cf_{CR}$ be the subgroup of $Z_3(\xt_n/StN_n)$
consisting of all the cycles with core-reduced cubical terms.

A cubical generator $x$ of $\C{3}$ and $cg(x)$ are  both said to be reduced 
if the partitioned graph $par(x)$
is isomorphic to one of  the nine partitioned graphs shown in figure \ref{reducedcores}.  
Assuming $x$ is reduced, we say that both $x$ and $cg(x)$ are of type Cn if $par(x)$ is isomorphic 
to the partitioned graph labeled Cn in Figure \ref{reducedcores}.
A 3-chain is said to be reduced if all its
cubical terms are.  Let $\cf_R$ be the subgroup of $Z_3(\xt_n/StN_n)$ consisting of all the reduced
cycles.

Let $\cf$ be the filtration
$$Z_3(\xt_n/StN_n) \supset \cf_{4,2} \supset \cf_{CR} \supset \cf_R \supset \cf_{C7} \supset 
\cf_{\theta P} \supset
\cf_{C5} \supset \cf_{C4}\supset
\cdots \supset \cf_{C1}$$
of $Z_3(\xt_n/StN_n)$ where $z\in \cf_{Cn}$ if $z\in\cf_R$ and 
none of the cubical terms of $z$ are of type Cm for $m\ge n$.
The subgroup $\cf_{\theta P}$ is defined just before lemma \ref{thetaP} below.
The lemmas in this section show that any element of a subgroup of the filtration $\cf$ 
is homologous in the space
$\xt_n/StN_n$ to an element of the next smaller subgroup of $\cf$.  Thus, any element of $H_3(\xt_n/StN_n)$
has a representative in $\cf_{C1}$.

\begin{lemma} \label{42}
Each $z\in Z_3(\xt_n/StN_n)$ is homologous to a cycle in $\cf_{4,2}$.
\end{lemma}

\begin{proof}
Let $x$ be any 
cubical term of a cycle $z\in Z_3(\xt_n/StN_n)$.
Let $G$ be the colored graph $cg(x)$ where the  edgesets are $\{a_1,a_2\}$, $\{b_1,b_2\}$ and $\{c_1,c_2\}$.
Let $v$ be a 4-valent, $2$-color
vertex of $G$ meeting, say, the edges $a_1$, $a_2$, $b_1$ and $b_2$.  Since $G$ has four vertices,
at least two of these edges form a monochromatic loop, say $a_1$ and $a_2$ do.
Then if $b_1$ and $b_2$ do not form a loop, $core(G)$  is as shown in figure \ref{AB}(a). 
And if $b_1$ and $b_2$ do form a loop, then both $c_1$ and $c_2$ meet the vertex of $G$ which 
is disjoint from both of the loops $a_1$,$a_2$ and $b_1$,$b_2$. The four possibilities are
shown in b(i) through b(iv) of figure \ref{AB}.
\begin{figure}[tb]
\centerline{\mbox{\includegraphics*{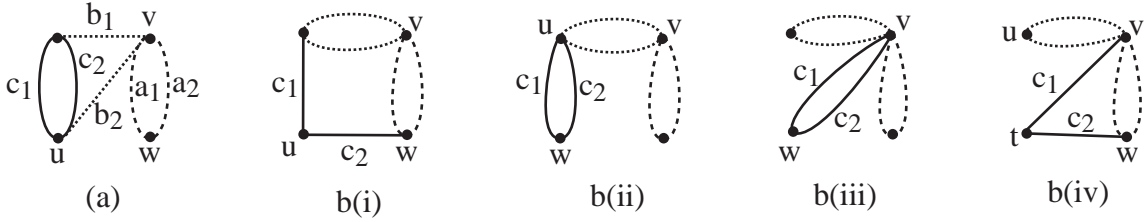}}}
\caption{All (unbased) cores with a 4-valent, 2-color vertex.}
\label{AB}
\end{figure}

In case the core of $G$ is as in figure \ref{AB}(a) and $u$ is the basepoint of $G$, let
$G'$ be a colored graph with five vertices with core as shown in figure \ref{CD}(a) and satisfying 
$G_b'=G$.  Orient $cell(G')$ so that the boundary $d(gen(G'))$ contains the term
$-x$. We say that $gen(G')$ is obtained from $x$ by separating the edges $a_1$ and $a_2$
of $G=cg(x)$ along the edges $b_1$ and $b_2$ of $G$.  Similar terminology is often used below.
Note that the two cubical terms
$G_{b_1}'$ and 
$G_{b_2}'$ of the chain $x + d(gen(G'))$ contain no 4-valent, 2-color vertices.
Such vertices can be eliminated in similar ways in the other cases by adding suitable boundaries.
For instance, in case the core of $G$ is as in figure \ref{AB}(a) with $w$ as the basepoint, add to
$x$ the boundary of the generator of $\C{4}$ obtained from $x$ by separating the edges $b_1$ and $b_2$ of $G$ 
along $a_1$ and $a_2$.  If $core(G)= b(i) u$, meaning that $core(G)$ is as in b(i) of figure \ref{AB}
with $u$ as the basepoint, or if $core(G)= b(i) w$ with a similar meaning, then separate, 
along the monochromatic loop through $v$ and
$w$,  the two edges of the other monochromatic loop of $G$. If $core(G)= b(ii) u$, then separate, along 
the monochromatic loop through $u$ and $v$, 
the two edges of the other monochromatic loop through $v$.  If $core(G)= b(ii) w$, first separate the 
two edges of the monochromatic loop through $u$ and $v$ along $c_1$ and $c_2$, and then use the first case above
twice. If $core(G)= b(iii) w$,  then first separate the edges $a_1$ and $a_2$ from $b_1$ and $b_2$ along
$c_1$ and $c_2$, and then use a previous case twice.  If $core(G)= b(iv) u$ or $w$, then separate, 
along the monochromatic loop through $*$, the two edges of the other monochromatic loop of 
$core(G)$.  Finally, if $core(G)= b(iv)t$, then use a graph with core as shown in figure \ref{CD}(b).
This covers all the cases, and so completes the proof.
\end{proof}
\begin{figure}[tb]
\centerline{\mbox{\includegraphics*{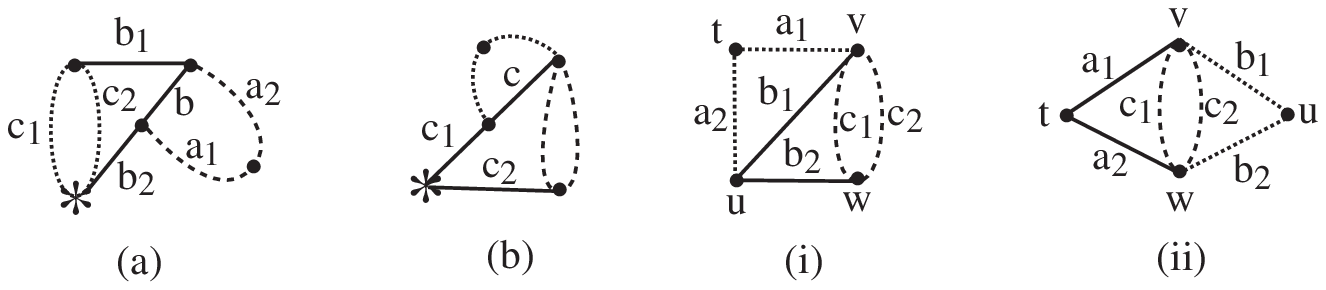}}}
\caption{}
\label{CD}
\end{figure}

\begin{lemma} \label{CR}
Each $z\in \cf_{4,2}$ is homologous to an element of $\cf_{CR}$.
\end{lemma}

\begin{proof}
	Let $x$ be a cubical term of $z$ such that $core(cg(x))$ has a 4-valent, 3-color vertex $v$.
Since $v$ meets all nonsingleton colors and is not a 4-valent, 2-color vertex, $v$ lies on exactly
one monochromatic loop.  We claim that $core(cg(x))$ is as shown in either (i) or (ii) of figure 
\ref{CD}.  Indeed, given that $a_1$, $b_1$, $c_1$ and $c_2$ meet $v$, of the six ways $a_2$ can 
be attached to the four vertices of $cg(x)$, all but the two ways shown are ruled out by the fact 
that  $cg(x)$ is colored and $v$ is not a 4-valent, 2-color vertex.  The same facts then imply that
$b_2$ must be attached as shown, or that, in figure \ref{CD}(ii), $b_2$ joins the vertices $t$ and
$u$.  This last possibility gives the same core as in (i).

Next we indicate, in all but one case, how the 4-valent, 3-color vertices of $cg(x)$ can be
eliminated by adding suitable boundaries to $x$, without introducing any 
4-valent, 2-color vertices.  Let $G=cg(x)$. If $core(G)=\ref{CD}(i)w$, meaning that the core of 
$G$ is the partitioned graph shown in (i) of figure \ref{CD} with $w$ as the basepoint, or if 
$core(G)=\ref{CD}(ii)w$ 
then, in both cases, separate $a_1$ and $b_1$ along $c_1$ and $c_2$.  If $core(G)=\ref{CD}(i)u$, use a graph 
with core as shown in figure \ref{F1}(a), and if $core(G)=\ref{CD}(ii)t$, first use a graph with core
as shown in figure \ref{F1}(b), and then use a previous case.
\begin{figure}[tb]
\centerline{\mbox{\includegraphics*{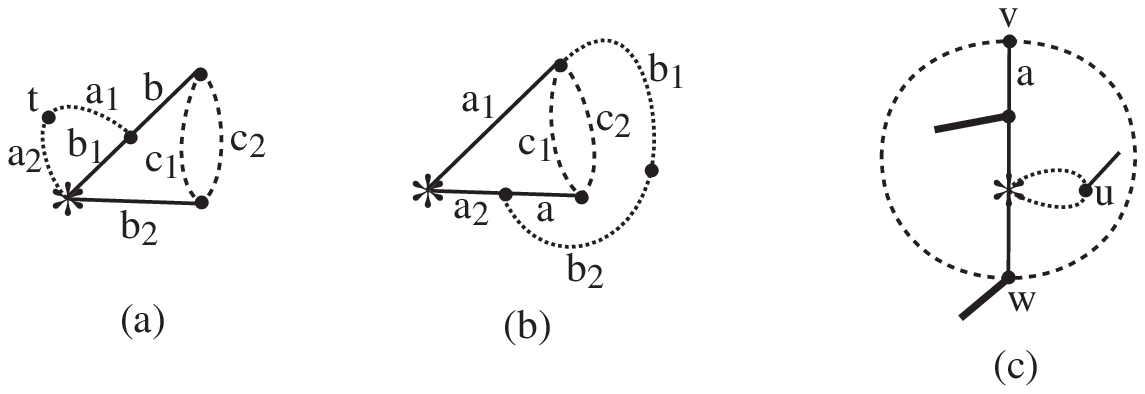}}}
\caption{}
\label{F1}
\end{figure}

 We may now assume that $G$ has no 
4-valent, 2-color vertices and if $G$ has a 4-valent, 3-color vertex, then $core(G)$ is as 
shown in figure \ref{CD}(i) with $t$ as the basepoint, that is, $core(G)$ is of type C7.
To complete the proof, we show that  $core(G)$ has type $Cn$ for some $n$ between 1 and 7.

First assume that $G$ has more than one monochromatic loop.  Then since $G$ has no
4-valent, 2-color vertices, these loops are either disjoint or meet at the basepoint.  If they're
disjoint, then two edges of the same color do not form a loop and these two edges are either
disjoint or meet at $*$.  In these cases $core(G)$ is either of type C6 or C4.  If the monochromatic loops
meet at $*$, then there are either two or three of them, and so $core(G)$ is of type C1, C2 or C3.

Now assume that $G$ has just one monochromatic loop consisting, say, of the edges $a_1$ and $a_2$,
and let $b_1$ and $b_2$ be two other edges of $G$ of the same color.  If $b_1$ and $b_2$ form a
disconnected subgraph of $G$, then $core(G)$  is of type C6. So assume $b_1 \cap b_2 \ne \varnothing$.
If $(b_1 \cap b_2) \cap (a_1 \cup a_2) \ne \varnothing$, then this intersection consists of the 
basepoint of $G$, and so $core(G)$ is of type C4.  If $(b_1 \cap b_2) \cap (a_1 \cup a_2) =\varnothing$,
then $(b_1 \cup b_2) \cap (a_1 \cup a_2)$ consists of either the one or two vertices of the loop
$a_1 \cup a_2$.  In both cases it follows from the assumptions in effect and the fact that each
non-basepoint vertex of $G$ is 2-valent, 1-color that $core(G)$ is of type C5 or C7.
\end{proof}

\begin{definition} \rm If a colored graph has a non-basepoint vertex $v$ such that the set of edges meeting 
$v$ has only one subset consisting of two edges of the same color, then, by separating the edges
meeting $v$, we mean separate, along these two edges of the same color, the other half-edges incident
with $v$.

By the $w$-vertex of a colored graph $G$ with type C2 core, we mean the vertex 
which corresponds to the vertex labeled $w$ of the graph in C2(a) of Figure \ref{reducedcores}. 
Similar conventions, using the labeling of the vertices in Figure \ref{reducedcores}, are used to designate the
other non-basepoint vertices of $G$, as well as the the vertices of colored graphs with other core types.
\end{definition}

\begin{lemma} \label{creduced}
Each $z\in \cf_{CR}$ is homologous to a cycle in $\cf_{R}$.
\end{lemma}

\begin{proof}
Let $x$ be any cubical term of $z\in \cf_{CR}$, and let $G=cg(x)$.  Each of the possible  seven core 
types for $G$ will be considered separately.

\noindent {\it Case C1.}
Assume first that $G$ has type C1 core.  If a non-basepoint vertex, say $v$, of $G$ has valence $>3$, then separate the edges of $G$ meeting $v$ to obtain the colored graph $G'$.  Choose the orientation
of the corresponding cell so that $x$ cancels with one of the terms of $d(gen(G'))$.
Then the colored graphs corresponding to the two cubical terms of 
$x+d(gen(G'))$ both have type C1 cores.  Also, the sum of the valences of the non-basepoint vertices 
of each of these is less than that of $G$.  So by separating edges repeatedly as needed, we may assume that
each cubical term of $z$ with type C1 core is reduced.

\vskip3pt
\noindent {\it Case C2.}
Assume now that $G$ has type C2 core.    By separating the edges incident with the $w$-vertex of $G$, 
and then possibly of resulting colored graphs with type C2 cores,
we may assume that the $w$-vertex has valence 3. 
Then by separating the edges meeting the $u$-vertex of $G$, we may assume that
both the $u$- and $w$-vertices of $G$ have valence 3.  Finally, let $gen(G')$ be obtained from $x$ by
separating the singleton edges meeting the $v$-vertex of $G$ from the nonsingleton edge which joins the
$u$- and $v$-vertices.  
Then one cubical term of $x+ d(gen(G')$ corresponds 
to a colored graph with type C1 core.  The other corresponds to one with type C2 core with $v$-vertex
having valence 3.  This and case C1 give a chain
$c$ such that $x+dc$ is reduced.  In fact, each cubical term has type C1 or C2.

\vskip3pt
\noindent {\it Case C3.}
Next assume that $G$ has type C3 core. By  repeatedly separating the edges incident with the $v$-vertex of $G$,
 we may assume that it has valence 3.  
Then separating edges incident with the $u$- and $w$-vertices
of $G$ and possibly resulting colored graphs with type C3 cores and the previous  cases give a  4-chain 
$c$ such
that $x + dc$ is  reduced.

\vskip3pt
\noindent {\it Case C4.}
If $G$ has type C4 core, then, as usual, we may assume that the $u$-vertex of $G$ has valence 3.
 Let 
$G'$ be the partitioned graph shown in (c) of figure \ref{F1}, colored and with dangling half-edges joined so
that $G_a' =G$. 
Orient the cell of $\xt_n/StN_n$ corresponding to $G'$ so that $d(gen(G'))$ contains
the term $-x$.  We say that $gen(G')$ is obtained from $x$ by moving the singleton edges incident with 
the $v$-vertex of $G$ toward the basepoint of $G$.  Similar terminology is often used below.  
The two colored graphs corresponding to the cubical terms of $x+d(gen(G'))$ have cores of type C2 and C4.
The valence of the $v$-vertex of the one with type C4 core is 3.  So moving the singleton edges
incident with its $w$-vertex toward the basepoint
and  case C2 give a chain $c$ such that $x+dc$ is reduced.

\vskip3pt
\noindent {\it Case C5.}
If $G$ has type C5 core then, as usual, we may assume that the $u$-vertex of $G$ has valence 3. 
Moving the singleton half-edges incident with 
the $v$-vertex of $G$ toward the $w$-vertex gives a chain $c'$ such that the two colored graphs 
corresponding to
the cubical terms of $x+dc'$ have cores of type C2 and C5. The valence of both the $u$- and $v$-vertices
of the one with type C5 core is 3.  So separating edges meeting its $w$-vertex and using case C2 gives a
chain $c$ such that $x+dc$ is reduced.

\vskip3pt
\noindent {\it Case C6.}
Assume that $G$ has type C6 core.  By moving the singleton edges of $G$ incident with the $u$-vertex
toward the basepoint, we may assume that the $u$-vertex has valence 3.  For this, a 4-valent, 2-color
vertex has to be eliminated as in the proof of \ref{42}, and case C5 is used. Then by
moving the singleton edges meeting the $v$-vertex of $G$ toward the $w$-vertex and using case 
C3, we may assume that the $v$-vertex of $G$ has valence 3. Finally, moving the singleton edges incident
with the $w$-vertex toward the basepoint along either edge and then using case C4 
gives a chain $c$ such that $x+dc$ is reduced.  The terms of $x+dc$ can be of all possible types except
for C7.

\vskip3pt
\noindent {\it Case C7.}
If $G$ has type C7 core, then first let $gen(G')$ be obtained from $x$  by moving the singleton edges
incident with the $w$-vertex of $G$ toward the $v$-vertex of $G$.  Then one cubical term of $x+d(gen(G'))$
corresponds to a colored graph with type C7 core with $w$-vertex of valence 3. The other has a 4-valent,
2-color vertex.  This can be eliminated as in the proof of \ref{42} at the expense of introducing a colored
graph with a type C5 core, and another with a 4-valent, 3-color vertex. Eliminating this vertex, in turn, as
done above, gives colored graphs with cores of type C3 and C5. So by previous cases we may assume that the
$w$-vertex of
$G$ has valence 3.

We may also assume that the $v$-vertex of $G$ also has valence 3 by moving the singleton edges incident with
the $v$-vertex toward the basepoint of $G$, and then using  case  C5.
Finally, by moving the singleton edges incident with the $u$-vertex of $G$ toward the basepoint, eliminating,
as in the proof of \ref{CR}, one of 
the 4-valent, 3-color vertices which arise, and then using  cases C4 and C5,
we obtain a chain $c$ such that $x+dc$ is reduced. The terms of $x+dc$ could be of any type other than C6.
\vskip3pt
This covers all the cases and so shows that each term of $z$ is homologous to some element of 
$\cf_{R}$.  So $z$ also is.
\end{proof}

We next discuss cell orientations of colored graphs and the effect the cellular boundary map 
has on these. Let $G$ be a StN-marked graph with edgesets $h_1,\ldots ,h_n$, and
let $cm(G)$ be the corresponding cell in $\widetilde X$. 
Also let $c(G)$ be
the cell in $\xt /StN$ corresponding to the unmarked colored graph $G$.  By the
discussion following Lemma \ref{triinj}, the cell $cm(G)$ is
homeomorphic to the product $\sigma^{\| h_1\|} \times \ldots \times \sigma^{\|
h_n\|}$ where $\sigma^{\| h_i\|}$ is the standard simplex of dimension $\|
h_i\| = |h_i| -1$.  Also, the vertices of each $\sigma^{\| h_i\|}$ correspond,
in a natural way, to the elements of the set $h_i$.  Thus a choice, up to an
even permutation, of an ordering of the edges in $h_i=\{e_i^j\}_{j=0}^{k_i}$
serves to orient $\sigma^{ \|h_i\|}$.  For each $i$, let $e_i^{\pi_i(0)},
\ldots , e_i^{\pi_i(k_i)}$ be such a choice (where $\pi_i$ is a bijection of
the set $\{1, \ldots , k_i\}$) and let $[e_i^{\pi_i(0)},\ldots
,e_i^{\pi_i(k_i)}]$ denote $\sigma^{\|h_i\|}$ with the resulting
orientation.  The cells $cm(G) \approx \sigma^{\|h_i\|} \times \ldots
\times \sigma^{\|h_n\|}$ and $c(G)$ are then oriented by the product
orientation.  Let $[e_1^{\pi_1(0)},
\ldots , e_1^{\pi_1(k_1)}] \times\ldots\times [e_n^{\pi_n(0)},
\ldots , e_n^{\pi_n(k_n)}]$ denote both the cell $c(G)$, so oriented, and also
the generator of the  cellular chain group of $\xt /StN$ it corresponds to.
With these conventions the cellular boundary map $d\colon \C{*} \to 
\C{{*-1}}$ is given by

\begin{eqnarray*}
  d\left([e_1^{\pi_1(0)}, \ldots, e_1^{\pi_1(k_1)}]
\times\cdots\times [e_n^{\pi_n(0)}, \ldots , e_n^{\pi_n(k_n)}] \right) = && \\
  \sum_{\buildrel {1\le i\le n, k_i\ne 0} \over {\scriptscriptstyle 0\le j\le
k_i}} (-1)^{k_1+\cdots +k_{i-1}} (-1)^j   [e_1^{\pi_1(0)},
\ldots , e_1^{\pi_1(k_1)}] \times\ldots\times & \!\!\![e_i^{\pi_i(0)},
\ldots,  \widehat{e_i^{\pi_i(j)}} ,\ldots, e_i^{\pi_i(k_i)}]
\times \cdots & \\ 
& \cdots \times  [e_n^{\pi_n(0)},
\ldots , e_n^{\pi_n(k_n)}] 
\end{eqnarray*}
 where the term containing $\widehat{e_i^{\pi_i(j)}}$ denotes the
generator of
$C_*({\xt /StN)}$ corresponding to the colored graph $G_{e_i^{\pi_i(j)}}$ with
orientation that inherited from $G$.  Also the exponent $k_i + \cdots
+k_{i-1}$ is $0$ in case $i=1$.  Just the orientation of the oriented cell 
$[e_1^{\pi_1(0)},
\ldots , e_1^{\pi_1(k_1)}] \times\ldots\times [e_n^{\pi_n(0)},
\ldots , e_n^{\pi_n(k_n)}]$ is denoted by  $[e_1^{\pi_1(0)},
\ldots , e_1^{\pi_1(k_1)}] \ldots [e_n^{\pi_n(0)},
\ldots , e_n^{\pi_n(k_n)}]$.   Singleton edges are often dropped from the
notation.  So, for instance, the orientation $[e_1^0],[e_2^0,e_2^1],[e_3^0]$
may be denoted by just $[e_2^0,e_2^1]$. 

\begin{figure}[tb]
\centerline{\mbox{\includegraphics*{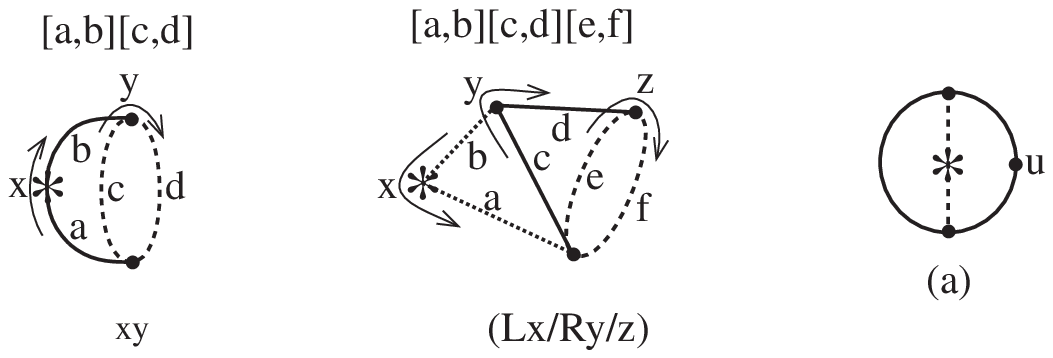}}}
\caption{}
\label{theta}
\end{figure}

In preparation for the proof of the next lemma, we introduce notation for a number of generators of 
$\C{*}$, and also define a homomorphism involving these which counts relations 
of the form $w_{ij}^L=w_{ij}^R$.
 Let $\theta_{xy}$ and $(Lx/Ry/z)$ be the generators of $\C{2}$ and $\C{3}$ indicated in figure 
\ref{theta}.  The three other generators of $\C{3}$, analogous to $(Lx/Ry/z)$, but with different 
orientations of the edgesets $x$ and $y$, are denoted by changing the
prefixal $L$ and $R$ accordingly.  
We say that a colored graph $G$ with three vertices is of type $\theta$ if the core of $G$ is 
equivalent to the colored graph $\theta_{xy}$ in Fig.\ \ref{theta} for some choice of colors $x$ and $y$.

\begin{figure}[tb]
\centerline{\mbox{\includegraphics*{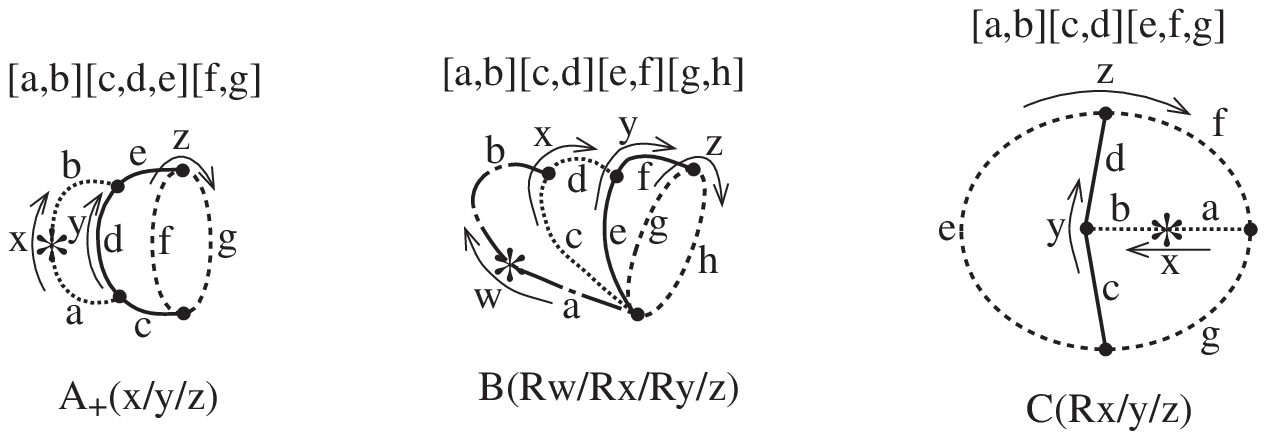}}}
\caption{}
\label{4theta}
\end{figure}

Let $A_+(x/y/z)$ be the generator of $\C{4}$ indicated in figure \ref{4theta}. 
  Another generator $A_-(x/y/z)$ is defined similarly, but 
with the orientation of the edgeset $y$ reversed.  Let $A'_+(x/y/z) = A_+(x/y/z)+c$  where $c$ is a chain
such that
\begin{equation} \label{A+}
d(A'_+(x/y/z)) \equiv (Lx/Ly/z) -(Rx/Ry/z)  \ \ \ \ (\text{mod} \ \cf_{C7}).
\end{equation}
The natural choice for $c$ is $gen(G)$ where $G$ is a colored graph as on the left of figure \ref{4theta},
 but with the partition
$\{a,b,d\}\{c,e\}\{f,g\}$ and with suitable cell orientation.  With this choice, the element of $\cf_{C7}$
by which the two equivalent elements in (\ref{A+}) differ
contains two cubical terms, both with type C6 cores.  Define $A'_-(x/y/z)$ similarly so that
\begin{equation} \label{A-}
d(A'_-(x/y/z)) \equiv (Lx/Ry/z) -(Rx/Ly/z) \ \ \ \ (\text{mod} \ \cf_{C7}).
\end{equation}

Let $B(Rw/Rx/Ry/z)$ be the generator indicated in figure \ref{4theta}.  The seven other generators obtained by changing the orientations of the edgesets
$w$, $x$ and $z$ are denoted in the same way, but with with suitable prefixes for $w$, $x$ and $z$. 
Applying lemmas \ref{42}, \ref{CR} and \ref{creduced} and their proofs to the six noncanceling faces of
$B(Pw/Qx/Sy/z)$, where
$P,Q,S\in \{L,R\}$, gives a chain $c$ such that, with $B'(Pw/Qx/Sy/z)=B(Pw/Qx/Sy/z) +c$, we have 
\begin{equation} \label{B}
d(B'(Pw/Qx/Sy/z)) \equiv -(Qx/Sy/z) +(Pw/Sy/z) - (Pw/Qx/z)\ \ \ \ (\text{mod} \ \cf_{C7}). 
\end{equation}

Let $C(Rx/y/z)$ be the generator indicated in  figure \ref{4theta},
 and let $C(Lx/y/z)$ be the same but with the
orientation of the edgeset $x$ reversed. 

For each $S\in \{L,R\}$, let $c_S$ be a generator of 
$\C{4}$ such that with $C'(Sx/y/z)=C(Sx/y/z)+ c_S$ we have
\begin{equation} \label{eqfive}
d(C'(Sx/y/z))\equiv (Sx/Ry/z)-(Sx/Ly/z)\ \ \  \ (\text{mod} \ \cf_{C7}).
\end{equation}
In both cases, $par(c_S)$ can be taken to be as in C6 of figure \ref{reducedcores}, but with the
basepoint a bivalent vertex on one of the monochromatic loops.
Define a homomorphism $f\colon \C{2} \to \oplus_{1\le i,j\le n} \mathbb{Z}_{i,j}$ by mapping
generators according to 
$$f(x)= \begin{cases}
0,   &\text{if $core(cg(x))$ is not of type $\theta$,}\\
1_{ij}, &\text{if $core(cg(x)) = \theta_{ij}$.}
\end{cases}$$
Here each $\mathbb{Z}_{i,j}$ is a copy of $\mathbb{Z}$ and $1_{i,j} \in \mathbb{Z}_{i,j}$ is the generator
corresponding to $+1\in \mathbb{Z}$.  Also, in the second line of the definition,
we are viewing  $core(cg(x))$ as a generator of $\C{2}$ with cell
orientation that inherited from $x$. Let $\Theta$ be the composite $f\circ d$ where $d\colon \C{3} \to 
\C{2}$ is the cellular boundary map.  Then $\Theta$ is trivial on simplicial generators since the cells
corresponding to these have no square faces.  Also, one can check that $\Theta$ is trivial on prismatic
generators.  The only reduced cubical terms on which $\Theta$ is nontrivial are those of type C7, and
for these  $\Theta(Sx/Ty/z)= 1_{xz} -1_{yz}$ for all four choices of $S,T\in \{L,R\}$.
\begin{lemma}
Each $z\in \cf_{R}$ is homologous to a cycle in $\cf_{C7}$.
\end{lemma}

\begin{proof}
Let $z_1$ be any type C7 cubical term of $z \in \cf_R$.  Since $\Theta(z)=0$, there is a cubical term
$z_2$ of $z$ different from $z_1$ such that some cancelation occurs among the terms of $\Theta(z_1)$ and
$\Theta(z_2)$.   The following three cases cover the various forms $z_2$ can have.  In all cases, 
we show that, by adding boundaries to $z$, it is possible to reduce the total number of type C7 terms of $z$,
while keeping $z\in \cf_R$.
\vskip3pt
\noindent {\it Case 1.}
Assume that $z_1 + z_2 = (Sx/Ty/z)+(Pw/Qx/z)$.  By (\ref{eqfive}) or (\ref{A+}) and (\ref{A-}), we may assume 
 that $S=Q$.  If $y\ne w$, then, using (\ref{B}), we see that  $z+
d(B'(Pw/Sx/Ty/z))$ is in $\cf_{R}$ and has one fewer type C7 term than $z$.  If $y=w$, we may
assume, again by (\ref{eqfive}), that $T=P$.   Then $z+d(B'(Lw/Sx/Ty/z) +B'(Lw/Ty/Sx/z)$ is in $\cf_R$ and has two
fewer type C7 terms than $z$.
\vskip3pt
\noindent {\it Case 2.} Assume $z_1 + z_2 = (Sx/Ty/z) - (Px/Qw/z)$.  By (\ref{A+}) and (\ref{A-}), 
we may assume that $S=P$.  If $y\ne w$, then $z+d(B'(Sx/Ty/Qw/))$ is in $\cf_R$ and has one fewer
type C7 term than $z$.  If $y=w$, and $z_1 + z_2\ne 0$, then in view of (\ref{eqfive}), both $z_1$ and $z_2$ can be
eliminated from $z$.
\vskip3pt
\noindent {\it Case 3.} The last case, in which $z_1 + z_2 = (Sx/Ty/z) - (Pw/Qy/z)$ is similar to the
previous ones.  First suppose that $x\ne w$. We may assume that $T=Q$ by
(\ref{eqfive}), and then add
$-d(B'(Sx/Pw/Ty/z)$ to $z$.  If
$x=w$, we may first assume that $S=P$ by (\ref{A+}) and (\ref{A-}), and then finish with the help of
(\ref{eqfive}).
\end{proof}

Both a generator $x$ of $\C{3}$ and the corresponding colored graph $cg(x)$ are referred to as 
$\theta$-prismatic if $core(cg(x))$ is as shown in (a) of figure \ref{theta}.  The $u$-vertex of a 
$\theta$-prismatic colored graph is the one corresponding to the vertex labeled $u$ in the figure.
The subgroup $\cf_{\theta P}$ of the filtration $\cf$ is defined by $z \in \cf_{\theta P}$ if $z\in \cf_{C7}$
and if $z$ contains no $\theta$-prismatic terms.

\begin{lemma} \label{thetaP}
Each $z\in \cf_{C7}$ is homologous to a cycle in $\cf_{\theta P}$.
\end{lemma}

\begin{figure}[tb]
\centerline{\mbox{\includegraphics*{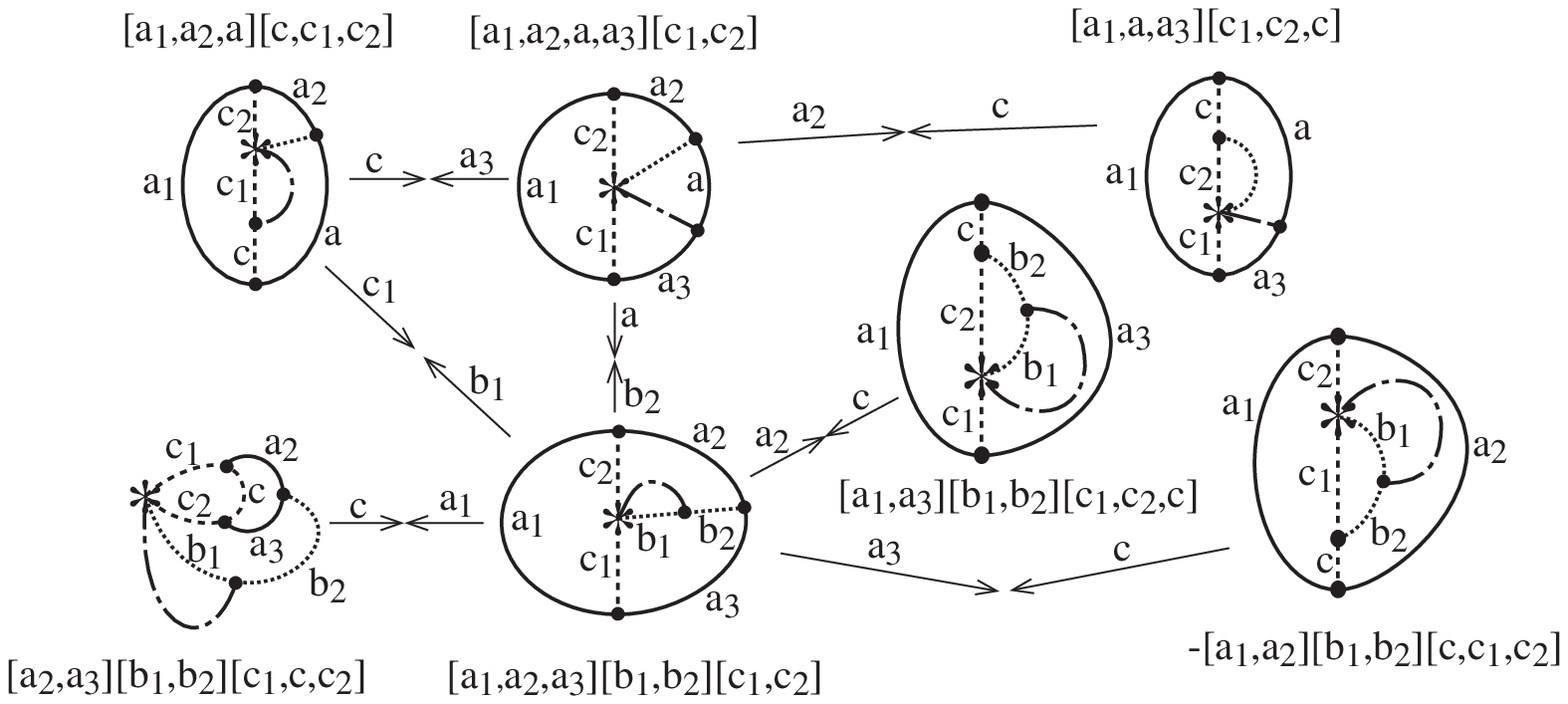}}}
\caption{}
\label{elimthetap}
\end{figure}

\begin{proof}
Any $\theta$-prismatic term of $z$ differs by boundaries from a sum of simplicial terms and prismatic
terms where each $\theta$-prismatic term in this sum corresponds to a graph with all nonbasepoint 
vertices of valence 3.  To see this, let $G$ be the colored graph corresponding to any 
$\theta$-prismatic generator.  By separating the singleton edges incident with the $u$-vertex of $G$,
and then doing this with any resulting $\theta$-prismatic colored graphs having $u$-vertices with valence $> 3$,
we may assume that the $u$-vertex of $G$ is trivalent.  Then by moving the singleton edges of $G$ 
incident with the other two non-basepoint vertices toward the basepoint, we may assume that each 
$\theta$-prismatic term of $z$ corresponds to a graph with all non-basepoint vertices trivalent.
Next we show how to eliminate such terms.

Let $c\in\C{4}$ be the sum of the seven generators, $g_1,\dots,g_7$ say, corresponding to the seven 
colored graphs and orientations shown in figure \ref{elimthetap}.  Any consistent choice of orientations
of the edgesets may be used.  Pairs of arrows pointing at each other indicate cancellation among
the terms $dg_i$ of $dc$.  For instance, if $G$ is the colored graph in the upper right of figure
\ref{elimthetap}, and $H$ is the colored graph to its left, then $G_c=H_{a_2}$, and the cell 
orientations are such that the two corresponding elements of $\C{3}$ cancel.  With the help of 
figure \ref{elimthetap}, it's routine to check that $dc$ is reduced with each cubical term of type 
C4 or C5, and that $dc$ has exactly one $\theta$-prismatic term coming from the colored 
graph $G_{c_2}$.  Since each non-basepoint vertex of $G_{c_2}$ is trivalent, all the $\theta$-prismatic terms
 of $z$ can be eliminated.
\end{proof}

\begin{figure}[tb]
\centerline{\mbox{\includegraphics*{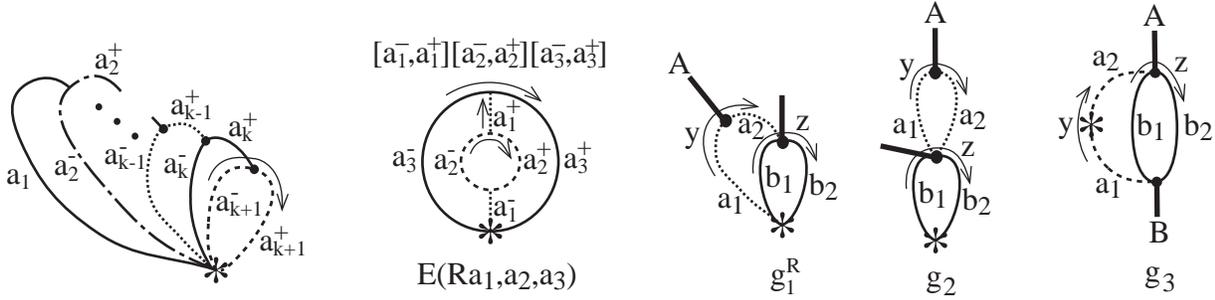}}}
\caption{Give each of $g_1^R$, $g_2$ and $g_3$ the  cell orientation $[a_1,a_2][b_1,b_2]$.}
\label{falls}
\end{figure}

Given $S_1,\ldots,S_k \in \{L,R\}$ and distinct positive integers $a_1,\ldots, a_{k+1}$, let 
$F_k(S_1a_1,\ldots,S_k a_k,a_{k+1})$, or just $F_k(S_ia_i,a_{k+1})$, be the colored graph with $k+1$ vertices
shown on the left of
 figure \ref{falls}.  Here the edgeset $h_{a_i}$ points away from the basepoint if $S_i=R$ and towards
it if $S_i=L$.  We also let $F_k(S_ia_i,a_{k+1})$ denote the generator of $\C{k}$ corresponding to the 
colored graph $F_k(S_ia_i,a_{k+1})$ with cell orientation $[a_2^-,a_2^+]\ldots [a_{k+1}^-,a_{k+1}^+]$.
Here we are assuming that $a_i^-$ meets the basepoint for $2\le i\le k$, and that the orientation of 
the loop which starts at the basepoint, goes along $a_{k+1}^-$, and then along $a_{k+1}^+$ 
back to the
basepoint, agrees with the 
orientation of the edgeset $\{a_{k+1}^-,a_{k+1}^+\}$.  A hat over a symbol in the following indicates,
as usual, omission of that symbol.

\begin{lemma} \label{0}
There is a chain $c\in\C{4}$ such that
\begin{align}
d(F_4(S_ia_i,a_5)&+c) \equiv -F_3(S_1a_1,\widehat{S_2a_2},\ldots,a_5) + \notag \\
&F_3(S_1a_1,S_2a_2,\widehat{S_3a_3},S_4a_4,a_5) -       
F_3(S_1a_1,\ldots,\widehat{S_4a_4},a_5) \ \ \ \ (\text{mod} \ \cf_{C3}).  \notag
\end{align}
\end{lemma}

\begin{proof}
We denote $F_3(S_1a_1,\ldots,\widehat{S_ia_i},\ldots,a_5)$ by just $F_3(\widehat{S_ia_i})$.
Assume that the nonsingleton edgesets of $F_4(S_ia_i,a_5)$ are $\{a_i^-,a_i^+\}$ as in figure \ref{falls}.
Let $d_{a_i^+}(F_4(S_ia_i,a_5))$ be the term of the boundary $d(F_4(S_ia_i,a_5))$ corresponding to 
the colored graph $F_4(S_ia_i,a_5)$ with the edge $a_i^+$ blown down.  Note that, for $i=2,3$ and $4$,
the colored graph $cg(d_{a_i^+}(F_4(S_ia_i,a_5)))$ contains a colored subgraph $\Gamma_i$ which, with
its inherited cell orientation, 
satisfies $gen(\Gamma_i) = (-1)^{i+1} F_3(\widehat{S_ia_i})$.

Next we define a chain $c=\sum_{i=1}^6 c_i \in \C{4}$ whose boundary contains terms corresponding to
the subgraphs $\Gamma_i$ and also terms which cancel the $d_{a_i^+}(F_4(S_ia_i,a_5))$ for $i=2,3$ and $4$.
Let $c_1$ be obtained from $d_{a_4^+}(F_4(S_ia_i,a_5))$ by separating the edges meeting the 4-valent  vertex of $cg(d_{a_4^+}(F_4(S_ia_i,a_5)))$. 
One cubical term of $d_{a_4^+}(F_4(S_ia_i,a_5))+dc_1$
is $-F_3(\widehat{S_4a_4})$. The other is of type C2(a).  Define $c_2$ similarly by separating the 
edges meeting the 4-valent vertex of $d_{a_3^+}(F_4(S_ia_i,a_5))$.  Then 
define $c_3$ by separating the edges meeting the 4-valent  vertex of the graph corresponding to the cubical
term of $dc_2$ different from $d_{a_3^+}(F_4(S_ia_i,a_5))$ and $F_3(\widehat{S_3a_3})$. 
 Define $c_4$, $c_5$ and 
$c_6$ similarly (by separating edges meeting 4-valent vertices) so as to isolate the subgraph $\Gamma_2$
of $d_{a_2^+}(F_4(S_ia_i,a_5))$ corresponding to $-F_3(\widehat{S_2a_2})$.  Its easy to check that
with $c=\sum_{i=1}^6 c_i$, the lemma holds.  No $\theta$-prismatic terms arise since each edgeset of each
colored graph considered contains an edge which meets the basepoint of $G$.
\end{proof}

\begin{lemma} \label{falll1}
For any choice of $S_1,S_2,S_3 \in \{L,R\}$ and distinct positive integers $a_1,\ldots,a_4$, there
is a chain $a\in \C{4}$ such that
\begin{equation}
da \equiv  F_3(S_1a_1,S_2a_2,S_3a_3,a_4)+ F_3(S_1a_1,S_3a_3,S_2a_2,a_4)\ \ \ \ (\text{mod} \ \cf_{C3}).
\notag
\end{equation}
\end{lemma}

\begin{proof}
Choose a positive integer $a_5$ distinct from each of $a_1,\ldots,a_4$.  Use lemma \ref{0} to pick
$c_1$ and $c_2 \in \C{4}$ such that, mod $\cf_{C3}$,
\begin{align}
d(F_4(&S_1a_1,S_2a_2,S_3a_3,S_5a_5,a_4)+c_1) \equiv  \notag \\ 
      &-F_3(S_1a_1,S_3a_3,S_5a_5,a_4)+F_3(S_1a_1,S_2a_2,S_5a_5,a_4)-F_3(S_1a_1,S_2a_2,S_3a_3,a_4),
   \notag
\end{align}
and 
\begin{align}
d(F_4(&S_1a_1,S_3a_3,S_2a_2,S_5a_5,a_4)+c_2) \equiv  \notag \\ 
      &-F_3(S_1a_1,S_2a_2,S_5a_5,a_4)+F_3(S_1a_1,S_3a_3,S_5a_5,a_4)-F_3(S_1a_1,S_3a_3,S_2a_2,a_4).
   \notag
\end{align}
Adding these two congruences completes the proof.
\end{proof}

Let $E(Ra_1,a_2,a_3)$ be the generator of $\C{3}$ indicated in figure \ref{falls}.  In case the 
edgeset $\{a_1^-,a_1^+\}$ has orientation opposite that shown, then the corresponding generator is
denoted by $E(La_1,a_2,a_3)$.

\begin{lemma} \label{falll2}
For any choice of $S_1, S_2 \in \{L,R\}$ and distinct positive integers $a_1,\ldots,a_4$, there
is a chain $a\in \C{4}$ such that, mod $\cf_{C5}$,
\begin{equation}
da \equiv  F_3(S_1a_1,S_2a_2,La_3,a_4)- F_3(S_1a_1,S_2a_2,Ra_3,a_4)-E(S_1a_1,a_3,a_4). \notag
\end{equation}

\end{lemma}

\begin{figure}[tb]
\centerline{\mbox{\includegraphics*{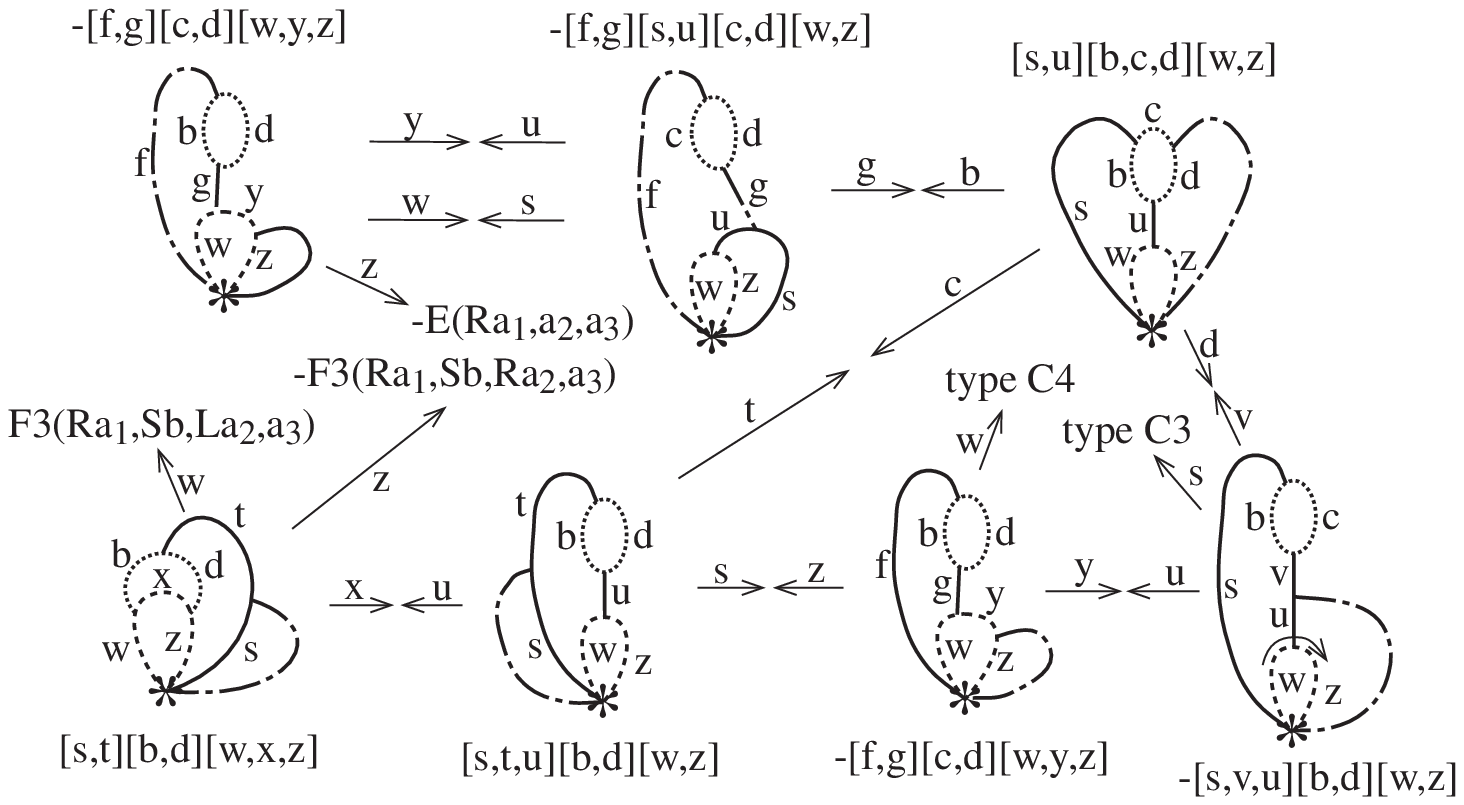}}}
\caption{}
\label{elimFEx}
\end{figure}

\begin{proof} 
Let $a\in\C{4}$ be the sum of the seven generators indicated in figure \ref{elimFEx}. As with figure
\ref{elimthetap}, pairs of arrows pointing at each other indicate canceling terms in $da$, provided
a consistent choice of orientations of edgesets is used.  For one such choice, the three specific generators
of $\C{3}$ on the left of the figure are terms of $da$.  Those terms of $da$ not accounted for in the figure
are all prismatic and none are $\theta$-prismatic.
\end{proof}

Let $g_1^R$, $g_2$ and $g_3$ be the generators of $\C{2}$ indicated in figure \ref{falls}.  Here 
$A=\{\alpha_i\}$ and $B=\{\beta_i\}$ are the sets consisting of  all the singleton half-edges incident with
the vertices  indicated in the figure.
Define $g_1^L$ just as $g_1^R$ is, except with the orientation of the edgeset $y$ changed.  
The notation for the $g_1^S$ and the $g_i$ should include some mention of the edgesets $y$ and $z$ and of
the sets $A$ and $B$, however we suppress this.
For lemma \ref{C5} below, we will count nontrivial commutator relations
using the homomorphism $f\colon \C{2}\to \oplus_{\overset{1\le x,y,z\le n}
{{\scriptscriptstyle S,T\in\{L,R\}}}}
\mathbb {Z}_{Sx,Ty,z}$ defined on generators by 
$$f(g)=
\begin{cases}
\sum_i1(\alpha_i,Sy,z),   &\text{if $g=g_1^S$,} \\
\sum_i(1(\alpha_i,Ry,z)-1(\alpha_i,Ly,z)),   &\text{if  $g=g_2$,}  \\
-\sum_i1(\alpha_i,Ry,z)-\sum_j1(\beta_i,Ly,z),   &\text{if  $g=g_3$,}  \\
0,  &\text{if $g\ne \pm g_1^S,g_2,g_3$.}
\end{cases}$$
Here each $\mathbb {Z}_{Sx,Ty,z}$ is a copy of $\mathbb {Z}$ with $1(Sx/Ty/z)\in \mathbb {Z}$ corresponding to 
$+1 \in \mathbb {Z}$.  

Let $F$ be the composite $f\circ d$.  Then $F$ is trivial on all simplicial generators,
and  on all prismatic generators except for $\theta$-prismatic ones.
Also, $F$ is trivial on all reduced cubical generators except for those of type
C5, C6 and C7.  So the only terms of any $z\in \cf_{\theta P}$ on which the homomorphism $F$ 
is nontrivial are type $C5$ and $C6$ cubical ones.  Moreover, for these we have
\begin{equation}
F(F_3(S_1a_1,S_2a_2,S_3a_3,a_4))=-1(S_1a_1,S_3a_3,a_4)+1(S_1a_1,S_2a_2,a_4),
\end{equation}
and
\begin{equation}\label{Fonex}
F(E(Sa_1,a_2,a_3))=1(Sa_1,Ra_2,a_3)-1(Sa_1,La_2,a_3).
\end{equation}
For $S\in \{L,R\}$, define $\bar {S}$  by $\{S,\bar{S}\}=\{L,R\}$.

\begin{lemma} \label{C5}
Each $z\in \cf_{\theta P}$ is homologous to a cycle in $\cf_{C5}$.
\end{lemma}

\begin{proof}
Write $z=\sum_{i=1}^n a_i$ where $n$ is minimal and each $a_i$ is a generator of $\C{3}$ corresponding
to some $3$-cell of $\xt_n/StN_n$.  Let $N_5(z)$ be the number of $a_i$ of type C5, and let
$N_6(a_i)$ be the number of type C6.  Also let $N_{5+6}(z)=N_5(z) +N_6(z)$.  We will show that if $N_{5+6}(z) >
0$, then $z$ is homologous to a cycle $z'\in \cf_{\theta P}$ with $N_{5+6}(z) < N_{5+6}(z')$.   

Since
$F(z)=0$, it follows from (\ref{Fonex}) that if $N_6(z)>0$, then $N_5(z)>0$ as well.  So we assume that
some $a_i$, say $a_1$, is of type C5 with $$F(a_1)=1(Sx,Ty,w)-1(Sx,Pz,w).$$ Again since $F(z)=0$, there is 
another $a_i$, say $a_2$, such that some cancelation occurs betweens the terms of $F(a_1)$ and
$F(a_2)$.  Two cases arise.
\vskip3pt
\noindent {\it Case 1.} Assume that $a_2$ is of type C5 with $F(a_2)=1(Sx,Qu,w)-1(Sx,Av,w)$.  
Then either $Qu=Pz$ or $Av=Ty$, say $Av=Ty$.  We may assume, by lemma \ref{falll1}, that $a_1=
F_3(Sx,Ty,Pz,w)$ and $a_2= -F_3(Sx,Ty,Qu,w)$.  If $u\ne z$, then, as in lemma \ref{0}, $z'=z+ 
d(F_4(Sx,Ty,Pz,Qu,w)+c)$ is in $\cf_{\theta P}$ and satisfies $N_{5+6}(z') = N_{5+6}(z) -1$.

If $u=z$ and $P\ne Q$, then lemma \ref{falll2} gives a boundary $b$ with $N_{5+6}(z+b) < N_{5+6}(z)$ 
and $z+b \in \cf_{\theta P}$.
\vskip3pt
\noindent {\it Case 2.} Assume that $a_2$ is of type C6 with either 
$$\text{$F(a_2)=1(Sx,\bar{T}y,w)-1(Sx,Ty,w)$ or $F(a_2)=1(Sx,Pz,w)-1(Sx,\bar{P}z,w)$}.$$
Then either $a_2=\pm E(Sx,y,w)$ or $a_2=\pm E(Sx,z,w)$.  In the first (resp.~second) case, we may assume,
by lemma \ref{falll1}, that $a_1=-F_3(Sx,Pz,Ty,w)$ (resp.~$F_3(Sx,Ty,Pz,w)$), and then lemma \ref{falll2}
gives a boundary $b$ with $z+b\in\cf_{\theta P}$ and $N_{5+6}(z+b) < N_{5+6}(z)$. 
\end{proof}

\begin{lemma} \label{C4}
Each $z\in \cf_{C5}$ is homologous to a cycle in $\cf_{C4}$.
\end{lemma}

\begin{figure}[tb]
\centerline{\mbox{\includegraphics*{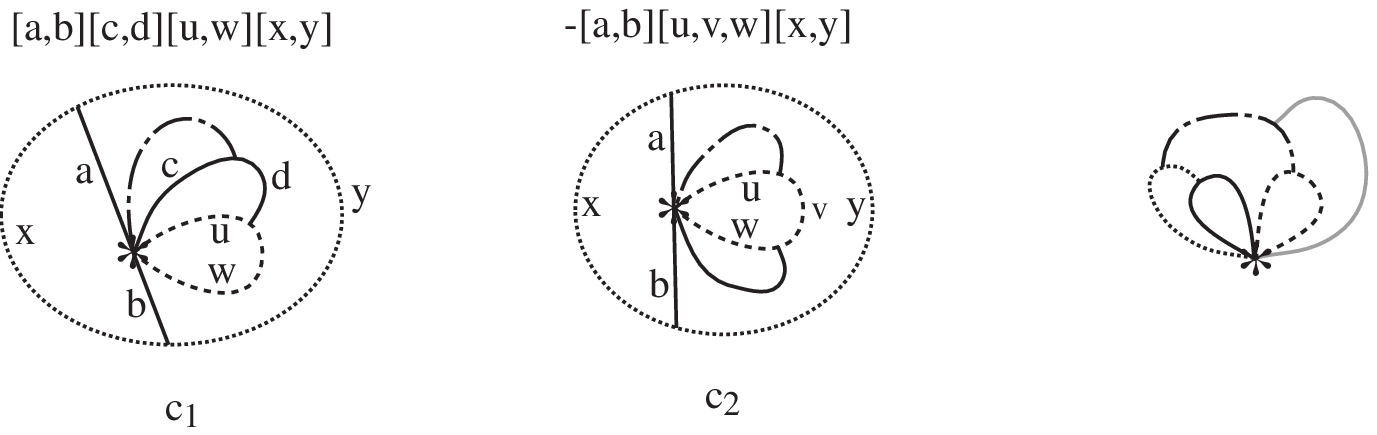}}}
\caption{}
\label{elimC4}
\end{figure}

\begin{proof}
Let $c_1$ and $c_2$ be the genertors indicated in figure \ref{elimC4}.  It's easy to check that $d(c_1+c_2)$
consists of two cubical terms of type C2(a), two prismatic terms (neither of which is $\theta$-prismatic),
and one cubical term of type C4.
\end{proof}

\begin{lemma} \label{C3}
Each $z\in \cf_{C4}$ is homologous to a cycle in $\cf_{C3}$.
\end{lemma}
\begin{proof}
Let $g\in \C{4}$ correspond to the colored graph shown on the right of  figure \ref{elimC4} with either cell
orientation. Then it's routine to check that there is a chain $c\in \C{4}$ consisting of 2,1,1-terms such that
$d(g+c)$ is in $\cf_{C4}$ and contains exactly one term of type C3.  Indeed, we can take $c=\sum_{i=1}^4 c_i$
where the colored graphs $cg(c_i)$, for $i=1$, 2 and 3, are defined by separating edges incident with the
4-valent vertices of the three colored graphs corresponding to the terms of $dg$ containing such vertices. 
Then just one  cubical term of $d(g+\sum_{i=1}^3 c_i)$ corresponds to a colored graph with a 4-valent vertex,
and $c_4$ is defined by separating edges incident with it.
\end{proof}

\begin{figure}[tb]
\centerline{\mbox{\includegraphics*{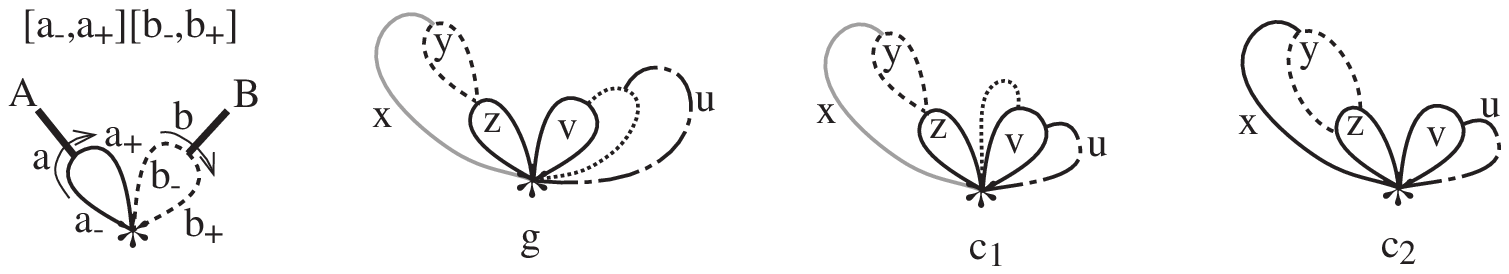}}}
\caption{}
\label{lastC}
\end{figure}

For the next lemma we will use a homomorphism 
$r_3\colon \C{2}\to \oplus_{\overset{1\le a,b,c,d\le n}
{{\scriptscriptstyle S,T\in\{L,R\}}}}
\mathbb {Z}_{\{Sa/b,Tc/d\}}$ which counts trivial commutator relations. Here each $\mathbb {Z}_{\{Sa/b,Tc/d\}}$
is a copy of $\mathbb {Z}$ with distinguished generator $1\{{Sa/b,Tc/d}\}= 1\{{Tc/d,Sa/b}\}$.
Let $g'\in \C{2}$ be the generator indicated on the left of figure \ref{lastC} where,
if $h_a=\{a_-,a_+\}$ and $h_b=\{b_-,b_+\}$, we are assuming that $a<b$, and where
 $A=\{\alpha_i\}$ and $B=\{\beta_j\}$ are the indicated sets of half-edges.  The homomorphism
$r_3$ is defined on generators by $$r_3(x)=
\begin{cases}
\sum_{i,j} 1\{\alpha_i /a,\beta_j /b\},  & \text{if $x=g'$} \\
0,        & \text{if $x\ne \pm g'$.}
\end{cases}
$$
Let $R_3$ be the composite $\cf_{C3} \hookrightarrow \C{3} \overset{d}{\to}\C{2} \overset{r_3}{\to}
\oplus\mathbb {Z}$.  Then $R_3$ is trivial on simplicial and prismatic terms as can be checked, and
is only nontrivial on generators of type C2.

	More preliminaries to the proof of the next lemma are needed.  View $\xt/StN$ as the direct limit
of the spaces $\xt_S/StN_S$ where $S$ varies over all finite subsets of $\mathbb {N}=\{1,2,\ldots\}$.
Given a generator $g$ of $\C{*}$, let $EC(g)$, the essential colors of $g$, be the subset of $\mathbb{N}$
of minimal cardinality such that 
there is an element $E(g)\in 
C_*(\xt_{EC(g)}/StN_{EC(g)})$ agreeing with $g$ in $\C{*}$.  
Given generators  $g_1$, $g_2 \in \C{*}$ with
$EC(g_1) \cap EC(g_2) = \varnothing$,  let $g_1 \vee g_2$ be the generator of $\C{*}$ 
corresponding to the 
one point union of $cg(E(g_1))$ and $cg(E(g_2))$ with the natural coloring and with cell orientation  that of $g_1$ followed by that of $g_2$.

Let $g=F_3(Sx,Aw,Ty,z) \vee F_1(Qu,v)$. Then there is 
a chain $c\in \C{4}$, obtained by separating edges as in the proof of lemma \ref{0},
 such that
\begin{equation} \label{first}
d(g+c) \equiv -(F_2(Sx,Ty,z) \vee F_1(Qu,v)) + F_2(Sx,Aw,z) \vee F_1(Qu,v)  \ \ \ \ (\text{mod}
\
\cf_{C2}).
\end{equation}
Now let $g$, $c_1$ and $c_2$ be the generators of $\C{4}$ indicated in figure \ref{lastC}.  Then, with  
appropriate cell orietations, we have 
\begin{equation}\label{second}
d(g+c_1+c_2) \equiv F_2(Sx,Ty,z) \vee F_1(Qu,v) - F_2(Sx,\bar {T}y,z) \vee F_1(Qu,v)  \ \ \ \
(\text{mod}
\ \cf_{C1}).
\end{equation}
Finally, let $g=F_2(Sx,Ty,z) \vee F_2(Qu,Aw,v)$. Then there are generators $c_1$ and 
$c_2$ in $\C{4}$ (obtained by separating edges of colored graphs corresponding to terms of $dg$) such that
\begin{equation} \label{third}
d(g+c_1+c_2)\equiv -F_1(Sx,z) \vee F_2(Qu,Aw,v) - F_2(Sx,Ty,z)\vee F_1(Qu,v)  \ \ \ \ (\text{mod $\cf_{C1}$).}
\end{equation}

\begin{lemma} \label{CC22}
Each $z\in \cf_{C3}$ is homologous to a cycle in $\cf_{C2}$.
\end{lemma}

\begin{proof} Let $z\in \cf_{C3}$.
We show that there is a 4-chain $x$ such that $z+dx$ is in $\cf_{C3}$ and contains fewer type C2 terms
than $z$.
Write $z=\sum_i z_i$.

\vskip3pt
\noindent {\it Case 1.}  Assume some $z_i$, say $z_1$, is of type C2(a) with $$z_1= \pm (F_2(Sx,Ty,z) \vee
F_1(Qu,v)).$$  Then $R_3(z_1) = \mp1\{Sx/z,Qu/v\}$.  So, since $R_3(z)=0$, there is some other $z_i$, say
$z_2$, with $R_3(z_2) = -R_3(z_1)$. Since $x\ne u$, either 
\begin{equation} \label{a} 
z_2= \mp (F_2(Sx,Aw,z) \vee F_1(Qu,v)) \end{equation} or 
\begin{equation} \label{b}
z_2 = \pm  (F_1(Sx,z) \vee  F_2(Qu,Aw,v))
\end{equation}  for some $A\in \{L,R\}$ and positive integer $w$.
If (\ref{a}) holds and $w\ne y$,  use the congruence (\ref{first}) to decrease the number of type C2 terms of
$z$ by two, while keeping $z \in \cf_{C3}$.
If $w=y$ and  $A\ne T$, use (\ref{second}) to do the same.
If (\ref{b}) holds, then by (\ref{first}),
we may assume that $w\ne y$. Then use (\ref{third}). 

\vskip3pt
\noindent {\it Case 2.}  If $z$ has a type C2(b) term, an argument like that of the previous case
can be used. The colored
graphs corresponding to the three generators denoted by $g$ in (\ref{first}), (\ref{second}) and (\ref{third})
 should be modified so that the two
singleton edges of each of these are replaced with a single singleton edge.
\end{proof}

The two cases in the previous proof could easily be combined, but then the wedge notation (which eliminates
the need for several figures) could not have been easily used.

\begin{lemma} \label{CC11}
Each $z\in \cf_{C2}$ is homologous to a cycle in $\cf_{C1}$.
\end{lemma}

\begin{proof}
If $x$ is a generator of $\C{3}$ of type C1, then $x =g \vee F_1(Sa,b)$ for some $g\in \C{2}$.  
Let $w$ be a color distinct from the five or six essential colors making up the set $E(x)$.  
Then there is a generator
$c$, obtained, as usual, by separating edges, such that $d(g \vee F_2(Sa,Sw,b) +c)=-x$.
\end{proof}
  
\section{${\mathbb R}^n$ in $X_k$} \label{RninXk}

	Let $k=2n-2$.  In this brief section, which will be used at the end of the next,
 we define a subgroup $\mathbb {Z}^k_*$ of
${\text {Aut}} (F_n)$ isomorphic to $\mathbb {Z} ^k$, the free abelian group of
rank
$k$, and also define a subcomplex $\mathbb  {R}^k_*$ of $X_n$.  The group
$\mathbb {Z}^k_*$  acts on
$\mathbb {R}^k_*$ in the usual way, and thus so does $\mathbb {Z}^k$, via an
isomorphism $\mathbb {Z}^k \approx \mathbb {Z}^k_*$ we specify. We also
give
$\mathbb {R}^k$ a standard triangulation and the usual 
$\mathbb {Z}^k$-action,
and then prove the following.

\begin{lemma} \label{textrninxn}
There is a $\mathbb {Z}^k$-equivariant simplicial homeomorphism
from $\mathbb {R}^k$ to $\mathbb  {R}^k_*$. Moreover, the classifying space
$\mathbb  {R}^k_* / \mathbb  {Z}^k_*$ 
for $\mathbb {Z}^k$ is a subcomplex of $X_n/\autfn = \widetilde X_n /
StN_n$.
\end{lemma}
 
Assume that $F_n$ is free on the set $\{x_1,\dots ,x_n\}$.
	The generators of $\mathbb  {Z}^k_*$ are, by definition, the $k$
automophisms $\alpha_i$ and $\beta_i$, for $1 \le i \le n-1$, defined by 
$\alpha_i(x_i)=x_nx_i$, $\beta_i(x_i)=x_ix_n$, and
$\alpha_i(x_j)=x_j=\beta_i(x_j)$ if  $i \ne j$. We identify $\mathbb  {Z}^k$
with $\mathbb  {Z}^k_*$ by letting $(m_1,\ldots ,m_k) \in \mathbb {Z}^k$
correspond to the automorphism $\alpha_1^{m_1} \beta_1^{-m_2} \ldots
\alpha_{n-1}^{m_{k-1}} \beta_{n-1}^{-m_{k}}$.

\begin{figure}[tb]
\centerline{\mbox{\includegraphics*{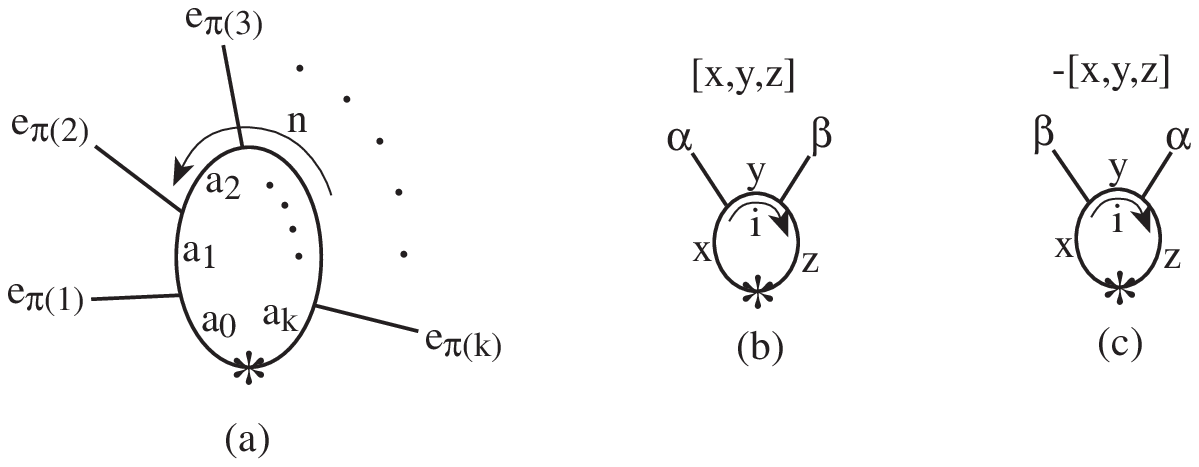}}}
\caption{}
\label{rninxn}
\end{figure}

A cell of $X_n$ is, by definition, in $\mathbb  {R}^k_*$ if it corresponds to
some blowdown of a colored graph of the sort shown in figure
\ref{rninxn}(a)  
(where
$\pi$ is any permutation of the set $\{1, \ldots,k\}$) with marking such 
that any (and hence every) $0$-cell 
in  the boundary of the cell
corresponds to an automorphism of $\mathbb 
{Z}^k_*$.  Thus $\mathbb  {R}^k_*$ consists of those cells whose
$0$-cells all correspond to elements of $\mathbb {Z}^k_*$.  

	The unit cube $[0,1]^k \subset \mathbb  {R}^k$ is triangulated by the $k!$
k-simplices $\sigma_{\pi(1)} \times \ldots \times \sigma_{\pi(k)}$ where
$\pi$ ranges over all permutations of the set $\{1, \ldots,k\}$ and
$\sigma_i=(0,\ldots,0,1,\ldots,1)$ is the (degenerate) k-simplex with $i$
zeros and $k+1-i$ ones.  (Here we are using the notation of, for instance, p.\ 244 of
Mac lane \cite{M}.)  Translation of this triangulation by the elements of
$\mathbb {Z}^k$ defines a $\mathbb {Z}^k$-invariant triangulation of all of
$\mathbb {R}^k$.

{\begin{proof}[proof of lemma \ref{textrninxn}]  Let $[a_0, \ldots ,a_k]_{\pi
,id}$ be the k-simplex of
$X_n$ corresponding to the colored graph shown in figure \ref{rninxn}(a) with
marking specified by $[a_0] = id \in  \autfn $.  Mapping, for each permutation
$\pi$, the k-simplex  $\sigma_{\pi(1)} \times \ldots \times \sigma_{\pi(k)}$ of
the cube
$[0,1]^k$ to $[a_0, \ldots ,a_k]_{\pi ,id}$ so that the ordering of the
vertices is preserved, and then extending equivariantly defines the desired
homeomorphism.
	
	Now assume that $g(a)=g(b)$ where $g\colon \mathbb {R}^k_*/\mathbb {Z}^k_*
\to
X_n/\autfn$ is the map induced by the inclusion $\mathbb  {R}^k_*\to
X_n$.  Then there are points $a',b'\in \mathbb  {R}^k_* \subset X_n$
representing $a$ and $b$ and some $f\in \autfn $ with $f\,a'=b'$.  Since $f$
acts cellularly, there are vertices $v$ and $w$ of $\mathbb  {R}^k_*$ with
$fv=w$.  Since $v$ and $w$ correspond to elements of $\mathbb {Z}^k_*$, so
does $f$.  Thus $g$ is injective.
\end{proof}

\begin{remark} \rm
\label{r}
While the number $n$ plays a special role in the definitions of
$\mathbb  {Z}^k_*$ and $\mathbb  {R}^k_*$, these could just as well be
defined using any of $1, \ldots , n-1$ in place of $n$, and the lemma would
still be true.  
\end{remark}	
 \subsection{} \label{basis}
 Let $\mathbb  {R}^k_{i*}$ 
denote the subcomplex of $\widetilde X_n / StN_n$ analogous to 
$\mathbb  {R}^k_* =\mathbb  {R}^k_{n*}$, but defined with $i$ in place of $n$.
We conclude this section with an explicit description of a basis for 
$H_2(\mathbb {R}^k_{i*}/\mathbb {Z}^k)$ which will be used in 
Lemmas \ref{P2} and \ref{P0}.
Order the half-edges of any colored graph by $Si<Tj$ if either $i<j$ or,
if $i=j$, $S=L$ and $T=R$.  Let $S_i(\alpha,\beta)$ (resp. $S_i(\beta,
  \alpha)$) be the generator of $C_2(\mathbb {R}^k_{i*}/\mathbb {Z}^k)$
indicated in figure \ref{rninxn}(b) (resp. figure \ref{rninxn}(c)) where we are
assuming that $\alpha < \beta$.  Denote by $1_i\{\alpha, \beta \}=1_i\{\beta
,\alpha \}$ the element of $H_2(\mathbb {R}^k_{i*}/\mathbb {Z}^k)$ which the 
cycle $S_i(\alpha,\beta) + S_i(\beta, \alpha)$ represents.  It follows from
the computation of $H_2(\mathbb {Z}^k)$ and lemma \ref{textrninxn} that the 
$1_i\{\alpha, \beta\}$ for all distinct half-edges $\alpha$ and $\beta$
different from $Li$ and $Ri$ form a basis for 
$H_2(\mathbb {R}^k_{i*}/\mathbb{Z}^k)$.

\section{Eliminating prismatic and simplicial terms}
We say that a colored graph $G$ has type $Pn$, for $n\ne 2$, if the core of $G$
is isomorphic to the core of the partitioned graph labeled
$Pn$ in figure \ref{prismcl}.  Type P2(a), (b) and (c) colored graphs are defined similarly, 
and a colored graph is of type P2 if it is of type P2(a), (b) or (c).

In order to eliminate all prismatic terms from $z$, we extend the
filtration $\cf$ of the group $Z_3(\xt_n/StN_n)$ of cycles to 
$$\cf_{C1}\supset \cf_{P5}\supset \cf_{P4}\supset\ldots \supset \cf_{P1}\supset \cf_{S}$$ 
where $z\in
\cf_{Pn}$ if none of the terms of $z$ correspond 
to colored graphs of type $Pm$, for $m\ge n$. 
(An additional technical requirement is placed on the elements of $\cf_{P4}$
in Def.\ \ref{reduced} below.)
\begin{figure}[tb]
\epsfbox{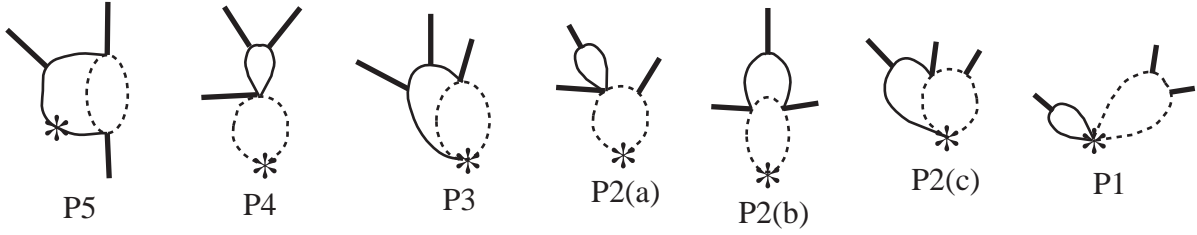}
\caption{Near classification of prismatic colored graphs.}
\label{prismcl}
\end{figure}
Also, $z\in \cf_S$ if none of the terms of $z$ are simplicial. 
Lemmas \ref{P5} to \ref{SS} below
  show that any cycle in a subgroup of the above filtration is
homologous in the space $\xt_n/ StN_n$ to a cycle in the next smaller subgroup of the filtration.   
Since $\cf_{C1} \subset \cf_{\theta P}$, it follows from the next lemma
 that the smallest subgroup $\cf_S$ of the filtration 
is the trivial group.
 Thus, every element of $\cf_{C1}$ is a boundary.

\begin{lemma} \label{prismclassification}
If $G$ is a prismatic colored graph, then either $G$ is of one of the types shown in figure \ref{prismcl}, 
or $G$ is $\theta$-prismatic.
\end{lemma}

 $\theta$-prismatic graphs were defined just before Lemma \ref{thetaP} above.

\begin{proof}
Let $G$ be any prismatic colored graph. Since $G$ has four vertices, the three edges of $G$ of the same
color form a connected subgraph of $G$. If they form a monochromatic loop, then, disregarding basepoints,
the core of $G$ is as shown in either P2(a) or P2(b) of figure \ref{prismcl}.  Otherwise the core is as
shown in P5 of figure \ref{prismcl}.  Each of the possible locations for the basepoint on these three cores
gives either one of the graphs in figure \ref{prismcl} or a $\theta$-prismatic graph.
\end{proof}

\begin{definition}\rm
	A colored graph $G$ with three vertices is of type $\theta$ 
(resp. x, y, w) if the core of $G$
is equivalent, for some choice of colors,  to the core of the 
graph labeled type~$\theta$ (resp.\ x, y, w) in figure \ref{typexy}.
The incidence number of a colored
graph of any of the three types $\theta$, $x$, or $y$ 
\begin{figure}[tb]
\centerline{\mbox{\includegraphics*{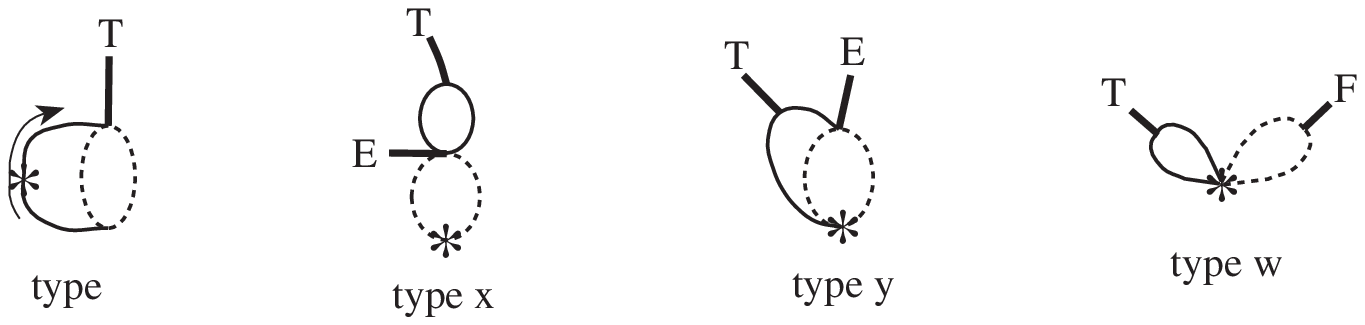}}}
\caption{}
\label{typexy}
\end{figure}
is the
cardinality of the set of half-edges corresponding to the set labeled $T$ in
the appropriate graph in the figure.  
The incidence number of a type $w$ colored graph is $|T| + |F|$, 
where $T$ and $F$ are the sets analogous to those in the figure.
The incidence number of a generator $g$ of $C_2(\widetilde X /StN)$ with
$cg(g)$ of type $x$, $y$, $w$ or $\theta$ is the incidence number of
$cg(g)$.

 \medskip
Let $G$ be a colored
graph of either type $P5$ or $\theta$.  The initial vertex of the first edge
of the (oriented) edgeset through the basepoint is called the initial vertex
of $G$.  We say that $G$ is reduced if its initial vertex is trivalent.  
\end{definition}

\begin{lemma} 
\label{P5}
Each $z\in \cf_{C1}$ is homologous in $\widetilde X _n/StN_n$ to some element of $\cf_{P5}$.
\end{lemma}
\begin{proof}

By
adding to $z$ generators of $C_4(\widetilde X_n/StN_n)$ corresponding to
(suitably oriented) colored graphs of the sort shown in figure
\ref{longth}(a) and (b),
\begin{figure}[tb]
\centerline{\mbox{\includegraphics*{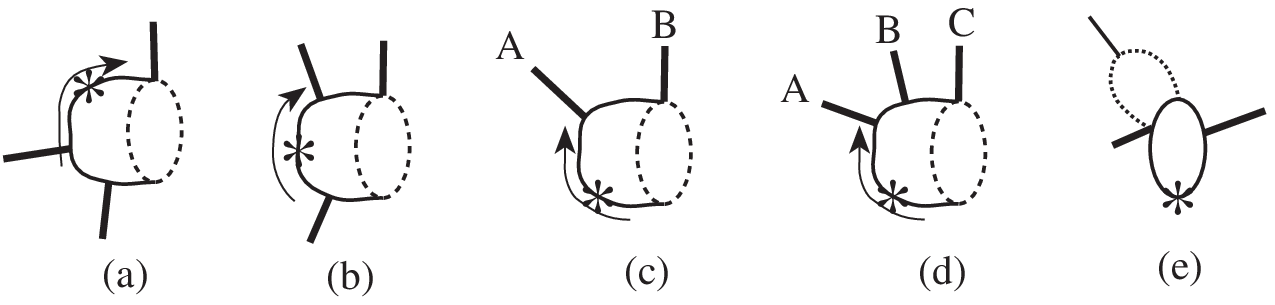}}}
\caption{}
\label{longth}
\end{figure}
we may assume that, for each type P5 term $z_i$ of
$z$,  $cg(z_i)$ is reduced.

Given a type $\theta$ colored graph $G$ which is not reduced, there is
exactly one reduced colored graph of type $P5$ having a blowdown equal to
$G$.  Thus, if $x$ is a term of $dz$ with $cg(x)$ of type $\theta$,
 then $cg(x)$ is reduced.  It follows that all the
type $P5$ terms of $z$ correspond  to colored graphs of the sort shown in
figure \ref{longth}(c) where $B$, but not $A$, may be empty. Such colored
graphs are said to have type $P5^+$. We denote by
$\theta_{pq}(A,B)$ the colored graph in figure \ref{longth}(c) and by
$\theta_{pq}(A,B,C)$ the one in figure \ref{longth}(d). Here
$q$ is the color of the edgesets shown as dotted and $p$ is the color of the edges in the other
edgeset with more than one edge.


Assuming $z\notin\cf_{P5}$,
let $m(z)$ be the maximal incidence number of all
the type $\theta$ faces of the terms of $z$. Choose two terms $z_1$ and $z_2$ of
$z$ such that $dz_1$  contains a term which realizes $m(z)$ and cancels
with a term of $dz_2$.  Then both
$z_1$ and $z_2$ correspond to type $P5^+$ colored graphs, and  if
$cg(z_i)=\theta_{pq}(A_i,B_i)$ for $i=1$ and $2$, then the set $A_1\cup B_1=A_2\cup 
B_2$ has cardinality $m(z)$. 

	If one of the  $B_i = \varnothing$, say $B_1$, then the boundary of 
the generator of $C_4(\xt_n/StN_n)$ corresponding to the suitably
oriented colored graph $\theta_{pq}(A_2,B_2,\varnothing)$ is $z_1 + z_2$ plus
two more terms. One of these is of type $P3$.  The $\theta$-faces of the
other, which is of type $P5^+$, all have incidence number $\le |B_2| < m(z)$.
Thus, if a $B_i=\varnothing$, the total number of terms of $z$ with a face
realizing $m(z)$ can be decreased.

	Now we assume that neither $B_i=\varnothing$ and show that this is still
true.  We first assume that $m(z) > 2$.
 Then if, say $|A_1| > 1$, any proper
subset $a$ of $A_1$ can effectively be deleted from $A_1$ and added to $B_1$
by adding to $z$ the boundary of the generator of $C_4(\xt /StN)$
corresponding to the (suitably oriented) colored graph $G$ shown in
figure \ref{ab}(a).
\begin{figure}[tb]
\centerline{\mbox{\includegraphics*{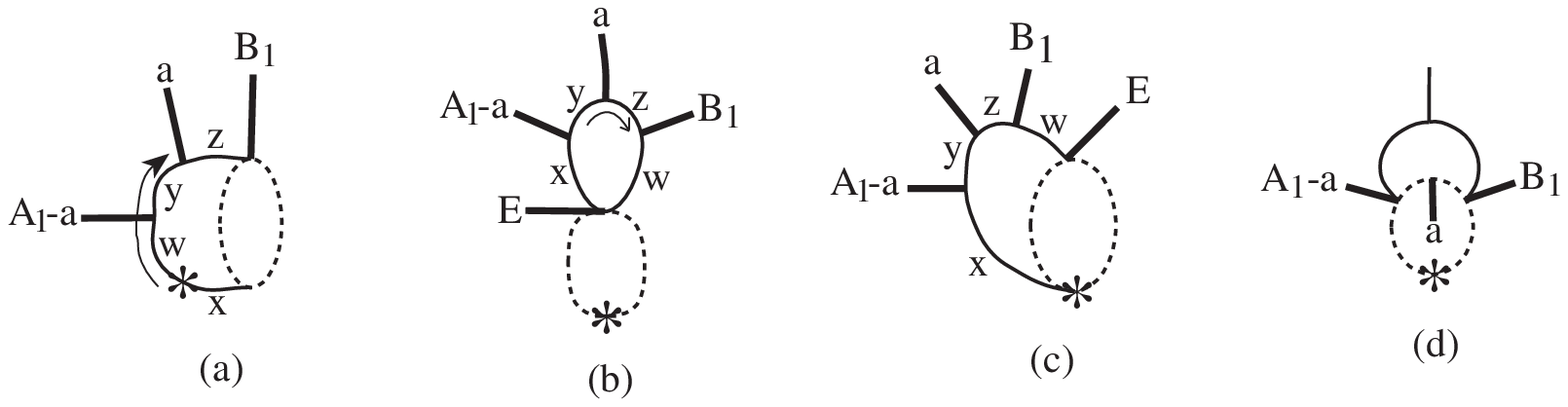}}}
\caption{}
\label{ab}
\end{figure}
This may also add to $z$ a type $P3$ term corresponding to $G_x$ and a type
$P5^+$ one corresponding to $G_w$.  But, again, the two
$\theta$-faces of this last term have incidence number $\le |a|+|B_1| < m(z)$.
  Similarly any proper subset of $B_1$ can be
moved from $B_1$ to $A_1$ without increasing $m(z)$.  By using either
one, two or three of these moves it's possible to eliminate $z_1$ and $z_2$
from $z$.  This reduces the total number of
terms of $z$ which have a type $\theta$ face realizing $m(z)$. Thus, if $m(z)
>2$, the number of such terms can be reduced to zero, thereby reducing $m(z)$.

	So now we assume that $m(z)=2$ and, by the above, that neither $B_i =
\varnothing$.  Then $A_1=B_2$ and $B_1=A_2$, and, with suitable orientations, 
$$z_1+z_2 + d\bigl( gen(\theta_{pq}(A_1,B_1,\varnothing)) -
gen(\theta_{pq}(B_1,A_1,\varnothing))\bigr)$$ has just two terms.  One is of
type $P3$. The $\theta$-faces of the other, which is of type $P5^+$, have
incidence number $\le 1$. So $m(z)$ can also be reduced if it equals $2$. 
This completes the proof since $m(z)$ cannot equal $1$.
\end{proof}

\begin{definition} \label{reduced}\rm
A type $P2$ colored graph $G$ is reduced if all its type $x$ and $y$ faces have incidence
number one.  All type P2 terms of elements of $\cf_{P4}$ are required to be reduced.
\end{definition}

\begin{lemma} Each $z\in \cf_{P5}$ is homologous to a cycle in $\cf_{P4}$. \label{P4}
\end{lemma}
\begin{proof}

	 If a type $P2$ colored graph $G$ is drawn as in P2(a), (b) or (c) of figure
\ref{prismcl}, then $G$ is reduced if the set of half-edges which are not incident with any dotted
edge has cardinality one. We may assme that all the $P2$ terms of $z$ are reduced by repeatedly
separating the half-edges of any of these sets which have cardinality $> 1$.

Assume $z\notin\cf_{P4}$. The argument that $z$ is homologous to a cycle in $\cf_{P4}$ is similar to
parts of the proof of lemma \ref{P5}. Let
$m(z)$ now be the maximal incidence number of all the type $x$ faces of the terms of $z$.  Let $z_1$ and
$z_2$ be two terms of $z$ having a common canceling face which realizes $m(z)$.  Since the $P2$ terms of
$z$ are reduced, both of the colored graphs
$cg(z_1)$ and $cg(z_2)$ are of type $P4$.  For $i=1$ and $2$, we denote by
$A_i$ the set of half-edges of $cg(z_i)$ corresponding to the set labeled $A$ of the graph
in the lower right of
figure \ref{4chain}, and by $B_i$ the one corresponding to $B$.  
Then $A_1 \cup B_1 = A_2\cup B_2$, and this set has cardinality $m(z)$.  If $m(z) > 2$ then the argument
in the proof of \ref{P5} showing $m(z)$ can be reduced is now applicable with $G$ the colored graph in figure
\ref{ab}(b) rather than \ref{ab}(a).  Note that all the type $x$ and $y$ faces of $G_x$ and $G_w$ have
incidence number $< m(z)$.

	So we now assume that $m(z) = 2$.  Then, using the above notation, we have
$A_1=B_2$ and $B_1=A_2$, and we assume that $A_1=\{\alpha\}$ and
$B_1=\{\beta\}$.  Figure \ref{4chain} indicates a $4$-chain $c$ with four terms
such that $z_1 + z_2 \pm dc \in \cf_{P1}$ for one of the choices of the signs.
\begin{figure}[tb]
\centerline{\mbox{\includegraphics*{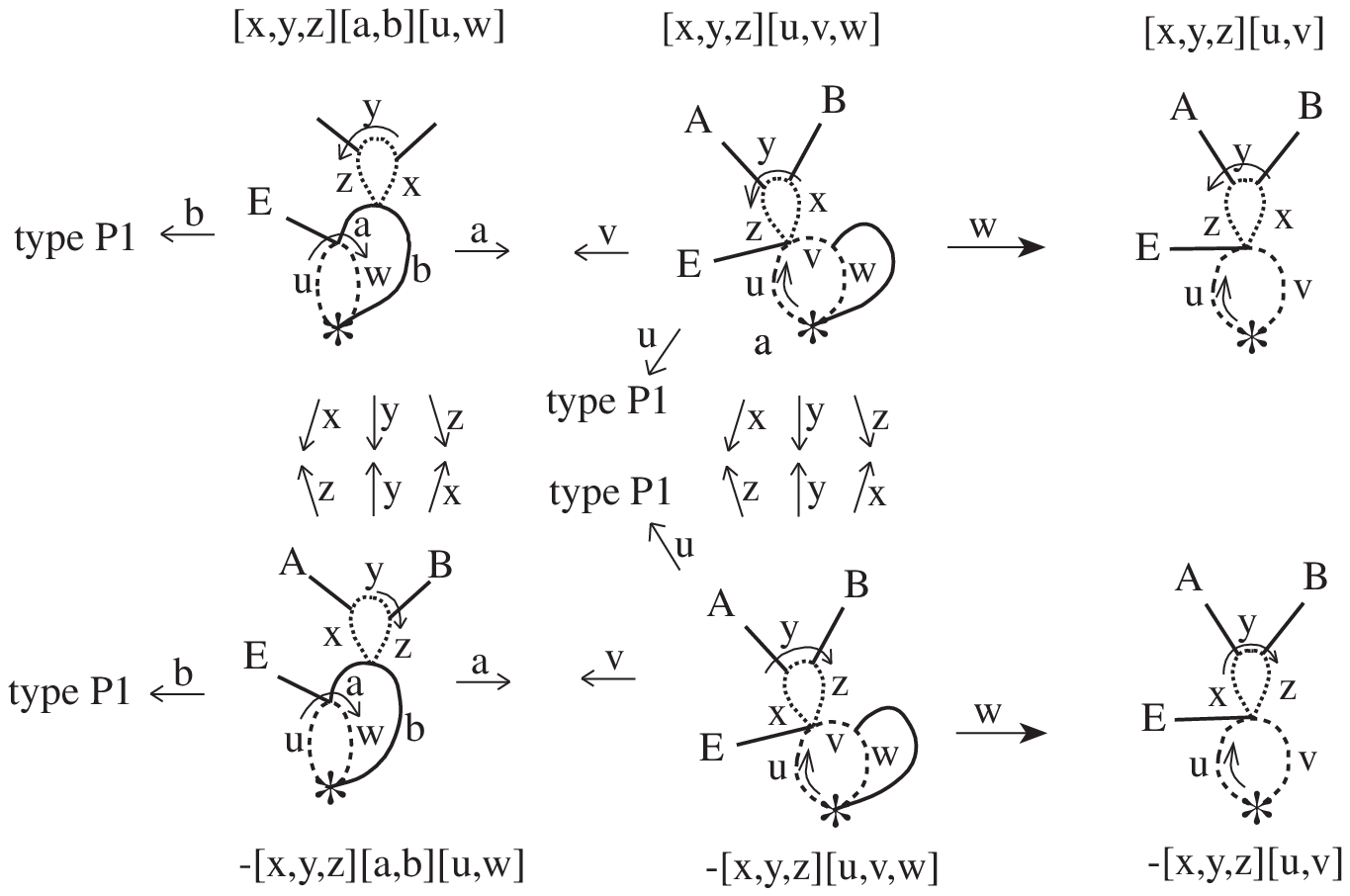}}}
\caption{}
\label{4chain}
\end{figure}
(Each graph in the figure is identical to the one above or below it except for the
orientation of the edgeset $\{x,y,z\}$.  The two colored graphs on the right
with the indicated orientations correspond to $\pm(z_1+z_2)$.) This completes the proof.
\end{proof} 

\begin{lemma} Each $z\in \cf_{P4}$ is homologous to a cycle in $\cf_{P3}$.
\label{P3}
\end{lemma}

\begin{proof}
	Assuming $z\notin \cf_{P3}$, let $m(z)$ now denote the maximal incidence
number of all the type y faces of the terms of $z$.  The argument of the previous lemma with
some changes shows that we may assume that $m(z) = 2$.  These changes include
using type y colored graphs in place of type x ones, using the graph in 
figure \ref{ab}(c) (with the edgeset $\{x,y,z,w\}$ suitably oriented)
instead of the one in figure \ref{ab}(b), and substituting type $P3$ for
$P4$.

	We will use a homomorphism $f\colon \C{2} \to \oplus_i \mathbb Z_i$ where there is one 
$\mathbb Z_i\approx \mathbb Z$
summand in the range corresponding to each colored graph shown in figure \ref{homf}(a). 
\begin{figure}[tb]
\centerline{\mbox{\includegraphics*{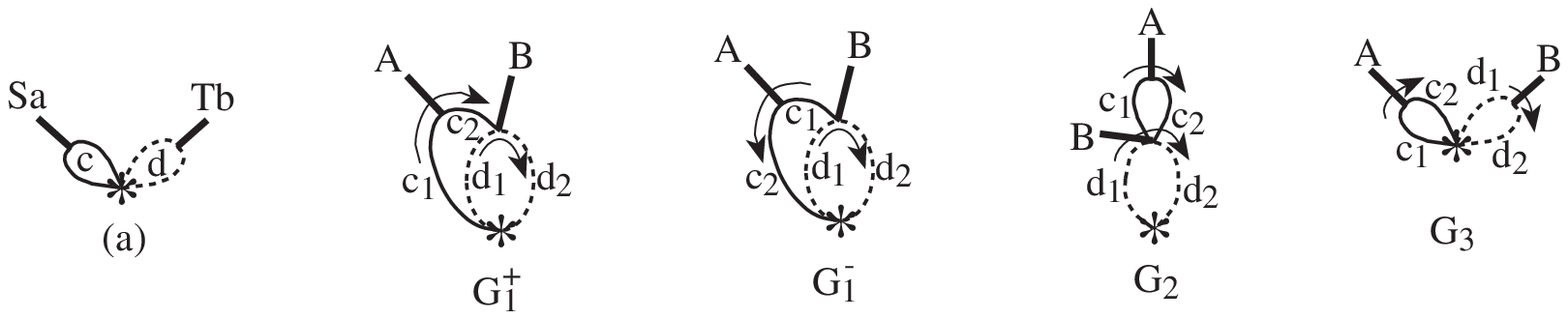}}}
\caption{}
\label{homf}
\end{figure}  
 Here $Sa$ and
$Tb$ are half-edges so that $S,T\in\{L,R\}$.  Also $a$, $c$ and $d$ are distinct colors, as are $b$, $c$
and
$d$. We let $1\{Sa/c,Tb/d\} =1\{Tb/d,Sa/c\}$ denote a generator of the $\mathbb Z$ summand
corresponding to the colored graph in figure
\ref{homf}(a).  Let $G_1^{\pm}$, $G_2$ and $G_3$ be the colored graphs indicated
in figure
\ref{homf} where we assume $A=\{\alpha_i\}$ and $B=\{\beta_i\}$.  Let $g_1^{\pm}=gen(G_1^{\pm})$ and 
$g_i=gen(G_i)$ be the corresponding generators of $\C{2}$ with the orientations of $G_1^{\pm}$ and the 
$G_i$ given by $[c_1,c_2][d_1,d_2]$ if $c<d$, and opposite this if $c>d$. Let $G^{\pm}$ be the two
colored graphs indicated in figure \ref{34cells} 
\begin{figure}[tb]
\centerline{\mbox{\includegraphics*{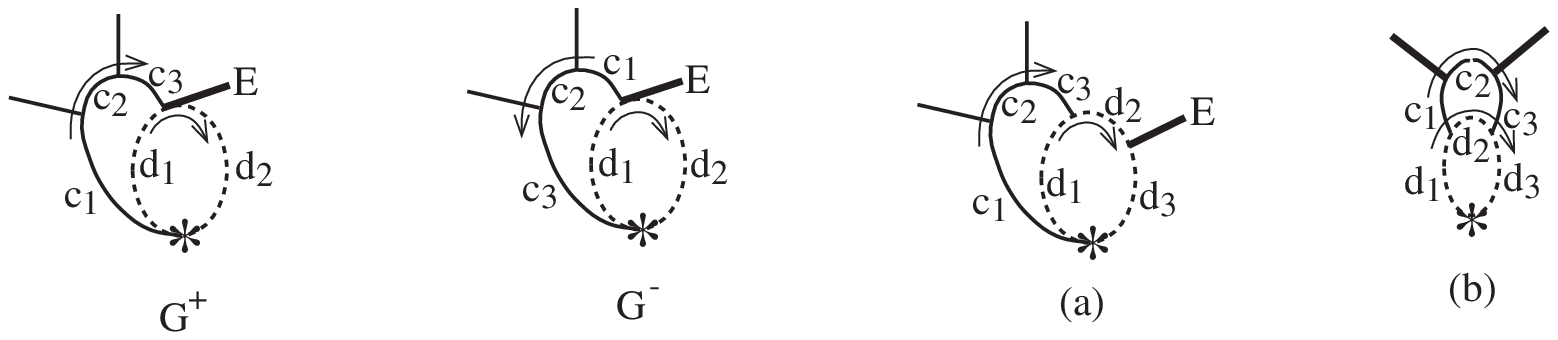}}}
\caption{}
\label{34cells}
\end{figure}  
where $\alpha$ and $\beta$ are single half-edges.
Let $g_{cd}^{\pm}(\alpha,\beta,E)=gen(G^{\pm})$ be the generators of $\C{3}$ corresponding
to $G^+$ and $G^-$, both with the orientation $[c_1,c_2,c_3][d_1,d_2]$.

The homomorphism $f$ is defined by 
$$ f(g_1^{\pm})=f(g_2)=f(g_3)= \left\{ \begin{array}{ll}
                                  \sum_{i,j} 1 \{ \alpha_i/c,\beta_j/d \}  &  \mbox{if
$B \ne \varnothing$}, \\
                                  0                                      & \mbox{if $B=\varnothing$},
                              \end{array} \right.  $$
and $f(g)=0$ for all other generators $g$ of $\C{2}$.
It's routine to check that $f(d(g_{cd}^{\pm}(\alpha,\beta,E)))=1\{\alpha/c,\beta/d\}$ where $d\colon
\C{3}\to \C{2}$ is the boundary homomorphism.  Also, if a generator $g$ of $\C{3}$ is simplicial or of
type
$P1$ or
$P2$, then
$f(d(g))=0$.

As usual let $z_1$ and $z_2$ be terms of $z$ with a common canceling face realizing $m(z)=2$.
The colored graph corresponding to this face is either $G_1^+$ or $G_1^-$, say it's $G_1^+$ (with a similar
argument if it's $G_1^-$).  Then, since the $P2$ terms of $z$ are reduced,
$$z_1 +z_2=g_{cd}^+(\alpha,\beta,E)-g_{cd}^+(\beta,\alpha,E)$$ for some $c,d\in {\bf n}$, half-edges
$\alpha$ and $\beta$, and set of half-edges $E$.  Let $g_4^+(\alpha ,\beta ,E)$ be the generator of 
$\C{4}$ corresponding to the colored graph shown in
figure \ref{34cells}(a) with  orientation $[c_1,c_2,c_3][d_1,d_2,d_3]$.  Then
$$z_1 +z_2 - \left(g_{cd}^+(\alpha,\beta,\varnothing)-g_{cd}^+(\beta,\alpha,\varnothing)\right)
+ d\left(g_4^+(\alpha,\beta,E)-g_4^+(\beta,\alpha,E)\right) \in \cf_{P2}$$  with  all the $P2$ terms
reduced.  So we assume that $E=\varnothing$.

  Since $f(d(z_1+z_2))=1\{\alpha/c,\beta /d\} - 
1\{\beta/c,\alpha /d\}\ne 0$ and $z$ is a cycle, $z$ must contain a term, $z_3$ say, distinct from
$z_1$ and $z_2$, such that $f(d(z_3))$ cancels with a term of $f(d(z_1 + z_2))$.  Let $z_4$ be a term of 
$z$ with a face canceling the one of $z_3$ which realizes $m(z)=2$.
  Then, for some set of half-edges $E'$, either
\begin{equation} \label{case1} z_3 + z_4 = -g_{cd}^+(\alpha,\beta,E')+g_{cd}^+(\beta,\alpha,E') 
\end{equation} or
\begin{equation}\label{case2}  z_3 + z_4 = -g_{cd}^-(\alpha,\beta,E')+g_{cd}^-(\beta,\alpha,E'). 
\end{equation}
In both cases we may assume, as above, that $E'=\varnothing$.  Then if (\ref{case1}) holds, $\sum_{i=1}^4 z_i
=0$.  To treat the other case, let $s_4(\alpha ,\beta )$ be the generator of $\C{4}$ corresponding to
the colored graph shown in figure \ref{34cells}(b) with  orientation $[c_1,c_2,c_3][d_1,d_2,d_3]$.  
  Then, if (\ref{case2}) holds, we have 
$$d\big(s_4(\alpha ,\beta ) - s_4(\beta,\alpha ) \pm c\big) - \sum_{i=1}^4 z_i \in \cf_{P2}$$
for a 4-chain $c$ like that indicated in figure \ref{4chain}.  Moreover, all the $P2$ terms are reduced.
 So in both cases the number of terms of $z$ realizing $m(z)=2$ can be reduced.  
\end{proof}

\begin{definition} \label{lowerincidencenum}\rm
The lower incidence number of a type x (resp.~y) colored graph
with set $E$ as shown in figure \ref{typexy} is $|E| + 2$ (resp.~$|E|+1$).
The lower incidence number of a generator $g$ of $C_2(\widetilde X /StN)$ with
$cg(g)$ of type $x$ or $y$ is the lower incidence number of
$cg(g)$.
\end{definition}

\begin{lemma} 
\label{P2}
Each $z\in \cf_{P3}$ is homologous to a cycle in $\cf_{P2}$.
\end{lemma}

\begin{proof} 
Terms of $z$ of type $P2$(a) can be exchanged for simplicial and reduced type $P2$(b) and (c) ones 
by using colored graphs as shown in figure \ref{longth}(e).
Thus, we assume that none of the terms of $z$
are of type $P2$(a).  

If $z$ contains any terms of type $P2$(b), let $m(z)$ be the maximal lower incidence
number of all the type $x$ faces of the terms of $z$.  We may assume that $m(z)=2$ by using the usual
argument, this time with colored graphs as, for instance, shown in figures \ref{ab}(d) and 
\ref{longth}(e).  Let $z_1$
and $z_2$ be two type $P2$(b) terms of $z$ with faces realizing $m(z)=2$. We may assume that the colored
graphs $cg(z_1)$ and
$cg(z_2)$ are identical except for the orientations of the edgesets with two edges. Then, for one choice
of the plus or minus sign, $$dy \pm (z_1 + z_2)$$ consists of just reduced type $P2(c)$ terms where $y$ 
is the sum
of the seven generators of $\C{4}$ indicated in figure \ref{4chain2} 
\begin{figure}[tb]
\centerline{\mbox{\includegraphics*{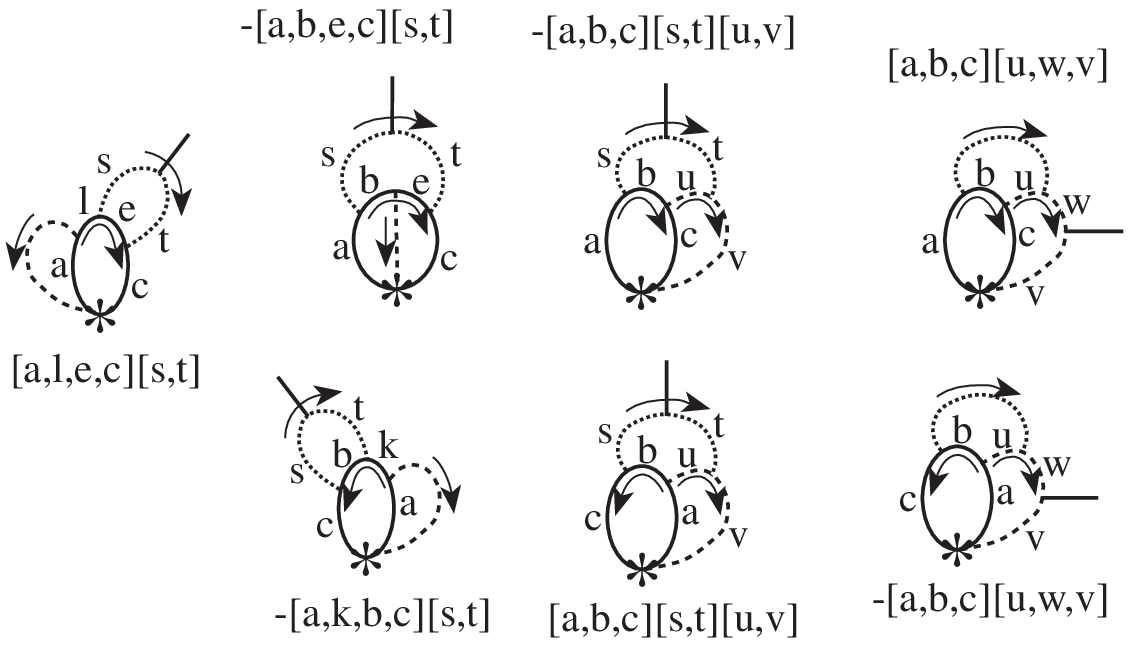}}}
\caption{}
\label{4chain2}
\end{figure}  
with  suitable choices of colors.
So we now assume that $z$ contains no terms of either type $P2$(a) or (b). 

The following notation will be used in the remainder of the proof.
Let $F_i^+(Ra,\alpha,\beta)$ (resp.$\ -F_i^-(Ra,\alpha,\beta)$) be the generator of $\C{4}$ 
corresponding to
the colored graph and  orientation in the upper left (resp.~lower left)
of figure \ref{3chain}.
Here we are assuming that $h_a=\{x,y\}$ and $h_i=\{s,t,u\}$.
Define $F_i^+(La,\alpha,\beta)$ (resp.$\ F_i^-(La,\alpha,\beta)$) as $F_i^+(Ra,\alpha,\beta)$ 
(resp.$\ F_i^-(Ra,\alpha,\beta)$) is defined, but switch the orientation of the edgeset $h_a$.

If $z$ contains any terms of 
type $P2$(c), then, as usual, we may assume that $m(z)=2$ where now $m(z)$ is 
the maximal lower incidence
number of all the type $y$ faces of the terms of $z$.  Write 
$$z= \rho + \sum_{i=1}^n \zeta _i $$
where $\rho$ contains all the type $P1$ terms of $z$, and $\zeta_i$ contains all the type $P2$(c) and 
simplicial terms corresponding to colored graphs whose edgesets with cardinality $> 2$ have color $i$.  
Choose a $q\in {\bf n}$ such that $\zeta_q$ contains a type $P2$(c) term. Write
$$ \zeta_q = s + \sum_j (z_j + z_j')$$
where $s$ contains all the simplicial terms of $\zeta_q$ and where, for each $j$, $z_j$ and $z_j'$ 
have a common canceling face realizing $m(z)$.
Since $m(z)=2$, we have, for each $j$,
\begin{equation}\label{star} z_j+ z_j' = \pm(F_q^+(S_ja_j,\alpha_j,\beta_j) - F_q^-(S_ja_j,\alpha_j,\beta_j))
\end{equation}
for some $S_j\in \{L,R\}$, $a_j\in {\bf n}$, and half-edges $\alpha_j$ and $\beta_j$, and where the
choice of the plus or minus sign depends on $j$.  The left part of figure \ref{3chain} 
\begin{figure}[tb]
\centerline{\mbox{\includegraphics*{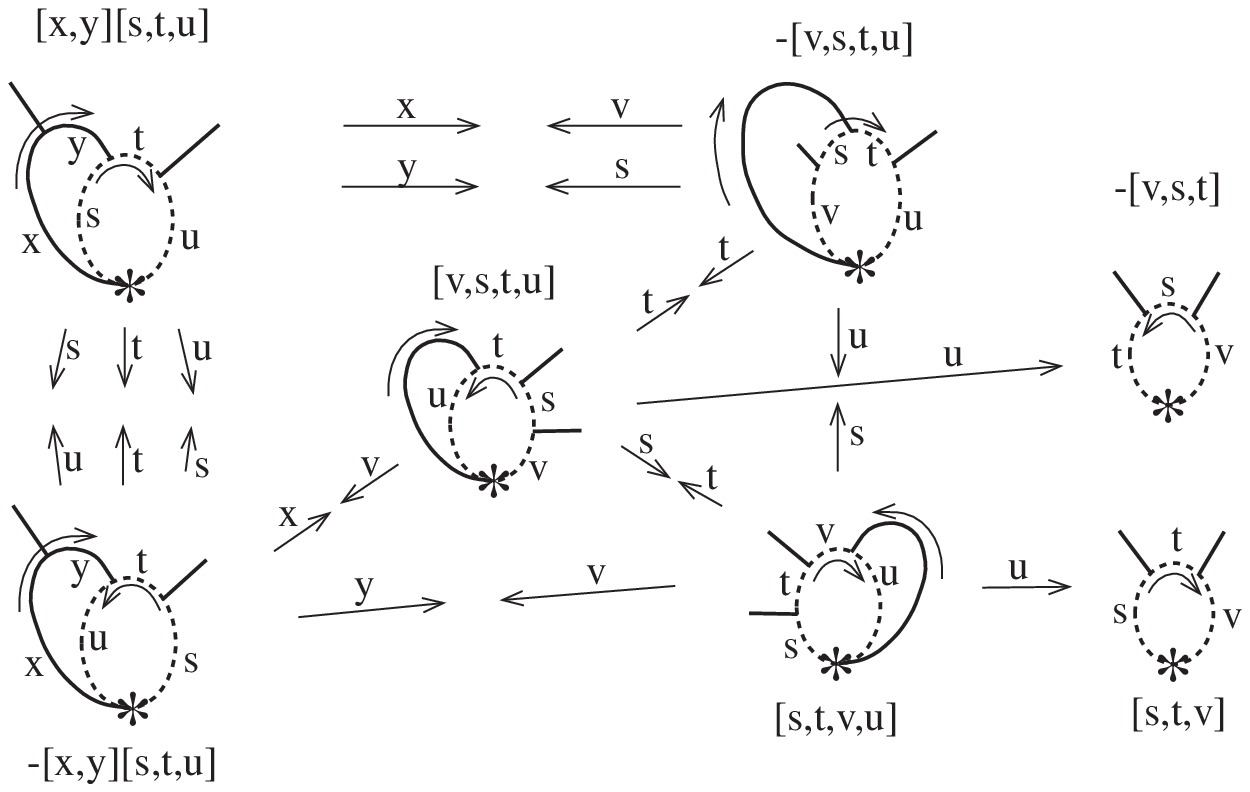}}}
\caption{}
\label{3chain}
\end{figure} 
indicates that each $d(z_j+z_j')$
is in the image of the homomorphism $C_2(\mathbb {R}^k_{q*} / \mathbb  {Z}^k)
\overset{i_*}{\to}\C{2}$. Here $\mathbb {R}^k_{q*} / \mathbb  {Z}^k$ is the subcomplex of 
$\widetilde {X}_n/StN_n$ defined in \S\ref{basis} and $i$ is the inclusion. So $d\zeta_q$ is also 
in the image of $i_*$.  By a similar argument, each $d\zeta_i$ is in the image of the map
$C_2(\mathbb {R}^k_{i*} / \mathbb  {Z}^k) \to \C{2}$ induced by inclusion.
Thus, since all the terms of $d\rho$ correspond to square 2-cells, each $\zeta_i$ 
is a cycle, as is $\rho$.

The 3-chain $s$ is itself in the image of 
$C_3(\mathbb {R}^k_{q*} / \mathbb  {Z}^k) \overset{i_*}{\to}\C{3}$. Thus $ds$, viewed as a 2-chain in
$\mathbb {R}^k_{q*}/ \mathbb {Z}^k$, represents the trivial element of $H_2(\mathbb {R}^k_{q*} / \mathbb 
{Z}^k)$. Since 
$\zeta_q$ is a cycle, $d\left(\sum_j(z_j+z_j')\right)$, viewed in the same way, also does.
But, with the notation of \S\ref{basis}, we have $d(z_j+z_j')=\pm 1_q\{\alpha_j,\beta_j\} \in H_2(\mathbb {R}^k_{q*} /
\mathbb {Z}^k)$, as figure
\ref{3chain} indicates. Thus, there are two values of $j$, say $j=1$ and $j=2$, such that
$d(z_1+z_1'+z_2+z_2')=0 \in  H_2(\mathbb {R}^k_{q*} / \mathbb {Z}^k)$. So, with the notation 
introduced in (\ref{star}), either 
\begin{equation}\label{case1a} \alpha_1=\alpha_2\ \ \text {and }\ \, \beta_1=\beta_2,  
\end{equation}
or
\begin{equation}\label{case2a} \alpha_1=\beta_2\ \ \text {and }\ \, \beta_1=\alpha_2.  
\end{equation}

To treat these two cases, we introduce two elements $c_1$ and $c_2$ of $\C{4}$. 
Assume that the finely dotted edges (resp.~other dotted edges, solid edges)
of the graphs in figure \ref{4chainP2c} have
color $a$ (resp.~$b$, $q$).
 Let $c_1$ be the sum
of the seven generators indicated in figure \ref{4chainP2c}, 
\begin{figure}[tb]
\centerline{\mbox{\includegraphics*{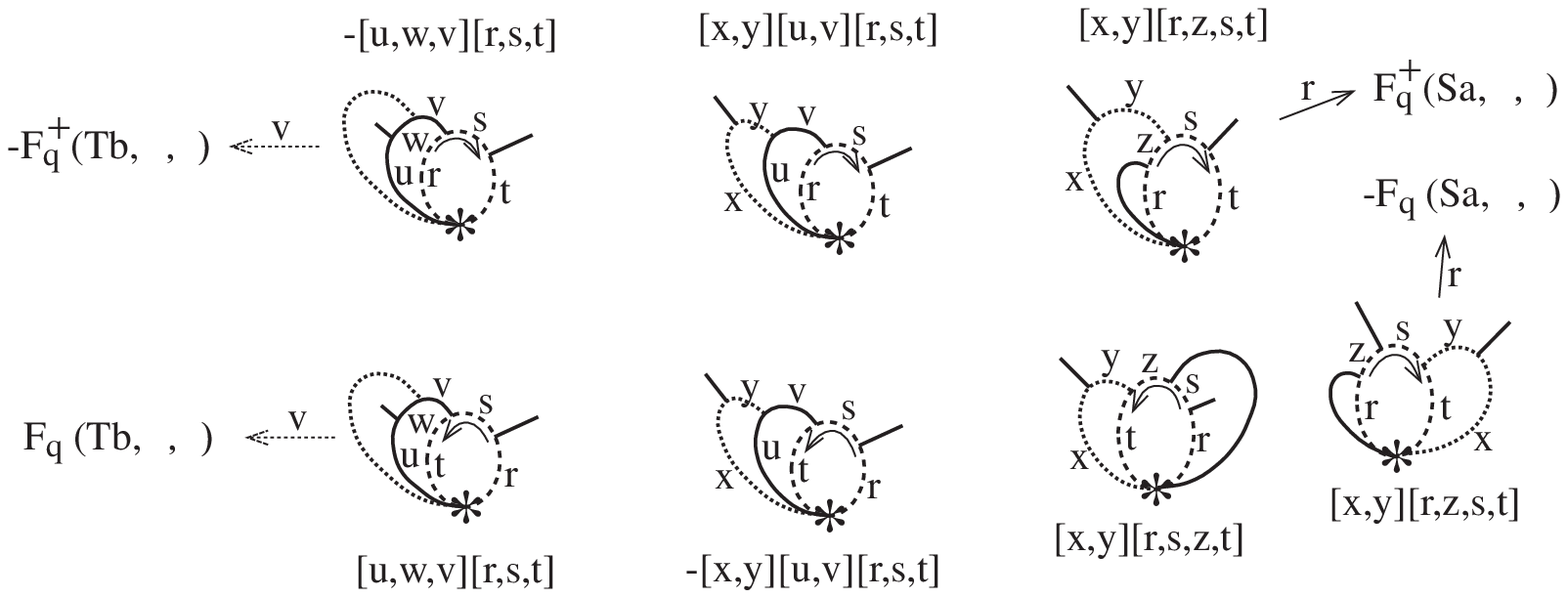}}}
\caption{}
\label{4chainP2c}
\end{figure} 
plus three more corresponding to the
three graphs on the right of the figure, but with the colors $a$ and $b$
switched, and also with the cell orientations of all three graphs switched. Then, for any consistent 
choice of orientations of the edgesets $a$ and $b$, we have 
\begin{equation}\label{c1}
dc_1-\bigl(F_q^+(Sa,\alpha,\beta)-F_q^-(Sa,\alpha,\beta)-F_q^+(Tb,\alpha,\beta)+F_q^-(Tb,\alpha,\beta)\bigr)
\in \cf_{P2}
\end{equation}
where $S,T\in\{L,R\}$.  

With the same coloring scheme,
let $c_2$ be the sum of the two generators of $\C{4}$ indicated in figure \ref{short4chain},
\begin{figure}[tb]
\centerline{\mbox{\includegraphics*{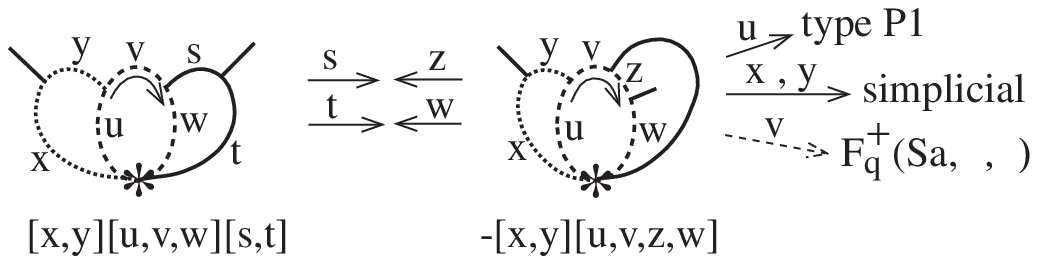}}}
\caption{}
\label{short4chain}
\end{figure} 
plus two more like these but with the orientations of the edgesets $q$ reversed, and also with
the cell orientations reversed, plus three
more like those on the right of  figure \ref{4chainP2c}, plus another five.  
These five more are similar to the five 2,4 terms of $c_2$ already described, but with both the colors
$a$ and $b$ and the half-edges $\alpha$ and $\beta$ interchanged.
Then, with cell
orientations chosen  properly,
\begin{equation}\label{c2}
dc_2-\bigl(F_q^+(Sa,\alpha,\beta)-F_q^-(Sa,\alpha,\beta)-F_q^+(Tb,\beta,\alpha)
+F_q^-(Tb,\beta,\alpha)\bigr) \in \cf_{P2}.
\end{equation}

Now, in view of (\ref{c2}), and (\ref{star}) with $j=1$ and $2$, we may assume that
(\ref{case1a}) holds, and not (\ref{case2a}).  Also, because of (\ref{c1}), we may assume 
that $a_1\ne a_2$.  Then adding to $z$ one of $\pm dc_1$ with $a=a_1$ and $b=a_2$ decreases 
by four the number of terms of $z$ with faces realizing $m(z)=2$.  This only adds type
P1 and simplicial terms to $z$, and so completes the proof.
\end{proof}

\begin{lemma} 
\label{P1}
Each $z\in \cf_{P2}$ is homologous to a cycle in $\cf_{P1}$.
\end{lemma}

\begin{proof}
Let $m(z)$ be the maximal incidence number of all the type $w$ faces of the terms of $z$.
By using colored graphs like the two on the left of figure \ref{elimP1}, we may assume that  
$m(z)=3$.
More specifically, assume $z_1$ and $z_2$ are two terms of $z$ with canceling type $w$ faces
realizing $m(z)> 3$.  If $cg(z_1)$ and $cg(z_2)$ have monochromatic loops of length two of different 
colors, then there is a generator $g$ with $cg(g)$ as shown in the lower  left of figure \ref{elimP1} such that
$m(z_1+z_2+dg)< m(z_1+z_2)$.  If the monochromatic loops have the same color then adding
one or two generators corresponding to the colored graph in the upper left achieves the same
result.

 Now the type P1 terms of $z$ can be paired off, with the elements of each pair
corresponding to colored graphs and orientations like those on the right of figure \ref{elimP1}.
The central portion of the figure shows that these are the boundary of a 4-chain.
\end{proof}
\begin{figure}[tb]
\centerline{\mbox{\includegraphics*{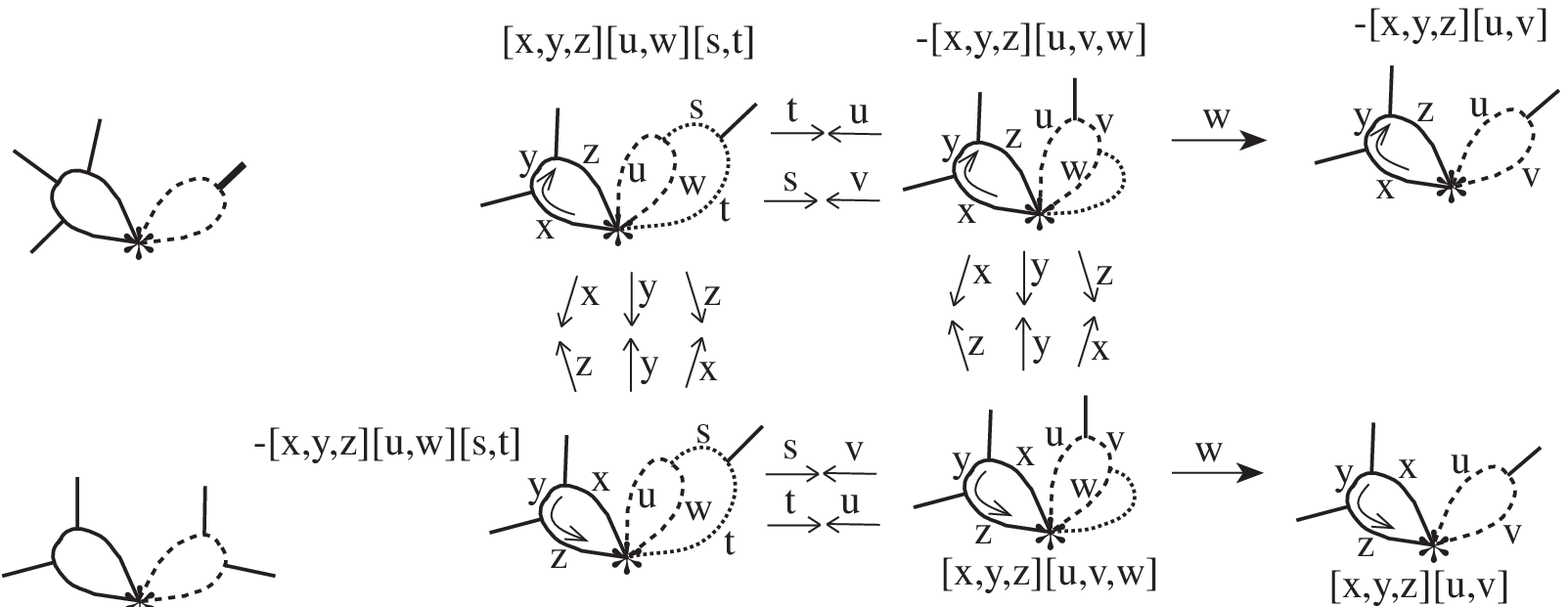}}}
\caption{}
\label{elimP1}
\end{figure} 
\begin{lemma} \label{SS}
\label{P0}
Each $z\in \cf_{P1}$ is homologous to $0\in \C{3}$.
\end{lemma}

\begin{proof}
Write 
$z= \sum_{i=1}^n \zeta _i $
where  each $\zeta_i$ is a cycle in the subcomplex 
$\mathbb{R}^k_{i*} / \mathbb {Z}^k$ of $\xt_n /StN_n$.  Each basis element of 
$H_2(\mathbb{R}^k_{i*} / \mathbb {Z}^k)$ described in \S \ref{basis} is trivial in 
$H_2(\widetilde X _n / StN_n)$ as figure \ref{3chain} indicates.   Thus, each $\zeta_i$,
and so also $z$, is trivial in $H_2(\widetilde X _n / StN_n)$.
\end{proof}
 
 Given the remarks preceeding lemmas \ref{42} and \ref{prismclassification}, this completes 
 the proof of Theorem \ref{X/S}.

\section{Computation of  $H_2(\xt_n/StN_n)$}\label{H2xmodstn}

This section is devoted to the proof of the following.

\begin{theorem} \label{h2q} The second homology group $H_2(\xt_n/StN_n)$ is free abelain of rank $n$ provided $n\ge 5$.
\end{theorem}

 We assume throughout this section that $n\ge 5$, although this is only needed in Lemma \ref{ne5}.
 Let ${\mathbb Z}^n$ be the free abelain group of rank $n$ with basis $\{e_1,\dots ,e_n\}$.
We will define a homomorphism $\Theta'\colon C_2(\xt_n/StN_n)\to {\mathbb Z}^n$ which roughly
counts how many terms of any given chain correspond to a colored graph which contains one 
of the colored 
graphs $\theta_{ij}$ as a subgraph.
(The colored graph $\theta_{xy}$  is shown on the left of Figure \ref{theta} in 
\S\ref{ECT}.)  The  lemmas in this section will show that $\Theta'$ induces an isomorphism
$\Theta\colon H_2(\xt_2/StN_n)\to {\mathbb Z}^n$, thereby proving the theorem.

The homomorphism $\Theta'$ is defined on oriented generators $x$ of  $C_2(\xt_n/StN_n)$ by

$$\Theta'(x)= \begin{cases}
0,   &\text{if $core(cg(x))$ is not a   $\theta_{ij}$,}\\
e_j, &\text{if $core(cg(x)) = \theta_{ij}$ with the same cell orientations.}
\end{cases}$$
  In the second line of the definition,
we are assuming the cell orientation of $\theta_{ij}$ is as shown in the upper left of Figure \ref{theta}.


\begin{lemma} The homomorphism $\Theta'$ vanishes on the image of the cellular boundary map 
$\partial \colon C_3(\xt_n/StN_n) \to C_2(\xt_n/StN_n)$.
\end{lemma}
\begin{proof} If $G$ is a colored graph with  $gen (G)\in C_3(\xt_n/StN_n)$ and no blowdown of
$G$ has core one of the $\theta_{ij}$, then clearly $\Theta'\circ \partial (gen(G))=0$.  So we will be concerned
with those $G$ having an edge $x$
such that the core of the blowdown $G_x$ is a $ \theta_{ij}$.  
Given such $G$ and $G_x$,
let $v$ be the vertex of $G_x$ which the edge $x$ is mapped to by the blowdown $G\to G_x$. 
Observe that 
the colored graph $G$ is determined, up to the color of the edge $x$, by
the colored graph $G_x$, the vertex $v$, and a partition into two 
nonempty sets $P_1$ and $P_2$ of the half-edges of $G_x$ meeting $v$.
(In case  $v$ is the basepoint of $G_x$, the basepoint of $G$ must also be specified as one of the two 
vertices of $x$.  This, however, will not be necessary for our purposes.)

Now we fix $G_x$ to be any colored graph with core some $\theta_{ij}$ and use the observation in the 
previous paragraph to determine all possibilities for $G$.

As one of essentially two possibilities for the vertex $v$, we 
first assume that it is either of the two non-basepoint vertices of $G_x$.
Let $a$ be the edge of the subgraph $\theta_{ij}$ which joins the basepoint $*$ to $v$, and let $c$ and
$d$ be the edges forming the monocromatic loop of $\theta_{ij}$.

As a first possibility for the partition, assume that $\{a,c,d\}\subseteq P_1$, where we are not distinguishing
between edges and half-edges incident with $v$.  Since $G$ is colored, the color of the edge $x$ of $G$ must be
different  from the two of $\theta_{ij}$, and so $x$ is a nonsingleton edge of $G$.
Let $x'$ be the edge of $G_x$ having the same color as $x$.  Then, since $G$ is colored, both $x$ and $x'$ 
are incident with the vertex of $G$ not meeting $a$, $c$ and $d$.  For each of the resulting
three possibilities
for the core of $G$, the map $\Theta'$ 
takes the corresponding generator to 0 since the subgraph $\theta_{ij}$
occurs in two faces of the boundary with opposite orientations.
The same is true, except in one case, for the other possibilities for $G$.
  We will not point this out in the rest of the proof
for each $G$, but will comment on the exceptional case.

As another possibility for the partition, assume $\{a\}\subseteq P_1$, $\{c,d\}\subseteq P_2$.  If $G$ is prismatic,
the color of $x$ must be the color of $a$.
If $G$ is cubical, with edges $x$ and $x'$ of the same color, then, in order that $G$ be colored,
one vertex of $x'$ must be incident with with the non-basepoint vertex of $a$.  Of the two possibilities
for the other vertex of $x'$, one gives a 
colored graph with $x$ and $x'$ forming a loop, the other 
the colored graph shown in the middle of Figure \ref{theta} in \S\ref{ECT}.

The latter possibility is the exceptional one in that a $\theta_{ij}$ occurs in four faces of the cell corresponding
to $G$.  Two of these faces are identical, but have opposite orrientations, so cancel in the boundary.
Of the other two, one contains $\theta_{ik}$, the other  $\theta_{jk}$ (for suitable distinct $i$, $j$ and $k$).
 These have opposite orrientations,
and so $\Theta'$ applied to their sum is 0. 

As the last possibility for the partition, assume $\{c\}\subseteq P_1$, $\{a,d\}\subseteq P_2$, although we could just
as well assume that $\{d\}\subseteq P_1$, $\{a,c\}\subseteq P_2$. If $G$ is prismatic, the color of $x$ muat be that 
of $c$ and the core of $G$ is as shown in Figure \ref{theta}(a).  If $G$ is cubical
with $x$ and $x'$ of the same color, then $G$ is only colored if $x$ and $x'$ form a loop in $G$.

Now assume that $v$ is the basepoint of $G_x$.  Let $b$ be the fourth edge of $\theta_{ij}$ different
from the edges $a$, $c$ and $d$ specified above. 
Assume $\{a,b\}\subseteq P_1$.  Then $G$ must be cubical, and the edges 
$x$ and $x'$ of the same color of $G$ must both be incident with the vertex of $G$ disjoint from $a$ and $b$.

Finally assume $\{a\}\subseteq P_1$, $\{b\}\subseteq P_2$.  If $G$ is prismatic, 
then $x$ has the color of $a$ and $b$, giving a case already covered.  If $G$ is cubical, 
the only posssibility for the core of $G$ is the partitioned graph C6 as shown in Figure 1 
of \S\ref{ECT}.
\end{proof}

\begin{lemma} The homomorphism $\Theta\colon H_2(\xt_n/StN_n)\to {\mathbb Z}^n$ induced by $\Theta'$ is onto.
\end{lemma}

\begin{proof} In \S\ref{pictures} we will exhibit elements of $H_2(\xt_n)$ whose images in 
$H_2(\xt_n/StN_n)$
are mapped to the generators of ${\mathbb Z}^n$ by $\Theta$.  Alternatively,  the left portion of 
Figure \ref{genus2} indicates a cycle
\begin{figure}[tb]
\centerline{\mbox{\includegraphics*{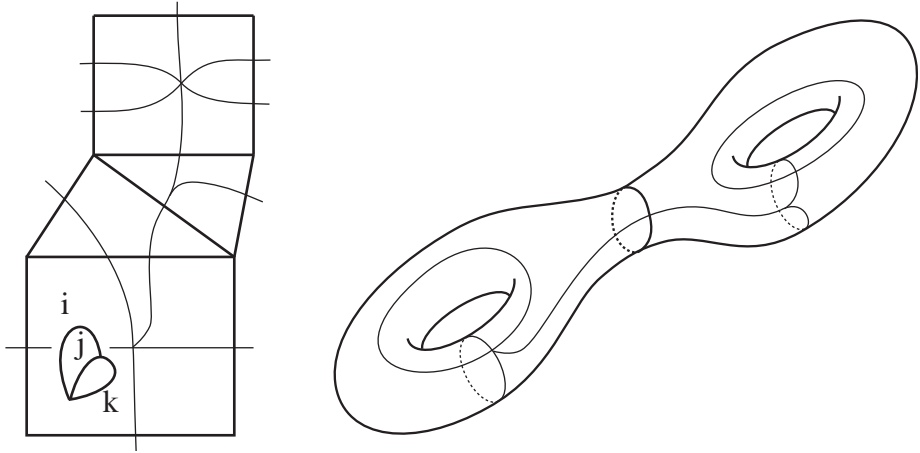}}}
\caption{Two views of a generator of $H_2(\xt_n/StN_n)$.}
\label{genus2}
\end{figure}
in $\xt_n/StN_n$ which $\Theta$ takes to the generator $e_k$ of ${\mathbb Z}^n$.
\end{proof}

This last cycle, which we denote by $t$, is 
carried by the image in $\xt_n/StN_n$ of a genus two surface.  The map of this surface to $\xt_n/StN_n$
is indicated on the right of the figure.
\begin{remark}\label{tt}  \rm  We note that $t$ may be viewed as an element of the homology of certain triangular 
subgroups of Aut$(F_n)$.  Specifically,
let $\sigma$ be any partial order of the set $\bf n$ with $i>k$, $j>k$, $i\ne  j$, 
and let $T$ be the corresponding triangular subgroup $T_n^\sigma$ as defined in \S\ref{tri}.
Consider the maps
 $$ C_2(ET/T)\buildrel i_*\over\longrightarrow C_2(\xt^{alg}_n/StN_n)\buildrel
 \mu_*\over\longleftarrow C_2(\xt_n/StN_n)$$ 
where $i_*$ is induced by the inclusion $ET\hookrightarrow \xt^{alg}$, and
 where $\mu_*$ is induced by the map $\tilde\mu_n\colon\xt\to \xt^{\rm alg}$ defined in \S\ref{xtoxalg}.    
 Choose $t^t\in C_2(ET/T)$ with $i_*(t^t)=\mu_*(t)$.
 We will see in Theorem \ref{last_thm}
 that $t^t$ represents a nontrivial element of $H_2(T)$ which is closely related to exotic
elements of $K_3(\mathbb Z)$. 
\end{remark}

A square (resp.\ triangular or simplicial) colored graph is a colored graph which corresponds to a
1,1-cell (resp.\ 0,2-cell).  Such cells were discussed near the beginning of \S\ref{ECT}.
Generators of $C_2(\xt_n/StN_n)$ and $C_2(\xt_n)$ are called square or triangular or simplicial
if the corresponding colored graphs are.

The following lemma is used to show that $\Theta$ is injective and will also be used in 
the proof of Lemma \ref{T}. 
We say that a chain $c$ in $C_2(\xt_n/StN_n)$ or $C_2(\xt_n)$
is reduced if, for each square term of $c$, all the non-basepoint
vertices of the corresponding colored graph have valence three.

\begin{lemma}\label{red} Let $G$ be any square colored graph with  $gen(G) \in C_2(\xt_n/StN_n)$.  Then there is a
reduced chain c in $C_2(\xt_n/StN_n)$ such that $\partial c=\partial(gen(G))$ and $c-gen(G)=0$ in $H_2(\xt_n/StN_n)$.
 If $G$ is marked as well as colored with $gen(G) \in C_2(\xt_n)$, then the same conclusion holds if
 the space $\xt_n/StN_n$ is replaced with $\xt_n$.
\end{lemma}

\begin{proof}
Disregarding basepoints, there are only two possibilities for the core of $G$.  Either it consists of two 
monocromatic loops sharing one vertex or of a single monocromatic loop with the two other edges meeting all
three vertices of $G$.  Considering basepoints, there are essentially four cases to consider, each of which is
straightforward.  The key point is that prismatic cells can be used to decrease the valence of selected vertices
of colored graphs corresponding to square cells, at the expense of introducing triangular cells.
\end{proof}

A similar lemma, but applicable to triangular rather than square graphs is given next.  It will be cited 
in the last lemma of this section and is also used in the proof of Lemma \ref{T}.

We will use the following notation.  Let $G$ be a simplicial colored graph.  Assume that the monocromatic edge
loop of $G$ of length $>1$ is $e_0, e_1,\ldots,e_k$ where $init(e_0)=term(e_k)=*$, the basepoint of $G$.
Let $A_i$ be the set of half edges of $G$, other than $e_{i-1}^+$ and $e_i^-$, which are incident with the
vertex 
$term(e_{i-1})=init(e_i)$. Then we write $[A_1|A_2|\ldots |A_k]_a$ for the oriented simplicial cell
$[e_o,e_1,\ldots,e_k]$ of $\xt_n/StN_n$, where $a$ is the color of the loop $e_0,\ldots,e_k$.  Often the subscript 
$a$ is dropped. 

If $G$ is also marked, then the corresponding oriented cell $[e_o,e_1,\ldots,e_k]$
of $\xt_n$ is also denoted by $[A_1|A_2|\ldots |A_k]_a$. 
This notation makes no mention of the marking, which is usually clear from context.
 
 With this bar notation, the boundary maps for both the space $\xt_n$ and $\xt_n/StN_n$ are 
$$d_0 [A_1|A_2|\ldots |A_k]= [A_2|\ldots |A_k], \quad \ d_k [A_1|A_2|\ldots |A_k]= (-1)^k[A_1|\ldots
|A_{k-1}],$$ and, if $0<i<k$,
$$d_i [A_1|A_2|\ldots |A_k]=(-1)^i [A_1|\ldots |A_{i-1}|A_i\cup A_{i+1}|A_{i+2}|\ldots |A_k].$$ 

By the cardinality of the cell $ [A_1|A_2|\ldots |A_k]$, we mean the maximum cardinality of the sets $A_i$.

\begin{lemma}\label{sred}
If $c$ is a simplicial 2-chain in $H_2(\xt_n/StN_n)$ or $H_2(\xt_n)$ with each 1-cell in the boundary of $c$ having
cardinality one, then $c$ is homologous rel boundary (in the sense of Lemma \ref{red})
to a  2-chain whose terms all correspond to simplicial cells of
cardinality two.
\end{lemma}

Geometrically, this says that $c$ is homologous to a chain consisting of pairs of simplicial 2-cells meeting 
along 1-cells of cardinality 2.  Each such pair forms a square corresponding to one of the  relations R1.

\begin{proof}[Proof of \ref{sred}]
Among all the edges of the 2-cells of $c$, let $e$ be one of maximum cardinality. Orient $e$ and assume
that $[A|B]$ and $[A'|B']$ are 2-cells corresponding to terms  of $c$ with
$d_1[A|B]=e=d_1[A'|B']$. Thus $[A\cup B]=[A'\cup B']=e$.  We will show that, if the 
cardinality of $e$ is $>2$, then the square
consisting of $-[A|B]$ and $[A'|B']$ is homologous to a simplicial 2-cycle whose cells all have cardinality
less than that of $e$. This will complete the proof.

As a first case, assume that $A'=A-a$, where $a\subsetneq A$, and $B'=B\cup a$. Then the 3-cell
$[A'|a|B]$ provides the desired homology since 
$$d[A'|a|B]=[a|B]-[A|B]+[A'|B']-[A'|a]$$ and the first and last terms have cardinality less than 
that of $e$.

For another case, assume that $A'=A\cup b -a$ and $B'=B\cup a-b$, where $a\subset A$ and 
$b\subset B$ both consist of one  element.
Since the set $A\cup B$ has cardinality $>2$, one of the containments  is proper, say $a\subsetneq A$.
Then the boundary of the 3-cycle $[A-a|a|B]-[A-a|b|B']$ is $-[A|B]+[A'|B']$ plus four terms corresponding 
to cells with cardinality less than that of $e$.

For the last case, assume that $A'=A\cup b -a$ and $B'=B\cup a-b$ as above, where now $a\subset A$ and
$b\subset B$  are nonempty subsets with $a\cup b$ containing  at least three elements.  Say $a=a_1\sqcup a_2$,
where both of the $a_i$ are nonempty.  Then the 3-chain 
$$[A-a_1|a_1|B]-[A-a_1|b|a_1\cup B-b]+[A'|a_2|B-b\cup a_1] $$
provides a homology between $-[A|B]+[A'|B']$ and six 2-cells, each with cardinality less than that of $e$.
\end{proof}

To show that $\Theta$ is injective we work with a filtration
$$ker\ \Theta'\cap Z_2(\xt_n/StN_n) \supset \cf_{R} \supset \cf_{\theta} \supset \cf_{F} \supset \cf_{Sq} \supset 
 0$$  by subgroups 
of $ker\ \Theta' \cap Z_2(\xt_n/StN_n)$, where $\cf_R$ consists of the reduced cycles in $ker\ \Theta'$.  
Other terms in the filtration will be defined as we go. 

 It follows from Lemma \ref{red} that all the chains in
$ker\ \Theta' \cap Z_2(\xt_n/StN_n)$ are homologous to reduced ones.  
The remaining lemmas in this section will show that the elements of each subgroup of the filtration
are homologous to ones of the next smaller subgroup of the filtration.  This will show that $\Theta$
is injective.

The subgroup $\cf_{\theta}$ is defined as all the elements of $\cf_R$  with no term
corresponding to a colored graph with core a $\theta_{ij}$.

\begin{lemma}  Each element of $\cf_{R}$ is homologous to one of  $\cf_{\theta}$. \end{lemma}
\begin{proof} Let $z\in \cf_R$. Then, since $\Theta'(z)=0$, if $x$ is a term of $z$ with $\Theta'(x)\ne 0$,
there is another term $x'$ of $z$ with $\Theta'(x+x')=0$.  Note that if $cg(x)=\pm \theta_{ik}$,
then  (assuming $x\ne -x'$)  $cg(x')=\mp\theta_{jk}$ for some $j\ne i$.

Now let $g_7$ be the generator of $C_2(\xt_n/StN_n)$ corresponding to the partitioned graph C7 in Figure 1,
oriented and colored so that, for a suitable sum $p$ of prismatic terms, the boundary $\partial(g_7+p)$ contains
terms which cancel $x+x'$ and the cycle $z+\partial(g_7+p)$ is reduced and has fewer terms on which $\Theta'$ is nonzero than $z$.  
A cycle homologous 
to $z$ with no such terms is obtained by, if necessary, adding more boundaries to $z$ like $\partial(g_7+p)$,
but with perhaps different colorings of the corresponding graphs.
\end{proof}

We say that a colored graph is of type F if it is square, reduced, and has a single monocromatic loop.
We us the notation $F_2(Si/Tj/k)$, or just $F(Si/Tj/k)$ for type F graphs as in \S\ref{ECT}  just before Lemma \ref{0}.

The group $\cf_F$ in the filtration is defined  as all chains in $\cf_{\theta}$ with no terms corresponding to a colored graph of type F.

\begin{lemma}  Each element of $\cf_{\theta}$ is homologous to one of $\cf_{F}$. \end{lemma}

\begin{proof}  Let $D$ be the free abelian group 
$\bigoplus_{1\le i,j\le n; i\ne j}(\mathbb Z^L_{ij}\oplus \mathbb Z^R_{ij}) $  where each $\mathbb Z^L_{ij}$
and $\mathbb Z^R_{ij}$ is a copy of $\mathbb Z$.  We let $e^L_{ij}\in \mathbb Z^L_{ij}$ and 
$e^R_{ij}\in \mathbb Z^R_{ij}$ be fixed generators.  Let $f\colon C_1(\xt_n/StN_n)\to D$ be the homomorphism
defined on oriented, reduced generators $x$ by 
$$f(x)= \begin{cases}
e^L_{ij},   &\text{if $x$ corresponds to $L_{ij}$,}\\
e^R_{ij},   &\text{if $x$ corresponds to $R_{ij}^{-1}.$}
\end{cases}$$ 
  On the remaining (unreduced) generators, we define $f$
 so that $f(\partial s)=0$ for all oriented triangular 
generators $s$ of $C_2(\xt_n/StN_n)$.

We also define a homomorphism $F\colon C_2(\xt_n/StN_n)\to D$ by $F(x)=f(\partial s)$ for each  
$x\in C_2(\xt_n/StN_n)$.  It is easy to check that $F$ vanishes on all square reduced generators 
in $\cf_{\theta}$ except for those of
type $F$, and, for these, $F(F_2(Si/Tj/k))=\pm e^S_{ik}$.

Now let $z\in \cf_{\theta}$.  Assume $z$ contains a term $c$ of type $F$.  Say $c=\pm gen(F(Si/Tk/j))$.
Then, since $F(z)=f(\partial z)=f(0)=0$, there is another term $c'=\mp gen(F(Si/Al/j)$ of $z$ such that  
$F(c')=-F(c)$.

If $k=l$, then there is a single prismatic term $p$ and a generator $g_6$ of $C_3(\xt_n/StN_n)$
corresponding to the partitioned 
graph C6 of Figure 1 suitably colored such that $z+\partial(g_6+p)$ is in $\cf_{\theta}$ and contains two
fewer type $F$ terms than does $z$.  Adding $\partial(g_6+p)$ to $z$ essentially trades $c$ and $c'$ for two 
triangular terms and one square one, not of type $F$.

If $k\ne l$, then there is a sum $p$ of prismatic terms such that $z+\partial(F_3(Si/Tk/Al/j)+p)$ is in 
$F_{\theta}$ and has two fewer type $F$ terms than $z$.  This time $c$ and $c'$ are exchanged for 
a number of triangular and square terms, where again the square terms are not of type $F$.
\end{proof}

Let $\cf_{Sq}$ be the elements of $\cf_{F}$ (or of $\cf_{R}$) with no square terms.

\begin{lemma} \label{ne5}
 Each element of $\cf_{F}$ is homologous to one of $\cf_{Sq}$. \end{lemma}

\begin{proof} There are two types of square colored graphs which terms of elements of $\cf_F$ may correspond to.
Both have two monocromatic loops, and one type has rank three, the other rank four.  Generators corresponding to 
both types have boundary zero.

Given a generator $x$ corresponding to either of these types of colored graphs, there is a chain $c$
consisting of one cubical and one prismatic term with boundary $x$. 
The cubical term corresponds to the partitioned graph
 C2(a) of Figure 1, if $cg(x)$ has rank 4, and to C2(b), if it has rank 3.
The prismatic term essentially isolates $x$. 
\end{proof}

\begin{remark} \label{rem} \rm The boundary of a lift to $\xt_n$ of one of  
the chains $c$ in the previous proof in the rank 4
case can be represented by a picture like that on the left of Figure \ref{lemT2}, 
but with the label of the bottom circle changed from $R_{ij}$ to $R_{kl}$.
All the generators corresponding to the crossings in the picture cancel in 
$C_2(\xt_n/StN_n)$ except for the one just above the middle of the picture.  This crossing can be 
taken to correspond to $x$.
\end{remark}

The proof of the Theorem \ref{h2q} is completed by the following.

\begin{lemma}  Each element of $\cf_{Sq}$ is homologous to $0\in C_2(\xt_n/StN_n)$. \end{lemma}

\begin{proof} Let $z\in \cf_{Sq}$.  Note that the terms of $z$ are all triangular.  By adding 
canceling pairs of more
triangular terms to $z$, we can write $z$ as a sum of chains whose boundaries all have cardinality 1.
So by Lemma \ref{sred}, we may assume that all the terms of $z$ have cardinality 2.  These can be paired off into 
cycles of two types--type R2 and type R3, where
 both colored graphs corresponding to the terms of a type Rn cycle have rank n.

Let $y$ be a cycle of type R2.  Then there is a
 3-chain $c=p_1+p_2 + S$ with  $\partial c= y$.
Both of the $p_i$ are prismatic terms and $S$ is a sum of three
simplicial terms.
 One $p_i$ corresponds to the (suitably colored) partitioned graph shown
in Figure \ref{lemT2}.  The other corresponds to the same colored graph, but with the orientation of the edgeset $j$
reversed.  Cell orientations are chosen so that contracting the edge $j_0$ in both graphs gives canceling
terms in $\partial c$.
The three simplicial terms effectively decrease the cardinality of two of the triangular terms in
$\partial(p_1+p_2)$ from 3 to 2, thereby isolating $y$ in $\partial c$.

If $y$ has type R3, then there is a 3-chain $c'$ with $\partial c'=y$.  The chain $c'$ is like $c$, but the
singleton edge of color $i$ in each of the five colored graphs corresponding to the terms of $c$ is replaced with two 
singleton edges. Each of these meets the basepoint of the relevant graph and one of the vertices which $i$ meets.
\end{proof}

As in Remark \ref{rem}, the picture in figure \ref{lemT2} represents a lift of a cycle of type R2 to $\xt_n$. 
The crossing just above the center of the picture is the only one that does not cancel in $C_2(\xt_n/StN_n)$.

%






%

%

\newcommand{\vt}{\mbox{$\widetilde V$}}

\section{The  tree complex $T(G)$}
Fix a graph $G$ (equipped with neither a partition nor a basepoint). After some preliminary definitions, we define the 3-skeleton of a $CW$-complex
$T(G)$.  The  2-skeleton of $T(G)$ was implicitly defined in \cite{HT}, where 
it  was shown that $T(G)$ is 1-connected.  
We show that $T(G)$ is 2-connected and use related facts in the proof of Lemma \ref{trivtree}, 
which in turn is used in the proof of Theorem \ref{thminW}.

A graph $G'$ is said to be a quotient of $G$ if it is obtained from $G$ by collapsing each component of an acyclic subgraph of $G$ to a point. 
If $G'$ is a quotient of $G$, with corresponding quotient map $q\colon G\to G'$, we take 
$E(G')=\{\,e\in E(G): q(e) \rm\  is \ not\ a\ point \,\}$,
with $q$ inducing the identity on $E(G) \cap E(G')$.
  Two quotients $G'$ and $G''$ of $G$ are equal if and only if $E(G')=E(G'')$.
Thus collapsing
distinct acyclic subgraphs (with no components consisting of just a vertex) gives, by definition,
distinct quotients, even if these are homeomorphic.


The partitioned quotients of $G$ form a poset $\cp\cq (G)$, where
$G'\le G''$  if $G'$ is a blowdown of $G''$
 (in which case the partition of  $G'$  is obtained by deleting the edges from the edgesets of 
$G''$ which are collapsed in the blowdown.)
  Maximal elements of $\cp\cq (G)$ correspond to the various partitions of $G$.
 If $G' \in \cp\cq (G)$, let 
 $\cp\cq_{\le G'}$ be the full subposet of $\cp\cq(G)$
consisting of all
$H\in \cp\cq(G)$ with $H\le G'$.  The analog of lemma \ref{A} with $\cp\cq_{\le G'}$ in place of 
 ${\rm b}G$ follows 
immediately from the definition of $\cp\cq(G)$.  So the analog of Prop.\ \ref{altdef} also holds.  That is,
the geometric realization $|\cp\cq(G)|$ of $\cp\cq(G)$ has a CW-structure with one $k$-cell for each
partitioned quotient of $G$ with $k+1$ vertices. 
Specifically, if $G'$ is such a partitioned quotient, the geometric realization $|\cp\cq_{\le G'}|$ of 
$\cp\cq_{\le G'}$ forms a k-cell of $|\cp\cq(G)|$.

Note that the 0-cells of $|\cp\cq(G)|$ correspond to the maximal trees of $G$.  Let $T_i$ be the maximal
tree corresponding to the 0-cell $v_i$.  A 1-cell of $|\cp\cq(G)|$ joins two 0-cells $v_i$ and $v_j$
if the number of edges in $T_i\cap T_j$ is one less than the number of edges in any maximal tree
of $G$.  Thus the 1-cells of $|\cp\cq(G)|$ correspond to elementary moves as defined on p.~225 of 
\cite{HT}.  

\subsection{The definition of $T(G)$} \label{defTG}
The 2-skeleton of $T(G)$ is defined by attaching a class of 2-cells, called $w$-triangles, to $|\cp\cq(G)|$.
(In the language of pseudo-isotopy theory, these $w$-triangles correspond to dovetail singularities.
The 2-cells of $|\cp\cq(G)|$ correspond to the other codimension two singularities.)  
For each quotient $G'$ of $G$ with two vertices, 
there is one $w$-triangle corresponding to each subgraph of $G'$
of the sort shown at the top of figure \ref{wcells}(i).  
In other words, there is one $w$-triangle corresponding to each triple $\{T_1,T_2,T_3\}$ of distinct
maximal trees of $G$ with $T_1\cap T_2 \cap T_3$ consisting of one fewer edge than $T_1$.
The lower part of figure \ref{wcells}(i)
 indicates the attaching map
of the $w$-triangle corresponding to the subgraph above it.  The edge labeled $a \wedge b$, for example, maps
to  the 1-cell of $|\cp\cq(G)|$ corresponding to $G'$ with the partition $\{a,b\}$.  The map is such that 
the vertex on the left maps to the 0-cell corresponding to the maximal tree of $G$ which contains $a$.

\begin{figure}[tb]
\centerline{\mbox{\includegraphics*{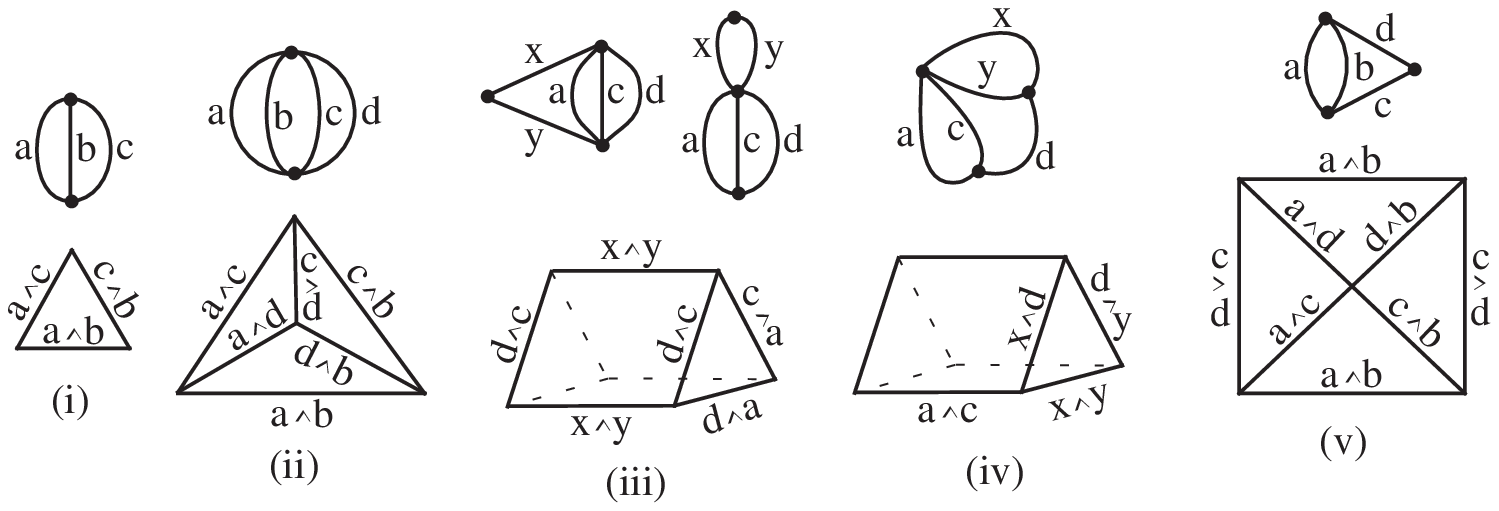}}}
\caption{Cells of $T(G)$ added to $|\cp\cq(G)|$.}
\label{wcells}
\end{figure}
To complete the definition of $T(G)$, three types of 3-cells, called $w$-prisms, $w$-tetrahedra and 
pillows, are attached to the 2-skeleton of $T(G)$.  
For each subgraph of the type shown in figure \ref{wcells}(ii) of each quotient of $G$ with 
two vertices,
or, equivalently, for each quadruple $\{T_1,T_2,T_3,T_4\}$ of distinct
maximal trees of $G$ with $T_1\cap \ldots \cap T_4$ consisting of one fewer edge than $T_1$,
there is one $w$-tetrahedron attached as 
indicated in figure \ref{wcells}(ii).
 Specifically, the four faces of the tetrahedron in (ii) are attached to the $w$-triangles corresponding 
to the four subgraphs consisting of three edges of the graph in (ii).

Now let $G'$ be any quotient of $G$ with three vertices. 
 There is one 
$w$-prism as shown in figure \ref{wcells}(iii) corresponding to each of the two subgraphs of $G'$ of 
the types shown at the top of (iii).
The two triangular faces of the $w$-prism in (iii) attach to $w$-triangles.  The three rectangular faces
correspond to partitions of the graph in (iii), and so attach to cells of $|\cp\cq(G)|$.
  There are two
$w$-prisms corresponding to each subgraph of $G'$ of the type shown in figure \ref{wcells}(iv).  
One of these two 
$w$-prisms is as in (iv), the other is as in (iii).  There is one pillow
as shown in (v) corresponding to each subgraph of $G'$ of the type shown in (v).  Two of the faces of this
pillow attach to 
$w$-triangles.  The other three faces attach to the 2-cells of $|\cp\cq(G)|$ corresponding to the graph
in (v) 
with the three  partitions $\{a,c,d\}$, $\{b,c,d\}$ and $\{\{a,b\},\{c,d\}\}$. 
This completes the definition of $T(G)$.
 
 As examples, let $G_1,\dots , G_6$ be the graphs from left to right in Figure \ref{wcells}.
 Then for $i=1,2,4$ and $6$, the complex $T(G_i)$ is the one just below $G_i$ in the figure. 
 $T(G_3)$ consists of the $w$-prism in (iii), with one $w$-tetrahedron attached at each end,
and with  three pillows then attached to this.

It is not difficult to see that $T(G)$ is connected.  This is the content of lemma 1.3 of \cite{HT}.
Lemma 1.5 of the same paper implies that $T(G)$ is simply connected.  For relations (1) and (3) of that
lemma correspond to 2-cells of $T(G)$ associated with partitioned graphs, and relation (2) to $w$-triangles.

\subsection{Simplification of cell-to-cell maps.}
A map of CW-complexes is said to be cell-to-cell if each closed cell of the domain is mapped
homeomorphically to a closed cell of the range.  The $CW$-complex with 
underlying space the 2-sphere $S^2$ and $CW$-structure $\cc$ is denoted by $(S^2,\cc)$.

Let $f\colon (S^2,\cc)\to T(G)$ be a cell-to-cell map. The following four types of operations will be used 
in the proof of lemma \ref{2conn} below to 
progressively simplify $f$.  

\begin{enumerate} 
\item[(a)] A deformation of $f$ to $f'\colon (S^2,\cc')\to T(G)$ which is 
 supported on the interior of a closed disk $D$ equal to the union of 2-cells
of $\cc$ and such that $f(D) \cup f'(D)$ is the boundary of a 3-cell or chain $\tau$ of $T(G)$.  Such deformations
arise below when
the image of $D$ is deformed, rel its boundary, through $\tau$.  By subdividing $D$, we may assume 
that $f'$ is cell-to-cell.
\end{enumerate}

Assume that the 2-cells $\sigma$ and $\tau$ of $\cc$ share at least an edge and that $f$ maps $\sigma$ and 
$\tau$ to the same 2-cell of $T(G)$.  As a subspace of $S^2$, $\sigma \cup \tau$ is either a disk, 
an annulus, or all of $S^2$ (since the boundary of each 2-cell of $T(G)$ consists of at most four 1-cells).

\begin{enumerate}
\item[(b)] If $\sigma \cup \tau$ is a disk $D$, a deformation of $f$ to $f'\colon (S^2,\cc')\to T(G)$ 
which is supported on the interior of $D$ and 
which has the effect that $f'(D)=f(\partial D)$.  By increasing the support
of the deformation slightly to a  neighborhood of $D$, we may assume that $f'$ is
cell-to-cell, where $\cc'$ is obtained from $\cc$ by deleting the interior of $D$, and then identifying points
in
$\partial D$ with the same image under $f$. \vskip2pt
\item[(c)] In case $\sigma \cup \tau$ is an annulus $A$, delete the interior of $A$ from $S^2$ and 
identify the points of the remaining edges of $\sigma$ and $\tau$ which have the same image under $f$.
The resulting quotient space $Q$ consists of two 2-spheres, and $f$ induces a cell-to-cell map $f/\!\!\sim$ from
$Q$ to $T(G)$. The operation we will use amounts to replacing $f$ with the two cell-to-cell maps
defined by restricting $f/\!\!\sim$ to each of these 2-spheres. \vskip2pt
\item[(d)] If $\sigma \cup \tau$ is all of $S^2$, any deformation of $f$ to a constant map.
\end{enumerate}

\begin{lemma} \label{2conn}
Let $f\colon (S^2,\cc)\to T(G)$ be a cell-to-cell map.  Then some sequence of the above four operations
transforms $f$ to a constant map or, if operation (c) is used, to several constant maps.
\end{lemma}

It follows that $T(G)$ is 2-connected, cf.~ lemma \ref{added}, but we will not use this fact.

Before proving \ref{2conn}, some terms and lemmas are given.
For each closed edge $e_i$ of $T(G)$, pick a point $m_i$, called the midpoint of $e_i$, in the interior of
$e_i$. The components of $e_i-m_i$ are called half-edges.  If $v$ is a vertex of $T(G)$, let $T(v)$ 
be the maximal tree of $G$ associated with $v$.  Let $e$ be an edge of $T(G)$ with vertices $a$ and $b$.
If $e_a$ is the halfedge of $e$ containing $a$, let $E_G(e_a)$ be the edge $T(a)-T(b)$ of the graph $G$.
This defines a function $E_G$, called the edge function, from the halfedges of the complex $T(G)$ to the
edges of the graph $G$.  A halfedge $h$ with $E_G(h)=\beta$ is called a $\beta$-halfedge of $T(G)$.

Now let $h$ be any halfedge of $T(G)$.  Assume that $h$ is contained in the 2-cell $\sigma$ of $T(G)$, 
and that $E_G(h)=\beta$.  Then, it follows from the definition of $T(G)$ that
 $\sigma$ contains exactly one other $\beta$-halfedge, $h'$ say.  The two $\beta$-halfedges $h$ and $h'$
meet in a vertex of $\sigma$ in case $\sigma$ is a $w$-triangle.  For the other 2-cells $\sigma$ of $T(G)$,
they are disjoint, and there is a unique edge of $\sigma$ which meets both $h$ and $h'$.  This unique
edge is said to be $\beta$-parallel for $\sigma$.  An edge of $T(G)$ is said to be $\beta$-parallel if it is 
$\beta$-parallel for at least one of the 2-cells of $T(G)$ containing it.

Assume now that $f\colon (S^2,\cc)\to T(G)$ is a cell-to-cell map. Use $f$ to
pullback the definitions of midpoints,
halfedges, the edge function $E(G)$, $\beta$-halfedges and $\beta$-parallel edges to the $CW$-structure $\cc$
on $S^2$.

For each of the following three lemmas, which will be used in the proof of lemma \ref{2conn},
 we assume that $\sigma$ and $\tau$ are 2-cells of $\cc$.

\begin{lemma} \label{1}
If $\sigma \cap \tau$ contains an edge which is $\beta$-parallel for $\sigma$ but not $\tau$,
then there is a deformation (consisting of type (a) operations above) of $f$ to a cell-to-cell
map $f'\colon (S^2,\cc')\to T(G)$ such that
\begin{enumerate}
\item the $CW$-structures $\cc$ and $\cc'$ agree except on the interior  of the disk $\sigma \cup \tau$,
as do the functions $f$ and $f'$,
\item the two edge functions associated with the two $CW$-structures $\cc$ and $\cc'$ have the same image, and
\item the $\beta$-parallel edges of $\cc'$ in the closed disk $\sigma \cup \tau$ are the two or three edges
of $\tau$ different from $\sigma \cap \tau$.
\end{enumerate}
\end{lemma}

\begin{proof}
Figure \ref{bp1} shows all possibilies for $\sigma$ and $\tau$, assuming these are as in the
statement of the lemma.  The six cases below treat these six possibilities, with case (j) treating
the configuratation shown in part (j) of the figure.
In each case, $G'$ denotes the
graph  obtained from $G$ by blowing down the edges common to all the maximal trees of $G$ associated with 
the vertices of $T(G)$ in the closed disk $f(\sigma \cup \tau)$.

\begin{figure}[tb]
\centerline{\mbox{\includegraphics*{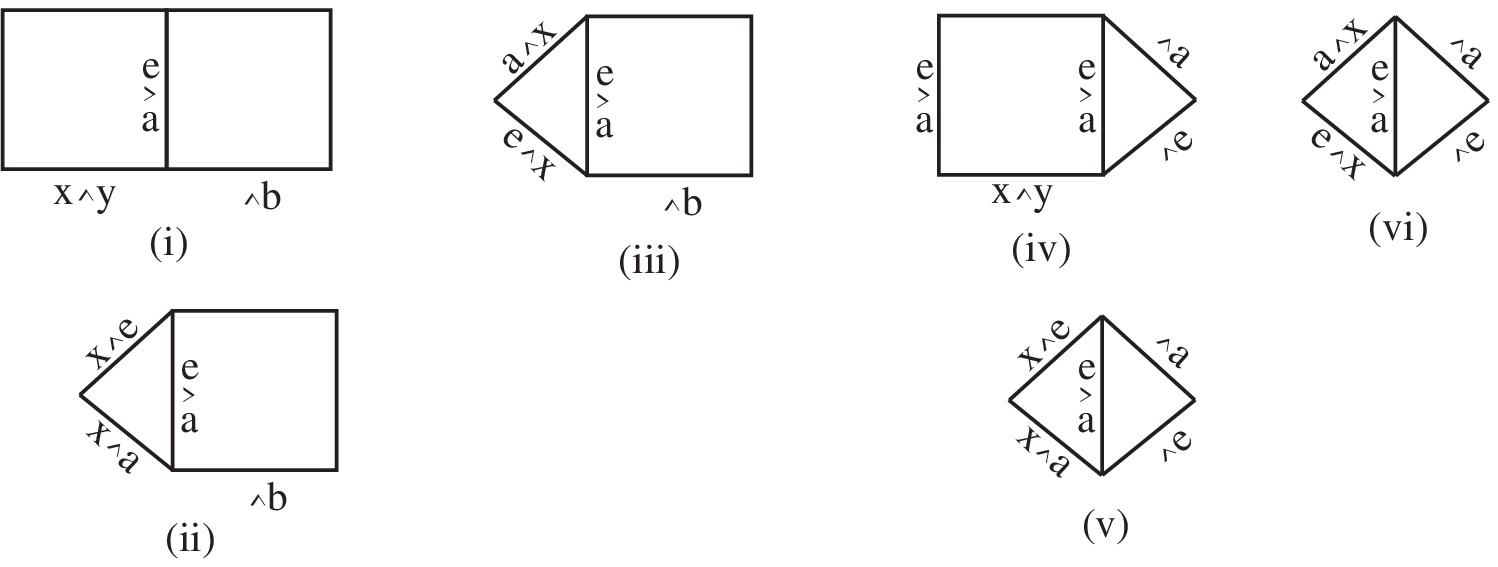}}}
\caption{}
\label{bp1}
\end{figure}
\subsubsection*{ Case (i).} First assume that $x=b$. Then we claim that $G'$ is partitioned by 
$\{\{a,e\},\{\beta,b,y\}\}$. 
(Singleton edges are not shown in the notation for partitions.)
Indeed, let $\omega$ be any cycle in $G'$. If $\omega_y$, the cycle in $G_y'$ corresponding to
 $\omega$, contains either a 
singleton edge of $G_y'$ or $\{a,e\}$, then $\omega$ contains either an edge of $G'$ different from
those in $\{a,e,\beta,b,y\}$ or $\omega$ contains $\{e,a\}$. Otherwise $\omega_y \supset\{\beta,b\}$.
Then $\omega_\beta \supset \{b,y\}$ since $G_\beta '$ is partitioned by $\{\{b,y\},\{a,e\}\}$. 
So $\{\beta ,b,y\} \subset \omega$, and $G'$ is partitioned as claimed. A deformation of type (a)
above through the corresponding 3-cell of $|\cp\cq(G)|$ completes case (i) provided $x=b$.

Now assume that $x\ne b$.  The following two observations will be used.
\begin{enumerate}
\item[(1)] No length two cycle of $G'$, except for possibly $\{b,x\}$, contains one element from two
of the three sets $\{a,e\}$, $\{x,y\}$ and $\{\beta,b\}$.
\item[(2)] No length three cycle of $G'$, except for possibly $\{a,b,x\}$ or $\{e,b,x\}$, 
contains one element from each of the 
 three sets $\{a,e\}$, $\{x,y\}$ and $\{\beta,b\}$.
\end{enumerate}

To prove (1), assume $\omega$ is a cycle of $G'$ of length two consisting of one edge from each of the sets
$\{a,e\}$ and  $\{x,y\}$.  Then $\omega_\beta$ would not be colored in $G_\beta$, a contradiction.
Similarly, $\omega$ does not consist of one element from each of $\{a,e\}$ and $\{\beta,b\}$
since $G_y$ is partitioned by these two sets.  If $\omega$ consists of one element from each of 
 $\{x,y\}$ and $\{\beta,b\}$ and $\beta \in \omega$, then $\omega_\beta$ is not colored in $G_\beta'$,
a contradiction. Similarly, $y \notin \omega$.  This leaves $\omega=\{b,x\}$ as the only possibility
and completes the proof of (1).

For the proof of (2), assume $\omega$ consists of one element from each of $\{a,e\}$, $\{x,y\}$ and
$\{\beta,b\}$.  If $\beta\in \omega$, then $\omega_\beta$ is not colored in $G_\beta'$,
and if $y\in \omega$, then $\omega_y$ is not colored in $G_y'$. So 
$\beta\notin \omega$ and $y\notin \omega$ which proves (2).

If $G'$ is partitioned by $\{\{a,e\},\{x,y\},\{\beta,b\}\}$, then $f(\sigma)$ and $f(\tau)$ are two
faces of the corresponding 3-cell of $|\cp\cq(G)|$. A deformation, rel the boundary of the disk 
$f(\sigma) \cup f(\tau)$, across this 3-cell finishes this case.

So assume that $G'$ is not partitioned by $\{\{a,e\},\{x,y\},\{\beta,b\}\}$, and let $\omega$ be any
cycle of $G'$ which is contained in $\{a,e,x,y,\beta,b\}$ and
which contains none of the three elements of $\{\{a,e\},\{x,y\},\{\beta,b\}\}$.
  Then the length of $\omega$ is either two or three, and so, by 
(1) and (2), $\omega$ equals either $\{b,x\}$, $\{a,b,x\}$ or $\{e,b,x\}$.  

First assume that 
$\omega=\{b,x\}$.  If $b$, $x$, $\beta$ and $y$ all meet the same vertex, then (1) implies that the 
other vertices of $b$, $\beta$ and $y$ are distinct.  But then of the 4 choose 2 ways that $e$ can 
be attached to the four vertices of $G'$, three are ruled out by (1) and the other three by (2).
It follows from this contradiction, and the fact that $\{\beta,y\}$ is not a cycle, that $x$, $y$ and $\beta$
form either a path or loop of length 3 in $G'$. If they form a loop, then it follows from (1) that the 
vertex of $G'$ not in this loop meets both $e$ and $a$.  And if they form a path and not a loop, 
then it follows from (1) and (2) that both $e$ and $a$ meet both endpoints of this path.  Either way,
$G'$ is partitioned by $\{\{e,a\}\{x,y,\beta\}\}$ and case (i) is completed (assuming $\omega= \{b,x\}$)
using two deformations of type (a) above through the two 3-cells of $T(G)$ shown in figure \ref{3tcells}(a).

\begin{figure}[tb]
\centerline{\mbox{\includegraphics*{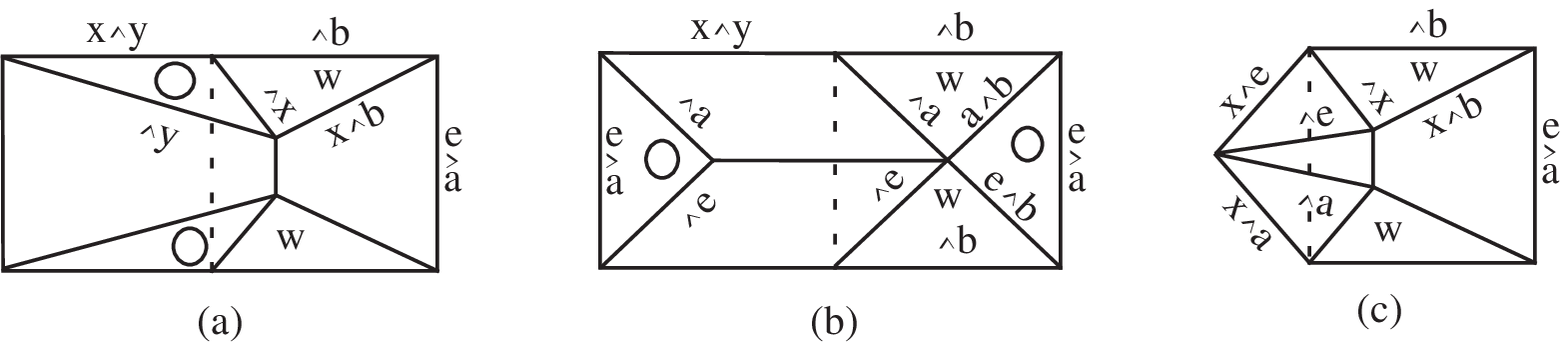}}}
\caption{}
\label{3tcells}
\end{figure}

Now assume $\omega=\{e,b,x\}$ (with a similar argument if $\omega=\{a,b,x\}$). By (1) and (2), both 
of the edges $y$ and $\beta$ meet the vertex of $G'$ which is not contained in the cycle $\omega$. 
Also 
$y$ and $\beta$ do not form a cycle of length two.  We consider three cases depending on which
of the three vertices of $\omega$ the edge $\beta$ meets.

First assume $\beta$ meets $x \cap e$.  Then $y$ meets $b\cap x$, because if $y$ met $b\cap e$ then
$\{y,\beta ,e\}$ would be a cycle of length three, contradicting (2). Now, by (1), either $a$ meets
both vertices of $e$, in which case a deformation as indicated in figure \ref{3tcells}(a) completes the
proof, or $a$ meets $y\cap \beta$ and $b\cap e$.  But the latter case cannot occur since it would
result in the cycle $\{a,b,y\}$ in contradiction to (2).

Now assume $\beta$ meets $b\cap e$. Then, by (2), $y$ meets $x\cap b$.  So, by (1), either $a$ meets both
vertices of $e$, or $a$ meets $\beta \cap y$ and $e\cap x$.  In the first case, use the deformation 
indicated by figure \ref{3tcells}(a), and in the second, use the one indicated by figure \ref{3tcells}(b).

Finally assume that $\beta$ meets $b\cap x$. Then, using (1) and (2), either $y$ meets $b \cap e$ and
$a$ meets both vertices of $e$, or $y$ meets $e\cap x$ in which case $a$ meets both vertices of $e$
or $a$ meets $b\cap e$ and $\beta \cap y$.  For all three possibilities, a deformation like that 
indicated by figure \ref{3tcells}(a) can be used to complete the proof.

\subsubsection*{Case (ii).} Assume first that $x\ne b$.  In case $\{b,\beta\}$ is a loop in $G'$, then, since
$G'$ is partitioned by  $\{\{a,e\},\{\beta,b\}\}$, the edges $e$ and $a$ both meet the vertex of $G'$
not in the loop $\{b,\beta\}$.  So does the edge $x$, since $\tau$ is a $w$-triangle corresponding
to the graph $G_\beta '$.  Thus there is a $w$-prism having $\sigma$ and $\tau$ as faces. This
3-cell of $T(G)$ can be used to complete the proof (assuming $\{b,\beta\}$ is a loop).

So assume $\{b,\beta\}$ is not a loop in $G'$.  Then $\{e,a\}$ is.  If $x$ meets both vertices of $e$, 
then, as just above, a single $w$-prism can be used to complete the proof.  If $x$ meets the vertex
$\beta \cap b$, and so also $e\cap b$, then a deformation as indicated in figure
\ref{3tcells}(c) completes the proof.

Now assume that $x=b$.  Then $\{b,\beta\}$ is not a loop in $G'$, so $\{e,a\}$ is, and the vertex 
$b\cap\beta$ is not contained in this loop.  Therefore there is a pillow in $T(G)$ having $\sigma$ and
$\tau$ as faces which can be used to complete the proof.

\subsubsection*{Case (iii).}
We will use the following.
\begin{enumerate}
\item[(1)] No cycle of length three in $G'$, except possibly for $\{a,e,b\}$, contains two edges from
$\{a,e,x\}$ and one from $\{\beta,b\}$.
\item[(2)] The only possible cycle of length two in $G'$  made up of edges from the set $\{a,e,x,\beta,b\}$
is $\{\beta,b\}$.
\end{enumerate}

To prove (1), assume $\omega$ is a cycle of length three with two edges in $\{a,e,x\}$ and one in
$\{\beta,b\}$.  If $\beta\in \omega$, then $\omega_\beta$ is not colored in $G'_\beta$. And if $x\in\omega$,
then $\omega_x$ is not colored in $G_x'$.  So $\beta\notin\omega$ and $x\notin \omega$.

For (2), let $\omega$ be a cycle of length two and assume that $\omega\subset\{a,e,x,\beta,b\}$ 
and $\omega\ne\{\beta,b\}$.  If $\beta \in \omega$, then $\omega_\beta$ is not colored in $G_\beta'$.
If $b\in \omega$, then $\omega_x$ is not colored in $G_x'$.  And if $\{b,\beta\} \cap \omega =\varnothing$,
then $\omega \subsetneq\{a,e,x\}$, and so is not colored in $G_\beta'$.  
These are all contradictions, so (2) is proved.

If $G'$ is partitioned by $\{\{a,e,x\},\{\beta,b\}\}$, then a deformation across the 
corresponding O-prism of 
$T(G)$ which contains $\sigma$ and $\tau$ as faces completes the proof. So assume $G'$ is not partitioned 
by $\{\{a,e,x\},\{\beta,b\}\}=\{h_i\}$, where, as usual, singleton $h_i$ are not shown in the set on the 
left of the equals sign.
Let $\omega$ be a loop in $G'$ which contains none of the
$h_i$. Then $\omega$ has length two or three, and so it follows from (1) and (2) that $\omega=\{a,e,b\}$.
Now (1) and (2) imply that both $\beta$ and $x$ meet the vertex of $G'$ not in $\omega$.
If $\beta$ meets $a\cap e$, then each of the three ways of attaching $x$ contradicts either (1) or (2).
So $\beta$ meets one vertex of $b$ and, by (1) and (2) again, $x$ meets the other.  Thus $G'$ is
partitioned by $\{a,e,x,\beta\}$, and figure \ref{3cells2}(d) shows how to complete the proof.

\begin{figure}[tb]
\centerline{\mbox{\includegraphics*{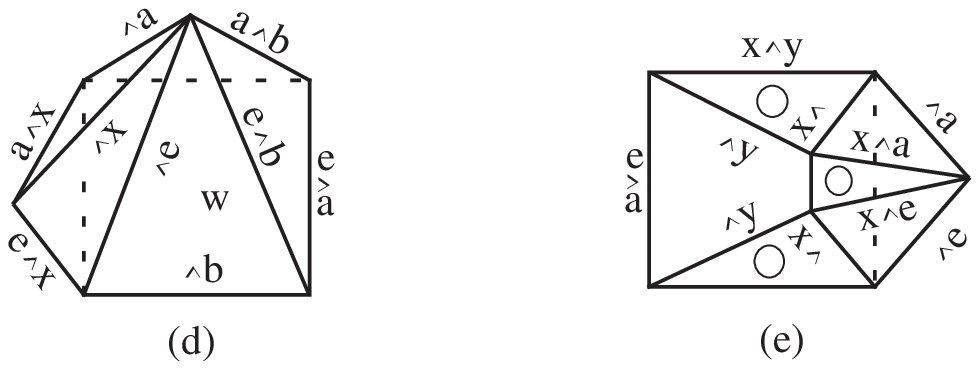}}}
\caption{}
\label{3cells2}
\end{figure}

\subsubsection*{Case (iv).} If $G'$ is partitioned by $\{\{x,y\},\{a,e,\beta\}\}$, then use the 
corresponding O-prism of 
$T(G)$ with $\sigma$ and $\tau$ as faces.  If not, then, just as in case (iii), $\{a,e,x\}$ and 
$\{x,y,\beta\}$ are loops in $G'$. So $G'$ is partitioned by $\{\{x,y,\beta\},\{a,e\}\}$ and figure
\ref{3cells2}(e) indicates how to complete the proof.

\subsubsection*{Case (v).} Since $G'$ is partitioned by $\{a,e,\beta\}$, the edges $a$, $e$ and $\beta$ 
form a loop in $G'$.  Because of the 2-cell $\tau$, the edge $x$ meets $a\cap e$ and one vertex of 
$\beta$. If $x$ meets $\beta\cap a$, then $G'$ is partitioned by $\{\{a,x\},\{\beta,e\}\}$ and the
pillow with $\sigma$, $\tau$ and the corresponding 2-cell as faces can be used to complete the proof.
If $x$ meets $\beta \cap e$, use the pillow with faces $\sigma$, $\tau$ and the 2-cell corresponding 
to $G'$ partitioned by $\{\{e,x\},\{\beta,a\}\}$.

\subsubsection*{Case (vi).} Let $\omega$ be a loop in $G'$ consisting of a subset of $\{a,e,x,\beta\}$.
Then $\{a,e,\beta\}\subset\omega_x$ since 
$G_x'$ is partitioned by $\{a,e,\beta\}$.  Similarly, $\{a,e,x\}\subset\omega_\beta$. Therefore $\omega=
\{a,e,x,\beta\}$ and so $G'$ is partitioned by $\{a,e,x,\beta\}$. Use the corresponding 3-cell of
$T(G)$ to complete the proof of this case and the lemma.
\end{proof}

\begin{lemma}\label{2}
Assume that $f(\sigma)$ and $f(\tau)$ are distinct 2-cells of  $T(G)$.  If $\sigma\cap\tau$ consists of
a single edge which is $\beta$-parallel for both $\sigma$ and $\tau$, then the same conclusions as in
the previous lemma hold, except that (3) of that lemma is replaced by: 
\begin{enumerate}
\item[($3'$)]
The $\beta$-halfedges of $\cc'$ in the closed disk $\sigma\cup\tau$ are the same as those of $\cc$
and are contained in two $w$-triangles of $\cc'$.
\end{enumerate}
\end{lemma}

\begin{proof}
We consider all possible configurations for $\sigma$ and $\tau$ in the following two cases.   
Define the graph $G'$ as at the start of the 
proof of lemma \ref{1}.
\subsubsection*{Case (i).} Assume $\sigma$ and $\tau$ are as in figure \ref{bp1}(i), but with $y=\beta$.
If $\{a,e\}$ is a loop in $G'$, then the edges $\beta$, $b$, and $x$ all meet the 
vertex of $G'$ not in the loop $\{a,e\}$. So there is a $w$-prism with with $f(\sigma)$ and $f(\tau)$
as faces. If $\{a,e\}$ is not a loop, then both $\{\beta,x\}$ and $\{\beta,b\}$ are, and so again 
$f(\sigma)$ and $f(\tau)$ are faces of some $w$-prism in $T(G)$. Deformations across these $w$-prisms
complete the proof in both cases.

\subsubsection*{Case (ii).} Assume $\sigma$ and $\tau$ are as in figure \ref{bp1}(iv), but with $y=\beta$.
Then, since $G'$ is partitioned by $\{a,e,\beta\}$,
the edges $a$, $e$ and $\beta$ form a loop in $G'$.
So $\{a,e\}$ is not a loop in $G'$, which implies that $\{x,\beta\}$ is, since $G'$ is also partitioned by
$\{\{a,e\},\{x,\beta\}\}$.  Use the pillow with $f(\sigma)$ and $f(\tau)$  as faces to complete the proof of 
this case and the lemma.
\end{proof}

Let $h$ be any halfedge of $(S^2,\cc)$. Say $E_G(h)=\beta$.  Let $\{\sigma_i\}$ be the set of 2-cells of $\cc$
which contain a $\beta$-halfedge.  For each $\sigma_i$, let $m_i$ and $m_i'$ be the midpoints
of the two edges of $\sigma_i$ which contain $\beta$-halfedges. Also let $a_i$ be an arc which joins
$m_i$ and $m_i'$ and which is otherwise contained in the interior of $\sigma_i$.  The union
of the $a_i$ form a collection  of disjoint circles in $S^2$, called $\beta$-circles.

\begin{lemma}\label{3}
Assume $c$ is a $\beta$-circle in $S^2$ contained in $\bigcup_{i=1}^n \sigma_i$ where each
$\sigma_i$ is a closed $w$-triangle of $\cc$ meeting $c$.  Then $f$ can be deformed to either 
a cell-to-cell or to a constant map $f'$ with the effect that 
\begin{enumerate}
\item $c$ is eliminated and no new $\beta$-circles are introduced and
\item if $E_G$ (resp.~ $E_G'$) is the edge function associated with $f$ (resp.~$f'$), 
the image of $E_G'$ is contained in the image of $E_G$.
\end{enumerate}
\end{lemma}

\begin{proof}
Assume $c\subset \bigcup_{i=1}^n \sigma_i$ as in the statement of the lemma.  We say that two of
the $\sigma_i$ are adjacent if they share at least one edge.  If $f$ takes two adjacent 
$\sigma_i$ to the same cell of $T(G)$, then $n$ can be decreased by two by using a deformation
of either type (b) of (d) above. If adjacent $\sigma_i$ map to different cells, then $n$ can be 
decreased by at least one with a type (a) deformation across the $w$-tetrahedron of $T(G)$
having the adjacent $\sigma_i$ as faces.
\end{proof}

\begin{proof}[Proof of lemma \ref{2conn}.]
Let $\beta$ be any element of the image of the edge function associated with the cell-to-cell map 
$f\colon (S^2,\cc)\to T(G)$, and let
$\{c_i\}$ be the set of $\beta$-circles in $S^2$.  Of the two disks in $S^2$ which each $c_j$ bounds, let $D_j$
be the one which contains all the $\beta$-parallel edges for the 2-simplices meeting $c_j$.
Let $D'$ be one of the $D_j$.  Among all the $D_j$ contained in $D'$, let $D$ be an innermost one.
If none of the $c_i$ are contained in the interior of $D$, we say that $D$ is strictly innermost.

First assume that $D$ is strictly innermost. Let $c$ be the $\beta$-circle which $D$ bounds.
Consider the set $\cb$ of $\beta$-parallel edges which are contained in 2-cells meeting $c$.
We assume, for now, that $\cb\ne \varnothing$.
If some $e\in \cb$ is contained in the intersection, $\sigma \cap \tau$ say, of two 2-cells of $\cc$
and $e$ is $\beta$-parallel for $\sigma$ but not
$\tau$, then an application of lemma \ref{1} decreases by one the total number of 2-cells contained in 
$D$.  Thus, since $D$ is strictly innermost,
 we may assume, by repeated applications of lemma \ref{1}, that, for each $e\in \cb$,
both 2-cells containing $e$ meet the $\beta$-circle $c$. 

So now, if $e\in \cb$, either lemma \ref{2} or one of the 
operations (b), (c) or (d) above can be applied to the 2-cells, $\sigma$ and $\tau$ say,
containing $e$. 
Applying the appropriate one decreases the cardinality of $\cb$ by one by eliminating $e$ from it.
 If lemma \ref{2} is applied, the $\beta$-circle $c$ is surgered into two $\beta$-circles
and the portion of these in the disk $\sigma\cup\tau$ is contained in two $w$-triangles. 
If (b) or (c) is applied, then $c$ is again transformed into two $\beta$-circles.
If (c) has been applied, these $\beta$-circles are in two different 2-spheres.
(The edges of the $CW$-structures of these 2-spheres are identified with those of $\cc$ in the natural way.)

By repeating the above procedure, once for each of the remaining edges in $\cb$, we may assume that $\cb
=\varnothing$. Then either $c$ has been eliminated or $c$ has been transformed
 into some number of $\beta$-circles.  In the
latter case, each of the $\beta$-circles is contained in a union of $w$-cells and so can be eliminated by lemma
\ref{3}.  This shows that  if $D$ is strictly innermost, then the number of $\beta$-circles can be decreased by
one, and the  image of the edge function is, if changed, made smaller. 
(In case operation (c) is used, then $f'$ has been 
transformed into some number of cell-to-cell maps $f_i\colon S^2\to T(G)$. Then ``number of $\beta$-circles"
means ``the sum of the number of $\beta$-circles in the domains of the $f_i$'' and ``image of the 
edge function'' means ``the union of the images of the edge functions of the $f_i$''.)

Assume now that $D$ is innermost, but not strictly so, and let $R$ denote the component of $S^2-
\bigcup c_i$ in $D$ having $c$ as limit points.  Define the nonempty set $\cb$ as above and, as above,
use lemma \ref{1} to decrease the number of 2-cells in $R$ until some edge of $\cb$ is contained in a
2-cell which meets a $\beta$-circle different from $c$. Then apply 
either operation (b) or lemma \ref{2} to the two 2-cells
containing this edge.  This process joins the two $\beta$-circles into one  and so again decreases
the number of $\beta$-circles by one without increasing the image of the edge function.

Thus all $\beta$-circles can be eliminated in such a way that the image of the edge 
function becomes a proper subset of what it was originally.
This completes the proof since $\beta$ was an arbitrary element of the image of the edge function.
\end{proof}

\section{The simplicial complex $\ct$}\label{SCT}

Let $\cs_T$ denote the sphere complex $S{\mathbb A}_n$ defined in $\S 2.4$ of 
\cite{HV}
together with, for each vertex $(G,\rho)$ of $S{\mathbb A}_n$, 
some choice of a maximal tree in $G$. 
In more detail, each vertex of the sphere complex $S{\mathbb A}_n$, which we denote by 
$\cs_n$, or just $\cs$, 
 is an equivalence class of
pointed marked graphs of rank $n$.  
Usually the adjective pointed is dropped in this context.
A maximal tree of any one representative of an equivalence class
serves to specify a maximal tree in every representative.  Such a compatible choice of maximal trees
is what is meant by a choice of a maximal tree in $G$.

We point out that there is a simplicial map ${\mathcal F}\colon \xt_n\to \cs_n$ which forgets partitions and
colorings.  It will be used in \S\ref{defP}.

If $v=(G,\rho)$ is a vertex of $\cs_T$ which is not a rose,
let $p(v)$ denote the edge of $\cs$ which joins $v$ to the marked rose $G/T$, where $T$ is the chosen 
maximal tree in $G$.  With the edges of $\cs$ oriented so that each points towards the vertex 
corresponding to the graph with fewer vertices, we write $term\ p(v)$ for the terminal vertex of 
$p(v)$.  If $v$ is a rose, then $p(v)$ is, by definition, just the vertex $v$, and we say, 
for convenience, that 
$p(v)$ terminates at $v$ and set $term\ p(v)=v$.  Let $\cc$ be the subcomplex of $\cs$ consisting of 
those simplices $\{v_0,\ldots,v_n\}$ with $term\ p(v_0)=\ldots =term\ p(v_n)$.

\begin{lemma}\label{con}
Each component of $\cc$ is contractible and contains exactly one marked rose.
\end{lemma}

\begin{proof}
Let $C_i$ be a component of $C$, and let $\{v_j\}$ be the set of vertices of $C_i$.  Since $C_i$ is 
connected all the $p(v_j)$ terminate at the same marked rose, $r$ say. 
Therefore $C_i$ is a cone on $r$ and so is contractible.
\end{proof}

Let $\ct$ be the simplicial complex obtained from $\cs$ by contracting the components of $\cc$ to 
distinct points.  Since $\cc$ contains all the vertices of $\cs$, $\ct$ has one vertex for each 
component $C_i$ of $C$, and $\{C_0,\ldots,C_n\}$ is a simplex of $\ct$ if and only if there is a simplex
$\{v_0,\ldots,v_n\}$ of $\cs$ with $v_i\in C_i$ for $i=0,\ldots,n$.  Equivalently, the vertices of $\ct$ correspond to 
the roses of $\cs$, and $\{R_0,\ldots,R_n\}$ is a simplex of $\ct$ if and only if there is a simplex
$\{v_0,\ldots,v_n\}$ of $\cs$ with $term\ p(v_i)=R_i$ for $i=0,\ldots   ,n$.  The simplicial quotient
map $q\colon \cs\to \ct$,  a homotopy equivalence by lemma \ref{con}, can be defined by taking 
each vertex $v$ of $\cs$ to the vertex $term\ p(v)$ of $\ct$.  A simplex $\tilde{\sigma}$ of 
$\cs$ is called a lift of the simplex $\sigma$ of $\ct$ if $q(\tilde{\sigma})=\sigma$ and if the 
dimensions of $\sigma$ and $\tilde{\sigma}$ are the same.  
Prop.\ \ref{prop} below shows that among all the graphs corresponding to the lifts of a given 
simplex, the one with the fewest number of vertices is unique.  Moreover, any such graph 
blows down to this unique lift.


Let $\sigma=\{R_0,\ldots,R_n\}$ be an $n$-simplex of $\ct$, and let $\tilde{\sigma}=G_n\to G_{n-1}\to 
\ldots\to G_0$ be a lift of $\sigma$ to $\cs$. (So the $R_i$ are marked roses, the $G_i$ are
marked graphs with markings compatible with the blowdowns $G_i\to G_{i-1}$, and
we may assume that
$q(G_i)=R_i$.)  For $i=0,\ldots,n-1$, let 
$T_i$ be the maximal tree of $G_n$ which is blowndown by the composite $G_n\to\ldots\to G_i\buildrel
{q}\over{\to} R_i$. Also let $T_n$ be the maximal tree of $G_n$ associated with the vertex $G_n$ of $\cs_T$.
Let $G_{\tilde{\sigma}}$ be the marked, pointed graph obtained from $G_n$ by collapsing
the edges of $\bigcap_{i=0}^nT_i$.

\begin{prop}\label{prop}
The vertex of the sphere complex $\cs$ corresponding to the marked graph $G_{\tilde{\sigma}}$ 
does not depend on the choice of
lift $\tilde{\sigma}$ of $\sigma$.
\end{prop}

The following definitions and conventions will be used in the proof.
Let $G$ be a graph with edges oriented and labeled.  
If $H$ is any subgraph of $G$, let $L(H)$ denote the set of labels of the edges of $H$,
and let $L(H)^{\pm}=L(H) \cup \{e^{-1}:e\in L(H)\}$. 
 
If $x\in L(H)$, then $init(x)$ and $term(x)$ denote the initial and terminal vertices,
respectively, of the 
edge $x$ labels.  The same edge as $x$, but with the opposite orientation,
is denoted by $x^{-1}$.
Thus $init(x)=term(x^{-1})$ and $term(x)=init(x^{-1})$.
An oriented edgepath in $G$ is 
specified by a word in the elements of $L(G)^{\pm}$ in the obvious way. 

Assume now that $G$ is
marked by $g\colon R_0\to G$ where the edges of $R_0$ are labeled by elements of the set $\{r_i\}$.
We often write
$g_E(r_i)$ for the (reduced) word corresponding to the edgeloop $g(r_i)$. 
Here we have identified a label of an edge with the edge itself.  We often do this.     
If $e$ is an edge of $G$, let $e\cdot g(r_i)$ be the number of occurrences of $e$ in the 
(reduced) word $g(r_i)$ plus the number of occurrences of $e^{-1}$ in the same.  
Let $e\cdot im(g)=\sum_i e\cdot g(r_i)$.  Also, if $e_i$ and $e_j$ are edges of $G$, let 
$e_ie_j\cdot im\ g=\sum_ke_i e_j\cdot g(r_k)$ where $e_i e_j\cdot g(r_k)$ is the number of occurrences
of the subword $e_i e_j$ in $g(r_k)$ plus the number of occurrences of $e_j^{-1}e_i^{-1}$ in the same.

\begin{proof}[Proof of \ref{prop}]
We first assume that $\sigma=\{R_1,R_2\}$ is a 1-simplex of $\ct$ with lifts $\sigt=G_1\to G_0$ 
and $\sigt'=G_1'\to G_0'$ in $\cs$. Let $G$ denote the marked graph $G_{\sigt}$ and let 
$S$ and $T$ be the maximal trees of $G$ which are the image under $G_1\to G$ of the maximal trees
$T_0$ and $T_1$ of $G_1$.  Define $G'$, $S'$ and $T'$ similarly.  Then there is a homotopy 
commutative diagram 

$$\xymatrix{&G \ar[dl]_S \ar[dr]^T \\
          R_2 & R_0\ar[u]_f \ar[d]^{f'} \ar[l]_h \ar[r]^g &R_1 \\
          &G' \ar[ul]^{S'} \ar[ur]_{T'}
}$$


\noindent
where $f$ and $f'$ are the markings of $G$ and $G'$, and where the arrows labeled $S$, $T$, $S'$ and $T'$
are blowdowns of the maximal tree of the label.

By viewing the diagram as four edges in $\cs$, and then acting on them with a suitable element of $\autfn$,
we may assume that $h$ is homotopic to a homeomorphism.  We often use this homeomorphism to identify
the edges of $R_2$ with those of $R_0$.

We say that an edge $e$ of $R_1$ is paired with  the edge $e'$ of $R_2$ if $e\cdot g(r_i)= e'\cdot h(r_i)$
for every edge $r_i$ of $R_0$.

\begin{lemma}\label{paired}
An edge $e$ of $R_1$ is paired with the edge $e'$ of $R_2$ if and only if $S(b)=e'$ and $T(b)=e$ 
for some edge $b$ of $G$.
\end{lemma}

\begin{proof}
Choose labels for the edges of $G$ and 
let $B=L(G)-L(S\cup T)$.  Label the edges of $R_1$ and $R_2$ so that $S(e)=e$ if $e$ is an edge
of $G$ not in the tree $T$, and  $T(e)=e$ if $e$ is not in the tree $S$. Then 
$$ L(R_2)=L(T)\cup B \ \ \  \text {and }\ \  L(R_1)=L(S)\cup B. $$  If $x\in L(R_2)\approx L(R_0)$, 
then $x\cdot im(h)=1$
since  $h$ is homotopic to a homeomorphism.  If $s\in L(S)\subset L(R_1)$, then $s\cdot im(g) \ge 2$ by 
the following lemma.  So an edge $e$ of $R_1$ can be paired with an edge $e'$ of $R_2$ only if $e\in B\subset
L(R_1)$.  Since $h$ is homotopic to a homeomorphism, it is easy to see that each $e\in B\subset
L(R_1)$ is paired with the unique edge of $B\subset L(R_2)$ with the same label as $e$.

The implication in the other direction is easy.
\end{proof}

\begin{lemma} \label{branch}
If $s_1\ldots s_n$ is an edgepath in the tree $S$ terminating at the  basepoint of $G$, then 
$$2\le s_1\cdot im\ g < s_2\cdot im\ g< \ldots <s_n\cdot im\ g.$$
\end{lemma}

\begin{proof}
Label the edges of $G$ and the $R_i$ as in the proof of lemma \ref{paired}.
If $t\in L(T)^{\pm}\subset L(R_2)^\pm \approx L(R_0)^\pm $, 
then, since $h$ is homotopic to a homeomorphism, $g_E(t)$ is a word in the letters 
$L(S)^{\pm}$, and is obtained from $f_E(t)$ by deleting the unique occurrence of $t$ in $f_E(t)$.  Also
if $b\in B$, then $g_E(b)=f_E(b)$.  So the lemma follows since each vertex of $G$ (other than the basepoint)
has valence $\ge 3$.
\end{proof}

\begin{lemma} \label{int}
If $ss'\cdot im\ g > 0$ where $s,s'\in L(S)^{\pm}$, then
\begin{align}
ss'\cdot im\ g=1 &\Leftrightarrow \text{there is a $t\in L(T)^\pm$ such that $sts'$ is an 
edge path in $G$, and}\notag \\
ss'\cdot im\ g\ge 2 &\Leftrightarrow \text{$ss'$ is an edgepath in $G$ where $init(s')=term(s) \ne *$.} \notag
\end{align}
\end{lemma}

\begin{proof}
If, in $G$, $term(s)\ne init(s')$, then both sides of the first 
of the two
equivalences in the statement of the lemma hold.  Otherwise, the same is true of the second equivalence
by \ref{branch} and a little more.
\end{proof}

Returning now to the proof of \ref{prop},
label the edges of $R_1$ and $R_2$ so that two edges have the same label if and only if they
are paired. Then label the edges of $G$ so that an edge $e$ of $G$ has the same label as an edge 
$r$ of $R_i$ if the blowdown $G\to R_i$ takes $e$ to $r$. Such a labeling is possible by lemma \ref{paired}.
Next, orient the edges of $G$ and the $R_i$
so that the blowdowns  $S$ and $T$ preserve orientations. Label and orient the 
edges of $G'$ similarly.  Let $B=L(R_1)\cap L(R_2)$.  Then we have $$L(G)=L(R_0)\cup L(R_1)=L(B)
\sqcup \left(L(R_1)-B\right) \sqcup \left(L(R_2)-B\right)=L(G').$$
  Thus a bijection from the edges of $G$ to those of $G'$
 is defined by taking an edge $e$ of $G$ to the edge
of $G'$ with the same label as
$e$.  We will show that this respects incidence relations and so defines a homeomorphism
$G\to G'$.

As a first step we show that $f_E$ is determined by $g_E$.  Let $t\in L(T)^\pm$.  Then $g_E(t)=
s_1\ldots s_n$ for some $s_i$ in $L(S)^\pm$.  Assume, for now, that $n>1$.  If
 $s_i s_{i+1} \cdot im\ g=1$ for some consequtive pair of letters $s_i$ and $s_{i+1}$ in $g_E(t)$,
then, by \ref{int} and the fact that $T$ has no cycles,
 $f_E(t)=s_1\ldots s_i\, t\, s_{i+1}\ldots s_n$. 
Otherwise, 
$f_E(t)$ equals either $s_1\ldots s_n\,t$ 
or $t\,s_1\ldots s_n$.  Which, is determined by $g$, according to \ref{branch}.

If $n=1$ so that $g_E(t)=s_1$, then $s_1\cdot im\ g\ge 2$ by \ref{branch}.
Thus there is an $e\in
L(R_0)=L(R_2)$ such that $s_1\cdot g(e) > 0$ and $e^{\pm 1} \ne t^{\pm1}$.  Assume $g_E(e)=s_1^{\pm1}$. If
$e\in L(T)^\pm$,  then the edges $e$ and $t$ would form a cycle in $T$. 
And if $e\in L(B)$, then $g_E(e)$  would
contain the letter $e$, another contradiction.  So the word $g_E(e)$ contains at least two letters, 
with either the first or last being $s_1$ or $s_1^{-1}$.  If the first is $s_1$ or the last is $s_1^{-1}$,
then $f_E(t) = s_1\,t$.  Otherwise, $f_E(t)=t\,s_1$.  Since $f_E(b)=g_E(b)$ if $b\in L(B)$,
we have shown how to compute $f_E$ from $g_E$.  Similar arguments apply just as well to $G'$ and 
$f_E'$, so we conclude that $f_E=f_E'$.

Now, for convenience, assume that the edges of the tree $S$ are oriented so that they point toward
the basepoint of $G$.  Let $x\in L(G)^\pm=L(G')^\pm$,  and let $y\in L(S)$.  Then
\begin{align}
\text{$*\ne term(x)=init(y)$ in $G$}  &\Leftrightarrow xy\cdot im\ f_E >0 \notag\\
                                      &\Leftrightarrow xy\cdot im\ f_E' >0 
\Leftrightarrow \text{$*\ne term(x)=init(y)$ in $G'$.} \notag 
\end{align}
For edges meeting the basepoint of $G$ and $G'$, we have
\begin{align}
\text{$init(x)=*$ in $G$} &\Leftrightarrow \text{$f_E(z)$ begins with $x$ for some $z\in L(R_0)^\pm$} \notag\\
                          &\Leftrightarrow \text{$f_E'(z)$ begins with $x$ for some $z\in L(R_0)^\pm$}
\Leftrightarrow \text{$init(x)=*$ in $G'$.} \notag
\end{align}
It follows that a homeomorphism $h\colon G\to G'$ is defined by identifying edges of $G$ and $G'$ 
with the same labels in an orientation preserving manner.  Since $f_E=f_E'$, we have $h\circ f=f'$.
Therefore $G$ and $G'$ are equivalent as marked graphs.

Now let $\sigma=\{R_1,\ldots,R_{n+1}\}$ be an $n$-simplex of $\ct$ with lift $\sigt=G_1\to\ldots \to G_{n+1}$
in $\cs$. Let $G=G_{\sigt}$ be the marked graph with maximal trees $T_i$ as defined above.  Also let 
$T_i$, for $i=1,\ldots ,n+1$, denote the image of $T_i$ under the blowdown $G_0\to G_{\sigt}$.  Thus the 
$T_i$ are maximal trees in $G$ with $\bigcap_i T_i$ containing no edges.  Each edge $G_i\to G_j$ 
of $\sigt$ gives rise to a commutative diagram
$$\xymatrix{ &G \ar[dl]_{T_i} \ar[d] \ar[dr]^{T_j} \\
            R_i& G/\,T_i\cap T_j \ar[l] \ar[r] &R_j \\
            &R_0 \ar[ul]^{g_i} \ar[u]_{f_{ij}} \ar[ur]_{g_j}
}$$
where all the arrows are blowdowns except for $g_i$, $g_j$, and $f_{ij}$ which are markings.  
These markings are induced by the marking $f\colon R_0\to G$ of $G$ which is not shown in the diagram.

Label the edges of $G$ and the $R_i$ as above (in the case where $\sigma$ was a 1-simplex)
so that an edge $e$ of $G$ has the same label as an edge $e'$ of $R_i$ if and only if
$T_i$ sends $e$ to $e'$.  Note that whether an edge of $R_i$ has the same label as an edge of $R_j$
can be determined from just $g_i$ and $g_j$ by precomposing $g_i$ with a suitable homotopy equivalence
of $R_0$ and then using \ref{paired}.
Also orient the 
edges of $G$ and the $R_i$ so that the maps $T_i$ preserve the orientations of the edges.  
  We will show that the function 
$f_E$ from $L(R_0)^\pm$ to the set of words in the elements of $L(G)^\pm$ is determined by the $g_i$.  
Since $L(G)=\bigcup_i
L(R_i)$, the proof can then be completed as above assuming, as we may, that 
 $g_0\colon R_0\to R_0$ is homotopic to a homeomorphism, and then using the tree $T_1$ of $G$
in place of the tree $S$ used in the above argument.

Since the edge $G_1\to G_i$ of $\sigt$ is a lift of the 1-simplex $\{R_1,R_i\}$ of $\ct$,
the functions $f_i^E\equiv f_{0i}^E$ (where $E$ has been written as a superscript rather than 
a subscript) are determined by the $g_i$.  Assume $f_i^E(r)=S_i r S_i'$ where $r\in L(R_0)$, and where
the $S_i$ and $S_i'$ are words in the letters $L(T_1-T_i)^\pm\subset L(T_1)^\pm$.  Since $\bigcup_iL(T_1-
T_i)=L(T_1§)$, the union of the letters in the words $S_i$ form an edgepath in the tree $T_1$ of $G$ from $*$
to $init(r)$.  Since $x\cdot im\ f= x\cdot im\ g_i$ if $x\in L(T_1-T_i)\subset L(G)$, the order of the 
letters in this edgepath can be determined using the $g_i$ and lemma \ref{branch}.  A similar argument 
using the $S_i'$ shows that the $g_i$ determine $f_E(r)$ and so also $f_E$.
\end{proof}

\subsection{The $CW$-complex $\ct^*$} \label{ct}
Because of the proposition just proved, we henceforth denote $G_{\tilde{\sigma}}$ by just $G_{\sigma}$.  Note
that the rule $\sigma \mapsto G_\sigma$ is functorial in the sense that it induces a simplicial map from the 
barycentric subdivision of $\ct$ to $\cs$.

If $G$ is a marked, pointed graph, let $T_M(G)$ denote the $CW$-complex $T(G)$ defined in section \ref{defTG},
but with, in addition, the graphs 
associated with the cells of $T_M(G)$ marked as blowdowns of $G$.  
The basepoint of $G$ does not figure into the definition of $T(G)$, but is needed for the definition
of $T_M(G)$ since markings necessarily preserve basepoints.
The map $T_M(G)\to T(G)$ which forgets markings is a cellular isomorphism,
and  $T_M$ becomes a functor if 
blowdowns $G\to H$ of marked, pointed  graphs are sent to the natural inclusions $T_M(H) \hookrightarrow T_M(G)$.
Composition of functors gives the functor $\sigma \mapsto T_M(G_\sigma)$ from the simplices of 
$\ct$ to $CW$-complexes.

We next use this last functor to define a cellular subdivision $\ct^*$ of the 2-skeleton of $\ct$.
(Only a subdivision of a finite subcomplex of $\ct$ will be needed in applications.)  Let $e$ be any
edge of $\ct$, and let $v_1$ and $v_2$ be the vertices of $e$.  Each $v_i$ corresponds to a marked rose
which is a blowdown of $G_e$, and so also corresponds to a vertex $v_i'$ of $T_M(G_e)$.  Choose a 
subdivision $e^*$ of $e$ so that there is a cellular isomorphism which takes $v_i$ to $v_i'$ and which
takes $e^*$ onto a shortest edgepath in $T_M(G_e)$ joining $v_1'$ with $v_2'$.  For example, if 
$G_e$ has just two vertices, then $e^*=e$. Subdivide each edge of $\ct$ in this way to obtain the 
1-skeleton of $\ct^*$. 
Each 1-cell of $\ct^*$  corresponds
to a marked partitioned graph with two vertices.

Let $\tau$ be any 2-simplex of $\ct$, and let $e_1$, $e_2$, and $e_3$ be the three faces of $\tau$.
Define a map $f\colon \partial \tau\to T_M(G_\tau)$ by piecing together the three composites
$e_i\to T_M(G_{e_i})\to T_M(G_\tau)$.  The first map in each of these was chosen above, and the 
second is provided by the functor $\sigma\mapsto T_M(G_\sigma)$.  If $f(\partial \tau)$ is a contractible
subspace of $T_M(G_\tau )$, then $\tau$ is called a degenerate simplex and is not subdivided.
Otherwise, let $\tau^*$ be a cellular subdivision of $\tau $ such that $f$ extends to a cell-to-cell
map $\tau^*\to T_M(G_\tau)$.  Subdivide each 2-simplex of $\ct$ in this way and let $\ct^*$ be the 
resulting $CW$-complex.  Also fix choices of the above cell-to-cell maps.  Then the 2-skeleton
of $\ct^*$ comes equipped with a map $G$ to the space $\bigsqcup_{\sigma\in\ct}T_M(G_\sigma)
/\!\sim$ where $x\sim y$ if there is a simplex $\tau$ of $\ct$ with a face $\sigma$ such that
the inclusion 
$T_M(G_\sigma)\hookrightarrow T_M(G_\tau)$ takes one of $x$ or $y$ to the other.  Degenerate simplices 
of $\ct$ are mapped by $G$ to the image of their boundaries.

\subsection{Standard $w$-crossings}\label{SWC}
If an element $w$ of the group $StN$ is either $w_{ij}^L$ or $w_{ij}^R$,
as defined in \S\ref{remarks}, then, according to \S1.4 of 
\cite{KN}, the conjugate
$w g^{-1}\, w^{-1}$ is a defining generator of the group $StN$, provided $g$ is.  
The corresponding relation in the
group $StN$ can be realized  (as consequences of defining relations for $StN$) in $\xt$ as the union of 
2-cells.  We will soon make specific choices for these 2-cells and call their 
union a standard $w$-crossing.  Each standard $w$-crossing is  
topologically a disk with boundary a loop corresponding to the relation the standard $w$-crossing 
represents.

Figure \ref{A1}(i) shows the standard $w$-crossing correponding to the relation 
$w_{ij}^L R_{ik} (w_{ij}^L)^{-1}=L_{jk}$, where any consistant markings of the graphs in the figure 
can be used. Of course, different markings correspond to different conjugates of the relation.
We mention that figure \ref{A1}(i) is essentially the pictorial version of the first computation 
in the proof of 1.3 in \cite{KN}.  Figure \ref{A1}(ii) is an abbreviation for the standard $w$-crossing
in figure \ref{A1}(i), and figure \ref{A1}(iii) is an even briefer abbreviation for the same.
Changing the orientation of the $i$ edgesets in \ref{A1}(i) gives the standard $w$-crossing in
\ref{A1}(iv), shown in abbreviated form.  Figure \ref{A1}(v) shows (the abbreviated form of) another 
standard $w$-crossing obtained by changing the orientations of the $j$ edgesets in \ref{A1}(i) and
then adding two triples of 2-cells representing the relation $w_{ij}^L = w_{ij}^R$.  These triples
are indicated by asterisks in \ref{A1}(v).  (Figure \ref{A1}(iv) represents the second computation
in the proof of 1.3 in \cite{KN}, and figure \ref{A1}(v) corresponds to a relation covered by 1.4,
but not 1.3, in \cite{KN}.

\begin{figure}[tb]
\centerline{\mbox{\includegraphics*{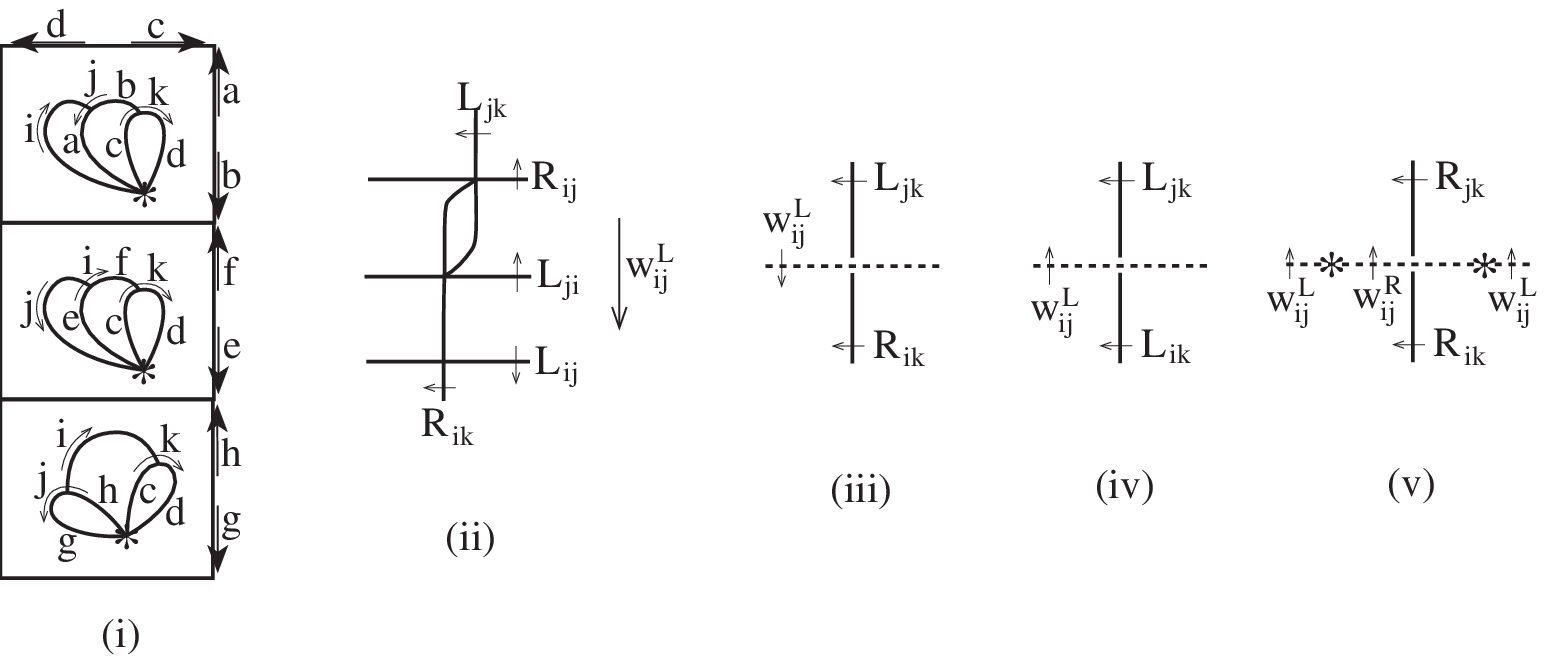}}}
\caption{}
\label{A1}
\end{figure}

Figure \ref{B1}(i) shows a second type of standard $w$-crossing.  It is superimposed on the faint 
outlines of the six corresponding 2-cells of $\xt$. The colored graphs corresponding to these 2-cells
should be clear.  Any consistant markings may be used.  Figure \ref{B1}(ii) is an abbreviation 
which we will use for the standard $w$-crossing of figure \ref{B1}(i).

\begin{figure}[tb]
\centerline{\mbox{\includegraphics*{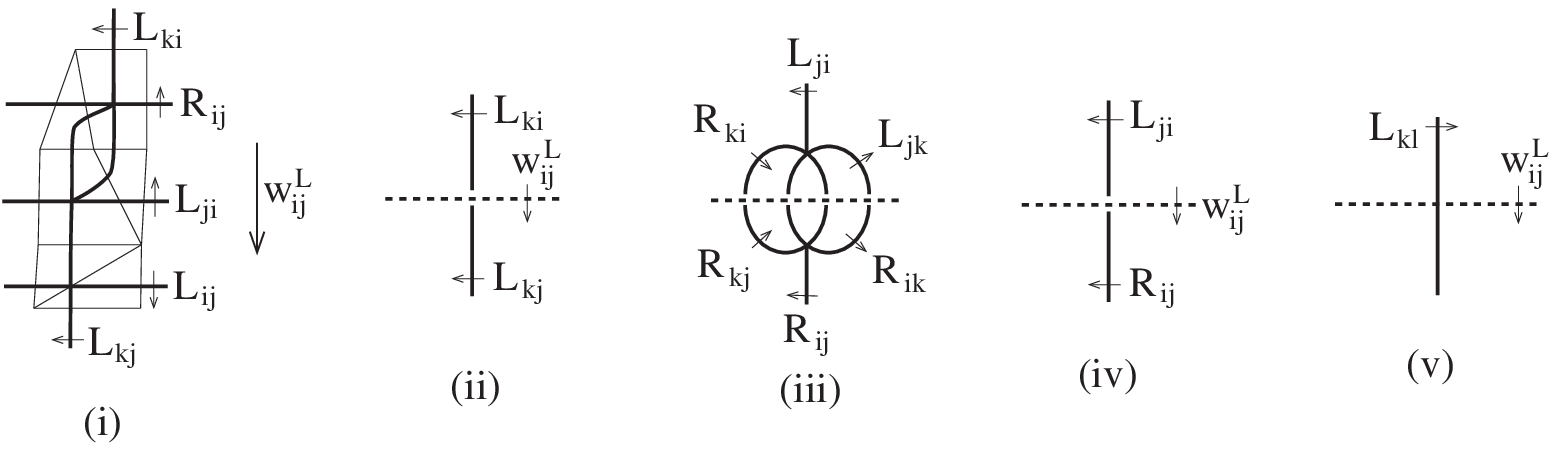}}}
\caption{}
\label{B1}
\end{figure}
Figure \ref{B1}(iii) indicates
a third type of standard $w$-crossing and figure \ref{B1}(iv) shows an abbreviation for it. 
This crossing will only be used in Figure \ref{swirl} below and discussion relevant to it. 
The use of the abbreviation is justified for our purposes
 since the choice of $k$ and the orientations of the $k$-edgesets
will not matter by Lemma \ref{2c}.
 There are four standard $w$-crossings in \ref{B1}(iii):
two are as in \ref{B1}(ii), but with the orientations of the edgesets changed, and two are as in 
\ref{A1}(iii).  

The last type of standard $w$-crossing are those involving four distinct subscripts 
and are abbreviated as in figure \ref{B1}(iv).  Each of these crossings consists of three 2-cells in 
$\xt$, where each 2-cell represents a trivial commutator relation.

The 2-cells in $\xt$  obtained from any one of the  above standard $w$-crossings by changing the
orientations of 
any combination of 
the $i$-, $j$-, or $k$-edgesets in all the cells of the crossing is also a standard 
$w$-crossing.

Let $C$ be any standard $w$-crossing, thought of as in abbreviated form, which contains no asterisks.
Then three other standard $w$-crossings are obtained from $C$ by adding asterisks on either or both sides of
the broken line in $C$, as in figure \ref{A1}(v). 

Let $\langle\!\langle A,x_k\rangle\!\rangle$ be any Whitehead element of the group $StN_n$, and let $w$
be any generator of $W_n$.  We next define a standard $w$-crossing $C$ corresponding to the relation
$w \langle\!\langle A,x_k\rangle\!\rangle w^{-1}=\langle\!\langle B,x_l\rangle\!\rangle$.
$C$ is the union of triangles in $\xt$ realizing the relations 
$\langle\!\langle A,x_k\rangle\!\rangle=\prod_{a\in A}\langle\!\langle a,x_k\rangle\!\rangle$ and
$\langle\!\langle B,x_l^{\pm 1}\rangle\!\rangle=\prod_{b\in B}\langle\!\langle b,x_l^{\pm 1}\rangle\!\rangle$
with the standard $w$-crossings involving $w$ and the $\langle\!\langle a,x_k\rangle\!\rangle$.
The order of the factors in the products is arbitrary except that if 
$\langle\!\langle a_p,x_k\rangle\!\rangle$ is the p$^{th}$ factor in
$\prod_{a\in A}\langle\!\langle a,x_k\rangle\!\rangle$, then $w_{ij}^L \langle\!\langle
a_p,x_k\rangle\!\rangle (w_{ij}^L)^{-1}$ must equal the p$^{th}$ factor in 
$\prod_{b\in B}\langle\!\langle b,x_l^{\pm 1}\rangle\!\rangle$.  Varying the order gives different 
standard $w$-crossings representing the same relations as $C$ does, as does varying 
the standard $w$-crossings in $C$.  

\subsection{The subcomplex $\widetilde{W}$ of $\xt$}\label{W}
For $n>4$, let $W_n$ be the group with generators $W_{ij}^L$ and $W_{ij}^R$, where $i$ and $j$ are
distinct integers between 1 and $n$, and the eleven relations \medskip
$$
\begin{array}{lll}
(1^L)\; &[W_{ij}^L,W_{kl}^L]  &(1^R)\quad [W_{ij}^R,W_{kl}^R]  \notag\\[3pt]
(2^L)\; &W_{ij}^L\, W_{ik}^L\, (W_{ij}^R)^{-1}\, W_{jk}^L &(2^R)\quad W_{ij}^R\, W_{ik}^R\,
(W_{ij}^L)^{-1}\, W_{jk}^R \notag\\[3pt]
(3^L)\; &W_{ij}^L\, W_{ki}^R\, (W_{ij}^L)^{-1}\, W_{kj}^L &(3^R)\quad W_{ij}^R\, W_{ki}^L\,
(W_{ij}^R)^{-1}\, W_{kj}^R \notag\\[3pt]
(4^L)\; &W_{ij}^R\, W_{jk}^L\, (W_{ij}^L)^{-1}\, (W_{ik}^L)^{-1}  &(4^R)\quad W_{ij}^L\, W_{jk}^R\,
(W_{ij}^R)^{-1}\, (W_{ik}^R)^{-1}  \notag\\[3pt] 
(5^L)\; &W_{ij}^L\, W_{kj}^L\, (W_{ij}^L)^{-1}\, (W_{ki}^L)^{-1}\quad &(5^R)\quad W_{ij}^R\, W_{kj}^R\,
(W_{ij}^R)^{-1}\, (W_{ki}^R)^{-1} \notag\\[3pt]
(6)     &W_{ij}^L=W_{ij}^R  & \notag
\end{array}$$

\noindent where $i$, $j$, $k$, and $l$ are distinct.  Note that each of the relations $(n^R)$ is obtained from
 $(n^L)$ by interchanging each occurence of $L$ and $R$, and so is a consequence of $(n^L)$ and
relation (6). 

\begin{lemma}
	There is an injective homomorphism $i\colon W_n\to StN_n$ induced by $i(W_{ij}^L)=w_{ij}^L$ and
$i(W_{ij}^R)=w_{ij}^R$.
\end{lemma}

\begin{proof} Figure \ref{2L} shows how the image under $i$ of the relations $(2^L)$
and $(3^L)$ are consequences of the defining relations of the group $StN_n$. 
Both 2-cycles in the figure consist of the union of three standard $w$-crossings, and, in the 
case of part $(2^L)$ of the figure, a triple of 2-cells realizing the relation $w_{ij}^R=w_{ij}^L$.
Figure \ref{wLeqwR} indicates, up to markings, these three 2-cells.

Inverting $x_i$ in part $(3^L)$ of the figure (by which we 
mean changing each occurance of $x_i$ to $x_i^{-1}$ and each occurance of $x_i^{-1}$ to $x_i$) shows 
how the image of the relation $(5^L)$ is a consequence of defining relations in $StN_n$.
Similarly, inverting  both $x_j$ and $x_k$, and then all of $x_i$, $x_j$ and $x_k$ shows that
relations $(5^R)$ and $(3^R)$, respectively,  map to relations in $StN_n$.  Also, inverting $x_i$,
then both of $x_j$ and $x_k$, and finally all of $x_i$, $x_j$, and $x_k$ in part $(2^L)$ of the 
figure shows
that relations $(4^L)$, $(4^R)$, and $(2^R)$, respectively, map to relations in $StN_n$.  Each of $(1^L)$ and
$(1^R)$ clearly map to relations, each of which are consequences of nine trivial commutator relations in the
definition of $StN_n$. Also, relation (6) maps to a defining relation for $StN_n$.  This shows that 
$i$ is well-defined.

That $i$ is injective follows, for instance, from Theorem 3.3 and
its proof in \cite{KN}.
\end{proof}

\begin{figure}[tb]
\centerline{\mbox{\includegraphics*{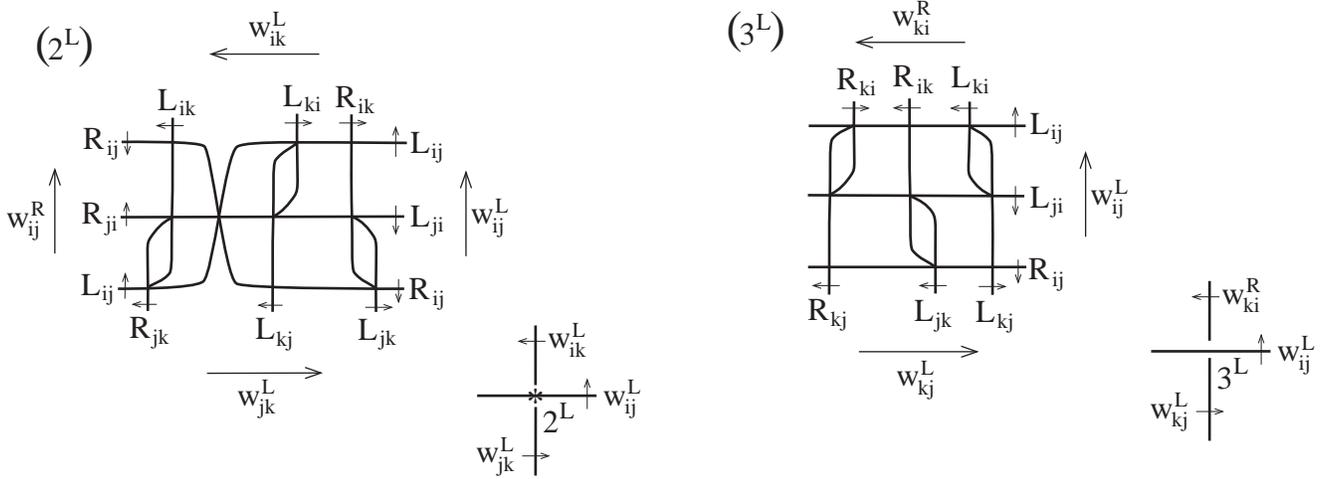}}}
\caption{Realizing the relations $(2^L)$ and $(3^L)$ in $\xt_n$.}
\label{2L}
\end{figure}


Let $BW_n^{(2)}$ be the $CW$-complex with a single vertex, with oriented edges corresponding to the 
generators of $W_n$, and with 2-cells (attached in the usual way)
 corresponding to the above defining relations for $W_n$.  The homomorphism $i$, both parts
of Fig.~\ref{2L} and the figures
 obtained from them by the inversions mentioned above, and Fig.~\ref{wLeqwR} serve to
define a cellular map $f\colon BW_n^{(2)} \to \xt_n/StN_n$.  Let $EW_n^{(2)}$ be the 
universal cover of $BW_n^{(2)}$, and let $\tilde{f}\colon EW_n^{(2)} \to \xt_n$ be any lift of $f$ 
which takes the vertices of $EW_n^{(2)}$ to those vertices of $\xt_n$ corresonding to the elements of the
subgroup $i(W_n)$ of $StN_n$.

Let $\widetilde{W}_n$ be the image of $\widetilde{f}$, and let $\widetilde{W}$ be the limit of the
$\widetilde{W}_n$ as $n\to \infty$.  The image under $\widetilde{f}$ of each 2-cell of $EW^{(2)}$
other than those corresonding to the relation (6) is called a standard $ww$-crossing. 
Such crossings  may, by definition, also contain asterisks.

Abbreviations for standard $ww$-crossings are also shown in Figure \ref{2L}.  The asterisk
in the abbreviation for the one corresponding to the relation ($2^L$) is often omitted, as are
the labels, ($2^L$) and ($3^L$) in the figure, of the relations they correspond to.

\section{Lemma T}\label{lemmaT}

We next describe a number of types of 2-cycles in $\xt_n$. 
Any cycle of one of these types is called a conjugate cycle.  
Such cycles arise in \S\ref{defP} as part of the ambiguity in passing from 2-cycles of $\ct$
to those of $\xt$. 

\begin{figure}[tb]
\centerline{\mbox{\includegraphics*{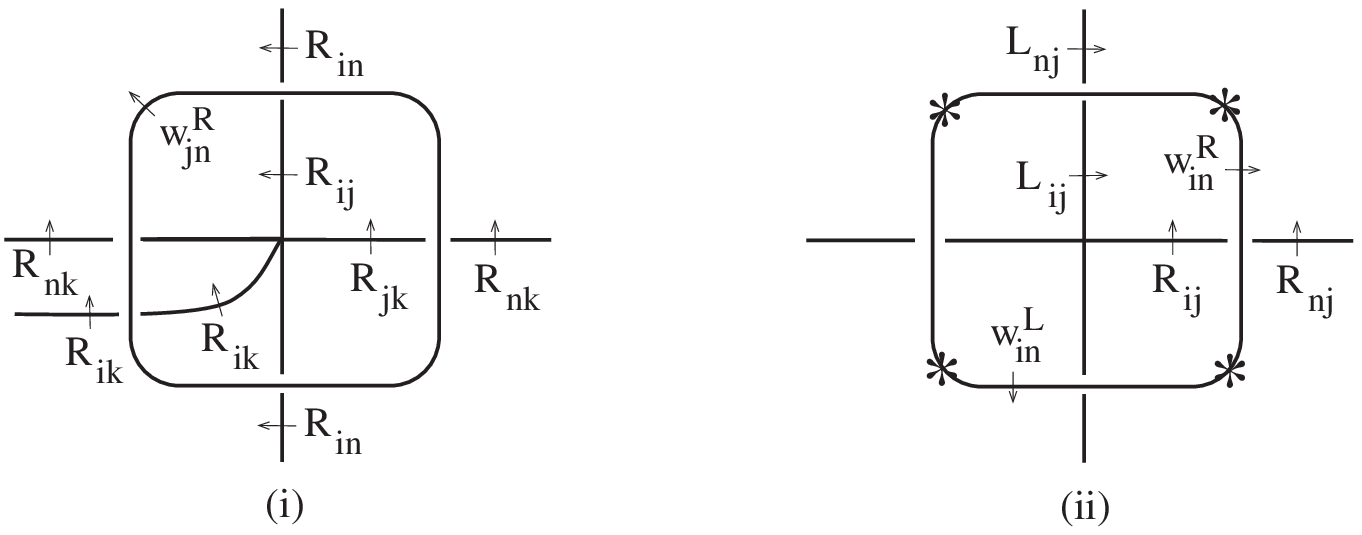}}}
\caption{}
\label{T1}
\end{figure}

Two examples of conjugate cycles are indicated in Figure \ref{T1} using the same conventions as in 
Figures \ref{A1}  and \ref{B1}.  In
general, there is, by definition, one type of conjugate cycle corresponding to each 
of the defining relations (R1)-(R3) in $\S$\ref{W} for the group $StN$.  Let $C$ be a 2-chain in $\xt$ 
realizing such a relation.  Each conjugate cycle of the corresponding type consists of standard 
$w$-crossings arranged along a circle surrounding $C$, as in the figure, together with a translate 
$C^w$ of $C$ on the opposite side of the circle.  The chain $C^w$ is only 
 implied by the loose ends in the figure which are understood to join in a point 
like that in the center of (i) and (ii).
The cells of $C$ and $C^w$ are called conjugate cells. We note that the marked colored graphs 
corresponding to the conjugate cells of $C$ are the same as those of $C^w$, except for the color of two edgesets of each graph.

 The circle in a conjugate cycle is called an $ij$-circle if crossing it in some direction corresonds 
to either  
$w_{ij}^L$ or $w_{ij}^R$.  An \{i,j\}-circle is either an
$ij$-circle or a $ji$-circle.  Some points on an $ij$-circle may  correspond to the relation
$w_{ij}^L=w_{ij}^R$, as, for example, in Figure \ref{T1}(ii). Such points are indicated by asterisks.

The relations (R1)-(R3) are realized by the reduced colored graphs with three vertices.
These are shown in Figure \ref{redsq}.  Thus conjugate cycles could have been defined in terms 
of these graphs.  We generalize the definition of conjugate cycles slightly so that there is one
type corresponding to each (not necessarily reduced) colored graph with three vertices.
So the  \{i,j\}-circle of a conjugate cycle may surround a chain corresponding to any number 
of the relations (R1)-(R3), rather than just one.

Figure \ref{T1}(i) indicates a conjugate cycle  corresponding to a relation of type (R2)
and  \ref{T1}(ii) to one of type (R1).  The right of
Figure \ref{ex} below shows another two conjugate cycles.  Both correspond to a relation 
of type (R3). These two types of conjugate cycles are the only ones which are nontrivial in $H_2(\xt)_{StN}$.

\begin{lemma}\label{T}  If $n\ge 5$, all conjugate cycles in $\xt_n$ are trivial as elements of $H_2(\xt_n)_{StN_n}$,
except for those which map to nontrivial elements of 
$H_2(\xt_n/StN_n)$.
\end{lemma}

\begin{proof}
We begin by showing that all conjugate cycles whose conjugate cells are reduced are trivial in
$H_2(\xt_n)_{StN_n}$, except for the noted exceptions.  

\begin{figure}[tb]
\centerline{\mbox{\includegraphics*{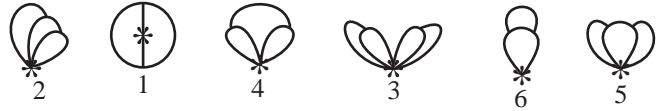}}}
\caption{All reduced colored graphs with three vertices.}
\label{redsq}
\end{figure}

We first assume that the conjugate cells correspond to the graph numbered 1 in the figure, where $k$ is the color of the two edges
forming a cycle and $i$ is the color of the other two edges.  If  the solid circle in the
corresponding conjugate cycle is labeled either  $w_{kn}$ or  $w_{nk}$, then the cycle is nontrivial in
$H_2(\xt_n)$ since the composite 
$H_2(\xt_n)\to H_2(\xt_n/StN_n) {\buildrel \Theta\over\to}  {\mathbb Z}^n$ takes this cycle
to $\pm(e_k\pm e_n)$. The map $\Theta$ was defined in \S\ref{H2xmodstn}.
So we assume that the subscripts are $i$ and $n$, and 
further assume that the conjugate cycle is as shown in Figure 3, since changing the orientations of the 
edgesets in this figure give all other possibilities.  

Let $G$ and $H$ be the colored graphs indicated in figure \ref{ignoreLR}  marked in any way such that
the marked roses $G_{b,d,f}$ and $H_{a,d,e}$ coincide.  With cell orientations as indicated in the figure,
let $g=gen(G)$ and $h=gen(H)$ be the corresponding 
generators of $H_3(\xt)$. Then there is a 3-chain $c$
in $\xt$ consisting of just additive prismatic terms such that $d(g+h+c)$ is the cycle indicated in 
figure \ref{ignoreLR}.  This cycle consists of six standard $w$-crossings, plus the terms corresponding to four relations
of type (iv) in the group $StN$.  The portion of the figure above the dotted line corresponds to terms
of $dg$, the portion below to those of $dh$.

\begin{figure}[tb]
\centerline{\mbox{\includegraphics*{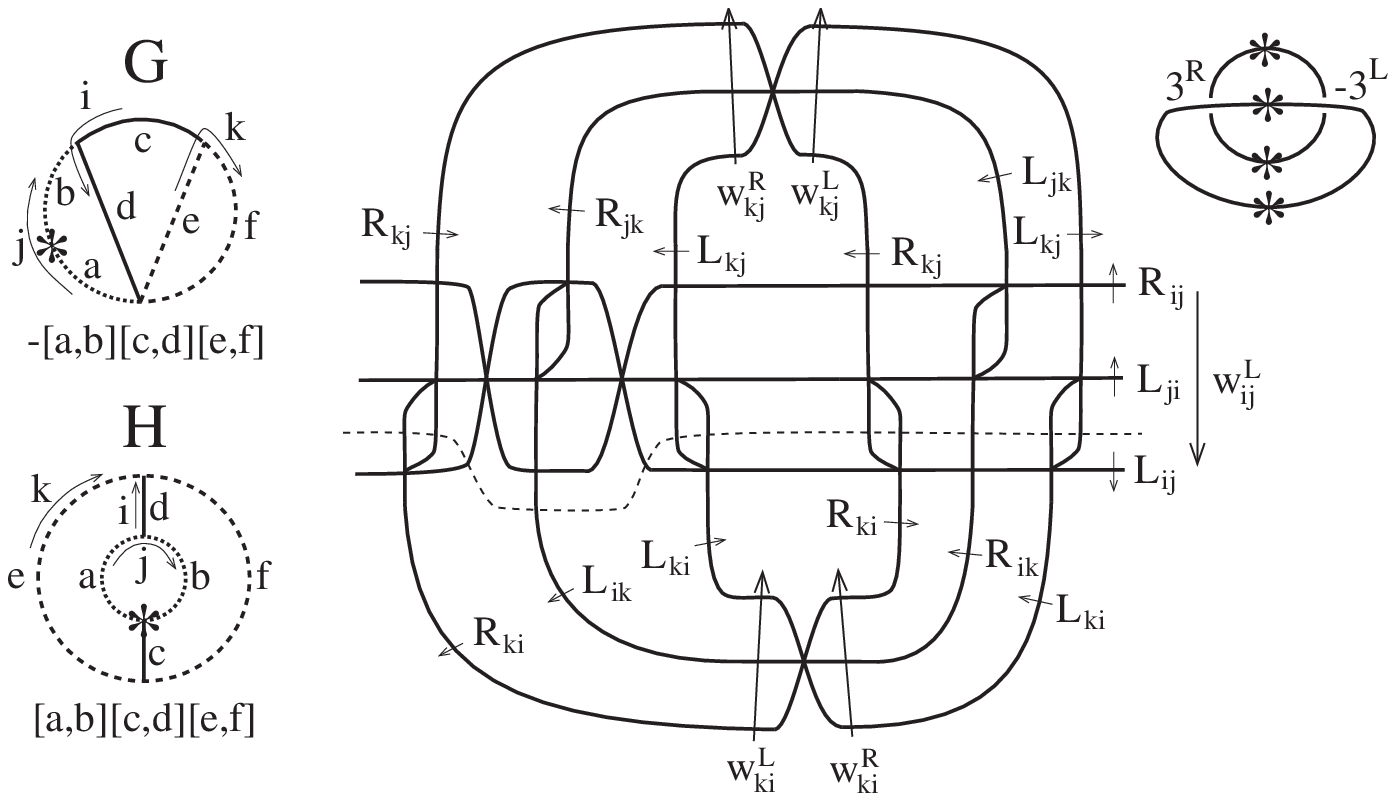}}}
\caption{}
\label{ignoreLR}
\end{figure}

By adding two $\theta$-prismatic terms to $g+h+c$, each of the two triples of cells nearest the dotted line
realizing the 
relation $w_{ij}^L=w_{ij}^R$ in figure \ref{ignoreLR} can effectively be moved around a standard $w$-crossing so as
to transform the picture in the center of figure \ref{ignoreLR} into the picture shown in the upper right
of this figure.  
Thus the 2-cycle indicated in this last
figure is a boundary. Changing the orientation of all of the $i$-edgesets in the above argument shows that the 2-cycle like
that  in the upper right of figure \ref{ignoreLR}, but with the relations $3^R$ and $-3^L$ replaced by $5^R$ and $-5^L$, is
also trivial in $H_2(\xt)$.

Now we assume that the conjugate cells of a conjugate cycle correspond to the graph numbered 2 in Figure
\ref{redsq}.  Let  $i$, $j$ and $k$ be the colors of the edges of the graph corresponding to one of the 
conjugate cells.  
There are then six cases depending on whether the $w$-crossings involve 
$w_{in}$, $w_{jn}$ or $w_{kn}$, or the same, but with the order of the subscripts switched.  
After, if necessary, moving asterisks across standard $w$-crossings
as in the previous case,
two of these six
cases are treated in Figure \ref{T2},   Each of the other cases are similar to one of these two. 
The six 3-cells on the left of the figure correspond to the colored graphs near them.  The markings of these
graphs can be any such that the indicated cancelations in the boundarys occur.  The boundary of these six
3-cells is then a conjugate cycle of the type represented by Figure \ref{T1}(i).  In a similar way,
the three 3-cells on the right show that conjugate cycles like those in Figure \ref{T1}(i), except with 
$w_{ni}$ in place of $w_{jn}$, are trivial as elements of $H_2(\xt_n)$.

\begin{figure}[tb]
\centerline{\mbox{\includegraphics*{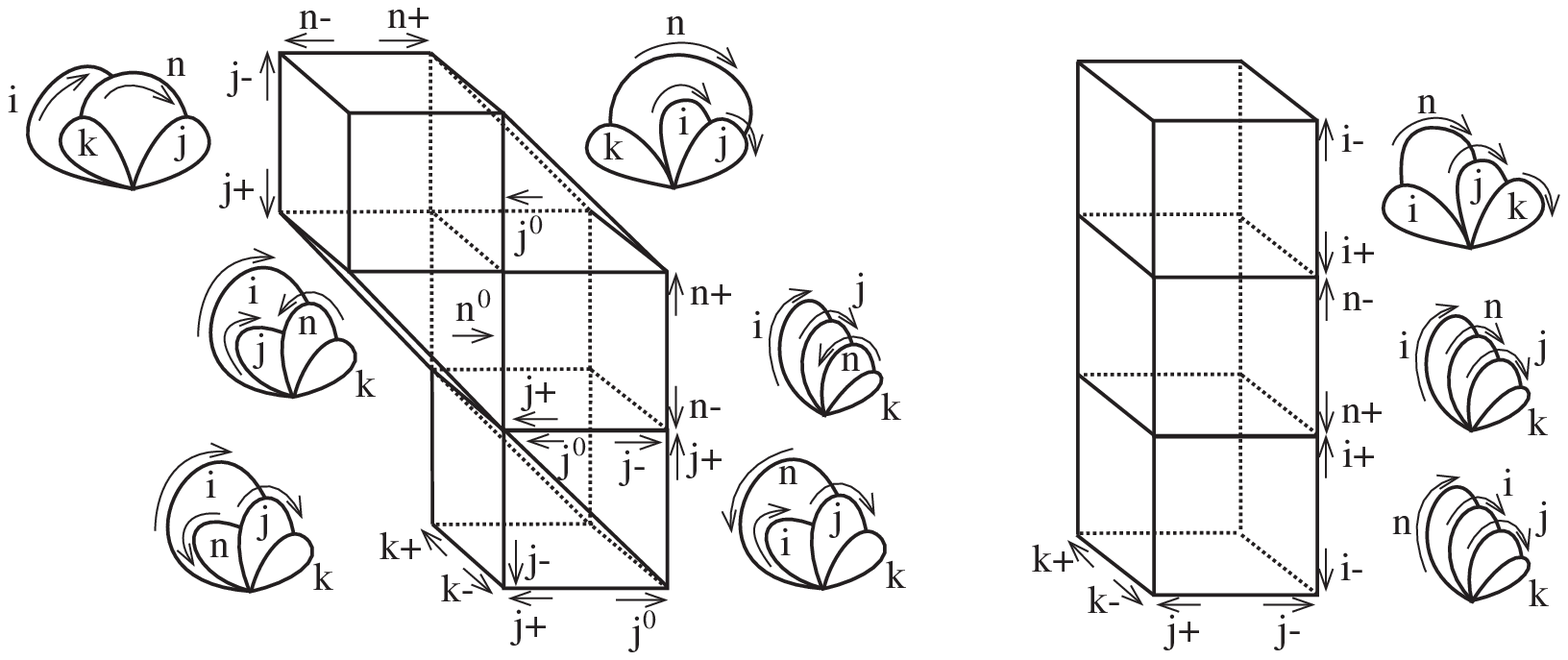}}}
\caption{}
\label{T2}
\end{figure}

Assume now that the conjugate cell corresponds to a colored graph with the underlying partitioned graph $G_1\vee G_2$ numbered 3 in Figure \ref{redsq}.  
If at least one of $i$ or $j$ of the $\{i,j\}$-circle of a correspoding conjugate cycle differs from 
the colors of $G_1\vee G_2$ , then
the conjugate cycle is the boundary of the 3-chain corresponding to the marked colored graphs making up a $w$-crossing, each wedged with $G_1$ or $G_2$.
If both $i$ and $j$ are colors of $G_1\vee G_2$, 
let $G_1'\vee G_2$ be the colored graph C2(a) of figure 1 in eliminating cubical terms section with 
the colors of the singleton edges and monocromatic loops of $G_1'$ the same as those of $G_1$.
Consider the 3-chain $C$ corresponding to $G_1'\vee G_2$ and the prismatic cell which makes all the 
cells of $\partial C$ reduced.
Then the sum of the conjugate cycles corresponding to the 2-cells
in $\partial C$  is itself a boundary.  But all of these conjugate cycles,  except for three,
have been considered in previous cases, so are boundaries.  One of the three exceptional cycles,
$x$ say, corresponds to $G_1\vee G_2$.  The other two, call them $y$ and $y'$, correspond to 
the same colored graph, but with opposite orientations.  Thus $y+y'=0$ in $H_2(\tilde X)_{StN}$.
So $x=0$  in $H_2(\tilde X)_{StN}$ since $x+y+y'=0$ in $H_2(\tilde X)$.

Assume now that the conjugate cycle corresponds to the graph numbered 4 in figure \ref{redsq}, where $i$ is the color of the singleton edge and $j$
and $k$ are the colors of the other two edgesets. If the conjugate cycle contains either a $\{j,n\}$- or
a $\{k,n\}$-circle, then previous cases essentially show that it is a boundary. If it contains an
$\{i,n\}$-circle, then replace each of the two conjugate cells with those five cells in  the boundary of
the 3-cell corresonding to the colored graph in Figure \ref {ignoreLR} containing two concentric circles
concentric circles which differ from the cell corresponding
to earmuffs.  Previous cases essentially show that the five conjugate cycles corresponding to these 
five cells are boundaries, which completes this case unless the conjugate cycle contains a $\{j,k\}$-
circle.  This last case is similar to the last one considered above involving $G_1\vee G_2$, but with
the graph  C2(b) in Figure \ref{reducedcores} in place of  $G_1'\vee G_2$ of C2(a).

Now let the conjugate cell correspond to the graph numbered 5, with the color of the monocromatic loop $k$, and $i$ and $j$ 
the other two colors. Cases like previous ones show that conjugate cycles containing $\{i,n\}$-
and $\{j,n\}$-circles are trivial. So we assume that the conjugate cycle $Z$ has an $nk$-circle, 
specifically one corresponding to $w_{nk}^L=R_{nk}^{-1} L_{kn}^{-1} L_{nk}$.
We show that $Z$ is the boundary of a 3-chain in $\xt_n$
 corresponding to 18 oriented 3-cells. Cells of these 18 which contain an edge corresponing to, 
respectively, $R_{nk}$, $L_{kn}$ and $L_{nk}$, are called, respectively, lower, middle and upper cells.
These three types of cells partition the 18 into three sets of six elements each.

The lower cells of $Z$ are in bijective correspondence with the permutations  of the
set $\{i,j,n\}$, and, accordingly, 
 triangulate  a cube $C$ in $\xt_n$ with four (oriented)
edges corresponding to each of $R_{ik}$, $R_{jk}$ and $R_{nk}$.  Specifically, the lower cell which 
corresponds to the permutation $p$ is
$[p(n)^+|p(j)^+|p(i)^+]_k$.  The edge path $[p(i)^+],  [p(j)^+],
[p(n)^+]$ consisting  of three of the edges of this cell is also an edgepath in the boundary of $C$, 
one of the six (of length three) from the vertex $init([p(i)^+])$ of $C$ to the 
diagonally opposite one $term([p(n)^+])$.

\begin{figure}[tb]
\centerline{\mbox{\includegraphics*{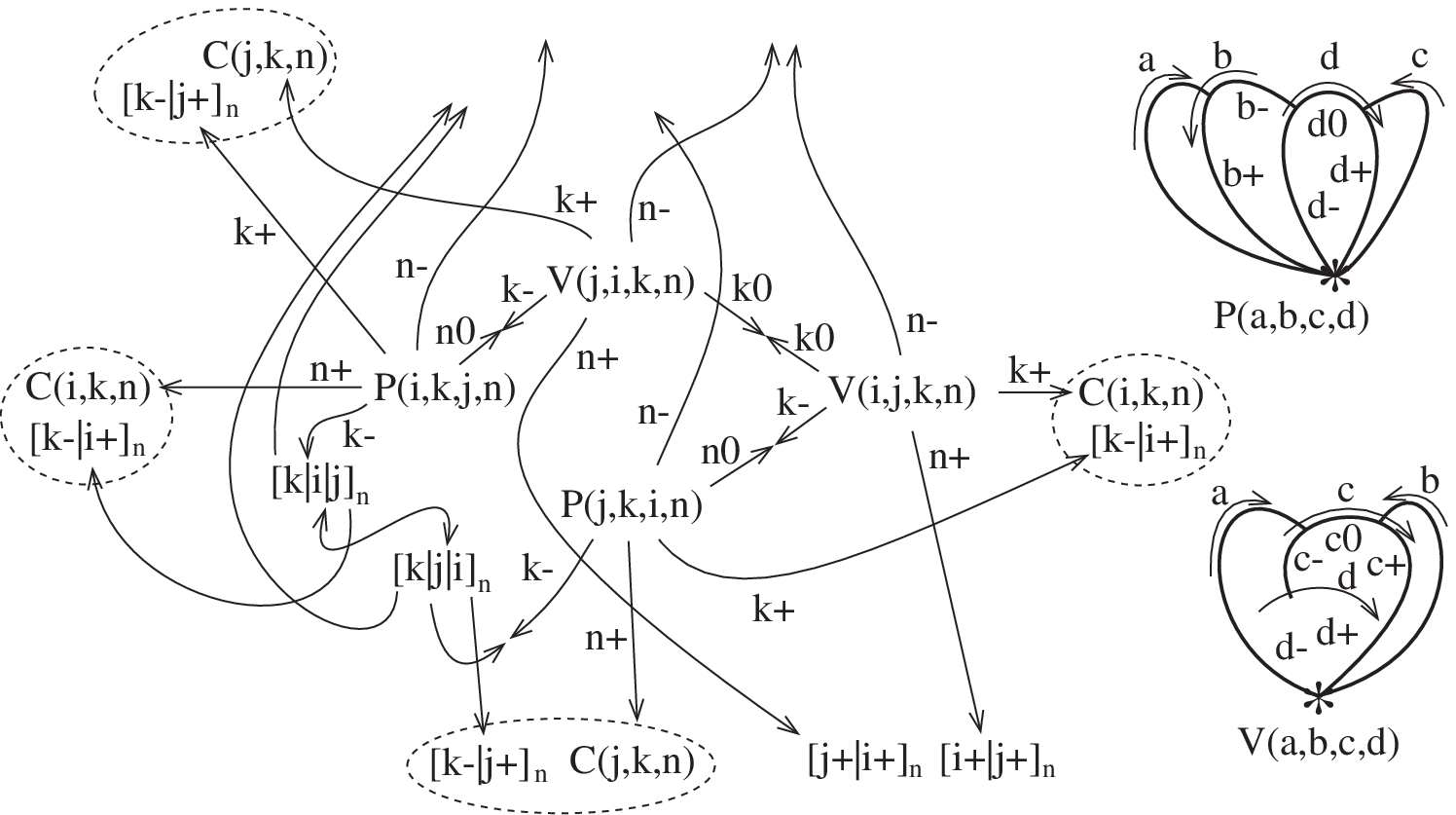}}}
\caption{}
\label{lemTpic}
\end{figure}

Let $P(a,b,c,d)$ and $V(a,b,c,d)$  be the colored graphs shown to the right of Figure \ref{lemTpic}.  The rest
of  the figure indicates the colored graphs corresponding to four of the middle cells.
The two other
middle cells, which are simplicial, are also indicated, using bar notation.
Markings and orientations are such that the indicated cancelations occur and such that two terms 
$[j^+|i^+]_n$ and $[i^+|j^+]_n$ to the lower right cancel with one face of the cube $C$.
Cells enclosed in dotted circles correspond to terms of the conjugate cycle $Z$.
This accounts for all the terms of the boundary of the middle cells except for the six
indicated by the upward pointing arrows at the top
center of the figure.

The upper cells correspond to the same colored graphs as do the middle cells, except that the 
colors $k$ and $n$ are switched, as are $i$ and $j$.  Orientations and markings are such that 
the six unaccounted for terms of the boundary of the middle cells cancel with the analogous six terms in the
boundary of the upper cells.  The  remaining terms of the boundary of the upper cells which do not cancel among
themselves are part of the conjugate cycle $Z$.

The case of conjugate cycles involving the graph numbered 5 and an $\{i,j\}$-circle is similar to the following case.

\begin{figure}[tb]
\centerline{\mbox{\includegraphics*{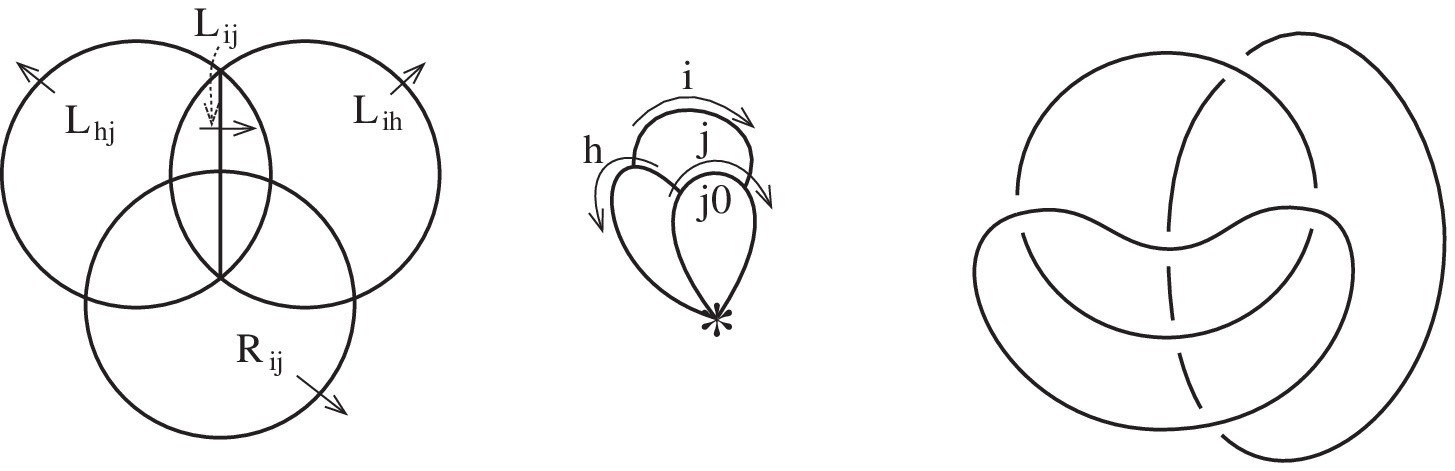}}}
\caption{}
\label{lemT2}
\end{figure}

Assume now that the conjugate cell is rank 2, where the monocromatic cycle has color $j$ and the other edge 
has color $i$.  This cell and an adjacent one
correspond to the crossing $x$ just above the middle of the picture on the
left of Figure \ref{lemT2}.
 Note that previous cases have dealt with conjugate cycles containing  conjugate cells
corresponding  to all the other crossings in the picture.  So, to complete this case, it's enough
to show that this
picture  represents a boundary.  Indeed it's  the boundary of a 3-chain $c$ with two prismatic terms:   
one corresponds to the
colored graph shown in the middle of Figure \ref{lemT2} and the other to the same colored graph, 
but with the orientation of the edgeset
$j$ reversed. Colorings and orientations  are chosen so that contracting the edge $j^0$ in both graphs gives
canceling terms in the boundary of the 3-cycle.
The chain $c$ also has three simplicial terms which isolate the crossing $x$ from the one just to its
right in the figure.

This completes the proof for conjugate cycles whose conjugate cells are reduced. 
By lemmas \ref{red}
and \ref{sred},
we may replace any nonreduced conjugate cell with a homologous reduced chain.  
Previous cases then apply.
 \end{proof}

The following  weaker version of Lemma \ref{T} will often be cited.

\begin{lemma}[Three Circles Lemma]\label{3circles}
Let  $n\ge 5$.
Let $P$ be a  picture in $\tilde W_n$ (or a translate of it in $\tilde X_n$) which
consists of three circles as shown on the right of Figure \ref{lemT2}
with each crossing corresponding to a defining relation
of the group $W_n$.  If the only asterisks of $P$ are on its single solid circle, then $P$ represents the trivial element of $H_2(\tilde X_n)_{StN_n}$.  If there are asterisks on any of the other edges of $P$, then $P$ may 
only represent the trivial element of the quotient $H_2(\tilde X_n)_{StN_n}/\langle O\!\!\!O \rangle $.

\end{lemma}

\section{Pictures and their Simplification}\label{pictures}

Theorem \ref{thminW}  below implies that each element of $H_2(\widetilde{X})_{StN}$ 
is represented by a 
2-cycle in $\wt$.  Such 2-cycles are studied in this section.  Figure \ref{ex} shows four examples 
of these, drawn using the notation introduced at the end of  \S\ref{W}. 
 Numbers in the figure near each crossing indicate which defining relation of \S\ref{W} 
the crossing corresponds to.
Following Igusa, representations of simple spherical cycles as in the figure are called pictures.  

Different 
pictures may represent the same cycle in $\wt$.  For instance, pictures which differ only as shown 
within the dotted circles in either (a) or (b)  of figure \ref{pict} are equal in $\wt$. 

Since $\wt$ is a
subcomplex of $\xt$, any picture gives a cycle in $\xt$.  Let $E$ be the subgroup of 
$H_2(\xt)_{StN}$ generated by the classes which the pictures 
 $C_{ij}^2$ and  $C_{ij}^4$ for $1\le i,j $ and $I_1^{-3,-4}$
 represent.   Notation is that of Figure \ref{ex}.  A subgroup $E_n$ of  $H_2(\xt_n)_{StN_n}$ is defined 
similarly, but with $ i,j \le n$.
The letter $E$ is used in the notation for these subgroups since, as we will see, the picture  $I_1^{-3,-4}$  maps to what has been 
called an exotic element of $K_3(\mathbb{Z})$.  
(Exotic elements are those not  in the image of the map
$\pi_3(B\Sigma^+) \to \pi_3(BGL(\mathbb{Z})^+)$.)
Let $\langle O\!\!\!O  \rangle$ be the subgroup of $E$ or $E_n$ generated by the classes of the pictures
$C_{ij}^2$ and $C_{ij}^4$ on the right of Figure \ref{ex}.
Most of this section is devoted to proving the following.

\begin{figure}[tb]
\centerline{\mbox{\includegraphics*{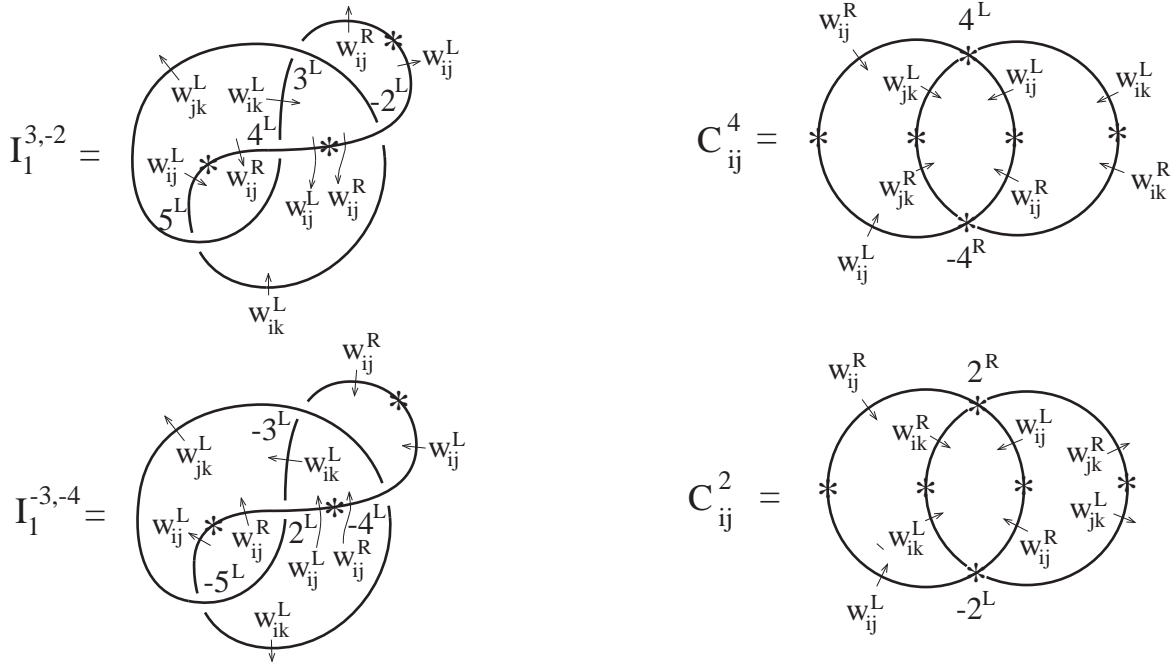}}}
\caption{Two preimages of an exotic element of $K_3(\mathbb{Z})$, and another pair of nontrivial elements of $H_2(\xt)_{StN}$.}
\label{ex}
\end{figure}

\begin{figure}[tb]
\centerline{\mbox{\includegraphics*{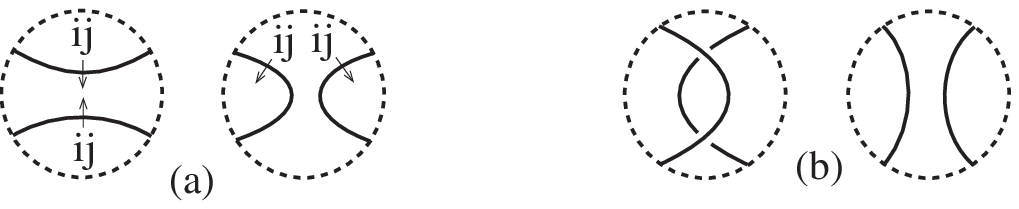}}}
\caption{}
\label{pict}
\end{figure}

\begin{theorem}\label{thm1}
Let $P$ be any picture in $\wt_n$ where $n \ge 5$.  Then the class of $P$ in $H_2(\xt_n)_{StN_n}$ 
is in the subgroup $E_n$.
\end{theorem}

This, along with theorem \ref{thminW} below, shows that  $ H_2(\xt_n)_{StN_n}=E_n$.

Numerous preparatory lemmas and some definitions proceed the proof proper of the theorem.
The first lemma we give is a convenience which shows that deleting the $L$ and $R$ superscripts and the
asterisks from any picture leaves a well-defined element of $H_2(\xt)_{StN}/\langle O\!\!\!O  \rangle$.

\begin{lemma} \label{noLR}
Changing the superscripts $L$ and $R$ and adding or deleting asterisks in any picture does not change the
class of the picture in $H_2(\xt_n)_{StN_n}/\langle O\!\!\!O  \rangle$ if $n\ge 5$.
\end{lemma}

This follows from the triviality of the conjugate cycles corresponding to the relation (R3) considered at the start of the proof of Lemma \ref{T}.



By forgetting some information,
any picture $P$ in $\wt_n$ may be thought of as a graph, also denoted by $P$, on $S^2$ with vertices,
all 4-valent, corresponding to the crossings 
  (and with some of the edges of the graph possibly containing
asterisks).  Thus we will speak of edges of a picture.  (Edges of pictures will be given 
a somewhat more general meaning after lemma \ref{sm} below.)
Because of lemma \ref{noLR}, transversely oriented edges of pictures are often labeled by just two letters,
with, for instance, $ij$ denoting either $w_{ij}^L$ or $w_{ij}^R$.  Of course, a picture drawn with this abbreviated notation
only determines a well-defined element of $H_2(\xt_n)_{StN_n}/\langle O\!\!\!O  \rangle$ and not
   a homology class in $\wt_n$, nor in $\xt_n$.
  
Figure \ref{swirls} shows a few pictures drawn with this abbreviated notation.  
The top and bottom portions of the picture on the left are examples of what we call swirls.
Swirls are of two types: knot-like and circular.  A knot-like (resp.~circular) swirl is, by definition, 
any 2-chain in $\xt$ with corresponding unlabeled partial picture as shown in 
   the top (resp.\ bottom) portion of the picture on the left of Figure \ref{swirls}.
There are eight types of knot-like swirls, the two shown in figure \ref{swirls} on the right, with three more obtained from each 
of these by changing the direction of each of $i$, $j$, and $k$.  By changing the direction of, for instance, $i$,
we mean changing the orientation of all the $i$-edgesets of the colored graphs corresponding to 
the terms of the cycle a given picture represents. 
For instance, the swirl at the bottom of the picture on the left of Figure \ref{swirls} is obtained from
the one in the upper right by changing the direction of $j$.
Another four knot-like swirls are shown in figure \ref{4swirls} below.

\begin{figure}[tb]
\centerline{\mbox{\includegraphics*{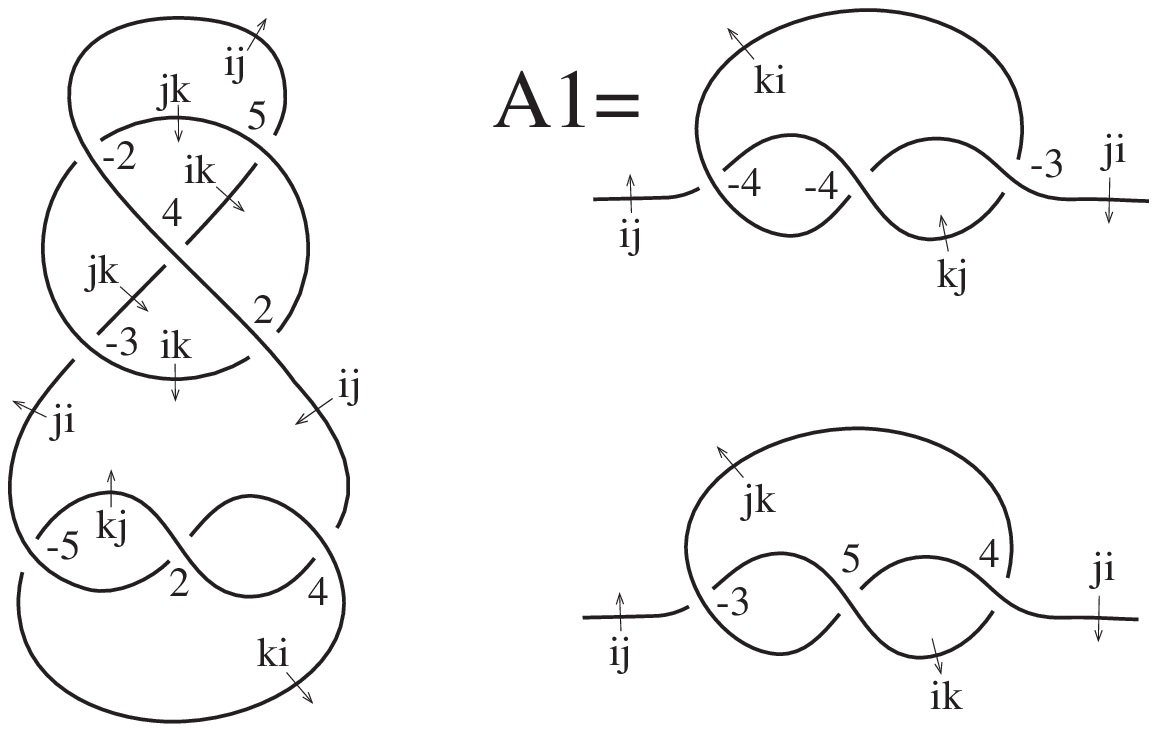}}}
\caption{}
\label{swirls}
\end{figure}

Swirls will often arise in the course of the proof of theorem \ref{thm1}.  In fact, the proof will show that,
by introducing swirls in any given picture $P$, the number of crossings in $P$ which are not part of swirls,
can always be reduced.  Thus $P$ is homologous in $\xt$ 
(modulo elements of $\langle O\!\!\!O  \rangle$) to a picture consisting solely of swirls.  Lemma
\ref{justsw} deals with pictures of this last sort.


Given a picture $P$, let $P_k$ be the subgraph of $P$ consisting of all the edges 
having label  $w_{ab}^L$ or $w_{ab}^R$ where either 
$a=k$ or $b=k$.  
Edges of $P_k$ are often referred to as $k$-edges.
Each vertex of $P_k$ is either bivalent or trivalent as a quick examination of the defining relations in 
\S\ref{W}  shows.
  Delete the bivalent 
vertices from $P_k$, and join the two edges meeting each of these into one edge.  
Denote by $P_k^*$ the result of this process.  Note that, as subsets of $S^2$, we have 
$P_k^*=P_k$.
In case $P_k$ has any 
components with only bivalent vertices, the corresponding components of $P_k^*$ are circles
with no vertices.  These circles are referred to as circular edges.  
The (possibly empty) complement of the circular edges in $P_k^*$
 has the structure of a graph with all vertices trivalent.  Any reference to  a vertex, edge, or cycle
of $P_k^*$ refers to one of this graph.
  
  Here is another application of lemma \ref{T}.
  
\begin{lemma}\label{2c}
Let $P$ be a picture in $\wt_n$. 
  Changing the transverse
orientations of all the $k$-edges in any one component of $P_k$ does not change the class of $P$ in 
$H_2(\xt_{n+2})_{StN_{n+2}}/\langle O\!\!\!O  \rangle$, nor does changing the subscript $k$  of all the $k$-edges in one of the components to $n+2$.
\end{lemma}

\begin{proof}
Write the cycle $P$ as $C_k+C_k'$ where the terms of $C_k$ are the generators of $C_2(\xt_n)$ corresponding to the 
crossings in one component of $P_k$ and the terms of $C_k'$ are those of the other components.  Let $C_{n+1}$ be the same 
chain as $C_k$, except with each appearance of the subscript $k$ in the corresponding  cycle changed to $n+1$. Then
$$P=(C_k-C_{n+1})+(C_{n+1}+C_k')$$  where both $C_k-C_{n+1}$ and $C_{n+1}+C_k'$ are cycles.
In $H_2(\xt_{{n+2}})_{StN_{{n+2}}}/ \langle O\!\!\!O  \rangle       $, we have 
$$C_k-C_{n+1}=   w_{n+2,k}^L\cdot (C_k-C_{n+1}) = C_{n+2}-C_{n+1},$$ 
where $C_{{n+2}}$ is the chain $C_k$ with $k$ changed to $n+2$ in the corresponding picture.
Thus $P=C_{n+2}+C_k'$, which 
proves the last assertion  of the lemma.  The first is proved similarly by replacing $w_{n+2,k}^L$ in the above equation with 
 $(w_{n+2,k}^L)^2$. 
\end{proof}

Let $P$ be a picture such that $P_k^*$ contains a circular edge $e$.  Let $P-e$ denote the picture
obtained from $P$ by deleting all the edges of $P$ which are contained in $e$.
The next lemma is often used and, for instance, shows, along with lemma \ref{T}, that swirls can be moved 
from one side of an edge of a picture to another.

\begin{lemma}[Broken Circles Lemma]\label{circle}
If $n\ge 5$, then $P$ and $P-e$ represent the same class in $H_2(\xt_n)_{StN_n}/E_n$.
\end{lemma}

The following proof actually shows that $P$ and $P-e$ represent the same class in 
  $H_2(\xt)_{StN}/\langle O\!\!\!O  \rangle$.  The weaker 
statement of the lemma is explained by remark \ref{remark} below.

\begin{proof}
Let $e$ be a circular edge as in the statement of the lemma.  Let $\alpha_1$ and $\alpha_2$ be the 
two loops in $S^2$ obtained by deforming $e$ slightly off itself in each of the two possible
directions.  Let $A_1$ and $A_2$ be thin annular neighborhoods of $\alpha_1$ and $\alpha_2$
disjoint from $e$, and let $A$ be the annulus consisting of $A_1\cup A_2$ together with the 
region between $\alpha_1$ and $\alpha_2$ which contains $e$.  Since all the vertices of $P_k$ contained
in $e$ are bivalent, we may assume that the two words in $W_{n-1}$ associated with the loops
$\alpha_1$ and $\alpha_2$ are the same.  Denote this one word, the trivial element of
  $W_n$, by $w$.  Pick a specific trivialization
of $w$ using the defining relations of the group $W_{n-1}$, that is, write $w$ as a product
of conjugates of these relations.  (For this, we need that the map $W_{n-1}\to W_n$ is injective, which 
  follows, for instance, from the proof of 3.3 in \cite{KN}.)
  Then use this trivialization (four times) to deform $P$ only
in the neighborhoods $A_1$ and $A_2$ to a picture $P'$ which is disjoint from $\alpha_1$ and 
$\alpha_2$.  Do this so that $A\cap P'$ admits an involution which preserves labels and transverse
arrows, which fixes $e$ 
  and which interchanges the two portions of 
$A\cap P'$ on each side of $e$.

Let $C$ be the component of $P'$ which contains the component $e$.  The involution pairs those 
vertices of $C$ which are not contained in $e$.  These pairs can be eliminated, one pair at a 
time, using  the three circles lemma \ref{3circles} with either      
a deformation like that in Fig.~\ref{EE} or, if the crossings  at the
paired vertices $v_1$ and $v_2$ in the figure are different, by using a deformation somewhat 
simpler than the one in the figure.  Both types of deformations preserve the 
$\mathbb{Z}_2$-symmetry of
$a\cap P'$ and do not change $P_k^*$.  Thus we may assume that $C$ consists of the circle $e$ together
with solid circles.  These solid circles can be eliminated (using lemma \ref{T} if necessary).  

\begin{figure}[tb]
\centerline{\mbox{\includegraphics*{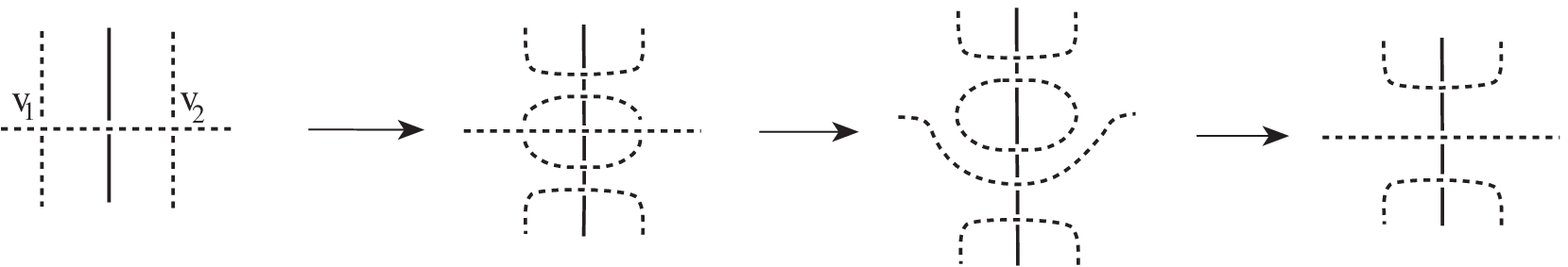}}}
\caption{}
\label{EE}
\end{figure}

So now we assume that $C$ consists of just the circle $e$.   Delete
$e$ from $P'$.  This does not change the class of $P'$ in the group $H_2(\xt)_W$ of coinvariants.  Next eliminate
the vertices in $A_1 \cup A_2$ of the remaining picture with a deformation which essentially reverses the
way in which they were introduced.  The resulting picture equals $P-e$ in $H_2(\xt_n)_{StN_n}/\langle O\!\!\!O  \rangle$.
\end{proof}

\begin{lemma}[Swirl Moving] \label{sm} Let $n\ge 5$. If two pictures agree except on neighborhoods
as shown in (a) and (b) or as in (c) and (d) of figure \ref{swirlmoving}, then the two pictures are equal in 
$H_2(\xt_n)_{StN_n}/\langle O\!\!\!O  \rangle$.
\end{lemma}

\begin{figure}
\centerline{\mbox{\includegraphics*{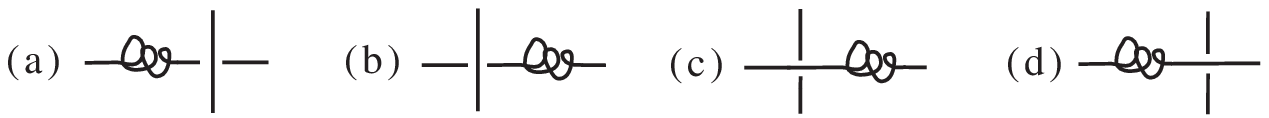}}}
\caption{}
\label{swirlmoving}
\end{figure}

\begin{proof} It follows from the broken circles lemma \ref{circle}  that two pictures which differ only as
shown in (c) and (d) are equal.  If they differ as shown in (a) and (b) they are equal by lemma \ref{T},
with the help of Lemma \ref{2c} if necesary.
\end{proof}

Since swirls can be moved ``across" edges, they will be,
to some extent, be ignored till lemma \ref{justsw} below, where pictures consisting solely of swirls 
are considered.  So the
collection of edges making up a swirl will often itself be referred to as an edge.
Pictures are often drawn accordingly without explicitly 
showing swirls (as in figure \ref{C}) which may be contained in their edges.  Or a swirl may be indicated as on the right of 
figure \ref{C2}(c).
Edges which contain no swirls may be refered to as swirl-free.

We next discuss conventions which are used in Figure \ref{C} and later figures.  Edges
in these labeled by $i$ and $j$ are $i$- and $j$-edges, respectively, of $P_k^*$, while dotted
edges are edges of $P$, but not $P_k^*$.  In pictures without arrows, any consistent choice of
these can be used.

\begin{lemma}\label{C1}
If two pictures differ only as shown in (a) and (b) or in (a) and (c) of Figure \ref{C}, then, after introducing a swirl
on one of the edges in the figure, the two pictures are equal in $H_2(\xt_n)_{StN_n}/\langle O\!\!\!O  \rangle$
if $n\ge 5$.
\end{lemma}

Note that an application of this lemma change $P_k$, but not the unlabeled graph $P_k^*$.

\begin{proof}
There are essentially eight cases depending on which of the four defining relations of the 
group $W$ the vertex in Fig.~\ref{C}(a) corresponds to and weather the two pictures differ as shown in (a) and
(b) or in (a) and (c). These cases correspond to the eight types of swirls mentioned above.   
For one case,
assume that a picture $P^a$ has a neighborhood like that shown in figure
\ref{C2}(a), while another picture $P^b$ is obtained from $P^a$ by changing this neighborhood
to the one shown in Figure \ref{C2}(b).  
By lemma \ref{sm}, we may assume that the edge of $P_k^*$ (completely) shown in \ref{C}(b) is swirl-free.
Then $P_a-P_b$ is the picture of \ref{C2}(c)
which is clearly trivial for a suitably chosen swirl.

For another case, assume the labels $ik$ in figure \ref{C2} are changed to $ki$.
Then the swirl appears on a vertical $k$-edge of figure \ref{C2}(b), but otherwise this
case is similar to the previous  one.  The remaining cases are similar to the two 
just considered.
\end{proof}

\begin{figure}[tb]
\centerline{\mbox{\includegraphics*{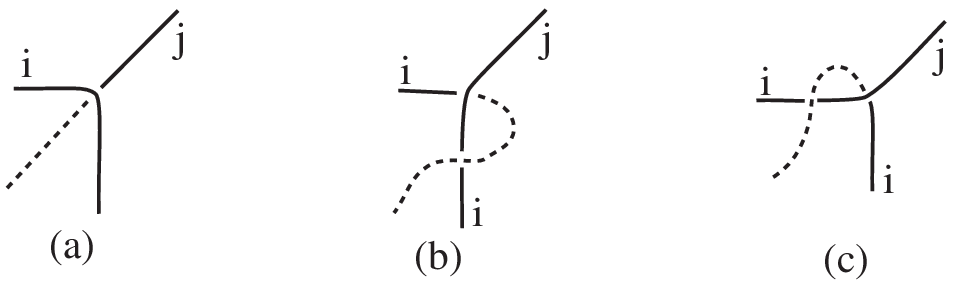}}}
\caption{}
\label{C}
\end{figure}

\begin{figure}[tb]
\centerline{\mbox{\includegraphics*{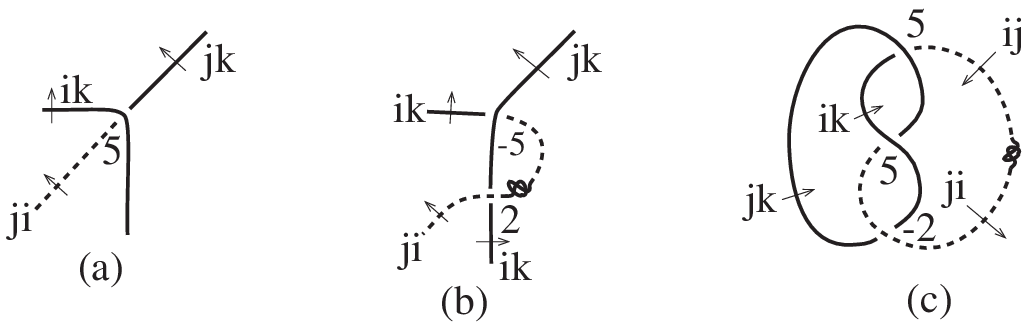}}}
\caption{}
\label{C2}
\end{figure}

Circular swirls will arise with applications of the following.

\begin{lemma}\label{LemC2}
Assume a picture $P$ has a neighborhood like that in Fig.~\ref{D}(a).  Then this 
neighborhood can be deformed to one like that in either (b) or (c) giving
a picture $P'$ equal to $P$ in $H_2(\xt)_{StN}/\langle O\!\!\!O  \rangle$, provided a suitable swirl is introduced on one
 of the $ij$-edges of $P'$.
\end{lemma}

\begin{figure}[tb]
\centerline{\mbox{\includegraphics*{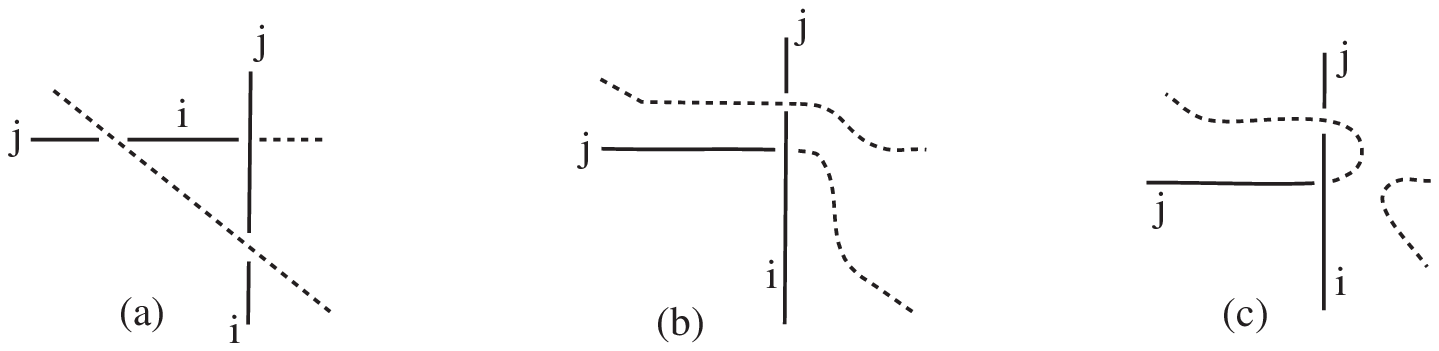}}}
\caption{}
\label{D}
\end{figure}

\begin{proof} 
There are numerous cases depending on the directions of the transverse arrows on the
edges and the labels of the edges.  For one case, assume two pictures differ by the 
neighborhoods shown in (a) and (b) of figure \ref{D2}.  
By the swirl moving lemma \ref{sm}, we may assume that the three edges forming the triangle in figure \ref{D2}(a) are swirl-free.
Then the difference of these
two pictures in $H_2(\xt)$ is represented by the picture in figure \ref{D2}(c).
This is trivial for a suitable choice of a circular swirl.  Other cases are similar.
\end{proof}

\begin{figure}[tb]
\centerline{\mbox{\includegraphics*{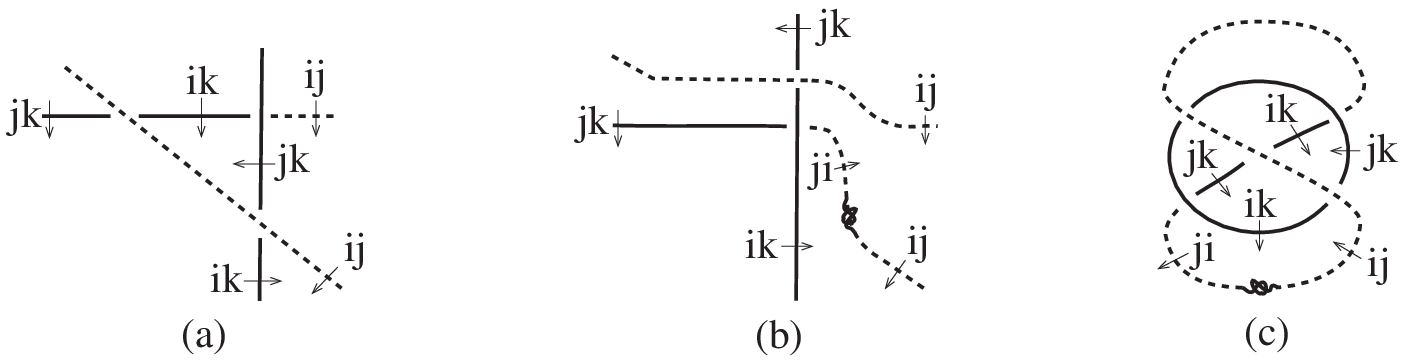}}}
\caption{}
\label{D2}
\end{figure}

An edgepath in a graph $G$ (which could contain circular edges) is a sequence $e_1,\ldots ,e_n$ of oriented (noncircular) edges of $G$ such that 
$term(e_i)=init(e_{i+1})$ for $i=1,\ldots ,n-1$.  An edgepath $e_1,\ldots ,e_n$ is a cycle if 
$term(e_n)=init(e_1)$, and a simple cycle if, in addition, the vertices $term(e_i), i=1,\ldots , n$ are
distinct.
 
The next lemma is used implicitly in the
definition  following it.

\begin{lemma}
Every cycle of $P_k^*$ consists of at least two edges.
\end{lemma}

\begin{proof} Suppose a single edge $e$ forms a cycle in $P_k^*$.  Then the loop in $S^2$ associated 
with $e$ can be deformed slightly to a loop $\omega$ such that the intersection $P_k\cap \omega$ is 
transverse and consists of a single point contained in that edge of $P_k$ which is incident with $init(e)=term(e)$ and which is
different from $e$.  But the element of the group $W$ associated with this loop cannot be trivial
since, for instance, the element of $\autfn$ it maps to does not take $x_k$ to $x_k$.  This contradiction
completes the proof.
\end{proof}

Let $P$ be a picture.
Fix an orietation of the 2-sphere containing $P_k$. This determines, at each vertex $v$ of $P_k^*$,
a cyclic ordering of the three edges of $P_k^*$ incident with $v$.  An edgepath $e_i,e_j$ is said to be
right-turning (at the vertex $term(e_i)=init(e_j)$) if $e_i,e_j,e$ is this cyclic ordering at the vertex 
$term(e_i)$, where $e$ is the third edge of $P_k$ incident with $term(e_i)$. A cycle $e_1,\ldots ,e_n$
in $P_k^*$ is right-turning if all the subedge-paths $e_i,e_{i+1}$ for $i=1,\ldots ,n-1$ as well as
$e_n,e_1$ are right-turning.  A left-turning cycle is one which is right-turning with respect to the 
orientation opposite to the fixed one.  A cycle is monotonic if it is either left-turning or 
right-turning.  A cycle in $P_k$ is monotonic if the corresponding cycle in $P_k^*$ is monotonic.

\begin{lemma}
Every monotonic cycle in $P_k$ is simple.
\end{lemma}

\begin{proof}
It suffices to show that monotonic cycles in $P_k^*$ are simple.  Suppose $\omega=e_1,\ldots ,e_n$
is a monotonic cycle in $P_k^*$, right-turning say, which is not simple.  Choose a simple subcycle
$\omega'$ of $\omega$ which is right-turning at all its vertices except for one, $v$ say. Then $\omega'$
can be deformed slightly to a loop in $S^2$ which meets $P_k$ transversely in a single point, namely
a point of that edge of $P_k$ incident with $v$ which is not part of $\omega'$.  The proof can now
be completed as was the previous one.
\end{proof}

Each monotonic cycle $\omega$ bounds two closed disks in $S^2$.  One of these contains the  component 
$C$ of $P_k$ which contains $\omega$.  The other has interior disjoint from $C$.  This latter disk is 
called a face of $C$, while the former is called a coface.  
The interiors of the faces of $C$ partition $S^2-C$.   
The edges of $P_k^*$ in the cycle forming the boundary
of  a face are often referred to as the edges of the face.  
A face with $n$ edges is called an $n$-gon.
We call 4-gons, 3-gons, and 2-gons, respectively,
squares, triangles, and bigons, respectively.

Fix a point $x\in S^2-P$.  Associate with each component $C$ of $P_k$ the coface of $C$
which does not contain $x$.  These cofaces  have disjoint boundaries, and so are partially ordered by
inclusion.  Components of $P_k$ which correspond to minimal or innermost such cofaces are called
innermost components of $P_k$.  

Let $C_k$ be an innermost component of $P_k$,
and let $C_k^*$ denote the component of $P_k^*$ corresponding to $P_k$.
 The face of $C_k$ containing $x$ is called   the backface of $C_k$.
The other faces are sometimes called frontfaces.
The interiors of frontfaces are disjoint from $P_k$, while the interior of each backface meets
$P_k$, unless $P_k$ is connected.

 A face of $P_k$ is, by definition, a clear face if the interior of the face does not meet $P$. 
We say that a frontface $F$ can be cleared if there is a deformation of $P$ which makes $F$ a clear face
without changing $P_k^*$.
\begin{lemma}[Clearing Lemma] \label{clearing}
Every front-face of an innermost component $C_k$ of a picture $P$ can be cleared.
\end{lemma}

\begin{proof}
Let $F$ be such a front-face. Use lemma 
\ref{C1} to deform $P$ so that $int(F)\cap P=\varnothing$ near each vertex of $C_k^*$ in $\partial F$.
Then vertices of $int(F)$ which are joined by a single edge of $P$ to a vertex of
$\partial F$ can be eliminatd, one at a time, using three circles lemma \ref{3circles} 
 or the broken circles lemma \ref{circle}.
Figure \ref{E} shows one way to do this for vertices in $int(F)$, which are joined by an 
unbroken edge to $\partial F$.  The deformation for the other type of vertex in
$int(F)$ is simpler.   

So now we assume that any remaining vertices in $int(F)$ are contained in components of $P$
different from the one of $\partial F$.  An edge of one of these components can
be deformed to meet $\partial F$ in two unbroken crossings.  Then the vertices of what was
this component can be eliminated from $int(F)$ as above.

By repeating the deormations described above as necessary, we may assume that all edges of 
$P$ in $int(F)$ meet two vertices  in $\partial F$.  These edges can be deformed outside
of $F$ using the three circles lemma and lemma \ref{LemC2}.
\end{proof}

\begin{figure}[tb]
\centerline{\mbox{\includegraphics*{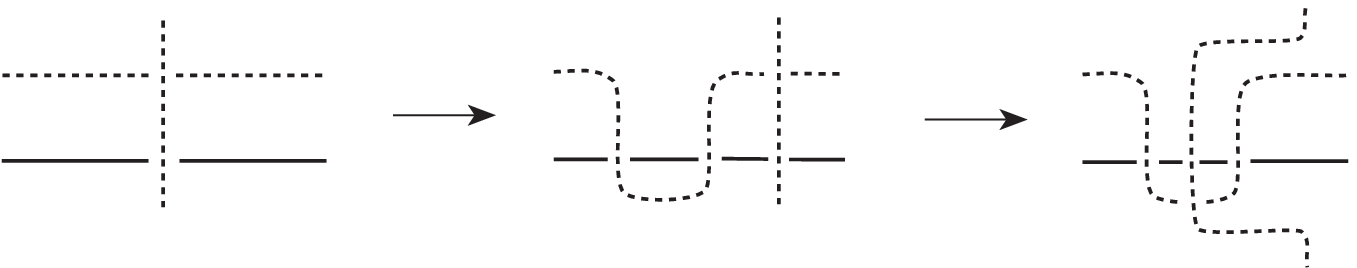}}}
\caption{}
\label{E}
\end{figure}

The clearing lemma just above often makes some of the following four lemmas applicable.

\begin{lemma}\label{one}
If $F$ is a clear face of $C_k$ with just one $l$-edge, for some fixed $l$,  then there
is a deformation of $P$ which decreases the number of vertices of $P_k^*$.
\end{lemma}

\begin{proof}
Assume $F$ is a clear face with just one $l$-edge. Use lemma \ref{C1}, if necessary,
to make the single $l$-edge of $F$ solid at both its vertices.  If $F$ is a bigon, it is 
then easy to deform $P$ in a neighborhood of $F$ so as to eliminate these two vertices.
This decreases the number of vertices of $P_k^*$ by two.

If $F$ is an $n$-gon, where $n>2$, then a deformation like that shown in figure \ref{F}
(or a similar deformation if the $j$-edge in the figure is broken at both of its vertices) 
produces a clear face with one $l$-edge and one fewer edge than $F$.  This deformation
changes $C_k^*$, but does not change the number of vertices of $C_k^*$.  Repeating such
deformations gives a clear bigon. The number of vertices of $P_k^*$ can then be reduced 
as above.
\end{proof}

\begin{figure}[tb]
\centerline{\mbox{\includegraphics*{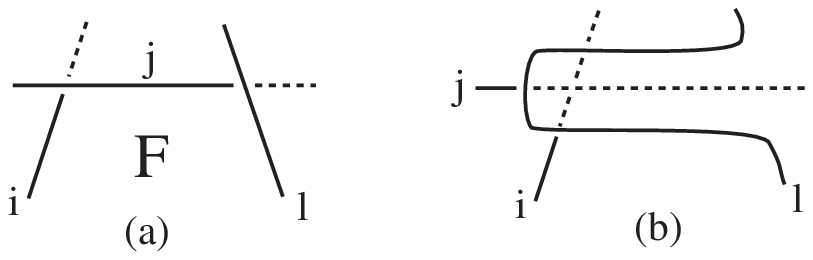}}}
\caption{}
\label{F}
\end{figure}

Two $l$-edges of a clear face $F$ are said to be consecutive if there is an edgepath in
$\partial F$ which joins the two edges and contains no other $l$-edges.

\begin{lemma}\label{two}
Assume $F$ is a clear face of $C_k$.  If two consecutive $l$-edges of $F$ have either
both inward- or  both outward-pointing transverse arrows, then there is a deformation of $P$
which decreases the number of vertices of $P_k^*$.
\end{lemma}

\begin{proof}
Let $e$ and $f$ be consecutive $l$-edges of $F$ which have either
both inward- or  both outward-pointing transverse arrows.  Surgery in $F$ as in Figure \ref{pict}(a)
using $e$ and $f$
produces two clear faces, at least one of which has just one $l$-edge.  So the previous 
lemma can be applied.
\end{proof}

\begin{lemma}\label{three}
If $F$ is a clear face with an odd number of $l$-edges, then there is a deformation of $P$
which decreases the number of vertices of $P_k^*$.
\end{lemma}

\begin{proof}
If the number of $l$-edges is odd and greater than one, then the arrows on all pairs of
consecutive $l$-edges do not alternate between inward- and outward-pointing.  So the previous
lemma applies.
\end{proof}

\begin{lemma}\label{four}
If $F$ is a clear face of $C_k$ with $i$-, $j$-, and $l$-edges, where $i,j,$ and $l$ are
distinct, then there is a deformation of $P$ which decreases the number of vertices of $P_k^*$.
\end{lemma}

\begin{proof}
If $\partial F$ contains three types of edges, then it has three consecutive edges of 
different types.  We may assume that these are arranged as in figure \ref{F}(a).
Then either $F$ has an odd number of $j$-edges or the deformation in figure \ref{F}
produces a clear face with an odd number of $j$-edges. Either way, the previous lemma
can be applied.
\end{proof}

While the clearing lemma \ref{clearing} will be used to provide clear faces to which lemmas 
\ref{one}--\ref{four} can be
applied, the clearing lemma itself can only be applied to frontfaces.  The next three lemmas
show that faces to  which lemmas \ref{one}--\ref{four}
can potentially be applied occur in numbers greater than
one.  Since $C_k^*$ has only one back-face, at least one of these must be a front-face.

\begin{lemma}\label{no!odd}
No component $C_k^*$ of $P_k^*$ can have exactly one face with an odd number of edges.
\end{lemma}

\begin{proof}
Let $F$ be a face of $C_k^*$ with an odd number of edges.  Since the interior of each edge 
of $C_k^*$ meets two faces we have
$$2(\text{number of edges of $C_k^*$})= \text{number of edges of $F$} +
\sum_{G\ne F}(\text{number of edges of $G$})$$
where the sum is over all faces $G$ of $C_k^*$ different from $F$.  It follows that some face
other than $F$ has an odd number of vertices.
\end{proof}

We mention that with slight modifications, the above proof shows that if all the faces of $C_k^*$
are $2n$-gons, then every cycle in $C_k^*$ consists of an even number of edges.

We say that a face $F$ of $C_k^*$ is surrounded by squares if each edge of $F$ is also the edge of 
a square (different from $F$, in case $F$ is itself a square).  We say that
$C_k^*$ is $n$-regular if $C_k^*$ has two disjoint $n$-gons $F_1$ and $F_2$ and another $n$ squares 
for faces,  if each edge of $C_k^*$ not in $F_1 \cup F_2$ meets both $F_1$ and $F_2$, and if whenever two 
of these edges meet consecutive vertices of $F_1$,  they also meet consecutive vertices of $F_2$.

If $C_k^*$ is $n$-regular then, using the above notation, both $F_1$ and $F_2$ are surrounded by 
squares.  Conversely, it can be shown that if $C_k^*$ has an $n$-gon which is surrounded by squares, then $C_k^*$ is
$n$-regular.  But we only need and prove the following partial converse.

\begin{lemma}\label{squares}
Assume all faces of $C_k^*$ are $2m$-gons for (possibly varying) $m>1$.  If some face of $C_k^*$
with $n$ edges is
surrounded by squares, then $C_k^*$ is $n$-regular.
\end{lemma}

\begin{proof}
Assume $F$ is a face of $C_k^*$ which is surrounded by squares.  Suppose an edge of $C_k^* -\partial F$
joins two vertices, say $v$ and $w$, of $F$.  Let $e_1,\ldots,e_p$ be a simple edge path in 
$\partial F$ from $v$ to $w$.  Since all faces of $C_k^*$ are $2m$-gons, all cycles in $C_k^*$ have
even length, and so $p$ is odd.  Note that $p\ne 1$ since none of the faces of $C_k^*$ are bigons by
assumption.   Since $e_1,e,e_p$ is a monotonic path meeting, but not contained in $F$, there must be
an edge not in $\partial F$ joining $init(e_1)$ with $term(e_p)$.  Here we are using the
orientations which make $e_1,e,e_p$ an edgepath.  Repeating this argument as necessary
shows that $e_m$, where $m=\frac12 (p+1)$, is the edge of a bigon.  Therefore no edge of 
$C_k^* -\partial F$ joins two vertices of $F$.

Now let $v$ be a vertex of $F$ and let $e$ be the unique edge of $C_k^* -\partial F$ incident with $v$.
Let $f(v)$ be the other vertex of $e$.  This defines a function $f$ from the vetices of $F$ to the 
vertices of $C_k^* -\partial F$.  Suppose $f$ is not injective.  Choose consecutive vertices
$v_1,\ldots,v_n$ of $F$ such that $f(v_1)=f(v_n)$ and such that  $f$ is injective  when restricted
to each of the sets $\{v_1,\ldots,v_{n-1}\}$ and $\{v_2,\ldots,v_n\}$. Note that $n>2$ since $C_k^*$ 
has no cycles of length three.  

If $n=3$, consider the square $S_1$ with vertices $f(v_1),v_1,v_2$, and
$f(v_2)$ and the square $S_2$ with vertices $f(v_3),v_3,v_2$, and $f(v_2)$.  Note that the single edge
joining $v_2$ to $f(v_2)$ is an edge of both squares.  Let $e_1$ be the edge of $S_1$ joining $f(v_2)$
and $f(v_1)$, and let $e_2$ be the edge of $S_2$ joining $f(v_2)$ and $f(v_3)=f(v_1)$.  If $e_1\ne e_2$,
then the vertex $f(v_1)$ has valence $>3$, a contradiction.  And if $e_1=e_2$, then $S_1$ and $S_2$ 
share two consecutive edges, another contradiction.  

So assume $n>3$.  But then each of the four 
distinct vertices $v_1,f(v_2),f(v_{n-1})$, and $v_n$ is joined by an edge to the trivalent vertex
$f(v_1)$, which is impossible.  This proves that the function $f$ is injective.  The proof of the lemma is
now easily completed.
\end{proof}

By a biggest face of $C_k^*$, we mean a face which has as many or more edges than all the other faces of
$C_k^*$.

\begin{lemma}\label{notsquares}
Assume all faces of $C_k^*$ are $2m$-gons, where $m>1$.  If $F$ is a biggest (nonsquare)
face of $F$ which is not 
surrounded by squares, then at least two edges of $F$ are not edges of squares.\end{lemma}

\begin{proof}
Suppose that $F$ is a biggest face of $C_k^*$ which is not surrounded by squares, and that all the edges of 
$F$, except one, are edges of squares. Let $c$ be the one exceptional edge of $F$, and let $G$ be the 
face of $C_k^*$ other than $F$ having $c$ as an edge.  Note that $F\cap G=c$.  The first paragraph of the 
previous proof can now be applied (provided the edgepath $e_1,\ldots,e_p$ in that proof is chosen so that
it does not include $c$) and shows that an edge of $C_k^* -\partial F$ cannot join two vertices of $F$.
The second paragraph also applies (provided the consecutive $v_1,\ldots,v_n$ are chosen so that both vertices
of $c$ are not among them) and shows that the function $f$ defined as above is injective.

Next we show that  the image of $f$ contains only two vertices of $G$.  Suppose there is an edge $e$ not
in $F\cup G$ which joins a vetex $v_F$ of $F$ to a vertex $v_G$ of $G$.  Let $e_1,\ldots,e_n$ be a 
simple edgepath in $\partial F \cup \partial G - int(c)$ from $v_F$ to $v_G$.  Here $int(c)$ is the 
interior of the edge $c$. Note that $n$ is odd
as remarked after the proof of Lemma \ref{no!odd}. 
 Since $e_1$ is an edge of $F$ different from $c$, the monotonic path $e_1,e,e_n$ is contained
in a square.  Let $e'$ be the fourth edge of this square.  Note that $e'$ is not in $F\cup G$ since faces 
don't share consecutive edges.  If $e'$ is incident with a vertex of $c$, then this vertex would have
valence $>3$, a contradiction.  Otherwise, repeat the above argument with $e'$ in place of $e$ and the
edgepath $e_2,\ldots, e_{n-1}$ in place of $e_1,\ldots,e_n$.  Continuing in this way produces a
contradiction, and proves that if $x$ and $y$ are the vertices of $c$, then the only vertices 
of $G$ in the image of $f$ are $f(x)$ and $f(y)$.

Now let $v_1,\ldots, v_n$ be the vertices of $F$, listed so that $v_i$ and $v_{i+1}$ are vertices
of the same edge and where $x=v_1$ and $y=v_n$.  Then the face of $C_k^*$  which contains the edge
joining $f(v_i)$ and $f(v_{i+1})$ for $i=1,\ldots, n-1$ also contain two edges of $G$ and therefore
 is bigger than $F$.  This contradiction completes the proof.
\end{proof}

Sometimes the following lemma is applicable when none of the previous four are.

\begin{lemma}\label{disconnect}
Assume two faces of an innermost component $C_k$ of a picture $P$ meet in a disconnected set.
Then there is a deformation of $P$ which 
divides  $C_k$ into two nonempty components, at least one of which is innermost.
Moreover, the number of vertices of $P_k^*$ is unchanged.
\end{lemma}

\begin{proof}
Assume  $F$ and $G$ are faces of $C_k$ having edges $e$ and $f$ in distinct components of
$F\cap G$.  At least one of $F$ or $G$ is a front-face, say $F$ is.  Clear $F$. Let $\alpha$
 be a loop in $S^2$ which meets each of the edges $e$ and $f$ transversely in one point,
and otherwise is contained in $int(E) \cup int(F)$.  Let $w$ be the word in the group $W$ corresponding
to $\alpha$.  Note that the two letters $a$ and $b$, say, in $w$ which correspond to the edges $e$ and $f$ may be assumed to be adjacent
and that these are the only letters in $w$ which have $k$ in their subscripts. 

Since $w$ represents the trivial element of the group $W$, it follows that the
subscripts of $a$ and $b$ are made up of just two letters.  For, if for instance $a=w_{ki}^L$ and
$b=w_{kj}^L$, then even on the level of the symmetric group, $k$ would not be fixed.
So $a=b^{-1}$ or if, if the subscripts of $a$ and $b$ are in opposite order, we may assume
 $a=b^{-1}$ by introducing a pair of canceling swirls on either $e$ or $f$.
  So surgery in $F$ as in figure \ref{pict} using $e$ and $f$ can be done.  
Its easy to check
that if the basepoint $*$ of $S^2$ is in $F \cup G$, then the surgery can be done so that  both of the
created components are innermost.  And if $*\in S^2 - (F\cup G)$, then one of these components is 
innermost.
\end{proof}

The next two lemmas are also used in the proof of Theorem \ref{thm1}.


\begin{lemma}\label{justsw}
If a picture $P$ in $\wt_n$ consists of disjoint circular edges, then $P$ represents an element of 
$E_n$ provided $n\ge 5$.
\end{lemma}

\begin{proof}
If there are no swirls on $P$, then $P$ represents the trivial element.  Suppose $P$ contains a circular
swirl. Let $P'$ be the subgraph of $P$ consisting of this swirl, and assume that $k$ is the third 
subscript of $P'$, with $i$ and $j$ the other two.  Then either $P_i'\equiv P_i\cap P'$ or 
$P_j'\equiv P_j\cap P'$ contains a clear bigon, say $P_i'$ does.  Use lemma \ref{C1} to eliminate this
bigon.  After this, $P_j'$ contains a clear bigon.  Eliminate it using lemma \ref{C1} again.  Then after 
moving the two swirls introduced with the two applications of  \ref{C1}, $P'$ consists of one
edge with three knot swirls.  Thus we assume that all the swirls of $P$ are knot swirls.  These
are even in number, and $P$ equals, in $H_2(\xt)$, a picture consisting of half as many disjoint  
circlular edges with two knot swirls on each. 

So we now assume that $P$ contains exactly two
knot swirls, $S_1$ and $S_2$. Let $ij$ and $ji$ be the labels of the two edges of $P$ which meet both 
$S_1$ and $S_2$.  The single edge of each swirl which meets both of these edges is called the long edge of 
the swirl.  By Lemma \ref{2c}, we may assume that the arrows on the long edges of $S_1$ and $S_2$
point inward, that is in the same direction as those in Figure \ref{4swirls}. Then the classification of 
swirls shows that, up to permuting the subscripts $i$ and $j$, each of the swirls $S_1$ and $S_2$  is one
of the four in Figure \ref{4swirls}.

\begin{figure}[tb]
\centerline{\mbox{\includegraphics*{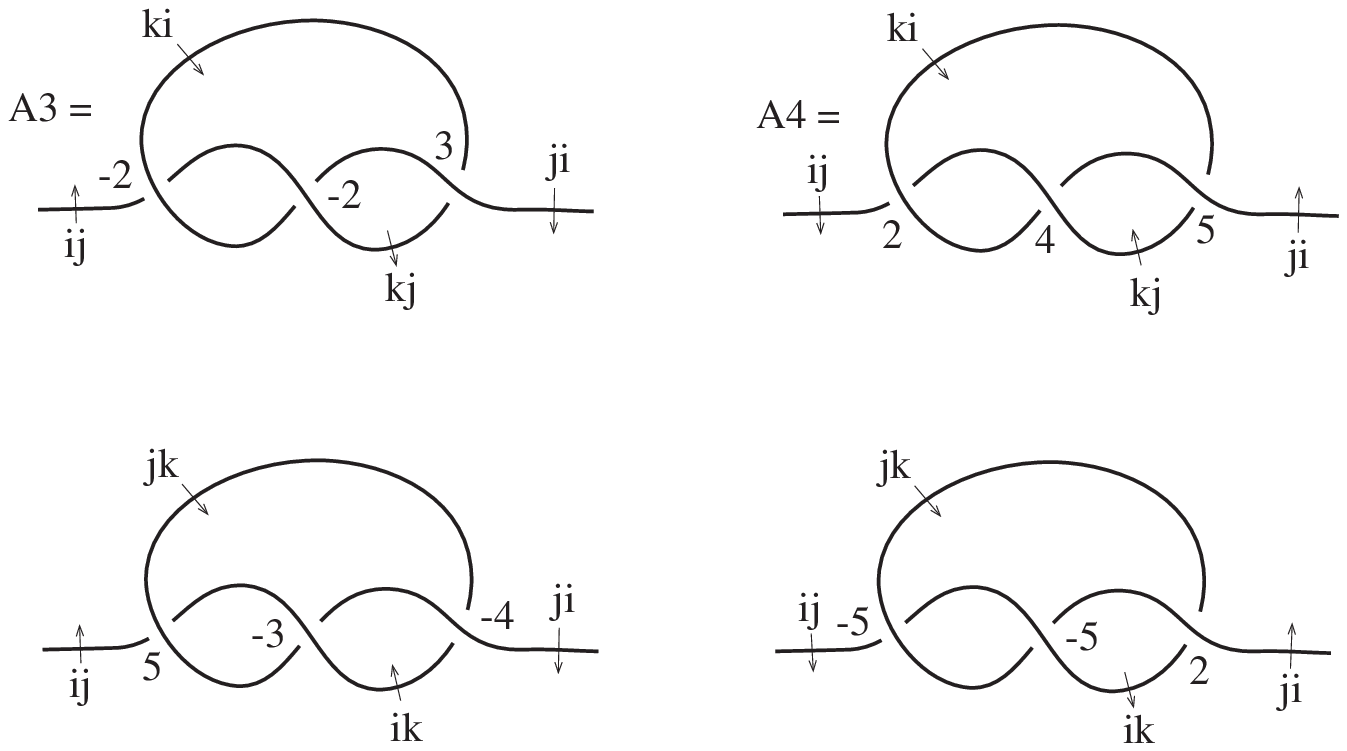}}}
\caption{}
\label{4swirls}
\end{figure}

As a first case, assume that $S_1$ and $S_2$ are the same swirl.  Then either they cancel trivially at the 
chain level giving $P=0$ in $H_2(\xt)_{StN}$, or the positions of the subscripts $i$ and $j$ in one of
$S_1$ or $S_2$ are different from those in the Figure.  In the latter case, the broken circles lemma 
\ref{circle}
can be used to effectively move a crossing of one of the swirls, say $S_1$, from one side of it to the other.
Then, after using Lemma \ref{2c} if necessary to change the direction of the arrows of the edges of $S_1$
with labels containing the third subscript, we have $P=0$ at the chain level giving $P\in E_n$.

For another case assume that $S_1$ and $S_2$ are either  the two swirls in the top row of Figure \ref{4swirls}
or the two in the bottom row.  By a previous case, we may assume that the $ij$ subscripts are just as in the 
figure, and, by Lemma \ref{2c}, that the third subscripts of both swirls are the same.  Two crossings of $P$ then
cancel leaving one of the generators of $E$ shown on the left of Figure \ref{ex}.

For the last case, assume that one of $S_1$ or $S_2$, say $S_1$, is in the top row of the figure and the other is
in the bottom.
By Lemma \ref{2c}, we may assume that the third subscript of $S_1$ and $S_2$ are different. 
Then the three circles lemma \ref{circle} can be used to move a crossing of $S_2$ 
from one side of this swirl to the other.  This changes the picture $P$ so that the position of the third subscript
of $S_2$ preceeds the other subscript in each of the labels in which the third subscript appears.  Previous cases
now apply to complete the proof.


\end{proof}


\begin{remark}\label{remark}\rm
Lemma \ref{circle}, the Broken Circles Lemma, which was stated and proved before the conventions that 
edges of pictures could contain swirls, is still true with this more general notion of edges.
The proof of \ref{circle}, appended as follows, shows this. 

At the start of the proof move all swirls near, but  not on $e$ away from $e$ outside of the annulus $A$.
The edge $e$ itself may contain swirls. Now, with the help of the swirl moving lemma, the proof of Lemma \ref{circle} goes through till the last paragraph.  So
assume the circle $C$ there  contains swirls.
  Deform $C$ into two circles, one of which contains all the swirls
on $C$ and bounds a disk with interior disjoint from the picture.  The other circle should contain
no swirls.  The resulting picture equals $P_1 + P_2$ in $H_2(\xt)$ where $P_1$ is represented by the
former circle and $P_2$ is one of the pictures considered in the proof of Lemma  \ref{circle}
where it was shown to equal $P_2-e$,
modulo an element of $\langle O\!\!\!O  \rangle$.
This completes the proof since $P_1$ represents an element of $E_n$  by lemma \ref{justsw}.
\end{remark}

\begin{proof}[Proof of Theorem \ref{thm1}]
Let $P$ be any picture $wt_n$.
Let $C_k$ be an innermost component of $P_k$.  We show that there
exists a deformation of $P$ which either
decreases the number of vertices of $P_k^*$ or
disconnects $C_k$, without changing the number of vertices of $P_k^*$, and
gives  
an innermost component whose vertices are a proper subset of those of $C_k^*$.
   This deformation does not enlarge the set of subscripts of those $w_{ij}$ which are labels of 
the edges of $P$.  Nor does it
 change the class of $P$
in $H_2(\xt_n)$, except possibly by an element of $E_n$. 
Since the number of vertices of $C_k^*$ is finite, iterating such deformations eventually gives 
one which reduces the number of vertices of $C_k^*$.  Thus all the vertices of
$C_k^*$ can be eliminated,  leaving
only circular edges.  These can be eliminated using lemma \ref{circle}, except for an
innermost circular edge with swirls.  So the theorem will follow from the previous lemma.

 If $C_k^*$ has a face with an odd number of 
vertices, then it has, by lemma \ref{no!odd}, at least two such faces. At least one of these must
be a front-face.  The Clearing Lemma \ref{clearing}, and then lemma \ref{three}, can be applied to
any such front face.  So we now assume that all faces of $C_k^*$ have an even number of edges.

If $C_k^*$ contains a bigon, then either $C_k^*$ consists of two vertices and three edges, or $C_k^*$ 
has two faces which meet in a disconnected set.  In the first case $C_k^*$ contains three
bigons, at least one of which is a front-face.  Thus the number of trivalent vertices can be reduced
by first applying the clearing lemma, and then using lemma \ref{one}.  In the second case, 
lemma \ref{disconnect} can be applied to disconnect $C_k$.

So now we assume that the number of edges of each face of $C_k^*$ is even and at least four. 

Choose a biggest face $F$ (as defined just before lemma \ref{notsquares}) of $C_k^*$.  Assume
that $F$ is surrounded by squares.  (The case in which $F$ is not surrounded by squares is considered
 in the last part of the proof.)
By lemma \ref{squares}, we may assume that $F$ is a front-face.  Clear $F$.  Then by lemma
 \ref{four}, we
may assume that $F$ has edges of just two types, say $i$- and $j$-edges.

If $F$ has more than four edges, then there must be two $i$-edges which have either both inward- or both 
outward-pointing arrows.  Surgery using these two edges produces a front-face which
meets what was the other largest face of $C_k$ in at least two components.  Thus $C_k$ can
be disconnected using lemma \ref{disconnect}.

So assume $F$ has four edges.  
We first will deform $P$ so that a disk neighborhood of $C_k$ is as shown in figure
\ref{cube}, with the back-face of $C_k^*$ corresponding to the unbounded face in the figure.
 We then will show that after adding a multiple of $exotic$ to $P$, the number of 3-valent
vertices of $P_k$ can be decreased.  

The face of $C_k^*$ which does not meet the back-face is called
the inner face of $C_k^*$ (because of the position of the corresponding face
in figure \ref{cube}).  Clear the inner face.  Then, if necessary, use lemma \ref{C1} to deform $P$
in neighborhoods of the four vertices of the inner face
so that the unlabeled graph $P$ agrees with the corresponding neighborhoods in the figure.
  Similarly, deform $P$ in small neighborhoods of the four vertices of the backface so 
that $P$ agrees with the corresponding neighborhoods in the figure.  

The four edges joining the inner and 
back-face of $C_k$ are called  diagonal edges.  Any crossing of $P$ on the interior of a diagonal
edge must be a broken $k$-crossing. These crossings can be pushed off the diagonal edges in the
direction of the backface using lemmas T and \ref{LemC2}.   This may introduce swirls and crossings on the edges 
of the backface,  but $C_k^*$ is unchanged.  Now the only crossings in $P$ on the interiors of the 
edges of $C_k^*$ are those on the edges of the backface.  These are all broken $k$-crossings.
So $P$ can be deformed, using the broken circles lemma \ref{circle},  so that,
 up to the  labels of the edges,  a disk containing 
$C_k$ meets $P$ in the graph shown in figure \ref{cube}.

\begin{figure}[tb]
\centerline{\mbox{\includegraphics*{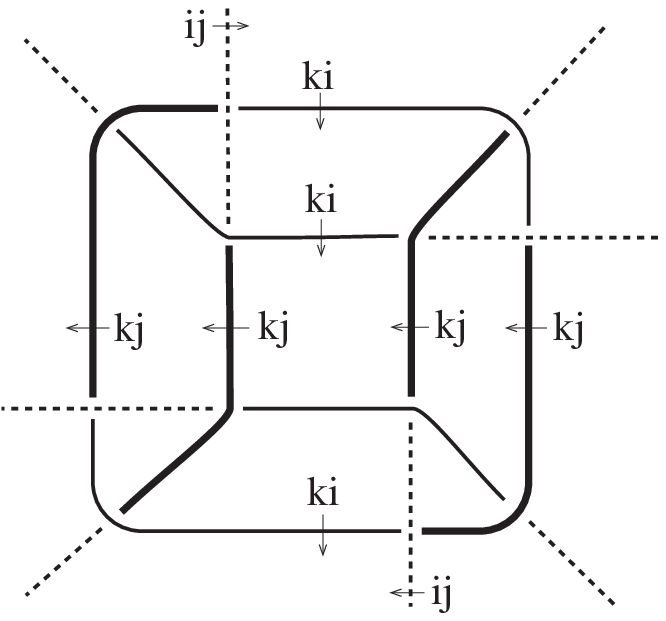}}}
\caption{}
\label{cube}
\end{figure}

We assume that the edges of the inner square are of only two types, say $i$- and $j$-edges, and
tranversely oriented as in the figure.  For otherwise, the number of vertices of $C_k^*$ can
be reduced as above using lemmas \ref{two} and \ref{three}.  
With this assumption, each diagonal edge of $C_k$ is either an $i$- or a 
$j$-edge, and the dotted edges in the figure are all $\{i,j\}$-edges.  So if the backface contains
an $l$-edge, where $l$ is different from $i$ and $j$, then all the edges of the back-face
are $l$-edges and any face meeting the backface, when cleared, contains just one $l$-edge.  
Because of lemma \ref{one},  we may assume, as above, that the edges of
the back-face are all either $i$- or $j$-edges with tranverse arrows as in the figure.  Otherwise surgery gives
bigons.


 So we now assume that a disk neighborhood of $C_k$ is  as in figure \ref{cube}. Let $v$ be any vertex
of the back-face and let $e$ be the $\{i,j\}$-edge of $P$ incident with $v$.  If the other vertex $w$
 of $e$ is part of a solid $ij$-crossing, then use lemma \ref{circle} to introduce a circle 
of edges (with all broken crossings) contained in the boundary of the disk neighborhood of $C_k$.
Do this so that the vertex introduced  on the edge $e$ cancels with the crossing at $w$.  
By repeating this procedure as necessary, we may assume that the vertex $w$ of $e$ is part 
of a solid 
 $ab$-crossing, where $\{a,b\}\not\subseteq \{i,j\}$.
 Now if $\{a,b\}\not\subseteq \{i,j,k\}$, with say $a\notin \{i,j,k\}$, then a deformation given by 
 lemma \ref{2c} which changes all the letters $k$ in figure \ref{cube} to $a$ eliminates $C_k$ with no other
 effect on $P_k$.  And if  $\{a,b\}\subseteq \{i,j,k\}$, then after possibly introducing a pair of canceling
 swirls, we may assume that  the vertex $w$ of $e$ is part of a solid
 $ki$- or $kj$-crossing.  Also since the boundary of a disk neighborhood of $C_k$
represents the word $w_{ij}^8$ (up to swirls), $w$ is not a vertex of $C_k$.  
 What's more, the above procedure may
increase the number of vertices of $C_k$, but not of $C_k^*$. 
 
 Now let $z$ be the cycle indicated in figure \ref{S} below. By adding $nz$ for suitable $n\in \{0,1,\dots ,7\}$
 to $P$, the graph $C_k$ can
effectively be rotated so that $w$ cancels with a vertex of the backface of $C_k$.  This
decreases the number of vertices of $P_k^*$, but may change $P$ by an element of the 
 group $E_n$, as the proof of Theorem  \ref{H2X} below shows.

Now assume that $F$ is a biggest face which is not surrounded by squares. Then by lemma 
\ref{notsquares}, two edges of $F$, say $e$ and $f$, are not edges of squares.  If $e$ and $f$
are both edges of a face $G$ different from $F$, then $F\cap G$
is not connected.  So, by lemma \ref{disconnect}, $C_k$ can be disconnected.

So assume $e$ and $f$ are edges of distinct faces $G$ and $H$, both different from $F$.  
At least one of $G$ or $H$ is a front-face, say $G$ is.  Clear $G$. By lemmata \ref{two}, \ref{three} and
\ref{four}, we may assume that $G$ is a $4n$-gon, where $n>1$, with arrows on consecutive edges of 
the same type pointing in opposite directions.

Assuming $e$ is an $i$-edge, let
 $e'$ be another $i$-edge of $G$ which points inward if $e$ does, and outward otherwise.
 Surgery in $G$ using
$e$ and $e'$ produces a face $F'$ of $C_k^*$ with two more edges than $F$ had.  It's
easy to check that $C_k$ is still innermost.  Now, if $F'$ is not surrounded by squares, repeat the
above argument with this new innermost component.  Since arbitrarily large faces cannot be produced
in this way, one of the previously considered cases must eventually occur.
\end{proof}

\section{The association $P$} \label{defP}
 
 Throughout this section stable versions of spaces defined above, such as $\xt$ and $\wt$, are used.
 However, all arguments go through with the corresponding unstable spaces, e.g.\ $\xt_n$ and $\wt_n$,
 in place of the stable ones provided $n\ge 5$. Theorem \ref{thminW} is used for these $n$. 
 
Let $f\colon (S^2,\cc) \to \ct^*$ be any cell-to-cell map. 
We next describe a way of associating with $f$ a 
2-cycle in $Z_2(X)$ along with a map $ S^2\to X$ representing it.  Both the cycle and map are denoted
by $p(f)$. Lifts of $p(f)$ from $X$ to $\xt$ will be denoted by $P(f)$.
This association is well-defined only in the sense that it induces a homomorphism from
$H_2(\ct)$ to the cokernel  of the composite $H_2(\widetilde{W})\to H_2(\xt)\to  H_2(\xt)_{StN}
$, as lemma \ref{Phom} below shows.
Of course, since $\ct$
is contractible, this is the trivial homomorphism. We will use this in conjunction with Theorem~\ref{thminW}
below in the computation of $H_2(\xt)_{StN}$.  In this section, we only use the fact that $\ct$ is 
simply-connected.

The map $P(f)$ is not a lift of $f$ up the map $F|\colon \xt \to \ct^*$ and, in general, no such lift exists.
Here $F|$ is the composite $\widetilde{X}\stackrel{\mathcal F}{\to} \mathcal S \stackrel{q}{\to} \mathcal T$ restricted to the 2-skeleton of 
$\widetilde{X}$ followed by the map $\mathcal T\to \mathcal T^*$.
From the point of view of pseudo-isotopy theory,  $p(f)$ is essentially obtained from $f$ by eliminating all
dovetail singularities by joining these in pairs.

Let $\{\sigma_i\}_{i=1}^n$ be the set of closed 2-cells of $\cc$ labeled so that, for $1<k\le n$,
$\bigcup_{i<k} \sigma_i$ is a disk and $\sigma_k\cap \cup_{i<k}\sigma _i$ contains at least one 
edge.  Such an ordering of the 2-cells is possible since each edge of $\ct^*$ has distinct vertices.
Orient each $\sigma_i$ so that, if $\sigma_i$ also denotes the generator of
the cellular chain group $C_2(S^2)$ corresponding to the oriented 2-cell $\sigma_i$, then
$\sum_{i=1}^n\sigma_i$ represents a generator of $H_2(S^2)$.

Any $\sigma_i$ which $f$ takes to a degenerate simplex of $\ct$, 
as defined in \S\ref{ct},
is itself referred to as degenerate.
Corresponding to each nondegenerate $\sigma_i$, we will define an element $\sigma_i'=\sigma_i^c +\sigma_i^w$
of the cellular chain group $C_2(X)$.  In case $\sigma_i$ is degenerate, $\sigma'_i=0$ by definition.
The cycle $P(f)$ will be defined as a lift of the cycle $\sum_{i=1}^n\sigma_i'$. 
We will also define sides of each $\sigma_i'$ as certain edgepaths in the boundary of the set
$\sigma_i'$, with one side corresponding to each 1-cell of the boundary of $\sigma_i$.
  One term of each side will be defined as the visible term of the side.  The other 
terms will be referred to as invisible.

We first define $\sigma_1'$ (or, if $\sigma_1'$ is degenerate, $\sigma_i'$, where $i$ is the least integer
such that $\sigma_i'$ is not degenerate).
If the 2-cell $f(\sigma_1)$ of $\ct^*$ corresponds to a marked 
partitioned graph, $G$ say, choose any coloring of $G$ consistent with the partition.
Let $\sigma_1'$ be the
2-cell of $X$ corresponding to the resulting marked colored graph.  The three or four sides of $\sigma_1'$
are defined as just the three or four edges of this 2-cell.

If $f(\sigma_1)$ is a $w$-triangle of $\ct^*$, let $G$ be the corresponding marked graph with edges
$a$, $b$, and $c$ such that $G/a$, $G/b$, and $G/c$ are the marked roses corresponding to the vertices
of $f(\sigma_1)$.  Choose a vertex of $f(\sigma_1)$, say the one corresponding to $G/a$, and an
edge of $f(\sigma_1)$ containing this vertex, say the one 
 corresponding to the partition $\{a,b\}$ of $G$.  (Different choices here
lead to exotic elements of $K_3(\mathbb Z)$.)  If more than three edges attach to the nonbasepoint vertex
of $G$, let $G^d$ be the marked graph shown in figure \ref{nodove}(i) with the requirement that $G^d/d=G$.
Note that $G^d$ depends on the choice of the vertex $G/a$ of $f(\sigma_1)$ made above since the edge
$a$ of $G$ is subdivided to obtain $G^d$.  We define $\sigma_1'$ to be the three 2-cells of $X$ 
indicated in figure \ref{nodove}(ii).  Here each 2-cell corresponds to the marked  graph $G^d$ partitioned 
as indicated.  Any consistant choice of colorings of the three partitioned graphs may be used.
One side of $\sigma_1'$, as indicated in (ii), is the sequence $b\wedge a, d\wedge b$, $b\wedge c$, $c\wedge d$
with $b\wedge a$ designated as the visible term.  The other two sides are the two one-term sequences
$a\wedge c$ and $c\wedge b$.  The visible term of each of these last two sequences is (necessarily)
defined as the only term of each.  If the choice of the edge made above had been $a\wedge c$ rather
than $a\wedge b$, then a figure like \ref{nodove}(ii), but reflected about the  vertical line bisecting it, 
would define $\sigma_1'$.
The sides of $\sigma_1'$ would be defined as the two one-term sequences $b\wedge a$ and $b\wedge c$, and
the four-term sequence $d\wedge b$, $b\wedge c$, $c\wedge d$, and $a\wedge c$, with $a \wedge c$ 
designated as the visible term of this last side.

\begin{figure}[tb]
\centerline{\mbox{\includegraphics*{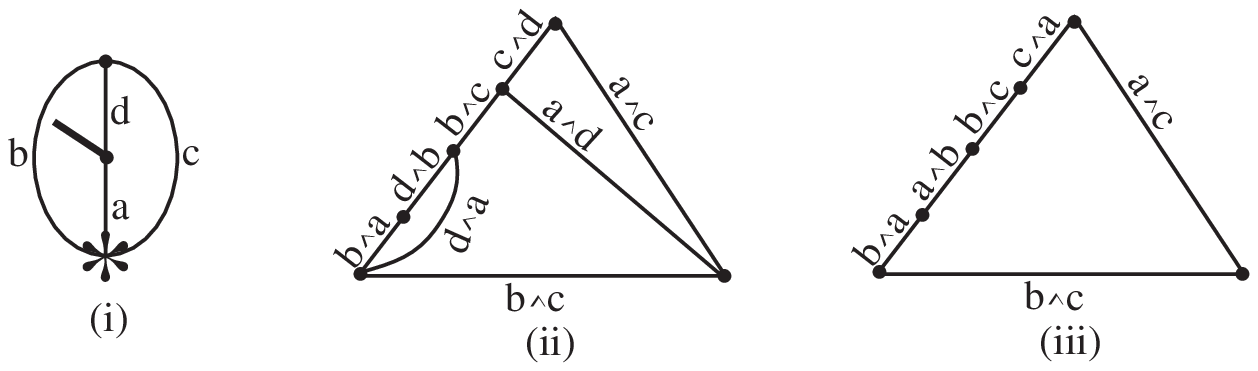}}}
\caption{}
\label{nodove}
\end{figure}
If $G$ has just three edges, then $\sigma_1'=0$, by definition. The sides of $\sigma_1'$ are still defined
though
and, given the first choice of edge and vertex of $f(\sigma_1)$ made above, are the sides of the triangle 
in figure
\ref{nodove}(iii).  Here each edge in $X$ of a side corresponds to the marked graph $G$ with the
indicated partition.  Any consistant choice of colorings may be used.  The  first term of the side $b\wedge a,
a\wedge b, b \wedge c, c\wedge a$ is defined as the  visible  term of that side.  This completes the definition
of $\sigma_1'$.  We set $\sigma_1^c=\sigma_1'$ and $\sigma_1^w=0$.

Now assume that $\sigma_i'$, along with its sides and one visible term per side, has been
defined for each $i<k$. Also assume that the following hold for  $i<k$:
\begin{enumerate}
\item[(a)] Forgetting colorings defines a bijection from the visible terms of the sides of $\sigma_i'$ to
the edges of $f(\sigma_i)$.
\item[(b)] If $i,j<k$, and  $\sigma_i$ and $\sigma_j$ share an edge, then the sides of $\sigma_i'$ and
$\sigma_j'$ corresponding to this edge are the same (up to a possible change in the orientation of the
edgepath constituting the side) and have the same visible term.
\item[(c)] Each of the two subsequences of each side of $\sigma_i'$ consisting of the invisible terms 
coming before and after the single visible one corresponds 
to a (possible empty) word in the generators of the group $W$.
\end{enumerate}

To define $\sigma_k'$, for $1< k<n$, so as to satisfy (a), (b), and (c), first assume that the edge
 $e$ of $\sigma_k$ is not contained in $\cup_{i<k}\sigma_i$. Let $s_j$ be the sides of the $\sigma_i'$,
for $i<k$, which correspond to the edges of $\sigma_k\cap \cup_{i<k}\sigma_i$.  If all the terms of $s_j$
are visible, then $\sigma_k'$ is defined as $\sigma_1'$ was, but using $\sigma_k$ in place of $\sigma_1$.
In case $f(\sigma_k)$ is a $w$-triangle, let $e$ be the chosen edge  of $\sigma_k$ so that the side of
$\sigma_k'$ corresponding to $e$ contains all the invisible terms of the sides of $\sigma_k'$.  Also if
$\tau\in \sigma_i$ is any edge of $\sigma_k\cap \cup_{i<k}\sigma_i$, then the colorings of the marked graphs 
associated with $\sigma_k'$ are defined to be those consistent with the coloring of the graph associated 
with the side of $\sigma_i'$ corresponding to $\tau$.

If any of the $s_j$ contain invisible terms, then $\sigma_k'$  is defined as $\sigma_k^c+\sigma_k^w$ where 
the term $\sigma_k^w$ corresponds to 2-cells of $X$ representing standard $w$-crossings.  The marked 
partitioned graph or graphs corresponding to $\sigma_k^c$ are, by definition, the same as those of $\sigma_k'$
as defined above in the case where the $s_j$ contain no invisible terms.  The side or sides of $\sigma_k'$
different from the $s_j$ contain one triple of terms corresponding to a generator of $W$ for each
such triple in the $s_j$.  In case $\sigma_k$ is a $w$-triangle, there is one additional such triple.
As usual, colorings are determined by those of graphs of any $\sigma_i'$ where $i<k$ and 
$\sigma_i \cap \sigma_k\ne\varnothing$. 

Figure \ref{lifts} shows how the cells of $\sigma_k^w$ and $\sigma_k^c$ are arranged in two examples 
and so should clarify the definition of $\sigma_k'$.  The outer square in (i) and the outer triangle 
in (ii) represent two possible boundaries of $\sigma_k$.   We assume, in each case, that only the bottom
edge of $\sigma_k$ is not one of the 
$s_j$.  Additional marks in the figure are superimposed on the $\sigma_k$ and indicate the cells
in $X$ corresponding to $\sigma_k$ as follows.  The dotted lines are $w$-lines.  Where these
cross the edges of $\sigma_k$ indicate the positions of the corresponding triples of terms in the
sides of $\sigma_k'$.  The solid lines which cross the edges of the $\sigma_k$ indicate the 
position of the visible terms in the sides of $\sigma_k$.  Crossings of these solid lines with the 
dotted $w$-lines indicate standard $w$-crossings.  Finally, the 2-cell of $\sigma_k^c$ is indicated
in each of (i) and (ii) by its edges, lightly drawn.  Thus (ii) treats the case where $f(\sigma_k)$
is a $w$-triangle.

\begin{figure}[tb]
\centerline{\mbox{\includegraphics*{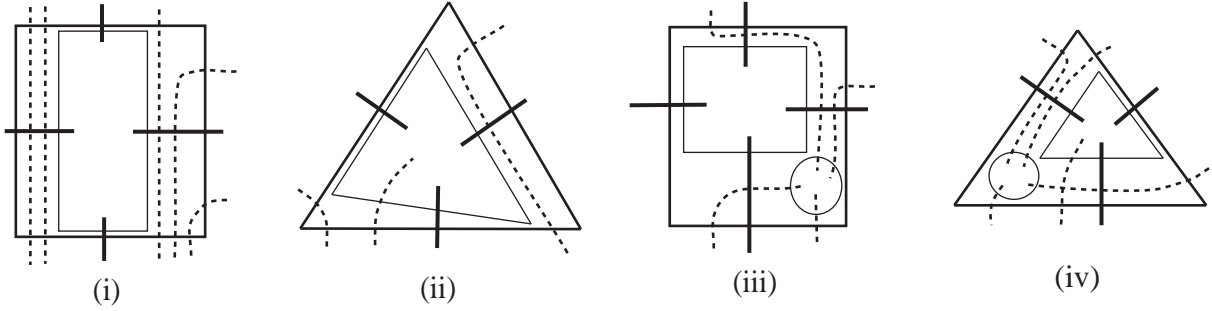}}}
\caption{In each figure the inner rectangle or triangle is $\sigma_i^c$.}
\label{lifts}
\end{figure}

Finally we define $\sigma_n'$.  The $s_j$, now defined as the sides in $\sigma_k\cap \cup_{i<n}\sigma_i$,
form a loop (representing a relation in  the 
group $\autfn$ for some $n$) in $X$ which will bound $\sigma_n'$. Again $\sigma_n'$ is defined as
$\sigma_n^c+\sigma_n^w$ with $\sigma_n^c$ and $\sigma_n^w$ as above, except that now $\sigma_n^w$ may also
contain 2-cells of
$X$ representing standard $ww$-crossings.  Parts (iii) and (iv) of figure \ref{lifts}, 
drawn with similar
conventions as used in parts (i) and (ii), should clarify this.  The circles in each indicate a loop in $X$
corresponding  to a word in the generators of the group $W$ which represent the trivial element of $W$.  The 
2-cells of $\sigma_n^w$ contained in these circles correspond to the defining relations of $W$ used
to trivialize this word.  This completes the definition of $p(f)$ except that we allow $\sigma_1^w$ 
to be nonzero provided the cells of $\sigma_1^w$ make up standard $w$-crossings.


Recall that $P(f)\colon S^2\to \xt$ is either lift of $p(f)\colon S^2\to X$ and also denotes the corresponding spherical 2-cycle in $\xt$.

Parts (i) through (iv) of figure \ref{lifts} suggest a natural subdivision of the $CW$-structure $\cc$ of
$S^2$ which makes the map $P(f)\colon S^2\to \xt$ cellular.  This subdivision has 2-cells, called core 2-cells,
correponding to the $\sigma_i^c$. The remaining 2-cells make up $w$- and $ww$-crossings.  The map $P(f)$ is
not, in general, cell-to-cell with respect to this subdivision since it takes  core 2-cells
$\sigma_i^c$ with $f(\sigma_i)$ correponding to graphs with three edges to the 1-skeleton of $X$.
But otherwise $P(f)$ restricts to a homeomorphism on each closed 2-cell of $\cc$.

The map $P(f)$ is sometimes specified by just giving the core cells in $S^2$ and the colorings of the 
corresponding  marked graphs, 
 along with solid lines 
and dotted lines as in (i) through (iv) of figure \ref{lifts}.  A solid line joins two sides of two core cells
if these two sides correspond to the same edge of $\cc$.  The dotted $w$-lines cross the solid lines and 
each other indicating $w$- and $ww$-crossings, and begin and end in core cells corresponding to
$w$-triangles.

\begin{definition}\label{spherical} \rm
A 2-cycle $z$ in a $CW$-complex $Y$ is called a simple spherical cycle if there is a cell-to-cell map
$f\colon (S^2,\cc)\to Y$ such that $z=\sum_i gen(f(\sigma_i))$.  Here
the $\sigma_i$ vary over the 2-cells of $\cc$ and the terms $gen(f(\sigma_i))$ are the generators of the 
cellular chain group $C_2(Y)$
corresponding to the 2-cells $f(\sigma_i)$, suitably oriented.  The 2-cycle $z$ is called a spherical
cycle if it is a sum of simple spherical cycles.
\end{definition}

The following technical lemma will be used often.

\begin{lemma}\label{added}
If $Y$ is any simply-connected simplicial complex, then, using the simplices of $Y$ as cells,
each element of $H_2(Y)$ is represented by a spherical cycle.
\end{lemma}

\begin{proof}
Let $z\in Z_2(Y)$ be a simplicial 2-cycle.  Assume $f'\colon S^2\to Y$ is any map representing
the homology class $[z]\in H_2(Y)$, meaning that $f'_*(1)=[z]$ for one of the 
generators $1$ of $H_2(S^2)$. Approximate $f'$ with a simplicial map $f\colon S^2\to Y$, and then
identify points in each simplex of $S^2$ which have the same image under $f$.  The resulting 
quotient space $S^2/\!\!\sim$ has a natural simplicial structure and, if $f'$ is not a constant map,
consists of 2-spheres and 
1-simplices which meet in discrete points. Moreover, the map $g\colon S^2/\!\!\sim \,\to Y$ induced by $f$
is  simplicial and is, in fact, cell-to-cell using the simplices of $S^2/\!\!\sim$ and $Y$ as cells. 
So $g$ restricted to each 2-sphere of $S^2/\!\!\sim$ is also cell-to-cell.  Thus, if  
the 2-simplices of $S^2/\!\!\sim$ are denoted by $e_i$ and
oriented so that the quotient map $q\colon S^2\to S^2/\!\!\sim$
preserves orientation (with $S^2$ oriented by the generator $1\in H_2(S^2)$ mentioned above), then 
the 2-cycle $z'\equiv \sum_i g(e_i)$ is spherical and $[z']= g_*(q_*(1))=f_*(1)=[z]$ in $H_2(Y)$.
\end{proof}

Now let $z\in Z_2(\ct)$ be any spherical 2-cycle.  Choose cell-to-cell maps 
$f_j\colon(S^2,\cc_j)\to \ct^*$ such that $z=\sum_{i,j} gen(f_j(\sigma_{ij}))$, where, for fixed $j$, the 
$\sigma_{ij}$ are the 2-cells of $\cc_j$.  Here the $gen(f_j(\sigma_{ij}))$ denote generators of 
$C_2(\ct^*)$ as in Def.~\ref{spherical}.  Going from $z$ to $f_j$ can be done by first using
Def.~\ref{spherical} to pick appropriate  cell-to-cell maps $S^2\to \ct$ and 
then subdividing the domains of these to get the $f_j$.  Let $P(z)=\sum_j P(f_j)$.

\begin{definition} \rm The group $H_2(\xt,\widetilde{W})_{StN}$   is the cokernel  of the composite $H_2(\widetilde{W})\to H_2(\xt)\to  H_2(\xt)_{StN}$.
\end{definition}

\begin{lemma}\label{BBB} The class of $P(z)$ in $H_2(\xt,\widetilde{W})_{StN}$
 is independent of the choices made in the
definition of $P(z)$.
\end{lemma}

Because of this, we often denote the class of $P(z)$ in $H_2(\xt,\widetilde{W})_{StN}$ by $P_H(z)$.
\begin{proof}[Proof of \ref{BBB}]
Let $z$ and $f_j$ be as above and let  $g_k\colon (S^2,\cd_k)\to \ct^*$ be another choice of cell-to-cell
maps such that $z=\sum_{k,l} g_k(\tau_{kl})$ where the $\tau_{kl}$ are the 2-cells of $\cd_k$.
Let $$\text{$P_1(z)=\sum_j P(f_j)= \sum_{i,j}(\sigma_{ij}^c+\sigma_{ij}^w)$ and
    $P_2(z)=\sum_k P(g_k)= \sum_{k,l}(\tau_{kl}^c+\tau_{kl}^w)$}.$$
Here the $\sigma$ and $\tau$ with superscripts and subscripts denote consistent lifts of the like denoted
cells in the definition of $p(f)$.
So, for instance, if the 2-cell
$g_k(\tau_{kl})$ of $\ct^*$
is not a $w$-triangle, then the marked colored graph associated with the 2-cell $\tau_{kl}^c$ is,
after forgetting colorings, the marked partitioned graph associated with $g_k(\sigma_{kl})$. 
Thus by using
lemma \ref{T} to surround core cells with dotted circles, we may assume that, if the 2-cells
$f_j(\sigma_{ij})$ and $g_k(\tau_{kl})$ of $\ct^*$  are the same and are not $w$-triangles, then the
corresponding core cells $\sigma_{ij}^c$ and $\tau_{kl}^c$ of $\xt$ are the same.  
Note that
applications of lemma \ref{T} may change the element of $H_2(\xt)_{StN}$ which $P_1(z)-P_2(z)$ represents, but not the element of $H_2(\xt,\widetilde{W})_{StN}$.
If $f_j(\sigma_{ij})=g_k(\tau_{kl})$ is a $w$-triangle, then  a deformation like
that indicated in figure \ref{swirl} may also be needed.  
One can show that the two pictures in this figure differ by the trivial element of $H_2(\xt)$ provided 
that the upper portion  of the dotted curve on the right is a suitably chosen swirl.
So we may now assume that $P_1(z)-P_2(z)$ consists entirely of cells making up swirls, and solid and dotted lines with 
standard $w$-crossings.  Following any solid line gives a loop in $X$ which can be eliminated 
(cf.~the proof of the broken circles lemma \ref{circle}).  Thus $P_1(z)-P_2(z)$ is homologous to a 2-cycle consisting of $w$-lines and standard $ww$-crossings.
Such a cycle represents the trivial element of $H_2(\xt,\widetilde{W})_{StN}$.
\end{proof}

\begin{figure}[tb]
\centerline{\mbox{\includegraphics*{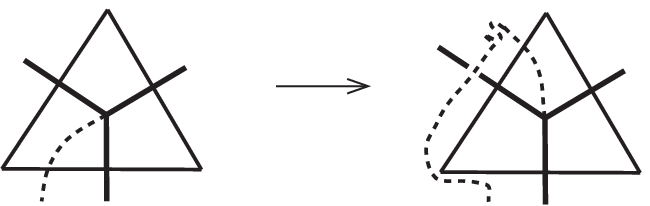}}}
\caption{}
\label{swirl}
\end{figure}

The definition of $P$ given above applys nearly verbatum to cell-to-cell maps $g\colon S^2\to T_M(G)$ 
(where $G$ is a fixed marked graph) and gives a 2-cycle in $\xt$ along with a map $S^2\to \xt$, both of 
  which we denote by $P'(g)$.
Since the map $G$ of section \ref{ct} takes degenerate simplices of the 2-skeleton of $\ct$ 
to the 1-skeleton of its  range,  it turns out to be convenient to also define $P(g)$ for maps
$g\colon S^2\to T_M(G)$ which restrict to cell-to-cell maps on the 1-skeleton of the domain and which
take each closed 2-cell either homeomorphically to a closed 2-cell or to the 1-skeleton.  Such maps
are referred to as almost cell-to-cell.  The only additional remark required to define $P'(g)$ for
an almost cell-to-cell map $g$ is that $\sigma_i'$ is defined to be 0 for each 2-cell $\sigma_i$
of $S^2$ on which $g$ is not injective.
Note that the previous lemma, with a simpler proof, holds with $P'(g)$ in place of $P(z)$.
We denote the class of $P'(g)$ in $H_2(\xt,\widetilde{W})_{StN}$ by $P_H'(g)$.

The rather lengthy proof of the next lemma is given at the end of this section.

\begin{lemma} \label{trivtree}
Let $f\colon S^2\to T_M(G)$ be an almost cell-to-cell map where $G$ is a fixed graph.  Then $P_H'(f)$ is the
trivial element of $H_2(\xt,\widetilde{W})_{StN}$ regardless of choices.
\end{lemma}

\begin{lemma} \label{Phom}
The association $P$ induces a homomorphism $P_H\colon H_2(\ct)\to H_2(\xt,\widetilde{W})_{StN}$.
\end{lemma}
\begin{proof}
Let $x\in H_2(\ct)$.
By lemma \ref{added}, there is a simplicial, spherical
2-cycle $z$ with $[z]= x$ in $H_2(\ct)$.
We show that setting $P_H(x)=P_H(z)$ gives a well-defined function 
$P_H\colon H_2(\ct)\to H_2(\xt,\widetilde{W})_{StN}$.

Assume $z'$ is another simplicial, spherical cycle with $[z']=x$.
Then $z-z'=\sum_i d\tau_i$ where the $\tau_i$ are generators of the group of simplicial 3-chains of
$\ct$.  (As usual, we also denote the corresponding 3-simplices by $\tau_i$.) 
Subdivide appropriate 2-simplices
so that  $z$, $z'$, and the $d\tau_i$ become spherical cycles in $\ct^*$.  Then, by 
 lemma \ref{BBB}, we have 
$P_H(z)=P_H(z')+\sum_iP_H(d\tau_i)$.  

With $\partial\tau_i$ denoting the topological boundary of the
3-simplex $\tau_i$, let $h_i\colon \partial\tau_i\hookrightarrow \ct$ be the inclucion map.  
Note that $P(d\tau_i)=P(h_i)$ in $H_2(\xt,\widetilde{W})_{StN}$
by the definition of $P$ on cycles and the cell-to-cell maps representing them.
Let
$g_i$ be the composite $\partial\tau_i\hookrightarrow \text{image of $h_i$}\to T_M(G_{\tau_i})$.  Here 
 the first arrow is
$h_i$ with range restricted, and the second is the restriction of the map $G$ (from the 2-skeleton of $\ct^*$
to $\bigsqcup_{\sigma\in\ct}T_M(G_\sigma) /\!\sim$) discussed at the end of section \ref{ct}.
  Then, with suitable
choices, we have $P(h_i)=P'(g_i)$ in $Z_2(\xt)$. So, by lemma \ref{BBB}, equality holds in
$H_2(\xt,\widetilde{W})_{StN}$ regardless of choices.
 Since $P'(g_i)=0$ by the previous lemma, each $P_H(d\tau_i)=0$, and so 
$P_H$ is well-defined.  It is also a homomorphism by
 lemma \ref{BBB}.
\end{proof}

The following theorem provides a crucial part of the computation of $H_2(\xt)$.
Let $F$ be the  composite $\xt{\buildrel {\mathcal F}\over \to}{\mathcal S} 
{\buildrel q \over  \to} {\mathcal T}$, where ${\mathcal F}$ is the forgetful function and $q$ the
quotient map, both defined near the beginning of \S\ref{SCT}. 

\begin{theorem} \label{thminW}
The composites $H_2(\xt)\buildrel{F_*}\over{\to} H_2(\ct)\buildrel{P_H}\over{\to}H_2(\xt,\widetilde{W})_{StN}$ and
$H_2(\xt)\to H_2(\xt)_{StN}\to H_2(\xt,\widetilde{W})_{StN}$ agree.

\end{theorem}

\begin{proof}
We first show that $H_2(\xt)$ is generated by classes which are represented by simple spherical cycles,
and then show that the two homomorphisms $H_2(\xt)\to H_2(\xt, \widetilde{W})_{StN}$ indicated 
in the statement of the theorem agree on such generators.

  Let $x\in H_2(\xt)$. Since $\xt$ is simply-connected, there is a map $S^2\to \xt$
representing $x$ with image contained in the 2-skeleton $\xt^{(2)}$ of $\xt$.  Make $\xt^{(2)}$ into a
simplicial complex
$\xt^{(2)}_*$ by adding one (of either) diagonal to each square 2-cell of $\xt^{(2)}$. Then, by lemma
\ref{added},
$x=\sum_i [x_i]$ where each $x_i$ is represented by a composite $S^2\buildrel{f_i}\over{\to}
\xt^{(2)}_*\hookrightarrow \xt$ with $f_i$ cell-to-cell.  

Fix a value of the subscript $i$.  Let $\{c_l\}$
be the set of 1-simplices which $f_i$ takes to the added diagonal edges of $\xt^{(2)}_*$. Write 
$\{c_l\}=\{d_j\}\sqcup\{e_k\}$ in such a way that there is a neighborhood of an interior point of
$c_l$ on which $f_i$ is injective if and only if $c_l$ is one of the $e_k$. For  each $j$, the union 
of the two 2-simplices $a_j$ and $b_j$ which meet along $d_j$ is a square. Subdivide each of $a_j$ and 
$b_j$ into two 2-simplices by adding the other diagonal $d_j'$ to this square.  Let $S^2_*$  be the 
resulting simplicial complex and let $m_j$ be the vertex of $d_j\cap d_j'$ of $S^2_*$. Define 
a simplicial map $f_i'\colon S^2_*\to \xt^{(2)}_*$ which agrees with $f_i$ on all the vertices of $S^2_*$
except the 
$m_j$.  Define $f'_i(m_j)$ 
to be the common image under $f_i$ of the two endpoints of $d_j'$.
Note that $f_i$ and $f_i'$ are homotopic. 

 Next delete the 1-simplices $e_k$ from $S^2_*$ so as 
to give a $CW$-complex $\cc$ with one square 2-cell corresponding to each $e_k$.  Note that $f_i'$, 
defined on $\cc$ in the obvious way, followed by the map $\xt^{(2)}_*\to \xt^{(2)}$, takes each of these
square cells homeomorphically to a square 2-cell of $\xt$. Thus, identifying points in cells of $\cc$
with the same image under $f_i'$ shows, as in the proof of \ref{added}, 
that $[x_i]=\sum_j[z_j]$ where each $z_j$ is a
cycle represented by a cell-to-cell map $S^2\to \xt$.  It follows that $H_2(\xt)$ is generated by 
classes represented by simple spherical cycles.

For the remainder of the proof, assume that  $x\in H_2(\xt)$ equals $[z]$ where $z$ is represented by
the cell-to-cell map $f\colon S^2\to \xt$. Subdivide the cells of $S^2$ giving $S^2_*$ 
so that $f\colon S^2_*
\to \xt$ is cell-to-cell using the simplices of $\xt$ as cells in the range.  Then the composite 
${\mathcal F}\circ f\colon S^2_*\to \ct$ is  simplicial. 
 As usual, let $q'\colon S^2\to  S^2/\!\!\sim$ be the quotient map
defined by identifying points in each simplex of $S^2_*$ with the same image under 
 ${\mathcal F}\circ f$, and let $g\colon S^2_*/\!\!\sim \,
\to\xt$ be the map induced by $F\circ f$.  At most one of the two 1-simplices making up
 each 1-cell of $S^2$ is
contracted to a point  by $q'$. Also, as one can  check, $q'(e_i)$ is a disk for each 2-cell
  $e_i$ of $S^2$. It
follows that $S^2_*/\!\!\sim$ is homeomorphic to a 2-sphere. So $P(g)$ is defined and
$P(g)=P_H(F_*(x))$ since the cell-to-cell map
$g$ represents $F_*(x)$.

For each 2-cell $e_i$ of $S^2$, let $G_i^{par}$ be the marked partitioned
graph, with coloring forgotten,
corresponding to the 2-cell $f(e_i)$ of $\xt$, and let $G_i$ be the same, but with the partition also
forgotten.  Denote the 2-cell of $T_M(G_i)$ corresponding to $G_i^{par}$ by $cell(G_i^{par})$.
Orient this cell so that the map from the 2-cell  $f(e_i)$ to it, defined by forgetting coloring,
preserves orientation. Here we are assuming that the cell $f(e_i)$ is oriented so that $f$ preserves
orientation, with $S^2$, in turn, oriented by the generator 1 of $H_2(S^2)$ such that $f_*(1)=x$.
Let $gen(G_i^{par})$ be the generator of the cellular chain group $C_2(T_M(G_i))$ corresponding to the
oriented 2-cell $cell(G_i^{par})$.

Let $\tau_{ij}$ be the 2-cells of $T_M(G_i)$ in the image of the composite $h\circ F\circ f|_{e_i}$, 
where $h$ is the map, with domain and range suitably restricted, from $\ct$ to $\ct_*$ followed by
the map $G$ defined at the end of section \ref{ct}.  Let $gen(\tau_{ij})$ be the corresponding
generators of
$C_2(T_M(G_i))$ with orientations chosen,  
in the case where $F$ is not injective on any of the three or four 1-cells of the boundary of $f(e_i)$,
so that 
$y_i=gen(G_i^{par}) + \sum_j gen(\tau_{ij})$ is a cycle.  
In case $F$ is injective on any of the edges of $f(e_i)$, let 
$y_i=gen(G_i^{par}) + \sum_j gen(\tau_{ij}) + \sum_k w_i^k$, where the $w_i^k$ are generators corresponding
to $w$-triangles, with one for each 1-cell of $f(e_i)$ on which $F$ is injective, and again with orientations
such that $y_i$ is a cycle.

Each $y_i$ is represented by a nearly cell-to-cell map $f_i\colon S^2\to T_M(G_i)$, and we
define $P'(y_i)$ as $P(f_i)$.
It is not difficult to check that $\sum_i P'(y_i)= x - P(g)$ in $H_2(\xt, \widetilde{W})_{StN}$.  This follows,
for instance, by an argument like that used to prove lemma \ref{BBB}.  Since each $P'(y_i)=0$ by lemma
\ref{trivtree} and $P(g)=P_H(F_*(x))$, the theorem is proved.
\end{proof}

\begin{proof}[Proof of \ref{trivtree}]
We first prove the lemma in the special case where $f\colon \partial \tau\to T(G)$ is the cell-to-cell
inclusion of the boundary of a 3-cell $\tau$ of $T(G)$. If $\tau$  corresponds to a partitioned 
graph $\Gamma$, then the 2-cycle $P'(f)$ bounds a 3-cell in $\xt$, the one corresponding to $\Gamma$ with
marking and coloring those chosen in the definition of $P'$. So, in this case,
$P'(f)$ is trivial in $H_2(\xt)$ regardless of choices.

If $\tau$ corresponds to a $w$-tetrahedron, then the graph $\Gamma$ associated with $\tau$ has two vertices,
the basepoint $*$, and $v$ say.  Let $a$, $b$, $c$, and $d$ be the four edges making up the four
maximal trees corresponding to the vertices of $\tau$.  
With suitable choices $P(f)$ consists of the 2-cells, 
suitably oriented, shown in figure  \ref{many1}.  The graphs associated with these 2-cells are drawn near 
them.  The heavily drawn dangling edge in these graphs represent all the halfedges incident with $v$,
other than $a$, $b$, $c$, and $d$.

\begin{figure}[tb]
\centerline{\mbox{\includegraphics*{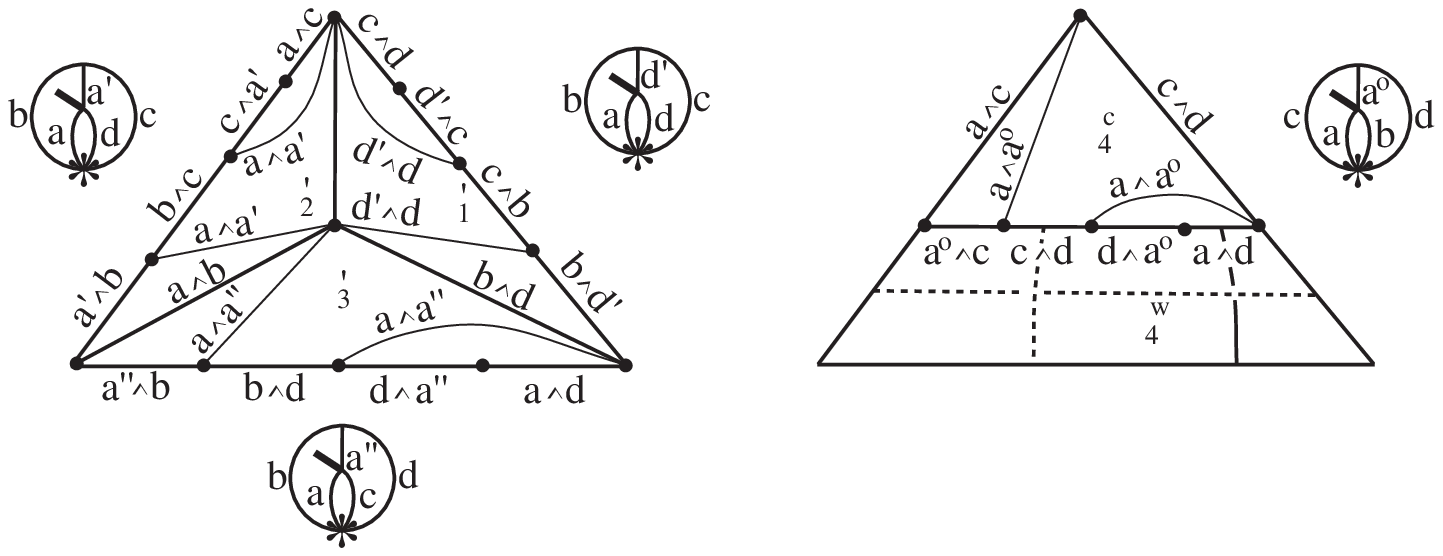}}}
\caption{}
\label{many1}
\end{figure}

We next show that we may assume that $\Gamma$ has only the four edges $a$, $b$, $c$, and $d$.
More specificaly, let $G'$ be any graph having as some blowdown a graph $\Gamma'$ with just two 
vertices and four edges, each edge of which meets both vertices.  Also let $\tau'$ be the 
$w$-tetrahedron of $T(G')$ corresponding to $\Gamma'$.  
With $f'\colon \partial \tau'\to T(G')$ the inclusion, we will show that, with suitable choices,
$P(f')=P(f)$ in $H_2(\xt)$.

Using 3-cells in $\xt$ corresponding to the graph $G$ in figure \ref{many2} on the right,
suitably marked and with various
colorings, one can see that the 2-chain $c_1$ in the figure,
 oriented so that $dc_i=d\sigma_1'$,
is homologous to $\sigma_1'$.  More specifically, $c_1-\sigma_1'$ is the boundary of a 3-chain with
six terms. These  correpond to the graph $G$ with each of the six partitions
$\{\{c,d,d'\},\{a,y\}\}$, $\{\{b,c\},\{d,d'\},\{a,y\}\}$, $\{\{b,d,d'\},\{a,y\}\}$, 
$\{\{a,d,y\},\{c,d'\}\}$, $\{\{a,d,y\},\{b,c\}\}$, and $\{\{a,d,y\},\{b,d'\}\}$.  Colorings are determined 
by the coloring at any vertex common to $\sigma_1'$ and $c_1$.  Similarly the 2-chain $c_2$ indicated in the
left of figure \ref{many2} is homologous to $\sigma_2'$ by way of the three 3-cells of $\xt$ corresponding to
the graph  on the left in the figure, suitably marked, with the three partitions  $\{\{a,a',c,x\}\}$, 
$\{\{a,a',x\},\{b,c\}\}$, and $\{\{a,a',b,x\}\}$. Also define 2-chains $c_3$ and $c_4$ similar 
to $c_2$ in such a way that $c_3$ and $c_4$ are homologous to $\sigma_3'$ and $\sigma_4^c$,
respectively, by way of graphs like those in Fig.~\ref{many1} corresponding to $\sigma_3'$ and $\sigma_4^c$,
but with the edges labeled $a$ subdivided by a vertex and with the dangling halfedges attached to this
vertex.  Then with choices like those giving figure \ref{many1}, the cycle 
$P(f')=\sum_{i=1}^4 c_i +\sigma_4^w$ is homologous to $P(f)$.

\begin{figure}[tb]
\centerline{\mbox{\includegraphics*{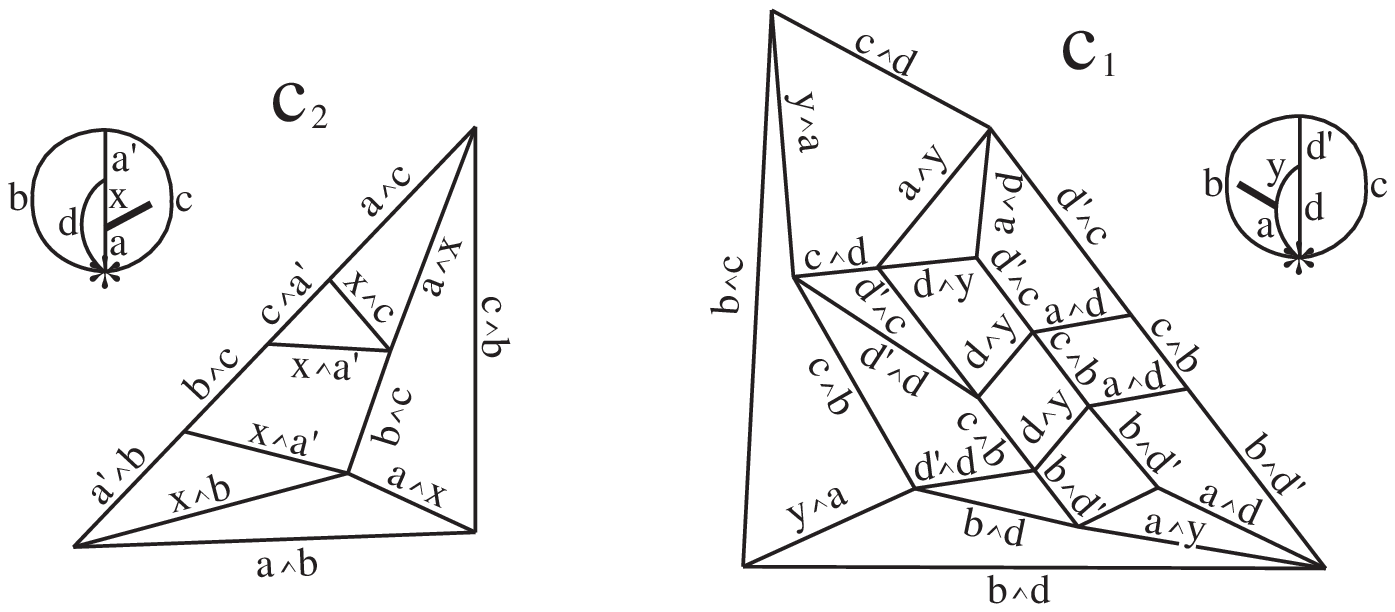}}}
\caption{}
\label{many2}
\end{figure}

So we now assume that $\Gamma$ has just four edges. Consider figure \ref{blocks} which shows three 3-cells
of $\xt$ corresponding to the colored graphs shown near them. The 3-disk these three cells form is 
denoted by $D$.  By the back face of $D$, we mean the three squares in the boundary of $D$ with one
edge $d\wedge a$ and the other edge one of $d'\wedge b$, $b\wedge c$, or $c\wedge d'$.  With suitable
markings, the nine faces in the boundary of $D$ other than the three back faces are homologous with
the standard $w$-crossing in $\sigma_4^w$.  Thus this $w$-crossing is also homologous to the three back faces 
of $D$. With the $w$-crossing represented in this last way, the portion of $P(f)$ other than
the $ww$-crossing in $\sigma_4^w$ is the negative of this $ww$-crossing.  Therefore, with the choices
above, $P'(f)=0$ in $H_2(\xt)$. So, by lemma \ref{BBB}, $P'(f)=0$ in $H_2(\xt,\widetilde{W})_{StN}$
regardless of the choices made.

\begin{figure}[tb]
\centerline{\mbox{\includegraphics*{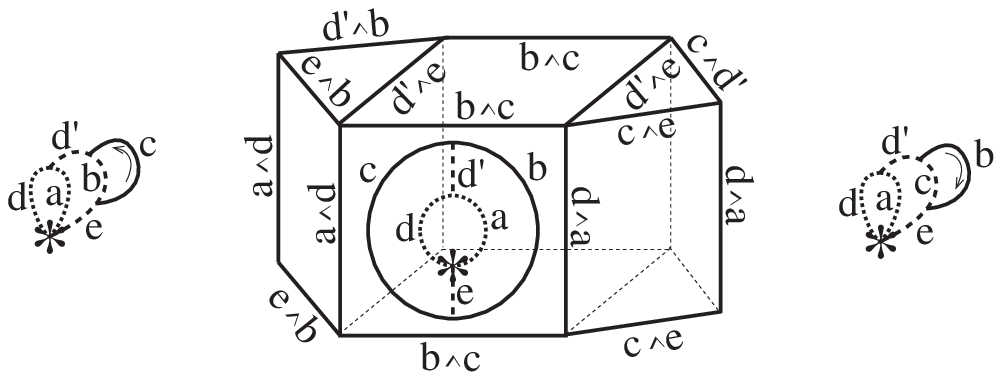}}}
\caption{}
\label{blocks}
\end{figure}

Now assume that $f\colon\partial\tau\to T(G)$ is the cell-to-cell inclusion of the boundary of 
a $w$-prism $\tau$ of $T(G)$.  Let $G'$ be the graph corresponding to $\tau$.  As a first case,
 assume $G'$ has a subgraph $\{a,c,d,x,y   \}$
as in figure \ref{graphs}(i) where the triangular faces of $\tau$ correspond to $G'/x$ and
$G'/y$.  
Possible additional edges of $G'$ are indicated by the dangling half edges in the figure.  There 
may be still more half edges incident with $*$.  But these remain so throughout the discussion, so
are not shown.
  Fig.~\ref{many3} illustrates $P(f)$ with certain choices.  Colored 
graphs are drawn within the cells they correspond to.  The three cells of $P(f)$ corresponding to
$G'/x$ are not shown, but we assume choices are such that these are the same as those corresponding 
to $G'/y$ except for the markings of the colored graphs.  Let $c$ be the cube of $\xt$ corresponding
to the colored graph in figure \ref{graphs}(iii) marked so that the vertex $v$ in figure 
\ref{many3}A is the vertex
of $c$ corresponding to the maximal tree consisting of the edges $b$, $d$, and $y$.  We will
show that $P(f)$ is homologous to the boundary of $c$. Let $p_d$, $p_y$, and $p_c$ be the O-prisms
shown in figure \ref{many3}B.  These correspond to the colored graphs drawn near them and are marked
so that the vertices $v$, $v'$, $w$, and $w'$  in this figure are the same as those with like labels
in figure \ref{many3}A. Let $p_x$ be the O-prism corresponding to the same colored graph as $p_y$, but 
marked so that the edge analogous to the one labeled $d\wedge c$ in $p_y$ joins $v'$ and $w'$.
Note that for $i=d,c,x,$ and $y$, the only square face of $p_i$ which is not contained in $P(f)$ is
the face of the cube $c$ corresponding to the graph in figure \ref{graphs}(iii) with the edge $i$ blown down.

Let $t_v$ and $t_w$ be the O-tetrahedra shown in figure \ref{many3}B, with $t_v$ the one incident with $v$.
Let $t_v'$ and $t_w'$ be the O-tetrahedra corresponding to the same colored graphs as $t_v$ and $t_w$,
respectively, but marked so that they are incident with $v'$ and $w'$, respectively.   Then, with suitable
orientations, we have 
$$ P(f)+d(p_d+p_y+p_c+q_y) +d(t_v+t_w+t_v'+t_w')=dc.$$
So, again by lemma \ref{BBB}, $P(f)=0$ in $H_2(\xt,\widetilde{W})_{StN}$ regardless of choices.

\begin{figure}[tb]
\centerline{\mbox{\includegraphics*{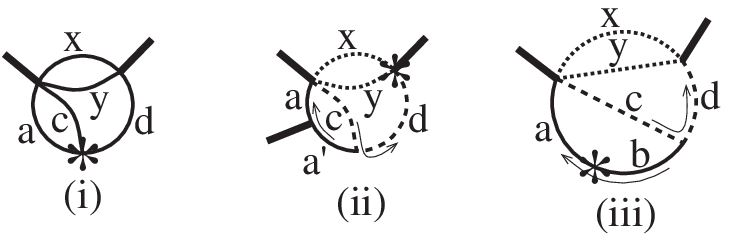}}}
\caption{}
\label{graphs}
\end{figure}

\begin{figure}[tb]
\centerline{\mbox{\includegraphics*{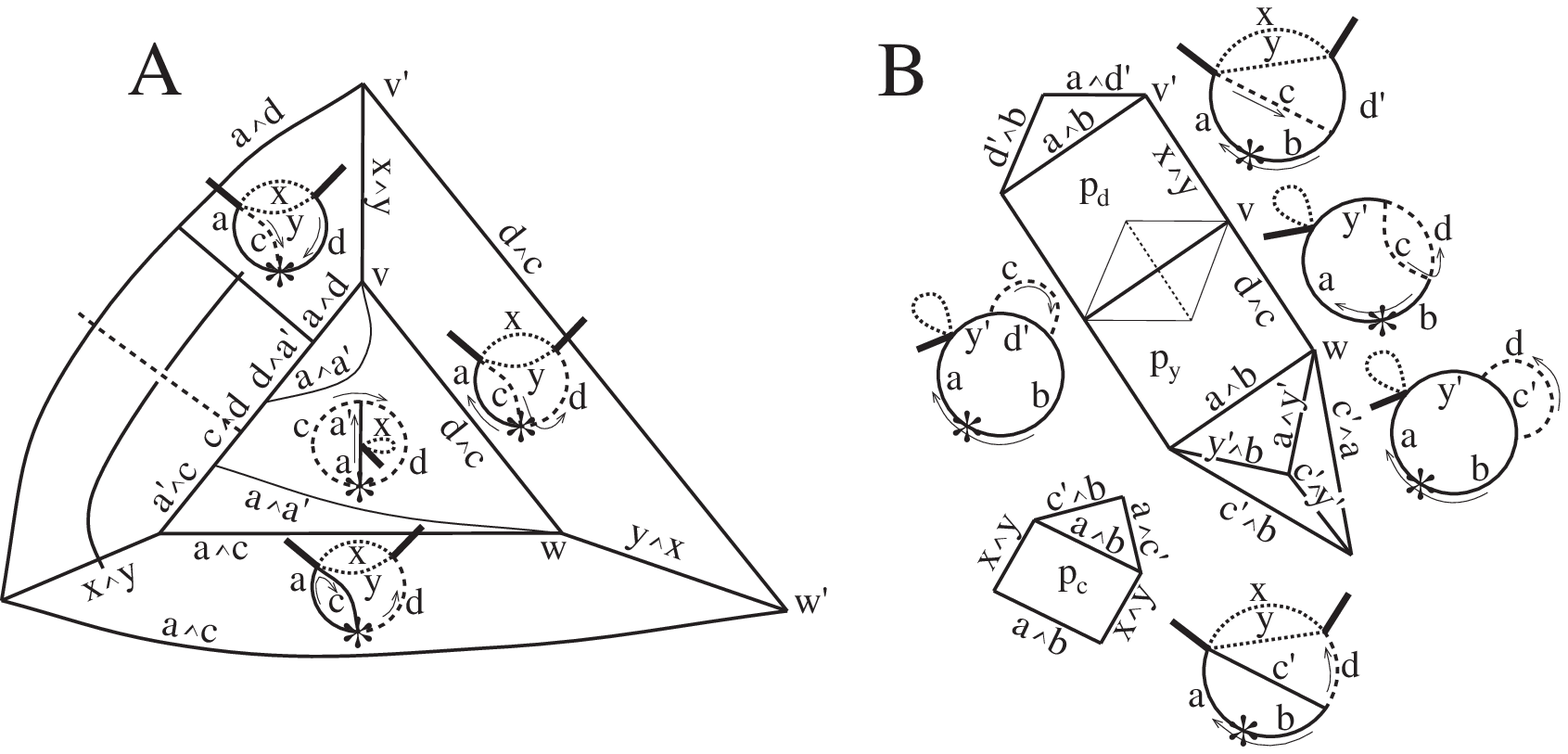}}}
\caption{}
\label{many3}
\end{figure}

For another case, assume $G'$ is as in figure \ref{graphs}(i), but with basepoint the vertex on the
the upper right.  With choices like those made in the previous case, $P(f)$ is as  in figure \ref{many3}A
except for the position of the basepoints and for the colored graph in  the middle of figure \ref{many3}A.
This should now be the colored graph $\Gamma/y$ where $\Gamma$ is the colored graph in figure
\ref{graphs}(ii) with suitable dangling halfedges. 
Let $c$ be the cube in $\xt$ corresponding to $\Gamma$, marked so that 
$\Gamma/\{a',d,y\}$ is the vertex of $P(f)$ analogous to the vertex $v$ in figure \ref{many3}A. We
show that $P(f)$ is homologous to $dc$.  Let $\gamma$ be the graph in figure \ref{graphs}(ii) and let
$p_d$ be the O-prism in $\xt$ corresponding to $\gamma$ with the partition $\{\{a,a',d\},\{x,y\}\}$
and marked so that $\gamma/\{a',d,y\}$ is the vertex analogous to $v$.  Also let $p_c$ be the 
O-prism corresponding to $\gamma$ partitioned by $\{\{a,a',c\},\{x,y\}\}$ and marked so that 
$\gamma/\{a',c,y\}$ is the vertex analogous to $w$. Then, in $H_2(\xt)$, we have, with suitable 
orientations, $P(f)+d(p_c+p_d)=dc$ which completes this case.

The case in which the basepoint of $G'$ is the upper left vertex of the graph in figure \ref{graphs}(i)
is similar to the case just considered unless the vertex incident with the edges $a,c$, and $d$ is incident 
with no other edges.  Then choices can be made so that $P(f)$ is trivial as a 2-cycle in $\xt$.  The 
remaining cases in which $\tau$ is a $w$-prism are those where the edges $a$, $c$, and $d$ are all 
incident with the same two vertices.  These cases are easier than those already considered and are
ommited.

Also ommited are the cases in which $\tau$ is a pillow.  These are also easier than those above, especially
if choices are made in defining $P(f)$ so that there is a single $w$-line with no crossings.  Though in
some cases this $w$-line will pass through the cells which realize the relation $w_{ij}^L=w_{ij}^R$.

Now assume that $f\colon (S^2,\cc) \to T(G)$ is any cell-to-cell map.  We will show, with the help
of frequent implicit applications of lemma \ref{BBB}, that if $f$ is transformed to $g$ using 
any one of the operations (a) through (d) above, then $P(f)=P(g)$ in $H_2(\xt,\widetilde{W})_{StN}$.
This, along with lemma \ref{2conn}, will complete the proof.  First assume that operation (a) is used
to transform $f$ to $g$, and let $D=\cup_i\, \sigma_i$ and $\tau$ be as in the statement of (a).
Also let $i\colon \partial\tau\to T(G)$ be the inclusion of the boundary of $\tau$.  By, for instance,
choosing the $\sigma_i$ to come first in the orderings of the 2-simplices of $\cc$ and $\partial\tau$ 
used to define $P(f)$ and $P(i)$, we may assume that 
these agree when restricted to $D$. Thus, since $P(f)-P(i)=P(g)$, 
 the first part of the proof implies that $P(f)=P(g)$.

Next assume that $g$ is obtained from $f$ by an operation of type either (b), (c), or (d), and
let $\sigma$ and $\tau$ be as in the statement of (b), (c), or (d), as appropriate. In the case of 
(b) or (d), we may assume, with suitable choices, that $\sigma'=-\tau'$.  Here, and in the next few
lines, the primes have the same meaning as in the definition of $P$ in section \ref{defP}.
So it follows that 
$P(f)=P(g)$.  In the case of (c), we may assume, after introducing $w$-crossings and then $ww$-crossings
as needed, that $\sigma'+\tau'$ is a cycle which, as in the proof of \ref{BBB}, represents the
trivial element of $H_2(\xt,\widetilde{W})_{StN}$.  So again, $P(f)=P(g)$.  Here, $P(g)=P(g_1)+P(g_2)$
where $g_1$ and  $g_2$ are the two cell-to-cell maps which $f$ is transformed to by (c).
\end{proof}

\section{The computation of $H_2(\xt_n)_{StN_n}$} \label{lastcomp}
 
Let $ {\mathbb Z}^n$ be the free abelian group of rank $n$.


\begin{theorem}  \label{H2X}
If $n\ge 5$, then $H_2(\xt_n)_{StN_n}\approx \mathbb{Z} _{24}\oplus {\mathbb Z}^n$.
\end{theorem}  

\begin{proof}
As remarked after theorem \ref{thm1}, we have $ H_2(\xt_n)_{StN_n}=E_n$.
We begin by finding 
 a convenient set of generators for this  group. 
The element $w_{14}^L\cdot I_1^{-3,-4}$ of $H_2(\xt_n)_{StN_n}$ is represented by the picture $I_1^{-3,-4}$
surrounded by a circle labeled $w_{14}^L$ and with transverse orientation pointing inward.  By the 
three circles lemma \ref{3circles}, this circle can be moved past all the crossings in the picture, except for 
those corresponding to the relation $w_{ij}^L=w_{ij}^R$. There are five of these, but because of their signs,
four cancel.  So,
 in $H_2(\xt_n)_{StN_n}$, we have
$$ I_1^{-3,-4} = w_{14}^L\cdot I_1^{-3,-4}=  I_4^{3,-2} -C_{14}^2.$$
Adding another surrounding circle gives
$$ I_1^{-3,-4} = (w_{14}^L)^2\cdot I_1^{-3,-4}= 
I_1^{3,-2} -(C_{14}^2+C_{14}^4).$$

\begin{figure}[tb]
\centerline{\mbox{\includegraphics*{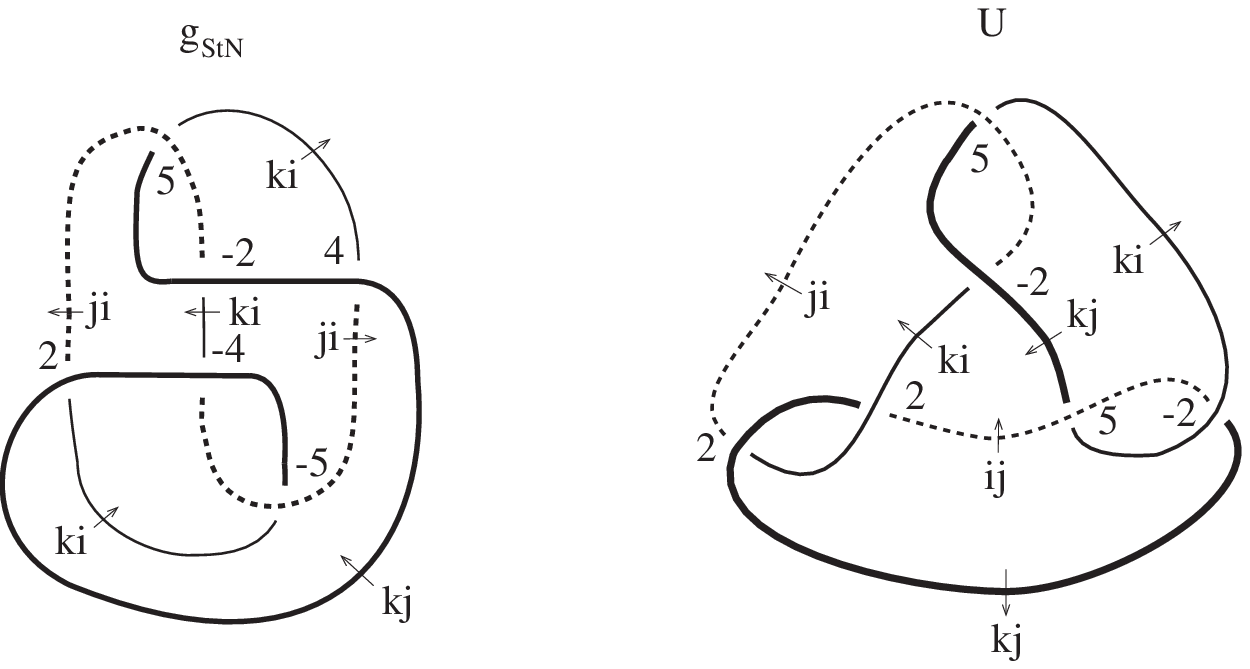}}}
\caption{}
\label{gstn}
\end{figure}

Now consider the picture $g_{StN}$ shown in Figure \ref{gstn}. 
Cutting it horizontally through the middle and adding the  crossing 3 to the top component 
and -3 to the bottom one gives the picture $-I_k^{-3,-4}+ w_{kj}^L\cdot (-I_k^{3,-2})$.
Thus
$$ g_{StN}=   -( I_k^{-3,-4} +I_k^{3,-2}) $$  in $H_2(\xt_n)_{StN_n}$.

It follows from the three displayed equations above and Theorem \ref{thm1} that the group $H_2(\xt_n)_{StN_n}$ 
is generated by the classes of
$I_i^{-3,-4}$, $i\in {\bf n}$ and $g_{StN}$.  The map $\Theta\colon H_2(\xt_n)_{StN_n}\to H_2(\xt_n/StN_n)$
takes the $n$ classes $I_i^{-3,-4}$, $i\in {\bf n}$ 
 onto the free abelian group $H_2(\xt_n/StN_n)$ of rank $n$.
Thus these generate a 
free abelian subgroup of $H_2(\xt_n)_{StN_n}$ of rank $n$.  The remaining generator $g_{StN}$ is in the kernel of  $\Theta$ and so generates it.
One can show that 
 $I_1^{-3,-4}$ and $I_1^{3,-2}$ are sent  by the map $H_2(\xt_n)_{StN_n}\to H_2(\vt)$
  to the generator  of $K_3(\mathbb Z) \approx \mathbb Z_{48}$ exhibited by Igusa in
Fig.\ 7.6 of   \cite{I}. Thus
  the class of $g_{StN}$ in $H_2(\xt_n)_{StN_n}$
has order at least 24.  (Alternatively, Hatcher observed in \cite{H} that Waldhausen's splitting theorem
  implies that the homology of the group $\rm{AUT}=\lim_{n\to \infty}$Aut$(F_n)$ contains that 
  of the symmetric group $\Sigma=\lim_{n\to \infty}\Sigma_n$ as a direct summand. Similarly for homotopy groups.  In particular, 
  $H_3(StN) =\pi_3(B{\rm StN}^+)=\pi_3(B\rm{AUT}^+)  $  contains 
  $\pi_3(B\Sigma^+)=\pi_3^S=\mathbb Z_{24}$ as a direct summand.
  As in the proof of Theorem \ref{thm} below,  the kernel of $\Theta$ is $H_3(StN)$.)
  So to conclude the proof we show that  the class of $g_{StN}$ in $H_2(\xt_n)_{StN_n}$ has order at most 24.

We will use the fact that $g_{StN}=U$ in $H_2(\xt_n)_{StN_n}$ where $U$ is as shown in Figure \ref{gstn}.
To see this, note that the upper left portions of both of these pictures are the same.  So they
cancel in the difference $g_{StN}-U$ leaving a circular edge  containing only the two swirls $-A4$, 
 shown in Figure \ref{4swirls} and $A2$, which is the knot-like swirl at the bottom of the picture at
 the right of Figure \ref{swirls}.
  This circle is trivial in 
 $H_2(\xt_n)_{StN_n}$ as in the first case of the proof Theorem \ref{justsw}.

\begin{figure}[tb]
\centerline{\mbox{\includegraphics*{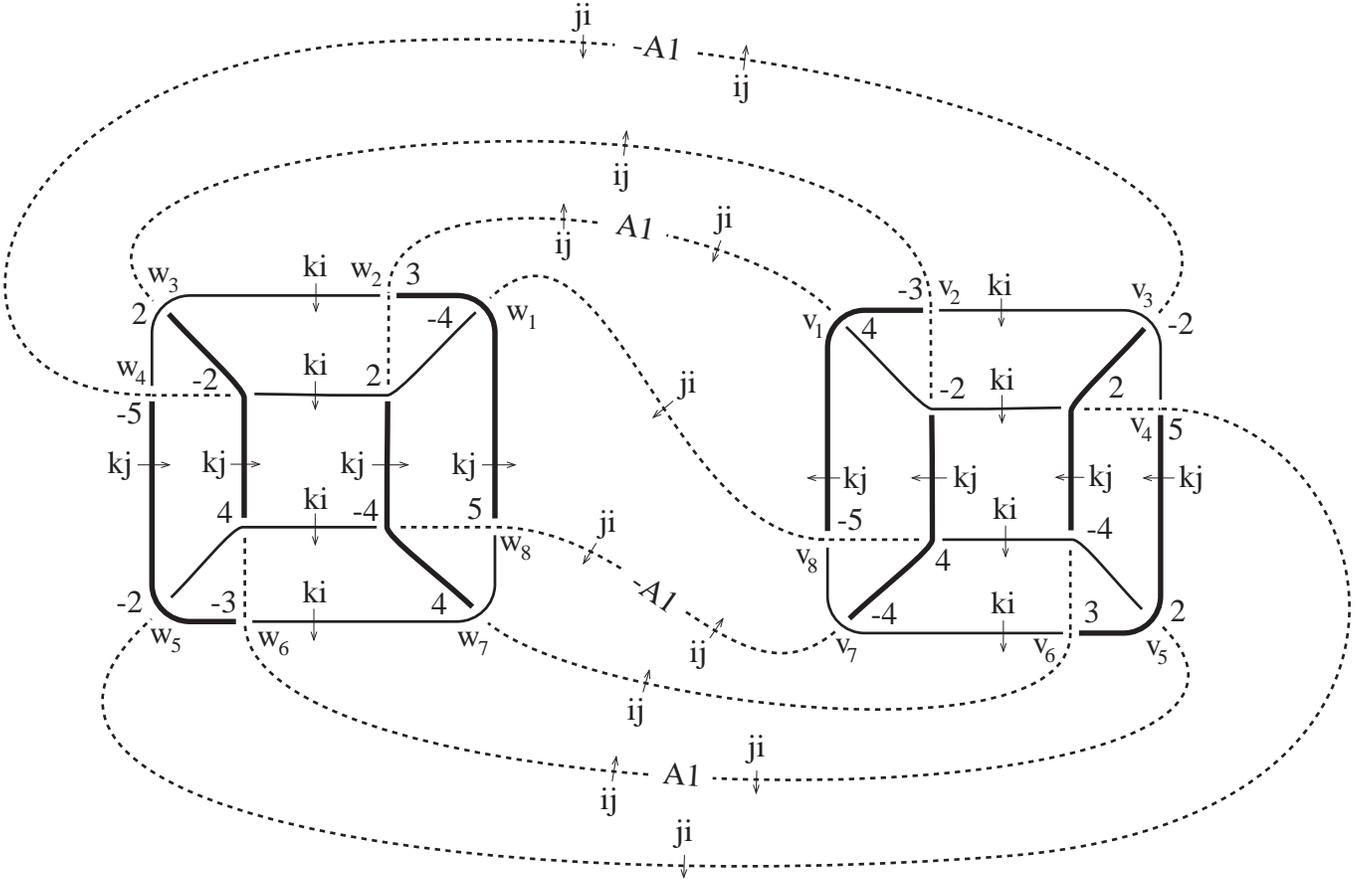}}}
\caption{The cycle S}
\label{S}
\end{figure}

Consider the picture $S$ shown in Figure \ref{S}.  The four swirls in the picture are all $A1$. Two have 
positive orientations and two negative ones as shown. We first show that $S=6g_{StN}$,
and then show that $4S=0$, which will complete the proof.
	
Note that the middle triangle of the picture $U$ in Figure \ref{gstn}, but with $i$ and $j$ permuted,
is the same as the upper right triangle  in $S$.  Similarly, the triangle in the lower right of $S$ is the
same as the inner one  of $U$ if the transverse orientations of the $i$-edges of $U$ are changed.
The other two triangles in the right-hand component of $S_k$ are the same as the inner triangle of two
other translates of $U$.  
 Since $\Theta(U)=0$, these four translates of $U$ are represented by pictures agreeing with $U$'s if 
 labels are forgotten.  So the class of $S-4U$ in $H_2(\xt_n)_{StN_n}$ is represented by a picture
 $(S-4U)'$ looking like the picture $S$ of Figure \ref{S} except all 12 crossings on the right are changed
 so that the broken edges are solid.  Using Lemma \ref{sm} each of the four swirls in  $(S-4U)'$ can
 be moved (in either  direction) past the first broken edge of  $(S-4U)'_k$ encountered.  
 Crossings then cancel in pairs leaving the picture $S-4U$ of Figure  \ref{Sm4U}.  Moving the swirls does not change
 the class of the picture in $H_2(\xt_n)_{StN_n}$ since all elements of $\langle O\!\!\!O  \rangle$ introduced 
 in moving the   swirls cancel.
\begin{figure}[tb]
\centerline{\mbox{\includegraphics*{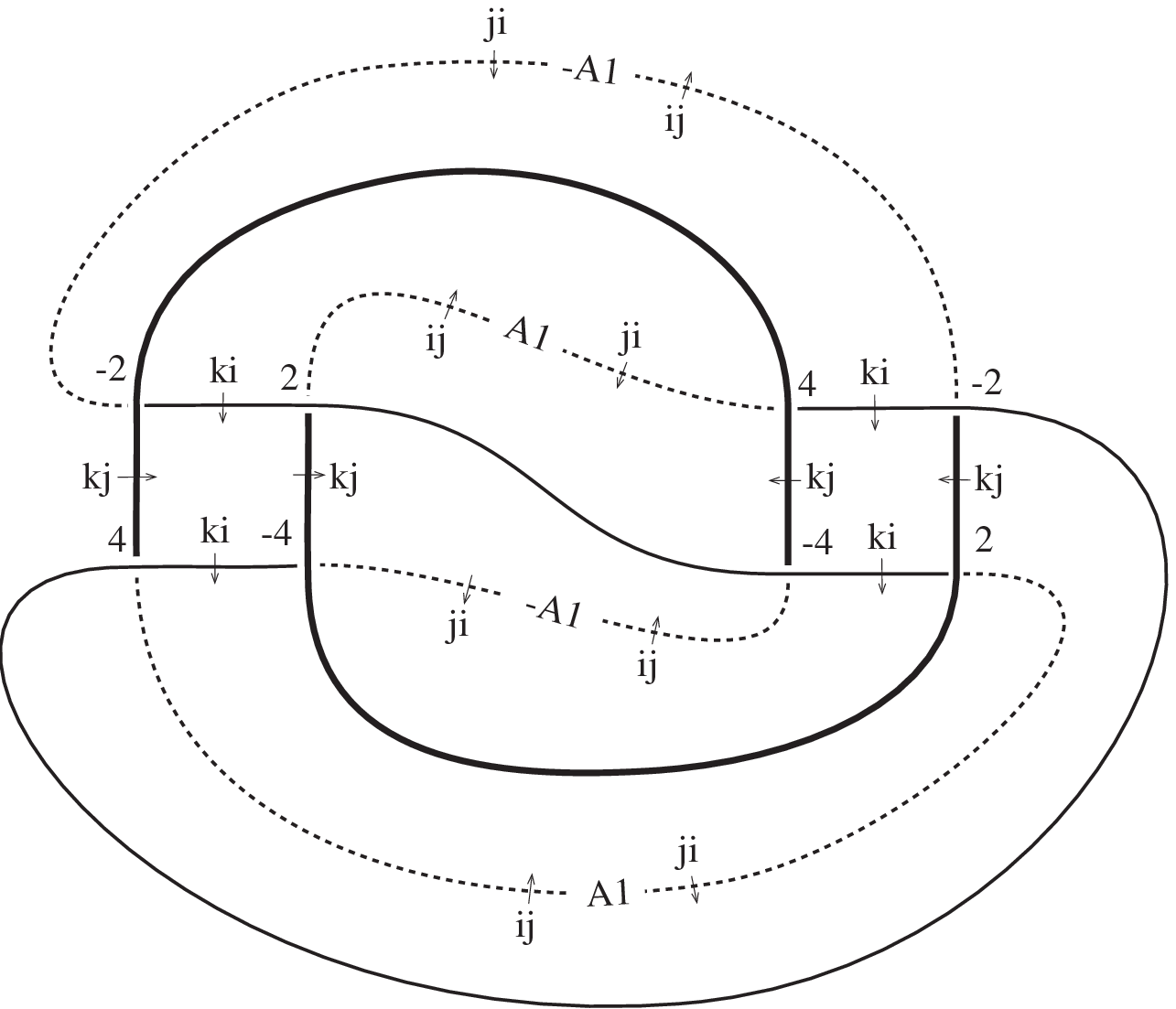}}}
\caption{The cycle $S-4U$}
\label{Sm4U}
\end{figure}
  
  \begin{figure}[tb]
\centerline{\mbox{\includegraphics*{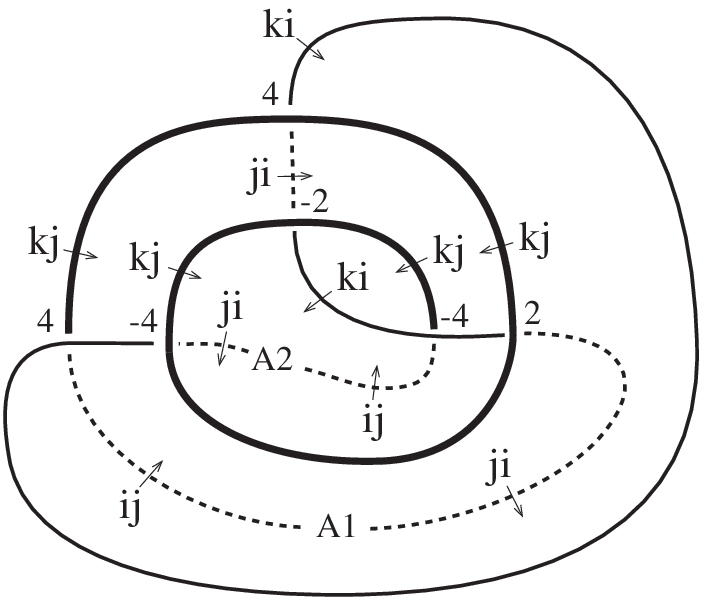}}}
\caption{The cycle $H_{A2}$}
\label{HA2}
\end{figure}
  
Now the lower half of the picture $H_{A2}$ in Figure \ref{HA2} agrees with the lower half of the picture $S-4U$,
except for the swirl $A2$.  Note that the two crossings labeled -2 and -4 which are part of one of the bigons in
$H_{A2}$ cancel with two of the crossings of the swirl $A2$ in $H_{A2}$.  Similarly, the two crossings
which are part of the other bigon cancel with crossings of $A1$, leaving the picture $-I_k^{3,-2}$.
Thus $H_{A2}=-I_k^{3,-2}$ at the chain level.

Let $H$ denote $H_{A2}$ with the swirl $A2$ changed to $-A1$ so as to agree with more of $S-4U$.
Then, in $H_2(\xt_n)_{StN_n}$, we have
$$H=H_{A2}-(A1+A2)=H_{A2}-I_k^{-2,3} =-(I_k^{3,-2}+I_k^{-3,-4})=g_{StN}=U.$$
The picture $(w_{il}^L)^2\cdot H$ is the same as $H$ except for the transverse orientations of its $i$-edges.
So a  portion of it agrees with the upper half of the picture $S-4U$, except that the two swirls of
$(w_{il}^L)^2\cdot H$ are $A4$ and $-A4$ rather than $A1$ and $-A1$.  But since the two swirls of $(w_{il}^L)^2\cdot H$
have opposite signs, they can be changed to $A1$ and $-A1$ without changing the  class of $(w_{il}^L)^2\cdot H$.
Let $H_{\Delta i}$ denote the picture $(w_{il}^L)^2\cdot H$ with its swirls changed in this way.  Thus, 
$H=H_{\Delta i}$ in $H_2(\xt_n)_{StN_n}$. Now $H+ H_{\Delta i}=S-4U$ at the chain level and so
$2U=2H=S-4U$ or $S=6U=6g_{StN}$ in $H_2(\xt_n)_{StN_n}$.

A picture representing $nS$ in $H_2(\xt_n)_{StN_n}$, where $n\in \mathbb Z$, looks, up to swirls, like the picture
$S$, except that 
if a dotted edge $e_i$ of $S$ joins the vertex $v_i$ to the vertex $w_{i+1}$, then in $nS$ there is an edge 
$e_i'$ joining  $v_i$ to $w_{i+n}$.
  Here $w_{i+n}=w_m$ if $i+n\equiv m \ \ (\text {mod}\ 8)$, and the swirls on $e_i'$ are those on the edges
$e_i,e_{i+1},\ldots,e_{i+n-1}$.
So, for example, $8S=0$ in $H_2(\xt_n)_{StN_n}$ even if $n=3$. (The four swirls 
on each of the eight edges which join the two components of $(8S)_k$ cancel in  pairs, 
  leaving the trivial cycle.)

To see that $4S=0$ if $n\ge 5$, use Lemma \ref{2c} to change the direction of the transverse orientations of all
of the $k$-edges in  one  of the two components of $(4S)_k$.  Then swirls again cancel leaving 
  the trivial cycle.
\end{proof}
%


The next theorem brings together a number of the main results.  It also incorporates 
Lee and Szczarba's result in \cite{LS} that $K_3({\mathbb Z})\approx {\mathbb Z}_{48}$.

\begin{theorem}\label{thm} If $n\ge 7$, there is a commutative diagram
$$
\xymatrix{
0\ar[r]  &\mathbb Z_{24} \ar[r]  \ar[d]^{\approx}    & \mathbb Z_{24}\oplus \mathbb Z^n \ar[r] \ar[d]^{\approx}    & \mathbb Z^n  \ar[r] \ar[d]^{\approx}_{\Theta^{-1}}     &   0\\
0\ar[r]  & H_3(StN_n)\ar[r]  \ar[d]^{\approx}    & H_2(\xt_n)_{StN_n}\ar[r] \ar[d]^f &H_2(\xt_n/StN_n)\ar[r] \ar[d] &   0\\
0\ar[r]  & \mathbb Z_{24} \ar[r]        & K_3(\mathbb Z) \ar[r]  & \mathbb Z_2     \ar[r]   & 0
}
$$ where $f$ is induced by abelianization.
All rows are exact, and four of the six vertical arrows are isomorphisms as indicated. The remaining two
are surjective. 
\end{theorem}
\begin{proof}
In the spectral sequence of the covering $\xt_n \to \xt_n   /StN_n  $ we have 
$E_{0,2}^2=E_{0,2}^3=H_2(\xt_n)$ and $E_{3,0}^2=E_{3,0}^3=H_3(StN_n)$ since 
$H_1(StN_n)=H_2(StN_n)=0$.  Also, since $H_3(\xt_n/StN_n)=0$, the differential 
$d^3:H_3(StN_n)\to H_2(\xt_n)$ is injective giving exactness of the middle row.
The other assertions follow from
Theorem   \ref{H2X} and its proof and Theorem \ref{h2q}. 
\end{proof}

 Exotic elements of $K_3(\mathbb Z)$   appear, 
in some sense, as nontrivial elements of the homology of triangular subgroups.
Using the notation of remark \ref{tt}, where this was discussed, we have

\begin{theorem}\label{last_thm}  The class
of the cycle $t^t$ in $H_2(T)$ is nontrivial.
\end{theorem} 

\begin{proof}  The map
$$\xymatrix{
\xt \ar[r]^{\tilde \mu} \ar[d]  &\xt^{alg}  \ar[d]  \\
\xt/StN \ar[r]   &\xt^{alg}/StN 
}$$
of covering spaces defined in \S \ref{xtoxalg} induces a map of the associated spectral sequences.  This in turn gives the 
 commutative diagram
 
 $$
\xymatrix{
0\ar[r]  & H_3(StN_n)\ar[r]  \ar[d]^=    & H_2(\xt_n)_{StN_n}\ar[r]^{\Theta} \ar[d]^{\tilde \mu_*} &H_2(\xt_n/StN_n)\ar[r] \ar[d]^{\bar \mu_*}     &   0\\
 & H_3(StN_n) \ar[r]^j \ar[dr]^d  &  H_2(\xt_n^{alg})_{StN_n}\ar[r]  \ar[d]^a & H_2(\xt_n^{alg}/StN_n)\ar[r] &   0\\
  &   &   K_3(\mathbb Z)& H_2(ET/T)\ar[u]^{i_*}
}
$$
with exact rows. By the choice of $t^t$ we have $i_*[t^t]=\bar \mu_*[t]$.  We claim this last class is
nontrivial.  Indeed, if $e\in  H_2(\xt_n)_{StN_n}$ satisfies $\Theta(e)=[t]$, then $a(\tilde\mu_*(e))$ is not
in the image of the injection $d$.  Thus $\tilde\mu_*(e)$ is not in the image of $j$, which gives $\bar\mu_*[t]\ne 0$.
So $[t^t]\ne 0$.

%
\end{proof}







\end{document}